# Generalization of Cycle Decompositions of Even Dimensional Hypercubes on $d-$Dimensional Toruses

Idael Martinez-Perez[*]

July 2023

### Abstract

We consider cycle decompositions of even, $2an$-dimensional hypercubes $Q_{2an}$, where $a \geq 3$ is odd and $n \geq 1$. Prior work done by Axenovich, Offner, and Tompkins focused on obtaining the existence of cycle decompositions for even-dimensional hypercubes using long cycles of a given form, leaving out cycles of shorter lengths and, in fact, cycles of even longer lengths than those obtained there, such as $C_{7 \cdot 2^{11}}$ in the case of $Q_{14}$. In this paper, we provide two novel methods for explicitly constructing cycle decompositions of virtually all possible cycle lengths, using cycles of a given form, on Cartesian products of cycles up to those known by the work of Axenovich, Offner, and Tompkins. In particular, we show that we can explicitly obtain cycle decompositions of even dimensional hypercubes $Q_{2an}$ for all lengths mentioned above while on the same Cartesian product of cycles. With this, the current understanding of cycle decompositions of even dimensional hypercubes is furthered constructively and is featured with some interesting consequences for when $a$ is a positive, even integer. Additionally, progress is made towards obtaining cycle decompositions using the longest admissible cycle lengths with the incorporation of a more explicit starting point from which such decompositions of $Q_{2an}$ can be studied further.

## 1 Introduction

Given a graph $G$, we denote the graph's set of vertices by $V(G)$ and those of edges by $E(G)$. The $n$-dimensional hypercube $Q_n$ is then the graph with $V(Q_n) = \{0,1\}^n$ and $E(Q_n)$ consisting of all vertex pairs differing in exactly one component. A subgraph $H$ of a graph $G$ is said to give an edge decomposition of $G$ if $G$ can be expressed as a pairwise edge-disjoint union of isomorphic copies of $H$. One class of subgraphs of major interest are cycles, precisely those subgraphs characterized for every one of their vertices being visited exactly once with two edges associated to each vertex. Those cycles that visit every vertex of the graph are then regarded as Hamiltonian cycles. Of equal importance are Cartesian products of graphs, specifically Hamiltonian cycles, used in defining structures known as tori on which we explicitly perform decompositions using cycles to ultimately obtain a cycle decomposition of the given torus and, more generally, that of a given hypercube expressed as a decomposition into isomorphic copies of the torus.

In the paper by Axenovich, Offner and Tompkins [4], a cycle decomposition method known as the *Wiggle decomposition* is introduced by which the existence of cycle decompositions for the longer cycle lengths on a given Cartesian product is established. However, to get cycle lengths that fall outside the windows we mention later on, the method would have to be applied on varying tori. Letting $n \in \mathbb{Z}^+$, $d \in \mathbb{Z}^{\geq 2}$ and $a = 2^{i_1} + \cdots + 2^{i_d}$ be odd with $i_1 > \cdots > i_d = 0$, we see that on a given Cartesian product representation of $Q_{2an}$, the $2an-$dimensional hypercube, we can only obtain cycle lengths of the form $C_{a \cdot 2^\alpha}$ for $2n-1$ different consecutive $\alpha$ by way of the Wiggle decomposition in the case of aiming for the longest cycles. Even then, there is no explicit definition for constructing cycle decompositions via the Wiggle decomposition for tori of dimension greater than two, where the dimension of a given torus is one more than the number of times we take the Cartesian product against cycles in defining the torus.

---

[*]Carnegie Mellon University, 5000 Forbes Ave, Pittsburgh, PA, 15213, **irmartin@andrew.cmu.edu**

In section 4, we introduce the *Square decomposition* method in its most general form to give explicit cycle decompositions for the shorter cycle lengths on a $d-$dimensional torus. In the case of hypercubes $Q_{2an}$, we get cycles of the form $C_{a \cdot 2^\alpha}$ for $1 \leq \alpha \leq 2n$. Observe that the cycle lengths we can obtain via the Square decomposition are independent of the dimension of the torus, meaning we can obtain all cycles lengths possible by the Square decomposition applied to the $a-$fold Cartesian product we introduce further below. However, as we will see later when detailing the second method, it will be worth using the machinery known as anchored products developed in [4] to keep the dimension of our resulting torus low to maximize the cycle lengths currently possible for the second method while having all cycle decompositions being done on the same torus using cycles of all the lengths we obtain via the two methods. By restricting our view to these lower-dimensional sub-structures, this makes it more practical to visualize and consider computationally in terms of tractability.

Following the above, we introduce in section 5 the *Lock-and-Key decomposition* method in its general form to give explicit cycle decompositions for all cycle lengths starting at the longest cycle length given by the Square decomposition to the longest given by the Wiggle decomposition on the anchored product as we mentioned previously. In the case of cycle decompositions for the hypercube $Q_{2an}$, we get cycles of the form $C_{a \cdot 2^\alpha}$ for all

$$2n \leq \alpha \leq 2an - (d-1) - \sum_{k=1}^{d} i_k$$

all while on the same torus, thus extending the cycle length range beyond the

$$2(a-1)n - (d-2) - \sum_{k=1}^{d} i_k \leq \alpha \leq 2an - (d-1) - \sum_{k=1}^{d} i_k$$

window obtained for the longest cycles constructively for $d=2$ and non-constructively for $d>2$ by the Wiggle decomposition while on the same torus. Given the dependence of the possible cycle lengths on the dimension $d$ of our torus, we find that keeping the dimension of the torus as low as possible via anchored products allows us to obtain more of the particularly long cycles. The Lock-and-Key decomposition works on a torus of any dimension, but some of the longer cycles directly obtainable on the anchored product are not immediately obtained on the $a$-fold $Q_{2n}$ Cartesian product due to the dimension being more dependent on $a$ by a linear factor as a consequence of the emphasis placed by the method on its symmetries along each dimension. This is to say that we can translate the longer cycles obtained on the anchored product to those on the $a$-fold $Q_{2n}$ Cartesian product, but in the process the symmetries are distributed at key points along a given cycle and so one would not have paths of length $a$ with one edge coming from each dimension if one were to parse a given cycle $C_{a \cdot 2^\alpha}$ by every $a$ edges, a property characteristic of the Lock-and-Key decomposition. Hence, the cycle lengths obtained by directly applying the Lock-and-Key decomposition on the $a$-fold $Q_{2n}$ Cartesian product can be understood as those lengths for which one can directly define a cycle decomposition with these symmetries.

One example of this in effect would be to work on the $a-$fold $Q_{2n}$ Cartesian product representation of $Q_{2an} = Q_{2n} \square \cdots \square Q_{2n} = C_{2^{2n}} \square \cdots \square C_{2^{2n}}$, in which case the cycle lengths are for $2n \leq \alpha \leq 2an - (a-1)$ for cycles of the form $C_{a \cdot 2^\alpha}$. So choosing an increasingly large odd $a \geq 3$ that requires the same number of powers of two in its binary representation as a smaller odd $a \geq 3$, we can observe that a considerable number of the longer cycles do not translate to cycles with the characteristic symmetries when working on the $a-$fold Cartesian product from before relative to those obtained by working on an anchored product to get the wider range of $\alpha$ we presented earlier. Thus, anchored products not only allow us to maximize the cycle lengths to get the longest possible in our torus of focus, but they also allow us to extend this reach to a larger class of hypercubes $Q_{2an}$ with odd $a \geq 3$ that have the same number of powers of two in their binary representation.

Combining the Square and Lock-and-Key decomposition methods, we get a set of necessary and sufficient conditions to decompose a hypercube $Q_{2an}$ with odd $a \geq 3$ and $n$ as above into cycles of the form $C_{a \cdot 2^\alpha}$ for all

$$1 \leq \alpha \leq 2an - (d-1) - \sum_{k=1}^{d} i_k$$

with constructive definitions to explicitly construct the cycle decompositions all while on the torus derived from an anchored product we demonstrate in this paper. We also would like to remark that said torus is also the largest sub-structure of $Q_{2an}$ that we can fully decompose, meaning it is the largest torus for which cycle decompositions for all $\alpha$ above is possible and cycles of greater length are not admitted in the case $d \geq 3$.

Given that the methods we present in this paper can be applied onto the $a-$fold Cartesian product of $Q_{2n}$'s, these methods provide a more practical, symmetric means by which to decompose said Cartesian product. Additionally, they offer a more concrete starting point to focus on in proving the existence/non-existence of cycle decompositions using longer cycles that cannot exist in our chosen torus but can in the $a-$fold Cartesian product representation of $Q_{2an}$, which are precisely all cycles of the form $C_{a \cdot 2^\alpha}$ for

$$2an - (d-1) - \sum_{k=1}^{d} i_k < \alpha \leq 2an - \lceil \log_2(a) \rceil$$

when $d \geq 3$ as we obtain cycle decompositions of $Q_{2an}$ for all possible cycle lengths when $d = 2$. Note that this is strictly through toruses as we have mentioned. There is another approach used in [4] that proves the existence of cycle decompositions using slightly longer cycles such as $C_{7 \cdot 2^{10}}$ in the case of $Q_{14}$, but it does not obtain the existence/non-existence of cycle decompositions using the longest cycle, namely $C_{7 \cdot 2^{11}}$. See the end of section 7 for a complete statement of the problem discussed above and a conjectured approach to resolve it given the constructive results obtained here.

We structure the remainder of the paper in the following way. Sections 2 and 3, respectively, further motivate the problem with background and remarks on past works pertaining to hypercube decompositions, and provide the definitions and results of the fundamental tools we use in establishing the two new cycle decomposition methods. Sections 4 and 5 present the Square and Lock-and-Key decomposition's edge set definitions, the sets of necessary-and-sufficient conditions under which they apply, and explicit applications of the methods to give cycle decompositions of a given torus. At the end of the latter section, we give a Corollary 8 of the Square Decomposition's Theorem 2 and the Lock-and-Key Decomposition's Theorem 5, focusing constructively on cycle decmpositions of a slightly different class of hypercubes that we implicitly obtain along the way that were also treated by Gibson and Offner in [7]. To conclude, sections 6 and 7 present the proofs of the Square Decomposition's Theorem 2 and Corollary 4 and the Lock-and-Key Decomposition's Theorem 5 and Corollary 7, all of which are stated in sections 4 and 5 respectively.

## 2  Background

With physical motivations for hypercube decompositions tracing back to processor allocation and parallel computing problems, Stout in [11] and Bass and Sudborough in [5] presented the first results involving hypercube packings and decompositions, and investigated hypercube decompositions via $k$-regular spanning subgraphs, respectively. Using particular trees, Stout [11], Honrak, Siran, and Wallis [8], Mollard and Ramras [9], and Wagner and Wild [13] showed that the hypercube can be decomposed by said trees. In the case of path decompositions of hypercubes, it was shown independently by Erde [6], and Anick and Ramras [2] that any odd-dimensional hypercube $Q_n$ can be decomposed into any path such that the length is at most $n$ and divides the number of edges in $Q_n$.

Spurring some of the first fundamental investigations of hypercube decompositions using cycles, Ringel [10] proved that when $n \geq 2$ is a power of two, $Q_n$ has a decomposition into Hamiltonian cycles and hence a Hamiltonian cycle decomposition. In seeking to broaden the above result to a larger class of hypercubes, Ringel asked whether $Q_n$ has a Hamiltonian cycle decomposition for all even $n \geq 2$. With a resolution in the affirmative, Ringel's question was treated implicitly by Aubert and Schneider in [3] and explicitly by Alspach, Bermond, and Sotteau in [1].

Ensuing the pivotal questions of Ringel and the complete extension of Hamiltonian cycle decompositions to all even hypercubes, Axenovich, Offner and Tompkins [4] showed in particular that even hypercubes $Q_n$ can be decomposed into long cycles whose length is of a given form. Similarly, Tapadia, Waphare, and Borse in [12] proved that $Q_n$ can be decomposed into short cycles whose length is a power of two. This result was then extended by Gibson and Offner in [7], where it is proven that $Q_n$ can be decomposed into cycles whose length is at least four and a power of two that divides $2^n$. Nevertheless, decompositions of hypercubes remain an open problem for trees, paths, and even in the case of cycles.

Furthering constructively the current understanding of cycle decompositions of Cartesian products of cycles, we present the Square and Lock-and-Key decomposition methods from which we get cycle decompositions of even-dimensional hypercubes from the shortest possible cycle lengths to the longest lengths obtained in [4] on anchored products, where the lengths are all of the same given form. In the process, we also get explicit cycle decompositions for $Q_n$ when $n$ is divisible by four using cycles of the lengths obtained in [7].

## 3  Definitions and Intermediate Results

From the definition of the $n$-dimensional hypercube presented earlier, we most notably get that $Q_n$ has $2^n$ vertices and $n2^{n-1}$ edges. For two graphs $G$ and $H$, the *Cartesian product* of $G$ and $H$, denoted $G \,\square\, H$, is the graph with vertex set $V(G \,\square\, H) = \{(j_1, j_2) : j_1 \in V(G), j_2 \in V(H)\}$ and edge set $E(G \,\square\, H) = \{(j_1, j_2)(j_1', j_2') : j_1 = j_1', j_2 j_2' \in E(H) \text{ or } j_2 = j_2', j_1 j_1' \in E(G)\}$. Further, taking $H_1, \ldots, H_k$ to be subgraphs of a graph $G$, we say that $\{H_1, \ldots, H_k\}$ defines a *decomposition* of $G$ if $G = H_1 \cup \cdots \cup H_k$ with $E(H_n) \cap E(H_m) = \emptyset$ if $n \neq m$ for $1 \leq n, m \leq k$ and $k \in \mathbb{Z}^+$. So in particular, if $H_1, \ldots, H_k$ are cycles, we say they define a *cycle decomposition* of $G$. In this paper, we use $C_b$ to denote a cycle of length $b$ for $b \in \mathbb{Z}^{\geq 2}$.

As in [4], we define the *anchored product* $(G_1, S_1) \boxplus (G_2, S_2)$ of two graphs $G_1$ and $G_2$ with $S_1 \subseteq V(G_1)$ and $S_2 \subseteq V(G_2)$ as the graph with vertex set

$$\{(j_1, j_2) : j_1 \in V(G_1), j_2 \in V(G_2), \text{ and } j_1 \in S_1 \text{ or } j_2 \in S_2\}$$

and edge set

$$\{(j_1, j_2)(j_1', j_2') : j_1 j_1' \in E(G_1),\ j_2 = j_2' \in S_2\} \cup \{(j_1, j_2)(j_1', j_2') : j_2 j_2' \in E(G_2),\ j_1 = j_1' \in S_1\}.$$

Note that if we take $S_1 = V(G_1)$ and $S_2 = V(G_2)$, we get that the anchored product as defined above is precisely $G_1 \square G_2$. The Cartesian product of two or more cycles is then what we call a *torus*. For a more detailed discussion of anchored products, see [4]. While every even dimensional hypercube $Q_{2mn}$ for $n, m \in \mathbb{Z}^+$ can be expressed as an $m$-fold Cartesian product of $C_{2^{2n}}$'s, we in turn end up working with an $m$-dimensional torus that becomes complex to fully visualize as we increase $m$. Nonetheless, the end-goal is to fully realize cycle decompositions on this torus for all cycle decompositions of all possible lengths that exist as this torus would be able to admit all such cycle decompositions.

To extract smaller toruses whose dimension is uniformly and comparatively low and invariant, we consider all hypercubes in the equivalence class $\{Q_{2an}\}$ with odd $a \geq 3$, where every member is characterized for having the same number of powers of two in their binary representation as that of a representative's $a$. From [1], we know that even-dimensional hypercubes $Q_n$ have a Hamiltonian cycle decomposition and we express this as $Q_n = \frac{n}{2} C_{2^n}$, where the coefficient represents the number of cycles in the cycle decomposition and $C_{2^n}$ is in particular a Hamiltonian cycle, since $|V(C_{2^n})| = 2^n = |V(Q_n)|$, and also a *spanning subgraph* as a consequence. By the above and the Cartesian product property $Q_{n+m} = Q_n \square Q_m$ for hypercubes with $n, m \in \mathbb{Z}^{\geq 0}$, we have for $Q_{2an}$ with $a = 2^{i_1} + 2^{i_2} + \cdots + 2^{i_d}$ for $i_1 > i_2 > \cdots > i_d = 0$ that

$$Q_{2an} = Q_{n2^{i_1}+1} \square \cdots \square Q_{n2^{i_d}+1} = n2^{i_1} C_{2^{n2^{i_1}+1}} \square \cdots \square n2^{i_d} C_{2^{n2^{i_d}+1}}.$$

Then, by Propositions 7 and 8 from [4] with the Hamiltonian cycles as our spanning subgraphs, we get the above can be expressed as

$$Q_{2an} = n \left( \prod_{k=1}^{d} 2^{i_k} \right) \left[ (C_{2^{n2^{i_1}+1}}, S_1) \boxplus (C_{2^{n2^{i_2}+1}}, S_2) \boxplus \cdots \boxplus (C_{2^{n2^{i_d}+1}}, S_d) \right],$$

where the coefficient is the number of isomorphic copies of the graph resulting from the anchored product and $S_k$ is the set of every $2^{i_k}$th vertex from $C_{2^{n2^{i_k}+1}}$ for all $1 \leq k \leq d$. The above anchored product can then be viewed as what is called a *subdivided torus*, where all vertices in $S_k$ have $2d$-many edges associated to them in the anchored product and all other vertices belonging to $C_{2^{n2^{i_k}+1}}$ not in $S_k$ are 2-regular. From this, we see that there are $2^{i_k}$ edges serving as subdivisions between any two adjacent vertices in $S_k$ that together form a vertex pair, and hence a subdivided edge along every copy of $C_{2^{n2^{i_k}+1}}$ in the anchored product for every $1 \leq k \leq d$. Focusing on the edges along a given copy of $C_{2^{n2^{i_k}+1}}$ implicitly defined by the vertices in the $S_k$, which serve as what are called *representing sets* in defining and identifying this particular copy of the subdivided torus, we can extract an *underlying torus* by making all edges in the subdivided torus *distance regular* (DR). In other words, we treat every resulting edge from using the vertex pairs generated by the representing sets as one edge regardless of the number of subdivisions it has in the subdivided torus. By doing this, we have an isomorphism between the subdivided torus and the underlying torus. Hence, we are interested in the cycle decompositions obtained on the underlying $d$-dimensional torus for $d \geq 2$ as that then gives us a cycle decomposition of the hypercube from which we derived the underlying torus, where all cycles in the cycle decomposition are of the same length.

Note that we say that an anchored cycle is distance regular if there are the same number of subdivisions between any two adjacent vertices from the representing set of the anchored cycle. Further note that, since there may be partitions in the edges of a subdivided torus along a given dimension if the representing set excludes some vertices of the cycle from which the anchored cycle is derived, we refer to edges as *non-partitioned edges* if the anchored cycle is precisely the cycle from which it is derived i.e. $(C_b, S) = C_b$ with $S = V(C_b)$. Viewing each cycle in the Cartesian product defining the underlying torus as defining a dimension of the torus and taking into account the link between the underlying torus and the subdivided torus, this gives rise to the following result in a more general setting:

**Proposition 1** *Let $d \geq 2$ and suppose an underlying torus $C_{|S_1|} \square \cdots \square C_{|S_d|}$ has a decomposition into cycles such that every cycle has $a_1$ edges from the dimension defined by $C_{|S_1|}$, $a_2$ edges from the dimension defined by $C_{|S_2|}$, and so forth up to $a_d$ edges from the dimension defined by $C_{|S_d|}$. If $(C_{y_1^*}, S_1)$, $(C_{y_2^*}, S_2)$, ..., $(C_{y_d^*}, S_d)$ are all anchored-DR cycles, where $S_1$ is the set of every $k_1^{th}$ vertex of $C_{y_1^*}$, $S_2$ is the set of every $k_2^{th}$ vertex of $C_{y_2^*}$ and so forth with $S_d$ being the set of every $k_d^{th}$ vertex from $C_{y_d^*}$, then $(C_{y_1^*}, S_1) \boxplus \cdots \boxplus (C_{y_d^*}, S_d)$ can be decomposed into cycles with*

$$\sum_{i=1}^{d} a_i k_i$$

*edges.*

**Proof:** Let $d \geq 2$ and fix $1 \leq i \leq d$. Then, every one of the $k_i$ edges along the $i$th dimension belonging to a given cycle in the decomposition is partitioned $a_i$ times. Summing the edges of the given cycle along each dimension $i$ for $i = 1, \ldots, d$ tells us that the cycle has $a_1 k_1 + \cdots + a_d k_d$ edges. ∎

## 4 Generalization of the $d-$Dimensional $\lambda-$Square Decomposition's Edge Set Definition:

**Theorem 2 ($d$-Dimensional $\lambda$-Square Decomposition Conditions)** *For integers $d \geq 2$, $\lambda > 0$, and $y_1, y_2, \ldots, y_d > 0$, if $2\lambda \mid y_1$, $2\lambda \mid y_2, \ldots,$ and $2\lambda \mid y_d$, then $C_{y_1} \square C_{y_2} \square \cdots \square C_{y_d}$ can be decomposed into*

$$\frac{1}{2\lambda} \prod_{i=1}^{d} y_i$$

*cycles of length $2d\lambda$. Note that every cycle here has the same number of edges coming from each cycle defining a dimension of the torus, and $\lambda$ corresponds to the individual number of edges of a "side." This is all without distinguishing between partitioned and non-partitioned edges.*

The corresponding edge set definition established by the General $d-$dimensional Square decomposition is the following:

$$E(C_{\ell,t,p_1,z,s_1,\ldots,s_{d-2}}) =$$

$$\bigcup_{\substack{m_1,m_2,\ldots,m_d=0 \\ M(m_1,m_2,\ldots,m_d)=1}}^{1} \{(j_1,j_2,j_3,j_4,j_5,\ldots,j_d)(j_1+(m_2-m_1),\ j_2+\chi,\ j_3+(-1)^r(m_1-m_3),\ j_4+(m_3-m_4),$$

$$j_5 + (m_4 - m_5),\ \ldots,\ j_d + (m_{d-1} - m_d)) \mid j_1 = t + (2\nu + r + m_1)\lambda + (m_2 - m_1)x_1,$$

$$j_2 = t + (2\mu + r + m_2 + p_1)\lambda + \chi x_2 + p_2 + \sum_{k=3}^{d-1} p_k,\ j_3 = s_1 + 2z\lambda + (-1)^r(m_3\lambda + (m_1 - m_3)x_3),$$

$$j_4 = s_2 + m_4\lambda + (m_3 - m_4)x_4,\ j_5 = s_3 + m_5\lambda + (m_4 - m_5)x_5,\ \ldots,$$

$$j_d = s_{d-2} + m_d\lambda + (m_{d-1} - m_d)x_d,\ 0 \leq x_1, x_2, \ldots, x_d \leq \lambda - 1\},$$

where

$$0 \leq \ell \leq \frac{y_1 y_2}{2\lambda^2} - 1, \ 0 \leq t \leq \lambda - 1, \ \eta(y, c) = \begin{cases} 1, & \text{if } y \in \mathbb{Z}^{>c} \\ 0, & \text{if } y \notin \mathbb{Z}^{>c} \end{cases} \text{ for } y, c \in \mathbb{Z}^{\geq 2}, \ 0 \leq p_1 \leq \eta(d, 2),$$

$$0 \leq z \leq \left(\frac{y_d}{2\lambda} - 1\right)\eta(d, 2), \ 0 \leq s_1 \leq (\lambda - 1)\eta(d, 2), \ 0 \leq s_k \leq (y_{k+1} - 1)\eta(d, k+1) \text{ for } 2 \leq k \leq d-2,$$

$$p_{k^*} = p_{k^*}(s_{k^*-1}) = s_{k^*-1} - \left\lfloor \frac{s_{k^*-1}}{\lambda} \right\rfloor \lambda \text{ for } 2 \leq k^* \leq d-1,$$

$$\chi = \chi(d, m_1, m_2, \ldots, m_d) = 1 - 2\left\lfloor \frac{1}{d}\sum_{k=1}^{d} m_k \right\rfloor \left\lceil \frac{1}{d}\sum_{k=1}^{d} m_k \right\rceil \text{ for } (d, m_1, m_2, \ldots, m_d) \in \mathbb{Z}^{\geq 2} \times \{0,1\}^d,$$

$$M = M(m_1, m_2, \ldots, m_d) \text{ with the logical expression}$$

$$q(m_1, m_2, \ldots, m_d) = (m_1 = m_2 = \cdots = m_d) \vee (m_2 \neq m_1 = m_3 = \cdots = m_d) \vee \bigvee_{k=2}^{d-1}(m_1 = m_2 = \cdots = m_k \neq m_{k+1} = \cdots = m_d)$$

$$\text{is defined as } M(m_1, m_2, \ldots, m_d) = \begin{cases} 1 & \text{if } q(m_1, m_2, \ldots, m_d) \text{ is true} \\ 0 & \text{if } q(m_1, m_2, \ldots, m_d) \text{ is false} \end{cases}, r = r(\ell) = \ell - 2\left\lfloor \frac{\ell}{2} \right\rfloor,$$

$$\nu = \nu(\ell) = \left\lfloor \frac{\ell}{2} \right\rfloor - \left\lfloor \frac{\ell \lambda}{y_1} \right\rfloor \frac{y_1}{2\lambda}, \text{ and } \mu = \mu(\ell) = \left\lfloor \frac{\ell}{y_1} \right\rfloor.$$

**Note:** For $d = 2$, the convention is $(j_1, j_2) \in (\mathbb{Z}/y_1\mathbb{Z}) \times (\mathbb{Z}/y_2\mathbb{Z})$.

For $d = 3$, the convention is $(j_1, j_2, j_3) \in (\mathbb{Z}/y_1\mathbb{Z}) \times (\mathbb{Z}/y_2\mathbb{Z}) \times (\mathbb{Z}/y_3\mathbb{Z})$.

For $d = 4$, the convention is $(j_1, j_2, j_3, j_4) \in (\mathbb{Z}/y_1\mathbb{Z}) \times (\mathbb{Z}/y_2\mathbb{Z}) \times (\mathbb{Z}/y_4\mathbb{Z}) \times (\mathbb{Z}/y_3\mathbb{Z})$.

For $d \geq 5$, the convention is $(j_1, j_2, \ldots, j_d) \in (\mathbb{Z}/y_1\mathbb{Z}) \times (\mathbb{Z}/y_2\mathbb{Z}) \times (\mathbb{Z}/y_d\mathbb{Z}) \times (\mathbb{Z}/y_3\mathbb{Z}) \times \cdots \times (\mathbb{Z}/y_{d-1}\mathbb{Z})$.

With the above, we wish to emphasize that the definition of the edge set above is to be read and used in the following way. For $d = 2$, one considers components $j_1$ and $j_2$ in the edge set definition to define the decomposition of $C_{y_1} \square C_{y_2}$ and disregards components $j_3$ through $j_d$. All parameters associated to these extraneous components can be sent to zero and disregarded. For $d = 3$, we consider the first three components $j_1, j_2$ and $j_3$ to define the cycles to decompose $C_{y_1} \square C_{y_2} \square C_{y_3}$ and disregard the components $j_4$ through $j_d$ presented. As before, all parameters associated to the extraneous components get sent to zero. Continuing in this fashion, we can obtain the corresponding definition for $d = 4$ and $d \geq 5$ to decompose $C_{y_1} \square C_{y_2} \square C_{y_4} \square C_{y_3}$ and $C_{y_1} \square C_{y_2} \square C_{y_d} \square C_{y_3} \square \cdots \square C_{y_{d-1}}$, respectively. In this case, we emphasize this configuration for the torus for better viewability given how large underlying tori of hypercube graphs can come to be. Of course, once the dimension $d$ is chosen, the above definition can modified to have the cycles defining the Cartesian product in any desirable order without altering the intended way of viewing the cycle decomposition given the symmetry of the cycles defined by the Square decomposition and the method's independence from the dimension of the torus.

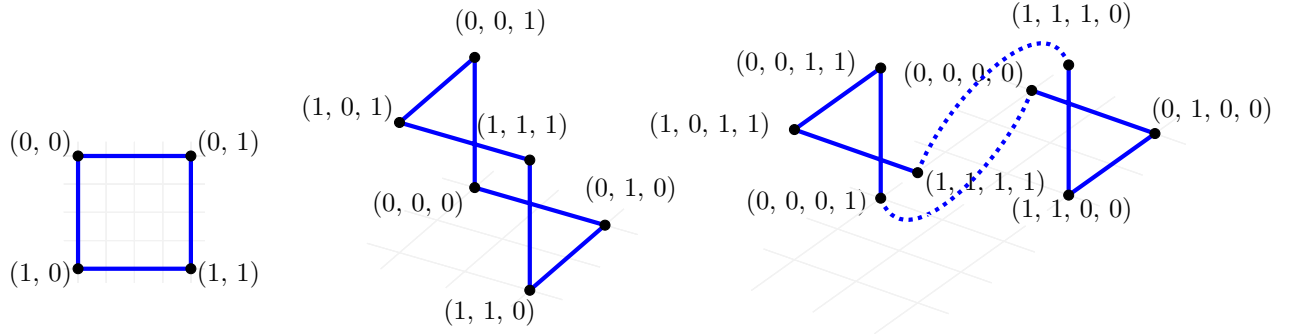

In the figures above, we show cycles defined by the Square decomposition for $d = 2$ (left), $d = 3$ (center) and $d = 4$ (right), where the $d$-tuples near every darkened vertex correspond to one of the $(m_1, \ldots, m_d)$ "moves" used by the Square decomposition's edge set definition to define a cycle for the cycle decomposition. The edge set definition starts from the $(0, \ldots, 0)$ $d$-tuple and follows the sequence of $d$-tuples reached by proceeding clockwise along the given cycles. Note that the $d$-tuples for edges along the third dimension are unique up to reflections along that dimension.

**Lemma 3** *Let $d \geq 2$ and $a = 2^{i_1} + 2^{i_2} + \cdots + 2^{i_d}$ be odd with integers $i_1 > i_2 > \cdots > i_d = 0$. Now let $(C_{y_1^*}, S_1)$ be such that $S_1$ has every $2^{i_1}$th vertex of $C_{y_1^*}$, $(C_{y_2^*}, S_2)$ with $S_2$ having every $2^{i_2}$th vertex of $C_{y_2^*}$ and so forth with $(C_{y_d^*}, S_d)$ having $S_d$ contain every vertex of $C_{y_d^*}$. Then, letting $\lambda = 2^{\alpha-1}$ for an integer $\alpha \geq 1$, if $2^{\alpha+i_1} \mid y_1^*$, $2^{\alpha+i_2} \mid y_2^*$, and so forth with $(i_d = 0)$ $2^\alpha \mid y_d^*$, it follows that $C_{a \cdot 2^\alpha}$ decomposes $(C_{y_1^*}, S_1) \boxplus (C_{y_2^*}, S_2) \boxplus \cdots \boxplus (C_{y_d^*}, S_d)$.*

**Proof:** Let $d \geq 2$, $a = 2^{i_1} + \cdots + 2^{i_d}$ be odd with integers $i_1 > \cdots > i_d = 0$, and $\lambda = 2^{\alpha-1} > 0$ with $\alpha \geq 1$. Given that the Square decomposition defines the cycles to have $2\lambda = 2^\alpha$ edges along each dimension in the underlying torus by Theorem 2, we see by Proposition 1 that $2^\alpha$ edges in the underlying torus along the $k$th dimension for $k = 1, \ldots, d$ corresponds to $2^{i_k} \cdot 2^\alpha$ edges in the subdivided torus. Hence, summing the edges along each dimension gives us that our cycle decomposition of the underlying torus using cycles of length $2d\lambda$ corresponds to a cycle decomposition of the given anchored product using cycles of length $a \cdot 2^\alpha$, where all cycles are of the same length as a consequence.

∎

**Corollary 4** *For $n \geq 1$, $d \geq 2$ and $a = 2^{i_1} + 2^{i_2} + \cdots + 2^{i_d}$ with integers $i_1 > i_2 > \cdots > i_d = 0$, a $d$-dimensional $2^{\alpha-1}$-square decomposition of $(C_{2^{n2^{i_1}+1}}, S_1) \boxplus (C_{2^{n2^{i_2}+1}}, S_2) \boxplus (C_{2^{2n}}, S_d) \boxplus \cdots \boxplus (C_{2^{n2^{i_{d-1}}+1}}, S_{d-1})$ can be constructed using $C_{a \cdot 2^\alpha}$ cycles, where $1 \leq \alpha \leq 2n$.*

**Note:** $2^{\alpha-1}$ *refers to the number of edges used on a given "side" without distinguishing between partitioned and non-partitioned edges. Further, observe that the form of the anchored product is dependent on $d \in \mathbb{Z}^{\geq 2}$ as presented earlier at the end of the General $d$-Dimensional Square Decomposition's edge set definition.*

Before presenting the second decomposition method, we provide applications of the Square decomposition to give cycle decompositions of two and three-dimensional toruses.

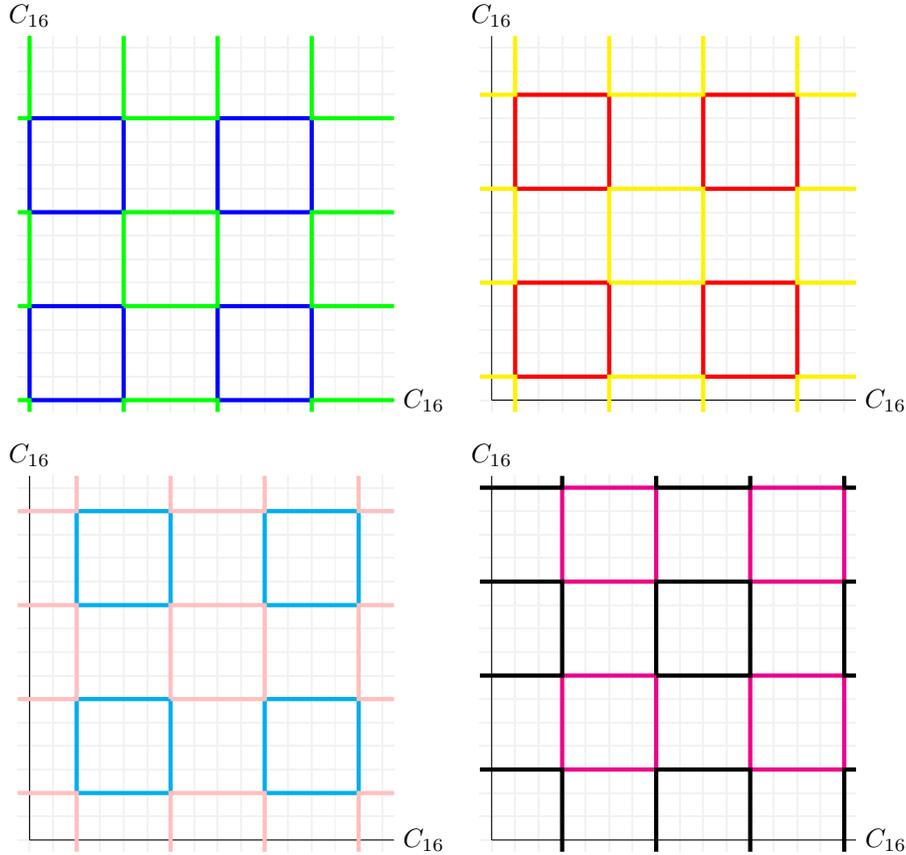

In the four figures above, we have the $C_{16}$ cycles defined by the Square decomposition's edge set definition on the torus $C_{16} \square C_{16}$ for the four translation sets admitted by this torus with $t = 0$ (top left), $t = 1$ (top right), $t = 2$ (bottom left), and $t = 3$ (bottom right).

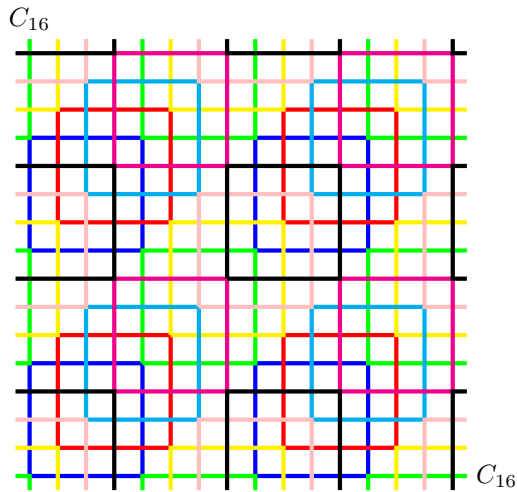

In the above, we show the Square decomposition of the torus $C_{16} \square C_{16}$ with $\lambda = 4$, where we have combined the $C_{16}$ cycles defined in each of the four translation sets shown previously. Replacing the cycles defining the vertical dimension with $C_{128}$'s and translating 7 copies of the above decomposition downwards by 16 edges more than the previous translation gives a cycle decomposition of $C_{128} \square C_{16}$, making this the cycle decomposition of an underlying torus of $Q_{12}$. Here, each vertical edge in the underlying torus corresponds to two edges in the subdivided torus and every horizontal edge corresponds to one edge in the subdivided torus. Hence, the above corresponds to a cycle decomposition of $Q_{12}$ using $C_{24}$'s.

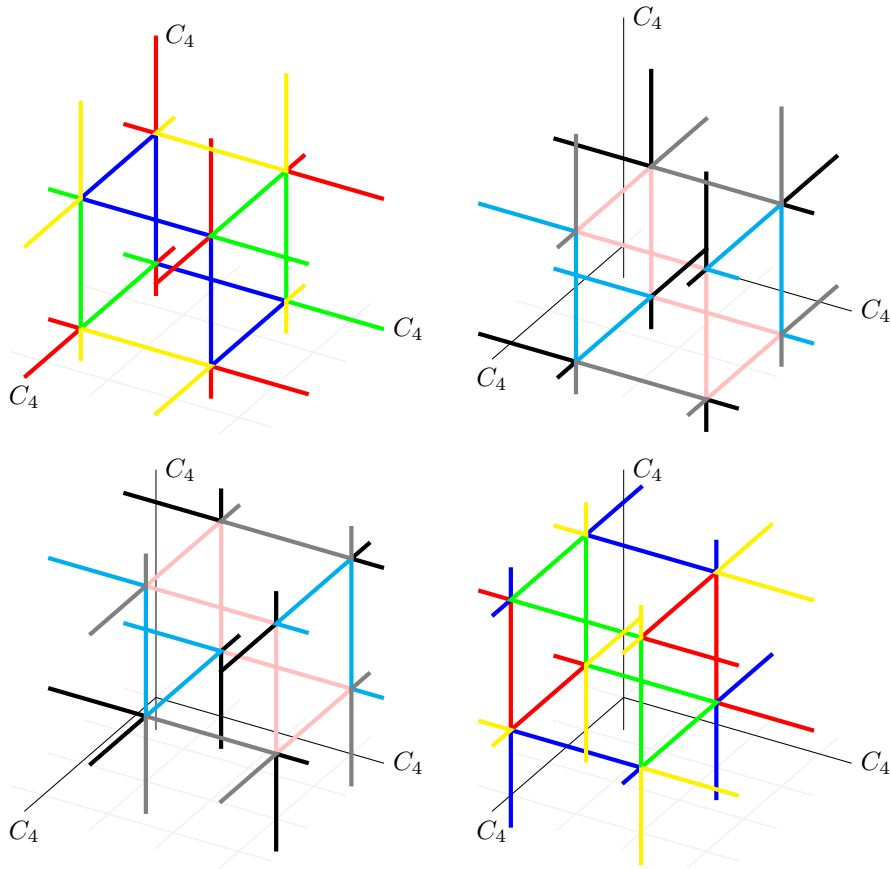

In the four figures above are the $C_{12}$ cycles defined by the Square decomposition for the three-dimensional torus $C_4 \square C_4 \square C_4$ with $\lambda = 2$, where each figure shows all cycles of a given translation set possible on a given level. In other words, there are two kinds of translations sets: those associated to $t$ and those associated to $s_1$'s remainder, $p_2$, when divided by $\lambda$. The higher-dimensional instances proceed analogously with regards to these translations, where we have chosen $j_2$ to be the main reference dimension with respect to which these translations are carried out. Replacing the cycles defining the first dimension by $C_{64}$'s and those defining the second dimension by $C_8$'s, we see that, by translating the above cycle decomposition along the modified dimensions analogously to the previous example, we get a cycle decomposition of an underlying torus $C_{64} \square C_8 \square C_4$ of $Q_{14}$. Mapping the edges from a given cycle to those in the subdivided torus with the corresponding number of subdivisions along each dimension, we get the decomposition below serves as a cycle decomposition of $Q_{14}$ using $C_{7 \cdot 2^2}$'s.

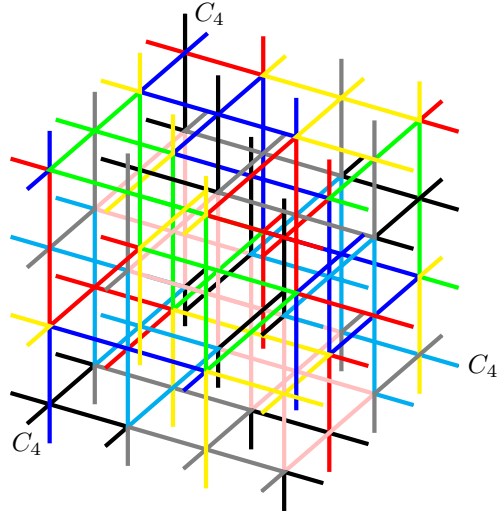

# 5 Generalization of the $d-$Dimensional Lock-and-Key Decomposition's Edge Set Definition:

**Theorem 5 ($d$-Dimensional Lock-and-Key Decomposition Conditions)** *Let $d \geq 2$, $z_1, z_2, \ldots,$ $z_d > 0$ all be even and $b_1, \ldots, b_{d-1} \geq 0$ all be odd or 0, where $b_1$ represents the number of key "bits." Further, let $A_k$ for $3 \leq k \leq d+1$ be as in the Lock-and-Key Decomposition's edge set definition presented below. For $d = 2$, let $4 \mid z_d$ and $A_3 \mid z_1$ and, for $d \geq 3$, let $4 \mid z_d$, $2A_{d+1} \mid z_{d-1}$, $\ldots$, $2A_4 \mid z_2$, and $A_3 \mid z_1$. Then, the tori $C_{z_1} \square C_{z_d}$ when $d = 2$, $C_{z_1} \square C_{z_d} \square C_{z_2}$ when $d = 3$, and $C_{z_1} \square C_{z_d} \square C_{z_2} \square \cdots \square C_{z_{d-1}}$ when $d \geq 4$ can all be decomposed into*

$$2 \prod_{k=0}^{d-2} \frac{z_{k+1}}{A_{k+3}}$$

*many cycles of length $\frac{dz_d}{2} \prod_{k=0}^{d-2} A_{k+3}$. Note that we are not distinguishing between the partitioned and non-partitioned edges, and that the number of edges coming from each of the $d$ cycles defining the torus is the same i.e. $\frac{z_d}{2} \prod_{k=0}^{d-2} A_{k+3}$ by the above convention.*

Below is the General $d-$dimensional Lock-and-Key decomposition's edge set definition on a torus of the form presented in Theorem 5 above:

$$E(C_{\ell,\gamma,t,p_1,s_1,p_2,s_2,\ldots,p_{d-3},s_{d-3},p_{d-2},s_{d-2}}) =$$

$$\bigcup_{\substack{m_1,m_2,\ldots,m_d=0 \\ M(m_1,m_2,\ldots,m_d)=1}}^{2A_1-1} \{(j_1,j_2,j_3,j_4,\ldots,j_{d-1},j_d)(j_1+(-1)^{(1-\eta_3)(1-A_2)\gamma}|m_2-m_1|,\ j_2+\chi,\ j_3+(-1)^{(1-A_2)}(m_1-m_3),\ j_4+$$
$$(-1)^{(1-A_2)}(m_3-m_4),\ldots,\ j_{d-1}+(-1)^{(1-A_2)}(m_{d-2}-m_{d-1}),\ j_d+(-1)^{(1-A_2)}(m_{d-1}-m_d)) \mid j_1 = (\ell+$$
$$(1-\eta_3)(1-A_2)\gamma)A_3 - (1-\eta_3)(1-A_2)\gamma + \eta_3(1-A_2) + (-1)^{(1-\eta_3)(1-A_2)\gamma}(m_1+2x_1+2\eta_3 m_d(1-m_1)),$$
$$j_2 = m_2 + \gamma + \eta_3(1-A_2)(s_1+1) + (4-2A_2)x_2,\ j_3 = p_1 + s_1 A_4 2^{\eta_{\alpha_2}} + (1-A_2) + (-1)^{(1-A_2)} m_3 + 2x_3,$$
$$j_4 = p_2 + s_2 A_5 2^{\eta_{\alpha_3}} + (1-A_2) + (-1)^{(1-A_2)} m_4 + 2x_4,\ \ldots,\ j_{d-1} = p_{d-3} + s_{d-3} A_d 2^{\eta_{\alpha_{d-2}}} + (1-A_2)$$
$$+(-1)^{(1-A_2)} m_{d-1} + 2x_{d-1},\ j_d = p_{d-2} + s_{d-2} A_{d+1} 2^{\eta_{\alpha_{d-1}}} + (1-A_2) + (-1)^{(1-A_2)} m_d + 2x_d,$$

$$0 \leq x_1 \leq \frac{A_3}{2} - 2, \ 0 \leq x_2 \leq z_d 2^{A_2-2} - 1, 0 \leq x_k \leq A_{k+1} - 1 \text{ for } 3 \leq k \leq d\}$$

$$\cup \bigcup_{\substack{m_1,m_2,\ldots,m_d=0 \\ M(m_1,m_2,\ldots,m_d)=1}}^{(2A_1-1)(2(1-A_2)-1)} \{(j_1,j_2,j_3,j_4,\ldots,j_{d-1},j_d)(j_1 + (-1)^{(1-\eta_3)\gamma+1}|\chi|, \ j_2 + (-1)^{r_1}((1-\eta_3)(m_2-m_1) + \eta_3(1-\eta_4)(m_1-m_3)$$
$$+ \eta_4(m_{d-1} - m_d)), \ j_3 + (m_2 - m_1), \ j_4 + (m_1 - m_3), j_5 + (m_3 - m_4), \ldots, j_{d-1} + (m_{d-3} - m_{d-2}), \ j_d+$$
$$(m_{d-2} - m_{d-1})) \mid j_1 = (\ell + 1 - (1-\eta_3)\gamma)A_3 - (1-\eta_3)\gamma + (-1)^{(1-\eta_3)\gamma+1}(2 + m_2 + 2m_1(1-m_2)$$
$$+ 2x_1 + \eta_3(2m_d(1-m_1) - r_1(A_3 - 1))), \ j_2 = 2 + \gamma + (1-\eta_3)m_1 + \eta_3(1 + s_1 + m_d - r_1) + 4x_2,$$
$$j_3 = s_1 + m_1 + r_1, \ j_4 = s_2 + m_3 + r_1, \ldots, \ j_{d-1} = s_{d-3} + m_{d-2} + r_1, \ j_d = s_{d-2} + m_{d-1} + r_1,$$
$$0 \leq x_1 \leq \frac{A_3}{2} - 2, \ 0 \leq x_2 \leq \frac{z_d}{4} - 1\}$$

$$\cup \bigcup_{\substack{m_1,m_2,\ldots,m_d=0 \\ M(m_1,m_2,\ldots,m_d)=1}}^{1} \{(j_1,j_2,j_3,j_4,\ldots,j_{d-1},j_d)(j_1 + (-1)^{((1-\eta_3)(1-A_2)+\eta_3(1-A_1))\gamma}((1-m_1)m_2 + ((1-\eta_3)(1-A_2) + \eta_3(1-A_1))\prod_{j=1}^{d} m_j + ((1-\eta_3)A_2 + \eta_3 A_1)(1-m_2)m_1), \ j_2 + \prod_{j=1}^{d}(1-m_j) + ((1-\eta_3)(1-A_2) + \eta_3(1-A_1))(1-m_2)m_1 + (-1)^{R_1+1}((1-\eta_3)A_2 + \eta_3 A_1)\prod_{j=1}^{d} m_j, \ j_3 + (-1)^{\eta_{\alpha_2}((1-R_1)R_2S_2+R_1r_2)+(1-A_1)\gamma}(m_1 - m_3), j_4+$$
$$(-1)^{\eta_{\alpha_3}((1-R_1)R_3S_3+R_1r_2)}(m_3 - m_4), \ldots, \ j_{d-1} + (-1)^{\eta_{\alpha_{d-2}}((1-R_1)R_{d-2}S_{d-2}+R_1r_2)}(m_{d-2} - m_{d-1}), \ j_d+$$
$$(-1)^{\eta_{\alpha_{d-1}}((1-R_1)R_{d-1}S_{d-1}+R_1r_2)}(m_{d-1} - m_d)) \mid j_1 = (\ell + 1 - ((1-\eta_3)(1-A_2) + \eta_3(1-A_1))\gamma)A_3 - ((1-\eta_3)(1-A_2) + \eta_3(1-A_1))\gamma + t + (-1)^{((1-\eta_3)(1-A_2)+\eta_3(1-A_1))\gamma}(-2 + m_1 + ((1-\eta_3)(1-A_2) + \eta_3(1-A_1))(1-m_2)m_1 + \eta_3(2m_d(1-m_1) - A_1(1-A_2)(A_3 - 2)(1-m_d))), \ j_2 = m_2 + \gamma + (4 - 2A_2)x_2 + ((1-\eta_3)(1-A_2) + \eta_3(1-A_1))(1-m_2)m_1 + 2((1-\eta_3)A_2 + \eta_3 A_1)(1-m_2)m_1 + \eta_3(A_1(1-A_2)s_1 + 2((1-m_1)m_d - \eta_{\alpha_2}(1-R_1)(m_1(1-m_2) + m_d(1-m_1)))), \ j_3 = p_1 + s_1 A_4 2^{\eta_{\alpha_2}+(1-A_1)} + ((-1)^{(1-A_1)\gamma} - 2\eta_{\alpha_2} r_2 \Psi_0)m_3+$$
$$2x_3 + 2\eta_4 \eta_{\alpha_3}(1-R_1)((1-r_2)X_0 - r_2\Psi_0)(1-m_3)m_d, \ j_4 = p_2 + s_2 A_5 2^{\eta_{\alpha_3}} + (1 - 2\eta_{\alpha_3} r_2 \Psi_1)m_4 + 2x_4+$$
$$2\eta_5 \eta_{\alpha_4}(1-R_1)((1-r_2)X_1 - r_2\Psi_1)(1-m_4)m_d, \ldots, \ j_{d-1} = p_{d-3} + s_{d-3}A_d 2^{\eta_{\alpha_{d-2}}} + (1 - 2\eta_{\alpha_{d-2}} r_2 \Psi_{d-4})m_{d-1} + 2x_{d-1} + 2\eta_d \eta_{\alpha_{d-1}}(1-R_1)((1-r_2)X_{d-4} - r_2\Psi_{d-4})(1-m_{d-1})m_d, \ j_d = p_{d-2} + s_{d-2}A_{d+1} 2^{\eta_{\alpha_{d-1}}}$$
$$+ (1 - 2\eta_{\alpha_{d-1}} r_2 \Psi_{d-3})m_d + 2x_d, \ 0 \leq x_2 \leq z_d 2^{A_2-2} - 1, \ 0 \leq x_k \leq A_{k+1} - 1 \text{ for } 3 \leq k \leq d\}$$

$$\cup \bigcup_{\substack{m_1,m_2,\ldots,m_d=0 \\ M(m_1,m_2,\ldots,m_d)=1}}^{2(1-A_2)-1} \{(j_1,j_2,j_3,j_4,j_5\ldots,j_{d-1},j_d)(j_1 + (-1)^{((1-\eta_3)+\eta_3(1-A_1))\gamma+1}\left(m_2(1-m_1) + \prod_{j=1}^{d} m_j\right), \ j_2+$$
$$(-1)^{A_1((1-\eta_3)m_1+\eta_3 m_d)}\left(\prod_{j=1}^{d}(1-m_j) + (1-A_1)(1-m_2)m_1 + A_1((1-\eta_3)(1-m_2)m_1 + \eta_3(1-\eta_4)(1-m_1)m_3 + \eta_4(1-m_{d-1})m_d)\right), \ j_3 + (-1)^{(1-A_1)\gamma}(1-A_1)(m_1 - m_3) + A_1(m_2 - m_3)m_1, \ j_4 + (1-A_1)($$
$$m_3 - m_4) + A_1(m_1 - m_4)m_3, \ j_5 + (1-A_1)(m_4 - m_5) + A_1(m_3 - m_5)m_4, \ldots, j_{d-1} + (1-A_1)(m_{d-2}$$
$$- m_{d-1}) + A_1(m_{d-3} - m_{d-1})m_{d-2}, \ j_d + (1-A_1)(m_{d-1} - m_d) + A_1(m_{d-2} - m_d)m_{d-1}) \mid j_1 = (\ell + 1$$
$$- ((1-\eta_3) + \eta_3(1-A_1))\gamma)A_3 - ((1-\eta_3) + \eta_3(1-A_1))\gamma + t + (-1)^{((1-\eta_3)+\eta_3(1-A_1))\gamma+1}(m_1$$

$$+m_1(1-m_2)+2\eta_3(1-m_1)m_d),\ j_2 = 2+m_1+\eta_3 A_1(s_1+(1-m_1)m_d+m_d)+\gamma+m_2(1-m_1)+4x_2$$
$$+2\eta_3(1-A_1)(1-m_1)m_d,\ j_3 = s_1 2^{(1-A_1)}+(-1)^{(1-A_1)\gamma}(1-A_1)m_3+A_1 m_1 m_3,\ j_4 = s_2+(1-A_1)m_4$$
$$+A_1 m_3 m_4,\ldots,\ j_{d-1} = s_{d-3}+(1-A_1)m_{d-1}+A_1 m_{d-2}m_{d-1},\ j_d = s_{d-2}+(1-A_1)m_d+A_1 m_{d-1}m_d,$$
$$0 \leq x_2 \leq \tfrac{z_d}{4}-1\},$$

where

$$\chi = \chi(d, m_1, m_2, \ldots, m_d) = 1 - 2^{\left\lfloor \frac{1}{d}\sum_{k=1}^{d} m_k \right\rfloor} \left\lceil \frac{1}{d}\sum_{k=1}^{d} m_k \right\rceil \text{ for } (d, m_1, m_2, \ldots, m_d) \in \mathbb{Z}^{\geq 2} \times \{0,1\}^d,$$

$$M = M(m_1, m_2, \ldots, m_d) \text{ with the logical expression}$$

$$q(m_1, m_2, \ldots, m_d) = (m_1 = m_2 = \cdots = m_d) \vee (m_2 \neq m_1 = m_3 = \cdots = m_d) \vee \bigvee_{k=2}^{d-1}(m_1 = m_2 = \cdots = m_k \neq m_{k+1} = \cdots = m_d)$$

is defined as $M(m_1, m_2, \ldots, m_d) = \begin{cases} 1 & \text{if } q(m_1, m_2, \ldots, m_d) \text{ is true} \\ 0 & \text{if } q(m_1, m_2, \ldots, m_d) \text{ is false} \end{cases}$, $\alpha^* = \dfrac{z_d}{2}\prod_{k=0}^{d-2} A_{k+3}$ with $A_{k+3}$ as defined below,

$\alpha_k = \alpha_{k-\max} = z_d \prod_{j=0}^{k-2} \dfrac{z_{j+1}}{2}$ for $1 \leq k \leq d$ denotes the maximum $\alpha^*$ possible on a $k$-dimensional Cartesian product of the form presented previously, $\eta(y,c) = \begin{cases} 1, & \text{if } y \in \mathbb{Z}^{>c} \\ 0, & \text{if } y \notin \mathbb{Z}^{>c} \end{cases}$ for $y, c \in \mathbb{Z}^{\geq 2}$, $\eta_k = \eta(d, k-1)$ for $3 \leq k \leq d$,

$\eta_{\alpha_k} = \eta\left(z_d(b_k+1)\prod_{j=0}^{k-2}\dfrac{z_{j+1}}{2}, \alpha_k\right)$ for $1 \leq k \leq d$ and with $b_d = 0$, $A_1 = \eta(z_d(b_1+1), \alpha_1)$, $A_2 = \eta(z_d(b_1+1)+1, \alpha_2)$,

$A_3 = (2(b_1+1))^{1-A_2} z_1^{A_2}$, $A_k = (b_{k-2}+1)^{(1-\eta_{\alpha_{k-1}})\eta_{\alpha_{k-2}}}\left(\dfrac{z_{k-2}}{2}\right)^{\eta_{\alpha_{k-1}}}$ for $4 \leq k \leq d+1$, $0 \leq \ell \leq \dfrac{z_1}{A_3}-1$, $0 \leq \gamma \leq 1$,

$r_1 = \left\lfloor\dfrac{2x_1+4}{A_3}\right\rfloor((1-\eta_3)(1-m_2)m_1+\eta_3(1-\eta_4)(1-m_1)m_3+\eta_4 m_d(1-m_{d-1}))$, $r_2 = x_2 - 2\left\lfloor\dfrac{x_2}{2}\right\rfloor$,

$0 \leq p_k \leq \left[2^{\eta_{\alpha_{k+1}}}-1\right]\eta_{k+2}$ for $1 \leq k \leq d-2$ and $0 \leq t \leq 2^{(1-A_1)\eta_3}-1$, $0 \leq s_1 \leq \left[\dfrac{z_2}{A_4 2^{\eta_{\alpha_2}+(1-A_1)}}-1\right]\eta_3$,

$0 \leq s_k \leq \left[\dfrac{z_{k+1}}{A_{k+3} 2^{\eta_{\alpha_{k+1}}}}-1\right]\eta_{k+2}$ for $2 \leq k \leq d-2$, $X_{k-2} = \prod_{j=3}^{k}\left\lfloor\dfrac{x_j+1}{A_{j+1}}\right\rfloor$ and $\Psi_{k-2} = \prod_{j=3}^{k}\left(1-\left\lceil\dfrac{x_j}{A_{j+1}}\right\rceil\right)$ for $2 \leq k \leq d$,

$R_1 = \eta_3((1-\eta_{\alpha_2})+\eta_{\alpha_2}((1-r_2)X_{d-2}+r_2\Psi_{d-2}))+(1-\eta_3)$, $R_k = (1-r_2)X_{k-2}+r_2\Psi_{k-2}$, and $S_k = (1-r_2)m_{k+1}+r_2(1-m_{k+1})$
for $2 \leq k \leq d-1$.

**Note 1:** For $d = 2$, the convention is $(j_1, j_2) \in (\mathbb{Z}/z_1\mathbb{Z}) \times (\mathbb{Z}/z_2\mathbb{Z})$.

For $d = 3$, the convention is $(j_1, j_2, j_3) \in (\mathbb{Z}/z_1\mathbb{Z}) \times (\mathbb{Z}/z_3\mathbb{Z}) \times (\mathbb{Z}/z_2\mathbb{Z})$.

For $d = 4$, the convention is $(j_1, j_2, j_3, j_4) \in (\mathbb{Z}/z_1\mathbb{Z}) \times (\mathbb{Z}/z_4\mathbb{Z}) \times (\mathbb{Z}/z_2\mathbb{Z}) \times (\mathbb{Z}/z_3\mathbb{Z})$.

For $d \geq 5$, the convention is $(j_1, j_2, \ldots, j_d) \in (\mathbb{Z}/z_1\mathbb{Z}) \times (\mathbb{Z}/z_d\mathbb{Z}) \times (\mathbb{Z}/z_2\mathbb{Z}) \times (\mathbb{Z}/z_3\mathbb{Z}) \times \cdots \times (\mathbb{Z}/z_{d-1}\mathbb{Z})$.

As in the Square decomposition, we emphasize that the definition of the edge set above is to be read and used in the following way. For $d = 2$, one considers components $j_1$ and $j_2$ in the edge set definition to define the decomposition of $C_{z_1} \square C_{z_2}$ and disregards components $j_3$ through $j_d$. All parameters associated to these extraneous components can be sent to zero and disregarded. For $d = 3$, we consider the first three components $j_1, j_2$ and $j_3$ to define the cycles to decompose $C_{z_1} \square C_{z_3} \square C_{z_2}$ and disregard the components $j_4$ through $j_d$ presented. As before, all parameters associated to the extraneous components get sent to zero. Continuing in this fashion, we can obtain the corresponding definition for $d = 4$ and $d \geq 5$ to decompose $C_{z_1} \square C_{z_4} \square C_{z_3} \square C_{z_2}$ and $C_{z_1} \square C_{z_d} \square C_{z_2} \square C_{z_3} \square \cdots \square C_{z_{d-1}}$, respectively. Unlike in the Square decomposition, $C_{z_d}$ needs to be the shortest cycle in the torus definition so that the decomposition using the shortest cycles defined by the Lock-and-Key decomposition agrees with that of the Square decomposition using the longest cycles it defines. Of course, one can swap cycles in the torus' definition provided the divisibility conditions are still satisfied and one is aware of how the modifications change the lower bound of the method given the construction's nature.

**Note 2:** Over each interval $\alpha_{k-1} < \alpha^* \leq \alpha_k$ for $2 \leq k \leq d$, we let $0 < b_{k-1} \leq \frac{z_{k-1}}{2} - 1$ be odd and set $b_j = 0$ for all $k \leq j \leq d - 1$. For all $1 \leq j < k - 1$ with $k$ as above, the General Lock-and-Key decomposition edge set definition has $b_j = \frac{z_j}{2} - 1$ for all $\alpha_{k-1} < \alpha^* \leq \alpha_d$. For $\alpha^* = \alpha_1$, $b_j = 0$ for all $1 \leq j \leq d - 1$. The idea is that we only vary $b_{k-1}$ from 0 if we seek to get cycle lengths that have $\alpha_{k-1} < \alpha^* \leq \alpha_k$ for $2 \leq k \leq d$ and leave all other $b_j$ at either 0 if $k - 1 < j \leq d - 1$ or $\frac{z_j}{2} - 1$ if $1 \leq j < k - 1$. Lastly, observe that the Lock-and-Key decomposition edge set definition essentially truncates $b_{k-1}$ at $\frac{z_{k-1}}{2} - 1$ in the definition of $A_{k+1}$ and leaves it as that for all $\alpha_k < \alpha^* \leq \alpha_d$.

Numbering the sets of the General Lock-and-Key decomposition edge set definition one through four from top to bottom, we state the role of each set during the intervals of $\alpha^*$ for which they are active and provide drawings of the foundational components that persist in the decomposition cycles for various cycle lengths, where one traces the given component starting from the $(0, \ldots, 0)$ $(m_1, \ldots, m_d)$ $d$-tuple:

**Set 1:** For all $\alpha_1 < \alpha^* \leq \alpha_d$ and $d \geq 2$, this set defines all keys that have the edges along $j_1$ defined by $j_1$ moves with positive orientation. The $(m_1, \ldots, m_d)$ $d$-tuples for each move used by the definition of set 1 to define part of a given cycle are labelled below for cycles on toruses of dimension $d = 2$ (left), $d = 3$ (center), and $d = 4$ (right):

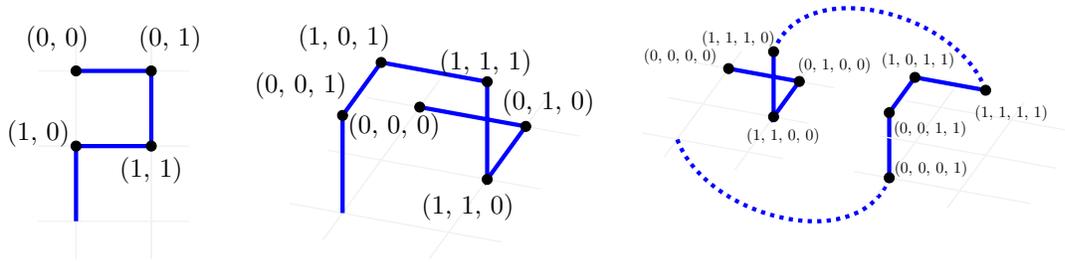

**Set 2:** For all $\alpha_1 < \alpha^* < \alpha_2$ and $d \geq 2$, this set defines all keys that have edges along $j_1$ defined by $j_1$ moves with negative orientation. The associated $d$-tuples for the edges of the components defined by set 2 are labelled below for cycles defined on toruses of dimension $d = 2$ (left), $d = 3$ (center), and $d = 4$ (right):

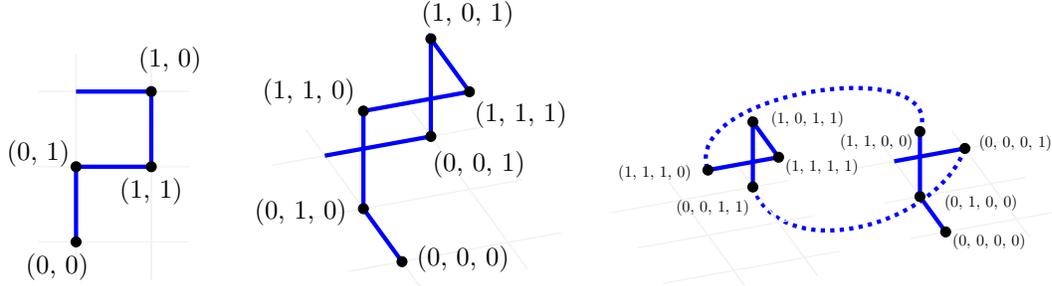

**Set 3:** When $d = 2$ and $\alpha_1 \leq \alpha^* < \alpha_2$, this set defines the left half of the block-like portion of a given cycle where the edges along $j_1$ are defined by $j_1$ moves with strictly positive orientation. If $\alpha^* = \alpha_2$ when $d = 2$, this set defines all bottom portions, including those that would otherwise be defined by set 4, by leveraging a new symmetry that becomes available when $\alpha^* = \alpha_2$. When $d \geq 3$ and $\alpha_1 \leq \alpha^* \leq \alpha_d$, it defines all states that have edges along $j_1$ defined by $j_1$ moves with positive orientation. In particular when $\alpha_2 \leq \alpha^* \leq \alpha_d$ with $d \geq 3$, set 3 defines all possible states that the "tails" take depending on which range $\alpha_{k-1} < \alpha^* \leq \alpha_k$ we find ourselves in for $3 \leq k \leq d$. These states include stair-casing moves, normal moves and column transitions, all of which we define at the start $\alpha^*$-case 4 of dimension case 2 in the proof of Proposition 15. For all $d \geq 2$, this set overall encapsulates $(6d-10)$-many different states total among all components. This is without taking into account all combinations of what type of $j_k$ move is being carried out based on $d$ and $\alpha^*$. Given the level of detail this set has to account for and it being active for all $\alpha_1 \leq \alpha^* \leq \alpha_d$ for every $d \in \mathbb{Z}^{\geq 2}$, it makes sense that this set is the longest of all four sets.

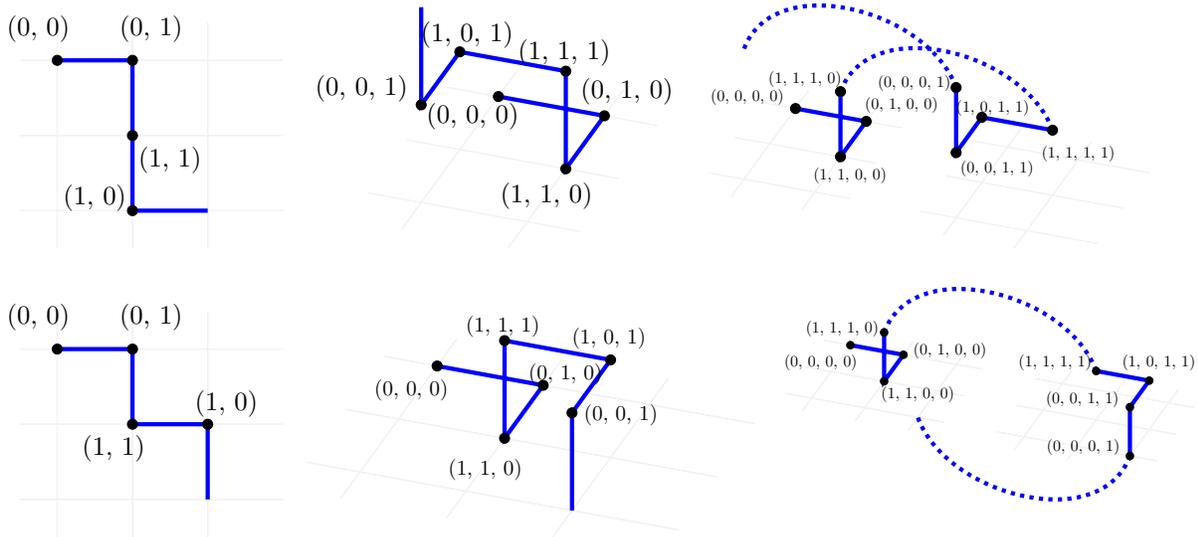

Above we exhibit some of the many components set 3 defines with different combinations of states as applicable with all having positive orientation overall. The left column corresponds to variations of column transitions when $d = 2$ and $\alpha_1 \leq \alpha^* \leq \alpha_2$. The center column's top figure shows stair-casing moves along all three dimensions while the bottom shows column-transition moves. These occur when $d = 3$ over $\alpha_2 < \alpha^* \leq \alpha_3$ and $\alpha_2 \leq \alpha^* \leq \alpha_3$, respectively. The right column's top figure shows stair-casing moves along the lowest three dimensions while the dashed arcs going along the fourth dimension correspond to normal moves. To have stair-casing moves along the fourth dimension for this figure, we would have the last dashed arc directed to a new copy of the same three-dimensional torus instead and this would take place when $d = 4$ and $\alpha_3 < \alpha^* \leq \alpha_4$. Lastly, the right column's bottom figure corresponds to column transition moves along all four dimensions. The right column states take place when $d = 4$ for $\alpha_2 < \alpha^* \leq \alpha_4$ and $\alpha_2 \leq \alpha^* \leq \alpha_4$, respectively.

**Set 4:** When $d = 2$ and $\alpha_1 \leq \alpha^* < \alpha_2$, this set defines the right half of the block-like portion of a given cycle, where the edges along $j_1$ are defined by $j_1$ moves with strictly negative orientation. When $d \geq 3$ and $\alpha_1 \leq \alpha^* < \alpha_2$, set 4 defines a slightly different version of the right half of the block-like section, where the second half is translated to the right in the case $\alpha_1 < \alpha^* < \alpha_2$ and remains the same for $\alpha^* = \alpha_1$. Below we show the component of a given cycle defined by set 4 on a torus of dimension $d = 2$ (left), $d = 3$ (center), and $d = 4$ (right). Note that the edge not drawn along $j_2$, the horizontal dimension, between the two sub-components when $d = 3$ and $d = 4$ as seen below comes from set 2 when $x_1$ is maximal as part of the transition from the $d = 2$ instance where it is used as part of the connected component defined by set 2.

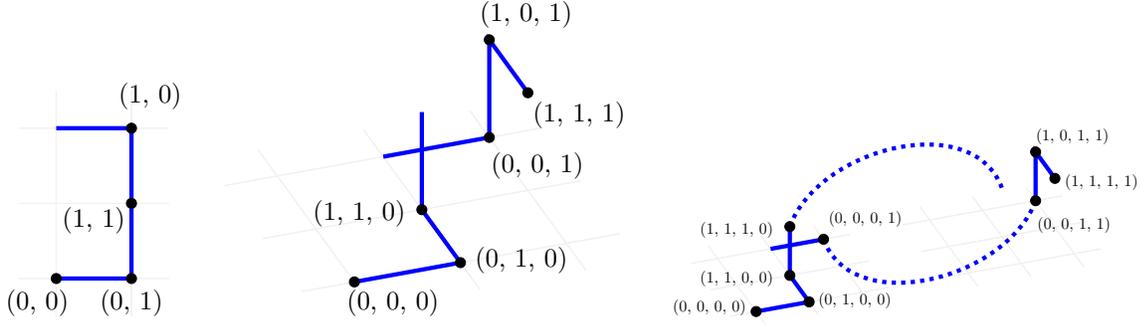

It is worth appreciating that by parsing the cycle definition as outlined above, partitioning the cycle into these sets with the corresponding roles to be carried out when applicable, we can most notably transition from the definition of the cycles for $\alpha_1 < \alpha^* < \alpha_2$ when $d = 2$ to that for the same range of $\alpha^*$ when $d \geq 3$ by simply performing translations of the fundamental pieces each set defines. Then, for $\alpha_2 \leq \alpha^* \leq \alpha_d$, most of the non-symmetric detail can be placed on set 3, where we have leveraged the new symmetries that become available for all $\alpha_2 \leq \alpha^* \leq \alpha_d$ to have set 1 define the remaining symmetric components, namely the key "bits".

**Remark:** The convention to use for the form of the Cartesian product based on $d \in \mathbb{Z}^{\geq 2}$ is presented above in Note 1. Note the following modifications that can be made to the parameters of the Lock-and-Key Decomposition's edge set definition presented earlier for easier usage in obtaining cycle decompositions of $Q_{2an}$, where $n \geq 1$ and $a = 2^{i_1} + \cdots + 2^{i_d}$ is odd with integers $i_1 > \cdots > i_d = 0$, for the established lengths assuming $z_j = 2^{n 2^{i_j+1} - i_j}$ for all $1 \leq j \leq d$:

$$\alpha = \log_2(\alpha^*), \ \alpha_k = \alpha_{k-\max} = \sum_{j=1}^{k-1}\left[n2^{i_j+1} - i_j\right] + 2n - (k-1) \ \text{ for } 1 \leq k \leq d \text{ denotes the maximum}$$

$\alpha$ possible on a $k$−dimensional Cartesian product of the form presented previously, $\eta_k = \eta(d, k-1)$ for $3 \leq k \leq d$, $\eta_{\alpha_k} = \eta(\alpha, \alpha_k)$ for $2 \leq k \leq d$, $A_1 = \eta(\alpha, \alpha_1)$, $A_2 = \eta(\alpha+1, \alpha_2)$,

$$A_3 = 2^{(1-A_2)(\alpha-\alpha_1) + A_2(\alpha_2-\alpha_1)+1}$$

with $b_1 = 2^{\alpha-\alpha_1} - 1$, and

$$A_k = 2^{(1-\eta_{\alpha_{k-1}})(\alpha-\alpha_{k-2})\eta_{\alpha_{k-2}} + \eta_{\alpha_{k-1}}(\alpha_{k-1}-\alpha_{k-2})}$$

with $b_{k-2} = 2^{\alpha-\alpha_{k-2}} - 1$ for $4 \leq k \leq d+1$.

Analogous modifications can also be made when all cycles $C_{z_j}$ defining the tori, of the forms presented in Note 1, have lengths greater than or equal to four that are powers of two and are not necessarily those on which we focus. We can take all as stated above and modify $\alpha_k$ for $1 \leq k \leq d$ by applying the logarithm base 2 to the definition of $\alpha_k$ in the Lock-and-Key decomposition's edge set definition to get

$$\alpha_k = \sum_{j=1}^{k-1}[\log_2(z_j)] + \log_2(z_d) - (k-1).$$

**Lemma 6** *Let $d \geq 2$ and $a = 2^{i_1} + 2^{i_2} + \cdots + 2^{i_d}$ be odd with integers $i_1 > i_2 > \cdots > i_d = 0$. Now let $(C_{y_1^*}, S_1)$ be such that $S_1$ has every $2^{i_1}$th vertex of $C_{y_1^*}$, $(C_{y_2^*}, S_2)$ has $S_2$ consisting of every $2^{i_2}$th vertex of $C_{y_2^*}$, and so forth with $(C_{y_d^*}, S_d)$ having $S_d$ contain every vertex of $C_{y_d^*}$, where each $y_j^* \geq 4$ for $1 \leq j \leq d$ is a power of two. Further, let $A_j$ for $3 \leq j \leq d+1$ and $\alpha_k$ for $1 \leq k \leq d$ be as in the above modification with $z_j = y_j^* 2^{-i_j}$ for $1 \leq j \leq d$ assumed as in the Lock-and-Key decomposition's edge set definition. If $4 \mid y_d^*$ and $A_3 2^{i_1} \mid y_1^*$ when $d = 2$ and $A_3 2^{i_1} \mid y_1^*$, $A_4 2^{i_2+1} \mid y_2^*$, $A_5 2^{i_3+1} \mid y_3^*$, $\ldots$, $A_{d+1} 2^{i_{d-1}+1} \mid y_{d-1}^*$ and $4 \mid y_d^*$ when $d \geq 3$, it follows that $C_{a \cdot 2^\alpha}$ decomposes $(C_{y_1^*}, S_1) \boxplus (C_{y_d^*}, S_d)$ when $d = 2$ and $(C_{y_1^*}, S_1) \boxplus (C_{y_d^*}, S_d) \boxplus \cdots \boxplus (C_{y_2^*}, S_2)$ when $d \geq 3$.*

**Proof:** Let $d \geq 2$, $a = 2^{i_1} + \cdots + 2^{i_d}$ be odd with integers $i_1 > \cdots > i_d = 0$, $\alpha \geq \alpha_1 = \log_2(z_d)$, $b_j = 2^{\alpha - \alpha_j} - 1$ for all $1 \leq j \leq d-1$, and $b_d = 0$. Given that the Lock-and-Key decomposition defines the cycles to have $(z_d/2) \prod_{k=0}^{d-2} A_{k+3} = (z_d/2) 2^{\alpha - (\alpha_1 - 1)} = 2^\alpha$ edges along each dimension in the underlying torus by Theorem 5, we see by Proposition 1 that $2^\alpha$ edges in the underlying torus along the $k$th dimension for $k = 1, \ldots, d$ corresponds to $2^{i_k} \cdot 2^\alpha$ edges in the subdivided torus. Hence, summing the edges along each dimension gives us that our cycle decomposition of the underlying torus using cycles of length $(dz_d/2) \prod_{k=1}^{d-2} A_{k+3}$ corresponds to a cycle decomposition of the given anchored product using cycles of length $a \cdot 2^\alpha$, where all cycles are of the same length as a consequence. ∎

The form of the anchored product above is emphasized as choosing $C_{y_d^*}$ as we do below to be the shortest cycle of all the $C_{y_j^*}$ allows us to obtain the longest cycle defined by Square decomposition as the shortest cycle defined by the Lock-and-Key decomposition while yielding cycle decompositions using cycles of all the lengths it covers.

**Corollary 7** *Let $n \geq 1$, $d \geq 2$, $a = 2^{i_1} + 2^{i_2} + \cdots + 2^{i_d}$ with integers $i_1 > i_2 > \cdots > i_d = 0$, $A_{k+2}$ with $b_k = 2^{\alpha - \alpha_k} - 1$ for $1 \leq k \leq d-1$ with $\alpha$ and $\alpha_k$ as in the remark following the Lock-and-Key edge set definition, and $y_j^* = 2^{n 2^{i_j}+1}$ for all $1 \leq j \leq d$. Then, a $d$-dimensional Lock-and-Key decomposition of $(C_{2^{n 2^{i_1}+1}}, S_1) \boxplus (C_{2^{2n}}, S_d) \boxplus (C_{2^{n 2^{i_2}+1}}, S_2) \boxplus \cdots \boxplus (C_{2^{n 2^{i_{d-1}}+1}}, S_{d-1})$ can be constructed using $C_{a \cdot 2^\alpha}$ cycles, where*

$$2n \leq \alpha \leq 2an - (d-1) - \sum_{k=1}^{d} i_k.$$

*Finally, observe that decomposing such anchored products is equivalent to decomposing a torus of the form presented below as*

$$(C_{2^{n 2^{i_1}+1}}, S_1) \boxplus (C_{2^{2n}}, S_d) \boxplus (C_{2^{n 2^{i_2}+1}}, S_2) \boxplus \cdots \boxplus (C_{2^{n 2^{i_{d-1}}+1}}, S_{d-1})$$
$$\cong C_{2^{n 2^{i_1}+1} - i_1} \square C_{2^{2n}} \square C_{2^{n 2^{i_2}+1} - i_2} \square \cdots \square C_{2^{n 2^{i_{d-1}}+1} - i_{d-1}}$$

*and the cycles in the decomposition of the above torus have the same number of edges along every dimension.*

Thus, from Corollaries 4 and 7 as applicable based on the desired cycle lengths, we can now see that $Q_{(2a)n}$ can be decomposed into

$$Q_{(2a)n} = Q_{2(2^{i_1}+2^{i_2}+\cdots+2^{i_d})n} = Q_{n2^{i_1}+1} \square Q_{n2^{i_2}+1} \square \cdots \square Q_{n2^{i_d}+1}$$
$$= n2^{i_1}C_{2^{n2^{i_1}+1}} \square n2^{i_2}C_{2^{n2^{i_2}+1}} \square \cdots \square n2^{i_d}C_{2^{n2^{i_d}+1}}$$
$$= n\left(\prod_{k=1}^{d} 2^{i_k}\right)\left[(C_{2^{n2^{i_1}+1}}, S_1) \boxplus (C_{2^{n2^{i_2}+1}}, S_2) \boxplus \cdots \boxplus (C_{2^{n2^{i_d}+1}}, S_d)\right]$$
$$= n2^{\sum_{k=1}^{d} i_k}\left[2^{n\sum_{k=1}^{d} 2^{i_k+1} - \sum_{k=1}^{d} i_k - \alpha}\right]C_{a\cdot 2^\alpha} = n2^{2an-\alpha}C_{a\cdot 2^\alpha}.$$

This is $n2^{2an-\alpha}$ cycles each of length $a \cdot 2^\alpha$ for $1 \leq \alpha \leq 2an - (d-1) - \sum_{k=1}^{d} i_k$.

As we will see in the proofs in the subsequent sections, our constructions for the Square and Lock-and-Key decompositions' edge set definition reveal a partition by which we can decompose the proofs of the edge-disjoint cycle decomposition components in the statements of Theorems 2 and 5. This is a consequence of the edge set definitions defining all edges with orientations that agree with those the edges would have if one were tracing out a cycle in the decomposition. In the case of the Lock-and-Key decomposition's edge set definition, we leverage the symmetries we referenced earlier so that we need only treat the $\alpha_1 < \alpha^* < \alpha_2$ for $d = 2$ and $d \geq 3$ separately using four major sets while $\alpha_2 \leq \alpha^* \leq \alpha_d$ can be adequately treated with major sets 1 and 3. Lastly, the presentation of the edge set definitions with explicit algebraic closed forms for each of the components will enable us to reformulate the edge-disjoint subsection proofs in terms of systems of equations and the number-types of their solutions.

We conclude this section with applications of the Lock-and-Key decomposition to yield cycle decompositions of Cartesian products of cycles, leading up to another corollary of Theorems 2 and 5.

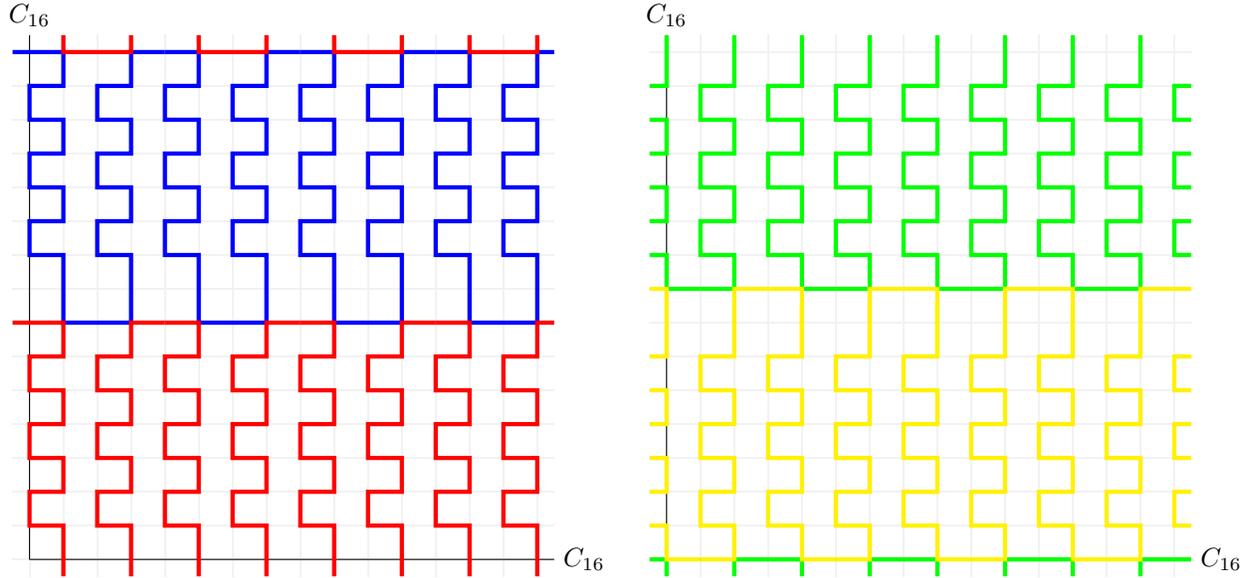

Here we have the $C_{128}$ cycles defined by the Lock-and-Key decomposition on the two dimensional torus $C_{16} \square C_{16}$ with $\gamma = 0$ (left) being "key" cycles and $\gamma = 1$ (right) being the complementing "lock" cycles. These two figures when combined yield the cycle decomposition below of $C_{16} \square C_{16}$ using $C_{128}$ cycles defined by the Lock-and-Key decomposition for all $\gamma \in \{0, 1\}$.

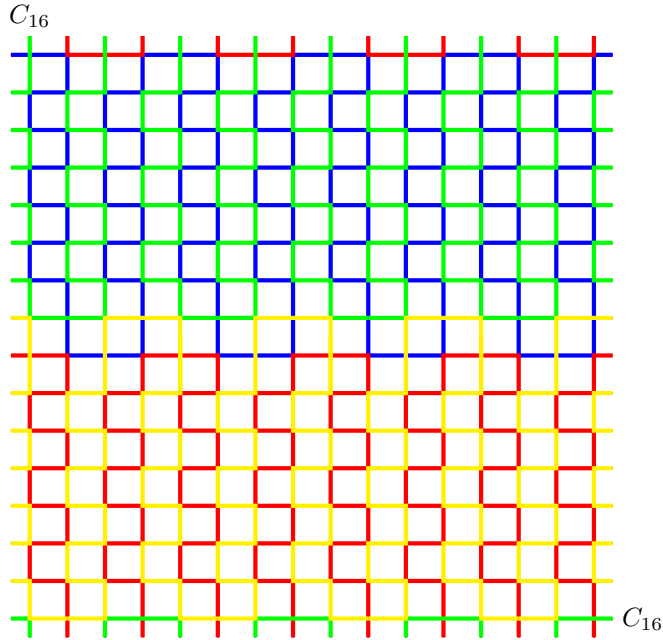

Replacing the cycles defining the vertical dimension of the torus above with $C_{2^{14}}$'s and performing $2^{10} - 1$ vertical translations of 16 edges each to the above cycle decomposition, we obtain a cycle decomposition of the underlying torus $C_{2^{14}} \square C_{16}$ of $Q_{20}$ using cycles that correspond to $C_{5 \cdot 2^6}$'s in the subdivided torus of $Q_{20}$. Since we enlarged our initial torus, note that we can construct cycle decompositions using longer cycles via the Lock-and-Key Decomposition by proceeding analogously as in our initial decomposition.

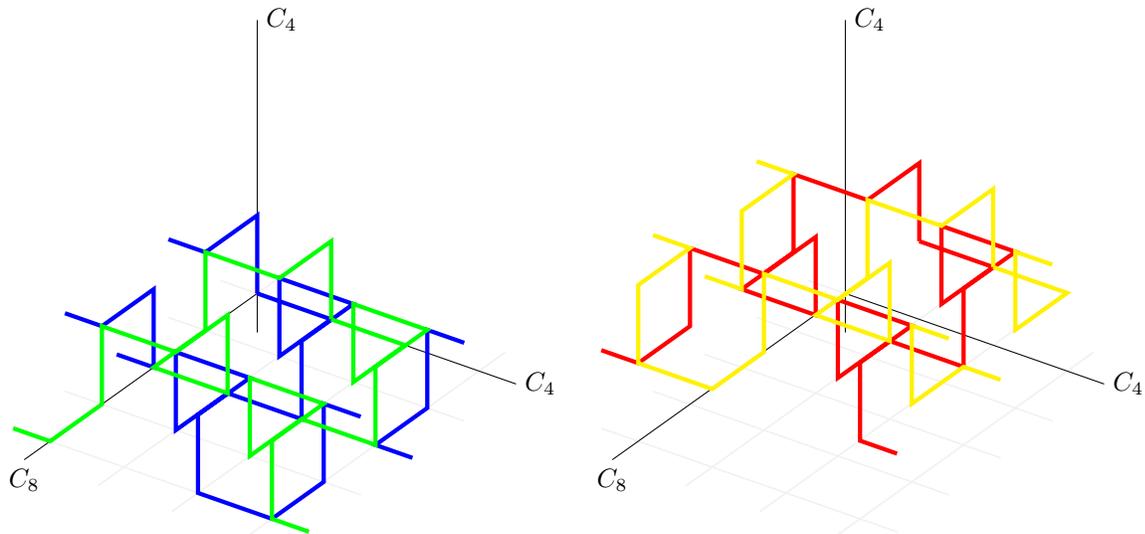

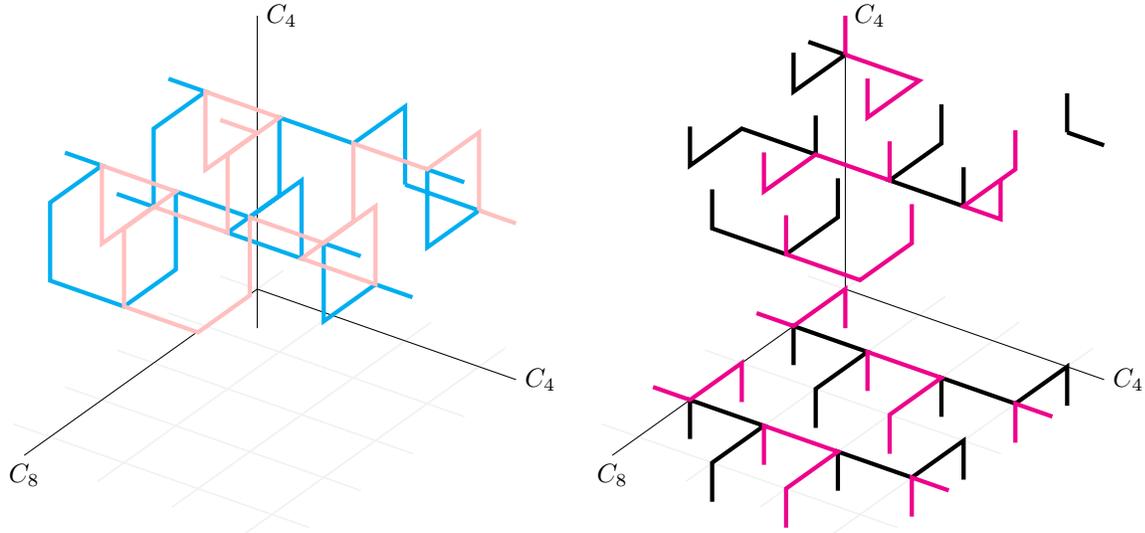

In the four figures above, we show the cycles defined by the Lock-and-Key Decomposition on a tract of the three-dimensional torus $C_8 \square C_4 \square C_4$ using $C_{24}$'s corresponding to the appropriate $(\ell, \gamma, s_1)$ pairs with $d = 3$ and all other parameters being zero. In this case, this corresponds to the interval of cycle lengths with $\alpha_1 < \alpha^* < \alpha_2$ as we are working with cycle lengths that require one level and hence have analogs in the two-dimensional torus $C_8 \square C_4$.

The left figure below gives a bird's-eye view of one of the cycles defined by the Lock-and-Key decomposition on the three-dimensional torus above. As we can see, the most notable difference between it and its two-dimensional analog is the absence of the lower block segment used in transitioning to the next column. Instead, we find that by translating the keys and performing the column transition in a manner more coherent with the definition of the keys, we are able to preserve the compatibility in symmetries along all required instances of a given cycle to leverage the translational symmetries that were present in the two-dimensional analog as well as those that arise in the higher-dimensional toruses. If one were to expand a cycle, defined by the Lock-and-Key decomposition of a two-dimensional torus such as $C_8 \square C_4$, along the third dimension to use in the decomposition of $C_8 \square C_4 \square C_4$, it would become apparent that the difference contains precisely the kind of translational symmetry that makes it incompatible with the symmetries of the keys. Conversely, we see that using the contracted cycle above to decompose $C_8 \square C_4$ would not be valid as we would not be able to introduce a cycle of the same length that in particular covers the square of edges formed at the center of the main key.

Hence, the two-dimensional definition of the Lock-and-Key decomposition is a special case of the Lock-and-Key decomposition relative to its definition for dimensions three and higher for $\alpha_1 \leq \alpha^* \leq \alpha_2$.

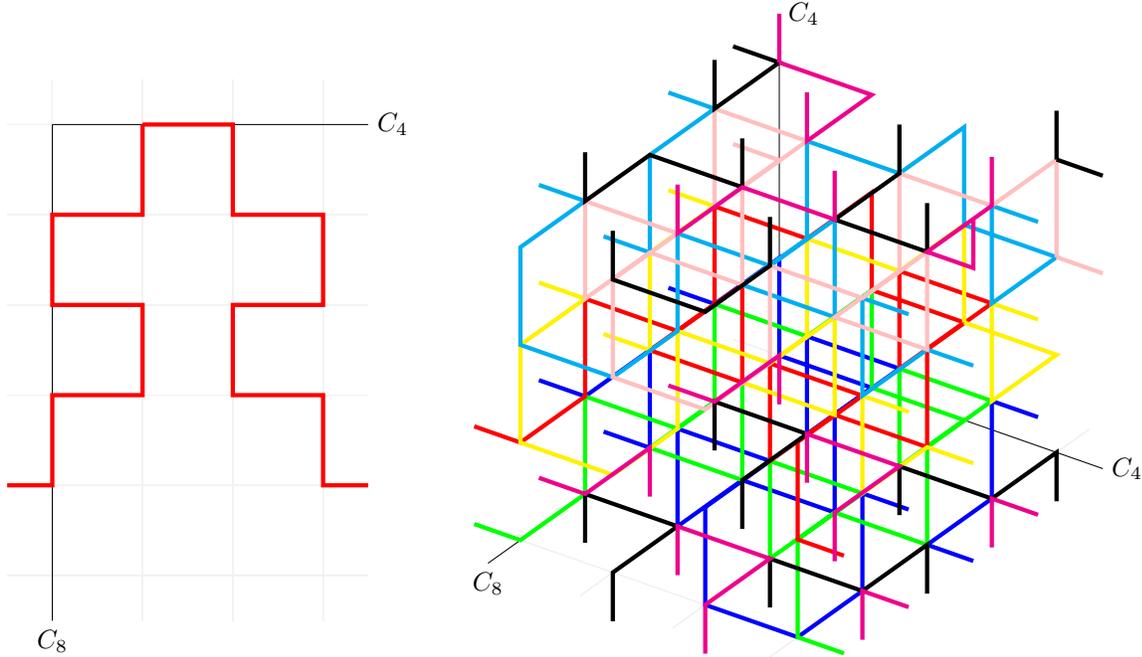

Combining the previous four figures of $C_{24}$ cycles defined by the Lock-and-Key Decomposition's edge set definition, we get the cycle decomposition on the right for a tract of the torus $C_8 \square C_4 \square C_4$. Observe that the two edges missing at each extreme of a given level are not present due to them belonging to other cycles we have not drawn that would fall in the other tract that would complete the torus. To get the complete decomposition for this torus, we simply need to translate all of the cycles above four units in the $j_1$ direction whose axis is defined by $C_8$.

However, we would like to remark that if our torus had instead been $C_4 \square C_4 \square C_4$, we can use the Lock-and-Key decomposition as defined for $\alpha^* = \alpha_2$. Alternatively, we can notice that the above definition defined for $\alpha_1 < \alpha^* < \alpha_2$ in fact works for $\alpha^* = \alpha_2$ as the first and last line along the $C_4$ in the $j_2$ direction would be the same. One would only need to partition the edges at the extremes in the $j_1$ direction as we did with the $j_2$ and $j_3$ edges found at the extremes of the torus.

While both yield the same cycle decomposition, the definition for $\alpha^* = \alpha_2$ is more compact requiring only two sets to define while using the observed alternative would require four sets as do the $\alpha_1 < \alpha^* < \alpha_2$. As can be seen, the definition for the Lock-and-Key decomposition we push forward for $\alpha_1 \leq \alpha^* < \alpha_2$ treats what we regard as special cases relative to the more coherent portion of the Lock-and-Key decomposition's definition that becomes apparent as we consider longer cycles with $\alpha_2 \leq \alpha^* \leq \alpha_d$ for all $d \in \mathbb{Z}^{\geq 2}$.

In the following figure, we show one $C_{128}$ cycle from the Lock-and-Key Decomposition of $Q_8$ on the torus $C_4 \square C_4 \square C_4 \square C_4$. To obtain the remaining 7 cycles, we can translate the cycle below to the right by one edge. Next, we can translate the two cycles up by one edge each, where we view the top left $C_4 \square C_4 \square C_4$ torus as the starting point of each cycle. Having four cycles at this point, the remaining four are obtained by starting the cycle decomposition from the top right $C_4 \square C_4 \square C_4$ torus to get the reference cycle that we can then translate as referenced previously to get the remaining three. In the case of the class of hypercubes of interest, we can consider the four-dimensional underlying torus $C_{2^{13}} \square C_{2^6} \square C_{2^3} \square C_{2^2}$ of $Q_{30}$ for example with the cycle decomposition using $C_{2^{23}}$'s following analogously to the figure below with an extended set of keys on twice as many levels as before. This would correspond to a cycle decomposition of $Q_{30}$ using $C_{15 \cdot 2^{21}}$-cycles.

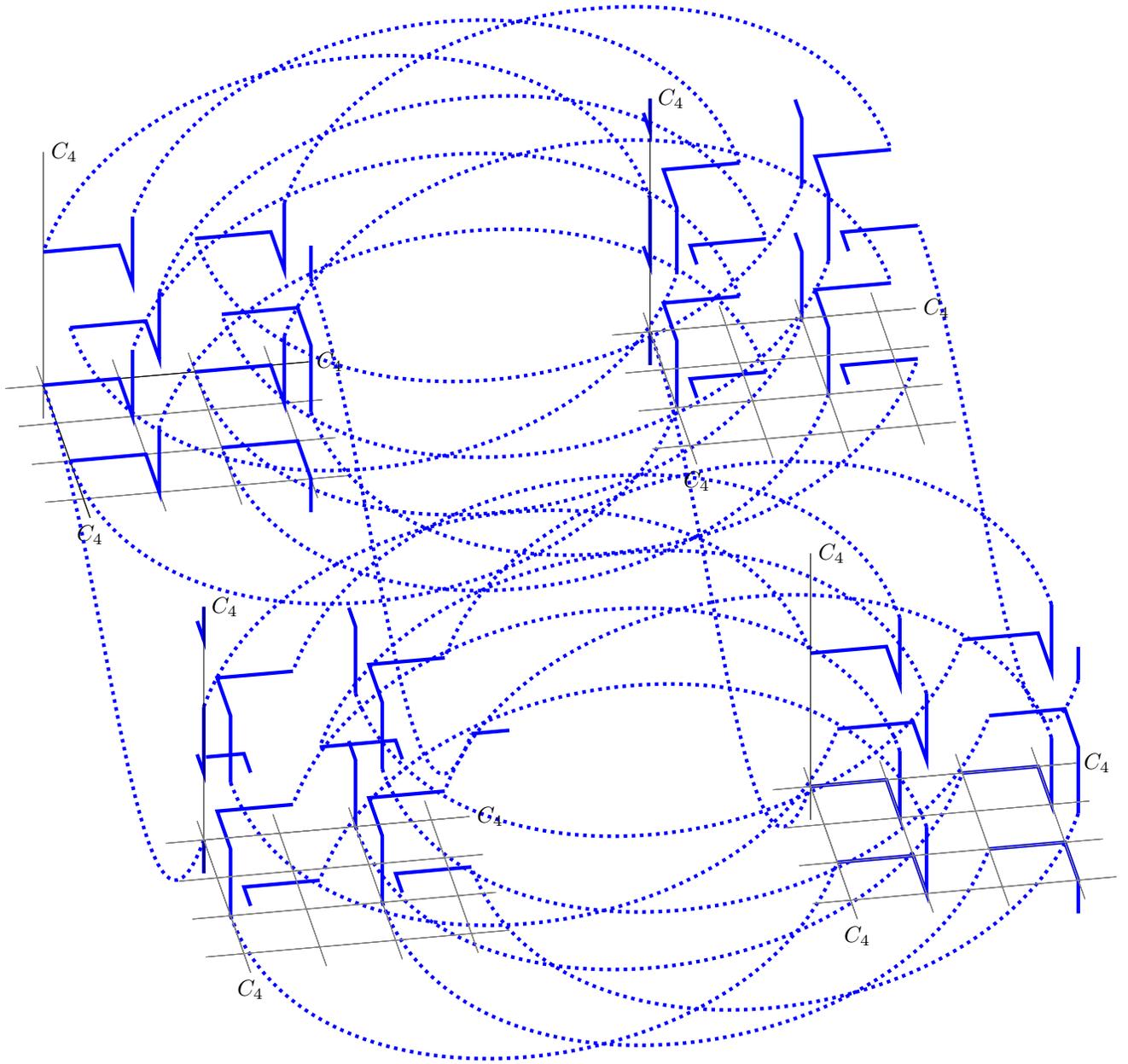

With the above being the smallest four-dimensional torus that can be decomposed using the Lock-and-Key Decomposition, we discover that this decomposition method gives rise to cycles decompositions into cycles of the form $C_{a\cdot 2^\alpha}$ for $2 \leq \alpha \leq 5$ with $a = 4$ and is hence not limited to odd $a$. In fact, this holds true for the Square decomposition as well, where for both methods we can have $a \geq 2$ be even and decompose the $a$-fold Cartesian product of $Q_{2n}$'s defining $Q_{2an}$. For said $a$, we can obtain cycle decompositions using cycles of the form $C_{2^\beta}$ for $\beta \geq 2$, which results in all admissible cycle decompositions in the case $a$ is a positive power of two, and additional cycle decompositions using cycles of said form otherwise to decompose hypercubes such as $Q_{220}$ and $Q_{420}$.

However, observe that the dimension of the torus is heavily dependent on $a$ if we consider the $a$-fold Cartesian product of $Q_{2n}$'s, causing the exclusion of some of the shortest and longest cycles that could have been used in the decomposition. This is due to how the dimension $a$ would predetermine the form of the cycle to be $C_{a\cdot 2^\alpha}$ for $\alpha \geq 1$ for both cycle decomposition methods, and the Lock-and-Key decomposition method's upperbound decaying linearly by $a$. To remedy this, recall from earlier that between the Square and Lock-and-Key decomposition methods we can fully decompose any two-dimensional torus of $Q_{2an}$ using cycles of all the lengths admissible by $Q_{2an}$. Hence, we have the following corollary of Theorems 2 and 5.

**Corollary 8** *Let $n \geq 1$ and $b \geq 2$ be even. Then, $Q_{2bn}$ can be decomposed into $bn2^{2bn-\beta}$-many cycles of length $2^\beta$ for all $2 \leq \beta \leq 2bn$.*

**Proof:** Let $d = 2$ and $\beta \geq 2$. Expressing $Q_{2bn}$ as $Q_{2bn} = Q_{bn} \square Q_{bn}$ with $b$ and $n$ as above, we find that we can immediately represent it as copies of a two-dimensional torus in the following way:

$$Q_{2bn} = Q_{bn} \square Q_{bn} = \frac{bn}{2} C_{2^{bn}} \square \frac{bn}{2} C_{2^{bn}} = \frac{bn}{2}(C_{2^{bn}} \square C_{2^{bn}}).$$

Letting $y_1 = 2^{bn}$, $y_2 = 2^{bn}$ and $\lambda = 2^{\beta-2} > 0$, the Square decomposition's Theorem 2 gives us that $C_{2^{bn}} \square C_{2^{bn}}$ can be decomposed into $2^{2bn+1-\beta}$-many cycles of length $2d\lambda = 2^\beta$ for $2 \leq \beta \leq bn+1$.

Now, taking $z_1 = 2^{bn}$ and $z_2 = 2^{bn}$, we see that since $b = 2c$ for some $c \geq 1$, $2^{bn} = 2^{2nc} = 4^{nc}$ and so it follows that $4 \mid 2^{bn}$. With $b_1 = 2^{\beta-(bn+1)} - 1$, we have $A_3 = 2^{(\beta-bn)(1-A_2)+A_2 bn}$ divides $2^{bn}$ for all $bn+1 \leq \beta \leq 2bn$. So by the Lock-and-Key decomposition's Theorem 5, we get that $C_{2^{bn}} \square C_{2^{bn}}$ can be decomposed into $2^{bn+1}/A_3 = 2^{2bn+1-\beta}$ cycles of length $dz_2 A_3/2 = 2^{bn} \cdot 2^{\beta-bn} = 2^\beta$ for $bn+1 \leq \beta \leq 2bn$ between the two cases $A_2 = 0$ and $A_2 = 1$.

Combining the above, we have that we can decompose the torus $C_{2^{bn}} \square C_{2^{bn}}$ into $2^{2bn+1-\beta}$-many cycles of length $2^\beta$ for all $2 \leq \beta \leq 2bn$. Hence, $Q_{2bn}$ can be decomposed into $bn/2 \cdot 2^{2bn+1-\beta} = bn2^{2bn-\beta}$-many $C_{2^\beta}$-cycles over the same range of $\beta$ as in the torus. ∎

## 6 Proofs of Theorem 2 and Corollary 4

We will prove that the proposed General Square Decomposition edge set definition indeed defines edge-disjoint cycles. Following this, we will show that all cycles are of the same length and in particular of length $a \cdot 2^\alpha$ for $1 \leq \alpha \leq 2n$ in the case of the cycle decompositions of $Q_{2an}$. Combining these results, we will have proven Theorem 2 and Corollary 4.

### 6.1 Edge-Disjoint Cycles Proofs:

Before we begin, we would like to define what we mean when we use the terms "critical start vertex", "intermediate vertices", and "critical end vertex" throughout the proofs.

Let's consider a move along $j_k$ for $1 \leq k \leq d$. **Critical start vertex** refers to the "corner" where a given move $(m_1, \ldots, m_d)$ begins. So $x_k = 0$ in the General Square Decomposition corresponds to the critical start vertex. **Intermediate vertices** refer to the vertices that are used along the way in defining the edges that end up in the edge set for a given move $(m_1, \ldots, m_d)$. Hence, given our previous definition and $0 \leq x_k \leq \lambda - 1$ by construction, intermediate vertices will be those vertices reached when $1 \leq x_k \leq \lambda - 1$, meaning we will require $\lambda \geq 2$. **Critical end vertex** refers to furthest end vertex (right vertex) in the last edge (vertex pair) that is included by the definition of the edge set for a given move $(m_1, \ldots, m_d)$. Since a given move along $j_k$ will have been fully performed, the end vertex in the $k$-th component will differ by $\lambda$ from the critical start vertex with the signage dependent on which of $(m_2 - m_1), \chi, (m_1 - m_3), (m_3 - m_4), \ldots, (m_{d-1} - m_d)$ applies to $j_k$ by its definition in the General Square Decomposition's edge set presented in Theorem 2.

Our main result for this section will be

**Theorem 9** *The General Square Decomposition's Edge set yields cycles that share an edge if and only if the two cycles are the same cycle.*

Now, let $C_1 = C_{\ell,t,p_1,z,s_1,\ldots,s_{d-2}}$ and $C_2 = C_{\ell_1,t_1,p_1^*,z_1,s_1^*,\ldots,s_{d-2}^*}$ with the parameters ranging over the intervals established in the general Square Decomposition's edge set definition. Note that any other parameters with no subscripts that are fixed for a given cycle, meaning they are identifiers of the cycle in some fashion, will have subscript 1 if they belong to $C_2$ and none if they belong to $C_1$. Further, if any of the parameters have a subscript prior to distinguishing between those belonging to $C_1$ and $C_2$, having no superscript will correspond to those parameters belonging to $C_1$ and a $*$ for a superscript will correspond to those parameters belonging to $C_2$.

Given the versatility of the edge set definition in that it performs a move $(m_1, \ldots, m_d)$ and its inverse when required, we must prove the following:

**Proposition 10** *Two cycles $C_1$ and $C_2$ cannot share an edge as a consequence of $C_2$ performing the inverse move to a given move $(m_1, \ldots, m_d)$ starting from the critical end vertex at which $C_1$ was left off by the end of the referenced move.*

**Proposition 11** *Two cycles $C_1$ and $C_2$ cannot share an edge as a consequence of $C_2$ performing the inverse move to a given move $(m_1, \ldots, m_d)$ starting from the intermediate vertices used in defining the edges belonging to $C_1$ as a consequence of the referenced move.*

**Proposition 12** *Two cycles $C_1$ and $C_2$ cannot share an edge as a consequence of $C_2$ performing the same move $(m_1, \ldots, m_d)$ performed by $C_1$ unless $C_2$ performs the move from the same critical start vertex corresponding to that move in $C_1$, which is if and only if $C_1$ and $C_2$ are the same cycle i.e. $(C_1 = C_2)$.*

With this, we will have shown that the cycles resulting from the General Square Decomposition edge set are edge-disjoint, meaning the edges in a cycle's edge set belong uniquely to that cycle.

Let $1 \leq \lambda \leq \lambda^*$, where

$$\lambda^* = \max\{\lambda \in \mathbb{Z}^+ : 2\lambda \mid y_j \text{ for all } 1 \leq j \leq d\}.$$

Observe that the set over which the maximum is being taken in the definition of $\lambda^*$ is non-empty as $\lambda = 1$ is in the set since each $y_j$ is even by assumption. Note that $(m_2 - m_1), \chi, (m_1 - m_3), (m_3 - m_4), \ldots, (m_{d-1} - m_d)$ are $\pm 1$ during $j_1, j_2, j_3, j_4, \ldots, j_d$ moves respectively, and 0 otherwise by the construction $M$ in Theorem 2. $M$ only admits $d$-tuples $(m_1, \ldots, m_d) \in \{0,1\}^d$ with the properties that make the above mappings the case for a given move. So only one of the mappings above will ever be $\pm 1$ for a given move and all others will be 0. Lastly, we will make use of the property $\lambda_j - \lambda_k = \text{sgn}(\lambda_j - \lambda_k)|\lambda_j - \lambda_k| = (-1)^{\lambda_k}|\lambda_j - \lambda_k|$ when $\lambda_j, \lambda_k \in \{0,1\}$ for $j, k \in \mathbb{Z}^+$ and hence

the case that $|\lambda_j - \lambda_k| \leq 1$. We will be using the properties above without further mention in the following proofs.

### 6.1.1 Proof of Proposition 10:

Let $1 \leq \lambda \leq \lambda^*$, $d \in \mathbb{Z}^{\geq 2}$, and for reference refer to Theorem 2 for specific bounds of each parameter, though we will bring some of them up in our arguments as necessary.

We would also like to make the reader aware that the form of the inverse move on $C_2$ to a move $(m_1, m_2, \ldots, m_d)$ on $C_1$ is $(1 - m_1, 1 - m_2, \ldots, 1 - m_d)$ in the case of moves along $j_k$ for $1 \leq k \leq d$ with $k \neq 3$. Given the alternating behavior for moves along $j_3$ as they depend on $r$ and $r_1$, the pairity of $\ell$ and $\ell_1$, respectively, it follows that inverse $j_3$ moves on $C_2$ with respect to a $j_3$ move $(m_1, m_2, \ldots, m_d)$ on $C_1$ is

$$(|r-r_1|m_1+(1-|r-r_1|)(1-m_1), |r-r_1|m_2+(1-|r-r_1|)(1-m_2), \ldots, |r-r_1|m_d+(1-|r-r_1|)(1-m_d)).$$

Note that we will be assuming for the sake of contradiction that two cycles share a vertex with the move configurations established below, meaning all components of their vertices agree, and we then show that such scenarios do not occur by the definition of the General Square Decomposition. In doing so, we will have shown that these cycles cannot share edges by way of these configurations. We now proceed with the proof of Proposition 10 by considering equalities at each component $j_k$ for $1 \leq k \leq d$ with $C_2$ configured to perform the inverse move to a given move $(m_1, \ldots, m_d)$ of $C_1$ along one of the components starting from the critical end vertex belonging to $C_1$ for that move:

**Major Case 1:** Looking at $j_1$, it follows that

$$t + (2\nu + r + m_1 + (m_2 - m_1))\lambda = t_1 + (2\nu_1 + r_1 + 1 - m_1)\lambda$$
$$\implies 2(\nu - \nu_1) + (r - r_1) + (m_1 + m_2) - 1 = \frac{t_1 - t}{\lambda}.$$

Given $2(\nu - \nu_1) + (r - r_1) + (m_1 + m_2) - 1 \in \mathbb{Z}$ and $|t - t_1| \leq \lambda - 1$ by construction, it must be the case $t = t_1$. Noting that $m_1 + m_2 = 1$ during moves in $j_1$, our equality in $j_1$ then implies

$$\nu - \nu_1 = \frac{r_1 - r}{2}.$$

Since $\nu - \nu_1 \in \mathbb{Z}$ and $|r - r_1| \leq 1$ by construction, we have $r = r_1$ and hence $\nu = \nu_1$.

<u>Dimension Case 1:</u> Let $d = 2$. Then our equation for $j_2$ becomes

$$t + (2\mu + r + m_2)\lambda = t_1 + (2\mu_1 + r_1 + 1 - m_2)\lambda$$
$$\implies (\mu - \mu_1) + m_2 = \frac{1}{2}.$$

This is a contradiction as $(\mu - \mu_1) + m_2 \in \mathbb{Z}$ for every value of $m_2 \in \{0, 1\}$ while $\frac{1}{2} \notin \mathbb{Z}$. So the above equality is not possible.

<u>Dimension Case 2:</u> Let $d \geq 3$. Looking at $j_3$, we see

$$s_1 + 2z\lambda + (-1)^r m_3 \lambda = s_1^* + 2z_1\lambda + (-1)^{r_1}(1 - m_3)\lambda$$

$$\implies 2(z - z_1) + (-1)^r m_3 - (-1)^{r_1}(1 - m_3) = \frac{s_1^* - s_1}{\lambda}.$$

Following from $2(z - z_1) + (-1)^r m_3 - (-1)^{r_1}(1 - m_3) \in \mathbb{Z}$ and $|s_1 - s_1^*| \leq \lambda - 1$ by construction, we conclude $s_1 = s_1^*$. So by the above and our result $r = r_1$ from before the dimension cases for this major case, we now have

$$(z - z_1) + (-1)^r m_3 = \frac{(-1)^r}{2}.$$

However, this is a contradiction as for every value of $m_3 \in \{0, 1\}$, it is the case that $(z - z_1) + (-1)^r m_3 \in \mathbb{Z}$ by the closure of the integers under addition while $\frac{(-1)^r}{2} \notin \mathbb{Z}$ with $r \in \{0, 1\}$. Hence, the above equality is not possible.

Thus, inverses in $C_2$ of $j_1$ moves cannot be performed starting from end critical vertices of the reference $j_1$ move in $C_1$.

**Major Case 2:** We now analyze moves along $j_2$:

<u>Dimension Case 1:</u> Let $d = 2$. Focusing on $j_2$ with the appropriate simplification by construction given $d = 2$, our initial assumption and the observation that during $j_2$ moves $2m_2 + \chi = 1$ give us

$$(t - t_1) + (2(\mu - \mu_1) + (r - r_1) + (2m_2 + \chi) - 1)\lambda = 0$$
$$\implies 2(\mu - \mu_1) + (r - r_1) = \frac{t_1 - t}{\lambda}.$$

Since $2(\mu - \mu_1) + (r - r_1) \in \mathbb{Z}$ and $|t - t_1| \leq \lambda - 1$ by construction, we see $t = t_1$. So our equation of $j_2$ now becomes

$$\mu - \mu_1 = \frac{r_1 - r}{2},$$

which implies $r = r_1$ and $\mu = \mu_1$ as $|r - r_1| \leq 1$ and $\mu - \mu_1 \in \mathbb{Z}$ by construction. Looking at $j_1$, we see by our results above

$$t + (2\nu + r + m_1)\lambda = t_1 + (2\nu_1 + r_1 + 1 - m_1)\lambda$$
$$\implies (\nu - \nu_1) + m_1 = \frac{1}{2}.$$

This is a contradiction as $(\nu - \nu_1) + m_1 \in \mathbb{Z}$ for every value of $m_1 \in \{0, 1\}$ by construction while $\frac{1}{2} \notin \mathbb{Z}$. Hence, the above equality is not possible.

<u>Dimension Case 2:</u> Let $d \geq 3$. Looking at $j_3$, it is the case that

$$s_1 + 2z\lambda + (-1)^r m_3 \lambda = s_1^* + 2z_1 \lambda + (-1)^{r_1}(1 - m_3)\lambda$$
$$\implies 2(z - z_1) + (-1)^r m_3 - (-1)^{r_1}(1 - m_3) = \frac{s_1^* - s_1}{\lambda}.$$

Noting that $2(z - z_1) + (-1)^r m_3 - (-1)^{r_1}(1 - m_3) \in \mathbb{Z}$ and $|s_1 - s_1^*| \leq \lambda - 1$ by construction, we find that $s_1 = s_1^*$. Hence, our equation for $j_3$ becomes

$$(z - z_1) + \frac{(-1)^r + (-1)^{r_1}}{2} m_3 = \frac{(-1)^r}{2}.$$

However, for every value of $r, r_1, m_3 \in \{0, 1\}$ by construction, $(z - z_1) + \frac{(-1)^r + (-1)^{r_1}}{2} m_3 \in \mathbb{Z}$ while $\frac{(-1)^r}{2} \notin \mathbb{Z}$. So we have a contradiction and conclude that the above equality is not possible.

Thus, inverses in $C_2$ of $j_2$ moves cannot be performed starting from end critical vertices of the reference $j_2$ move in $C_1$.

**Major Case 3:** Recall that inverses in $C_2$ of $j_3$ moves in $C_1$ are defined differently taking into account their relationship relative to each other in terms of $r$ and $r_1$. Further, given that our moves in this case are along $j_3$, we have $d \geq 3$. We now case on $|r - r_1|$:

<u>Case 1:</u> Let $|r - r_1| = 0$. Then, $r = r_1$ and our equation for $j_1$ then implies

$$t + (2\nu + r + m_1)\lambda = t_1 + (2\nu_1 + r_1 + |r - r_1|m_1 + (1 - |r - r_1|)(1 - m_1))\lambda$$

$$\implies 2(\nu - \nu_1) + (r - r_1) + (2m_1 - 1)(1 - |r - r_1|) = \frac{t_1 - t}{\lambda}.$$

Since $2(\nu - \nu_1) + (r - r_1) + (2m_1 - 1)(1 - |r - r_1|) \in \mathbb{Z}$ and $|t - t_1| \leq \lambda - 1$ by construction, we must have $t = t_1$. Applying our case assumption along with our result above, $j_1$'s equation implies

$$(\nu - \nu_1) + m_1 = \frac{1}{2},$$

but this is a contradiction as $(\nu - \nu_1) + m_1 \in \mathbb{Z}$ by construction while $\frac{1}{2} \notin \mathbb{Z}$. Hence, the above equality is not possible.

<u>Case 2:</u> Let $|r - r_1| = 1$. Focusing on $j_3$, we see

$$s_1 + 2z\lambda + (-1)^r(m_3\lambda + (m_1 - m_3)\lambda) = s_1^* + 2z_1\lambda + (-1)^{r_1}(|r - r_1|m_3 + (1 - |r - r_1|)(1 - m_3))\lambda$$

$$\implies 2(z - z_1) + (-1)^r m_1 - (-1)^{r_1}(|r - r_1|m_3 + (1 - |r - r_1|)(1 - m_3)) = \frac{s_1^* - s_1}{\lambda}.$$

Observing that $2(z - z_1) + (-1)^r m_1 - (-1)^{r_1}(|r - r_1|m_3 + (1 - |r - r_1|)(1 - m_3)) \in \mathbb{Z}$ and $|s_1 - s_1^*| \leq \lambda - 1$ by construction, we get that $s_1 = s_1^*$. Applying the result above, the case assumption and the observation that $m_1 + m_3 = 1$ during $j_3$ moves by construction, we get

$$z - z_1 = \frac{(-1)^{r_1}(m_1 + m_3)}{2} = \frac{(-1)^{r_1}}{2}.$$

However, this is a contradiction as for every value of $r_1 \in \{0, 1\}$, we see $z - z_1 \in \mathbb{Z}$ by the closure of the integers under addition while $\frac{(-1)^{r_1}}{2} \notin \mathbb{Z}$ for every value $r_1 \in \{0, 1\}$ by construction. Hence, the above equality is not possible.

Thus, inverses in $C_2$ of $j_3$ moves cannot be performed starting from end critical vertices of the reference $j_3$ move in $C_1$.

**Major Case 4:** For moves along $j_k$ for $4 \leq k \leq d$, we see $d \geq 4$ and that in $j_3$ it is the case that

$$s_1 + 2z\lambda + (-1)^r m_3 \lambda = s_1^* + 2z_1 \lambda + (-1)^{r_1}(1-m_3)\lambda$$
$$\implies 2(z-z_1) + (-1)^r m_3 - (-1)^{r_1}(1-m_3) = \frac{s_1^* - s_1}{\lambda}.$$

Given $2(z-z_1) + (-1)^r m_3 - (-1)^{r_1}(1-m_3) \in \mathbb{Z}$ and $|s_1 - s_1^*| \leq \lambda - 1$ by construction, we have that $s_1 = s_1^*$. So our simplified equation for $j_3$ implies

$$(z-z_1) + \frac{(-1)^r + (-1)^{r_1}}{2} m_3 = \frac{(-1)^{r_1}}{2}.$$

However, this is a contradiction as for every value of $r, r_1, m_3 \in \{0,1\}$, $(z-z_1) + \frac{(-1)^r + (-1)^{r_1}}{2} m_3 \in \mathbb{Z}$ by the closure of the integers under addition while $\frac{(-1)^{r_1}}{2} \notin \mathbb{Z}$. So the above equality is not possible.

Consequently, inverses in $C_2$ of $j_k$ moves, for $4 \leq k \leq d$ when $d \geq 4$, cannot be performed starting from end critical vertices of the reference $j_k$ move in $C_1$.

This concludes the proof of Proposition 10.

∎

Thus, we have the shown that our General Square Decomposition's edge set definition cannot produce a cycle $C_2$ that shares an edge with another by performing the inverse of a given move from the end vertex of the move belonging to $C_1$.

### 6.1.2 Proof of Proposition 11:

Let $2 \leq \lambda \leq \lambda^*$, $d \in \mathbb{Z}^{\geq 2}$, and for reference refer to Theorem 2 for specific bounds of each parameter, though we will bring some of them up in our arguments as necessary.

We would also like to make the reader aware that the form of the inverse move on $C_2$ to a move $(m_1, m_2, \ldots, m_d)$ on $C_1$ is $(1-m_1, 1-m_2, \ldots, 1-m_d)$ in the case of moves along $j_k$ for $1 \leq k \leq d$ with $k \neq 3$. Given the alternating behavior for moves along $j_3$ as they depend on $r$ and $r_1$, the pairity of $\ell$ and $\ell_1$, respectively, it follows that inverse $j_3$ moves on $C_2$ with respect to a $j_3$ move $(m_1, m_2, \ldots, m_d)$ on $C_1$ is

$$(|r-r_1|m_1 + (1-|r-r_1|)(1-m_1), |r-r_1|m_2 + (1-|r-r_1|)(1-m_2), \ldots, |r-r_1|m_d + (1-|r-r_1|)(1-m_d)).$$

Note that we will be assuming for the sake of contradiction that two cycles share a vertex with the move configurations established below, meaning all components of their vertices agree, and we then show that such scenarios do not occur by the definition of the General Square Decomposition. In doing so, we will have shown that these cycles cannot share edges by way of these configurations. We now proceed with the proof of Proposition 11 by considering equalities at each component $j_k$ for $1 \leq k \leq d$ with $C_2$ configured to perform the inverse move to a given move $(m_1, \ldots, m_d)$ of $C_1$ starting from an intermediate vertex belonging to $C_1$ for that move:

**Major Case 1:** Here we analyze moves along $j_1$:

<u>Dimension Case 1:</u> Let $d=2$ and observe that in $j_2$, its equation for this dimension case implies

$$t + (2\mu + r + m_2)\lambda = t_1 + (2\mu_1 + r_1 + (1-m_2))\lambda$$

$$\implies 2(\mu - \mu_1) + (r - r_1) + (2m_2 - 1) = \frac{t_1 - t}{\lambda}.$$

Given $2(\mu - \mu_1) + (r - r_1) + (2m_2 - 1) \in \mathbb{Z}$ and $|t - t_1| \leq \lambda - 1$ by construction, we have $t = t_1$. Now, by the above result, our $j_1$ equation by our initial assumptions is

$$t + (2\nu + r + m_1)\lambda + (m_2 - m_1)x_1 = t_1 + (2\nu_1 + r_1 + (1 - m_1))\lambda$$
$$\implies 2(\nu - \nu_1) + (r - r_1) + (2m_1 - 1) = \frac{(m_1 - m_2)x_1}{\lambda}.$$

Following from $2(\nu - \nu_1) + (r - r_1) + (2m_1 - 1) \in \mathbb{Z}$, $1 \leq x_1 \leq \lambda - 1$ with $\lambda \geq 2$ and $m_2 - m_1 = \pm 1$ given $m_1 + m_2 = 1$ during every $j_1$ move with $m_1, m_2 \in \{0, 1\}$ by construction, we see that $\frac{(m_1 - m_2)x_1}{\lambda} \notin \mathbb{Z}$. Hence, we have a contradiction as the above equality is not possible.

<u>Dimension Case 2:</u> Let $d \geq 3$. Focusing on $j_3$, we see

$$s_1 + 2z\lambda + (-1)^r m_3 \lambda = s_1^* + 2z_1 \lambda + (-1)^{r_1}(1 - m_3)\lambda$$
$$\implies 2(z - z_1) + (-1)^r m_3 - (-1)^{r_1}(1 - m_3) = \frac{s_1^* - s_1}{\lambda}.$$

Noting that $2(z - z_1) + (-1)^r m_3 - (-1)^{r_1}(1 - m_3) \in \mathbb{Z}$ and $|s_1 - s_1^*| \leq \lambda - 1$ by construction, we find that $s_1 = s_1^*$. So our simplified $j_3$ equation is

$$(z - z_1) + \frac{(-1)^r + (-1)^{r_1}}{2} m_3 = \frac{(-1)^{r_1}}{2}.$$

However, this is a contradiction as for every value of $m_3, r, r_1 \in \{0, 1\}$, it is the case that $(z - z_1) + \frac{(-1)^r + (-1)^{r_1}}{2} m_3 \in \mathbb{Z}$ by construction and the closure of the integers under addition while $\frac{(-1)^{r_1}}{2} \notin \mathbb{Z}$. Hence, the above equality is not possible.

Thus, inverses in $C_2$ of $j_1$ moves cannot be performed starting from any of the intermediate vertices of the reference $j_1$ move in $C_1$.

**Major Case 2:** We proceed to analyze $j_2$ moves under the assumptions presented at the start of the proof:

<u>Dimension Case 1:</u> Let $d = 2$ and observe that in $j_1$ it is the case that

$$t + (2\nu + r + m_1)\lambda = t_1 + (2\nu_1 + r_1 + 1 - m_1)\lambda$$
$$\implies 2(\nu - \nu_1) + (r - r_1) + (2m_1 - 1) = \frac{t_1 - t}{\lambda}.$$

Since $2(\nu - \nu_1) + (r - r_1) + (2m_1 - 1) \in \mathbb{Z}$ and $|t - t_1| \leq \lambda - 1$ by construction, we conclude $t = t_1$. Going to our equation in $j_2$ for this dimension case and applying the above result, it follows that

$$t + (2\mu + r + m_2)\lambda + \chi x_2 = t_1 + (2\mu_1 + r_1 + 1 - m_2)\lambda$$
$$\implies 2(\mu - \mu_1) + (r - r_1) + (2m_2 - 1) = -\frac{\chi x_2}{\lambda}.$$

Given $2(\mu - \mu_1) + (r - r_1) + (2m_2 - 1) \in \mathbb{Z}$, $1 \leq x_2 \leq \lambda - 1$ and $\chi = \pm 1$ during $j_2$ moves by construction, we observe that $-\frac{\chi x_2}{\lambda} \notin \mathbb{Z}$, giving us a contradiction as the above equality is not possible.

<u>Dimension Case 2:</u> Let $d \geq 3$. Then, looking at our equation in $j_3$, we see

$$s_1 + 2z\lambda + (-1)^r m_3 \lambda = s_1^* + 2z_1 \lambda + (-1)^{r_1}(1 - m_3)\lambda$$

$$\implies 2(z - z_1) + (-1)^r m_3 - (-1)^{r_1}(1 - m_3) = \frac{s_1^* - s_1}{\lambda}.$$

Since $2(z - z_1) + (-1)^r m_3 - (-1)^{r_1}(1 - m_3) \in \mathbb{Z}$ and $|s_1 - s_1^*| \leq \lambda - 1$ by construction, we get $s_1 = s_1^*$. So our simplified $j_3$ equation is

$$(z - z_1) + \frac{(-1)^r + (-1)^{r_1}}{2} m_3 = \frac{(-1)^{r_1}}{2}.$$

This is a contradiction as for all $r, r_1, m_3 \in \{0, 1\}$, $(z - z_1) + \frac{(-1)^r + (-1)^{r_1}}{2} m_3 \in \mathbb{Z}$ by construction while $\frac{(-1)^{r_1}}{2} \notin \mathbb{Z}$. Hence, the above equality is not possible.

Thus, inverses in $C_2$ of $j_2$ moves cannot be performed starting from any of the intermediate vertices of the reference $j_2$ move in $C_1$.

**Major Case 3:** Let $d \geq 3$ and observe that our equation in $j_1$ gives us

$$t + (2\nu + r + m_1)\lambda = t_1 + (2\nu_1 + r_1 + |r - r_1|m_1 + (1 - |r - r_1|)(1 - m_1))\lambda$$

$$\implies 2(\nu - \nu_1) + (r - r_1) + (2m_1 - 1)(1 - |r - r_1|) = \frac{t_1 - t}{\lambda}$$

Following from $2(\nu - \nu_1) + (r - r_1) + (2m_1 - 1)(1 - |r - r_1|) \in \mathbb{Z}$ and $|t - t_1| \leq \lambda - 1$ by construction, we see that $t = t_1$. We now case on $|r - r_1|$:

<u>Case 1:</u> Let $|r - r_1| = 0$. Then, $r = r_1$ and so we get

$$(\nu - \nu_1) + m_1 = \frac{1}{2}.$$

This is a contradiction as for every value of $m_1 \in \{0, 1\}$ by construction, $(\nu - \nu_1) + m_1 \in \mathbb{Z}$ while $\frac{1}{2} \notin \mathbb{Z}$. Hence, the equality above is not possible.

<u>Case 2:</u> Let $|r - r_1| = 1$. Then,

$$\nu - \nu_1 = \frac{(-1)^r}{2}.$$

This is a contradiction as $\nu - \nu_1 \in \mathbb{Z}$ by construction as the integers are closed under addition while for every $r \in \{0, 1\}$ by construction, $\frac{(-1)^r}{2} \notin \mathbb{Z}$. So the equality above is not possible.

Thus, inverses in $C_2$ of $j_3$ moves cannot be performed starting from any of the intermediate vertices of the reference $j_3$ move in $C_1$.

**Major Case 4:** For moves along $j_k$ for $4 \leq k \leq d$, we see $d \geq 4$ and that in $j_3$ it is the case that

$$s_1 + 2z\lambda + (-1)^r m_3 \lambda = s_1^* + 2z_1 \lambda + (-1)^{r_1}(1 - m_3)\lambda$$

$$\implies 2(z - z_1) + (-1)^r m_3 - (-1)^{r_1}(1 - m_3) = \frac{s_1^* - s_1}{\lambda}.$$

Since $2(z - z_1) + (-1)^r m_3 - (-1)^{r_1}(1 - m_3) \in \mathbb{Z}$ and $|s_1 - s_1^*| \leq \lambda - 1$ by construction, we deduce that $s_1 = s_1^*$. Hence, our equation for $j_3$ implies

$$(z - z_1) + \frac{(-1)^r + (-1)^{r_1}}{2} m_3 = \frac{(-1)^{r_1}}{2}.$$

This is a contradiction as for every value of $m_3, r, r_1 \in \{0, 1\}$ by construction, we get $(z - z_1) + \frac{(-1)^r + (-1)^{r_1}}{2} m_3 \in \mathbb{Z}$ by construction while $\frac{(-1)^{r_1}}{2} \notin \mathbb{Z}$.

Consequently, inverses in $C_2$ of $j_k$ moves, for $4 \leq k \leq d$ when $d \geq 4$, cannot be performed starting from any of the intermediate vertices of the reference $j_k$ move in $C_1$.

This concludes the proof of Proposition 11.

∎

**Remark** In the following proofs, note that $r$ and $r_1$ are the pairity of $\ell$ and $\ell_1$, respectively with $\ell$ and $\ell_1$ serving as parameters in the identity of $C_1$ and $C_2$ as defined prior to the statements of the propositions we have been proving. Now observe that given $\mu, \nu, r$, and $t$ agree between $C_1$ and $C_2$, it must be the case that $\ell = \ell_1$ as every other parameter $(p_1, z, s_1, \ldots, s_{d-2}$ as applicable) agreeing indicates where $C_1$ and $C_2$ are situated in the torus with respect to each other, and for every $(\nu, \mu)$ pair in a translation set with all else agreeing, there is only one even numbered cycle and only one odd numbered cycle by construction. So in particular, if we have $\nu, \mu, t$, and all other applicable (dimension-dependent) parameters $p_1, z, s_1, \ldots, s_{d-2}$ agreeing, then having $r = r_1$ implies that $\ell = \ell_1$.

More explicitly, knowing that $\mu = \mu_1$ on one hand gives us

$$\left\lfloor \frac{\ell}{y_1} \right\rfloor = \left\lfloor \frac{\ell_1}{y_1} \right\rfloor$$
$$\implies \left\lfloor \frac{\ell}{y_1} \right\rfloor - \left\lfloor \frac{\ell_1}{y_1} \right\rfloor = 0.$$

From this, we see that it must be the case $|\ell - \ell_1| < y_1$. On the other hand, $\nu = \nu_1$ and $r = r_1$ along with $\ell = 2\lfloor \frac{\ell}{2} \rfloor + r$ and $\ell_1 = 2\lfloor \frac{\ell_1}{2} \rfloor + r_1$ imply

$$\nu - \nu_1 = 0$$
$$\implies \left\lfloor \frac{\ell_1 \lambda}{y_1} \right\rfloor - \left\lfloor \frac{\ell \lambda}{y_1} \right\rfloor = \frac{2\lambda}{y_1} \left( \left\lfloor \frac{\ell_1}{2} \right\rfloor - \left\lfloor \frac{\ell}{2} \right\rfloor \right)$$
$$\implies \left\lfloor \frac{\ell_1 \lambda}{y_1} \right\rfloor - \left\lfloor \frac{\ell \lambda}{y_1} \right\rfloor = -\frac{\lambda}{y_1}(\ell - \ell_1).$$

Now, since $\lfloor \frac{\ell_1 \lambda}{y_1} \rfloor - \lfloor \frac{\ell \lambda}{y_1} \rfloor \in \mathbb{Z}$ by definition, $0 < \lambda < y_1$ by construction and $|\ell - \ell_1| < y_1$ by our previous argument, it must be the case that $\ell = \ell_1$. Otherwise $-\frac{\lambda}{y_1}(\ell - \ell_1) \notin \mathbb{Z}$, in which case we would get a contradiction as an integer cannot be a non-integral rational number.

### 6.1.3  Proof of Proposition 12:

Let $1 \leq \lambda \leq \lambda^*$, $d \in \mathbb{Z}^{\geq 2}$, and for reference refer to Theorem 2 for specific bounds of each parameter, though we will bring some of them up in our arguments as necessary.

Note that we will be assuming that $C_2$ and $C_1$ share a vertex with the configuration established below, meaning all components of the given vertex will agree while potentially differing in their parameter values. We will show that our edge set definition for the General Square Decomposition only allows this if two cycles begin the same move $(m_1, \ldots, m_d)$ from the same critical start vertex, which will then mean that the two cycles are indeed the same cycle as we will show. We now proceed with the proof of Proposition 12 by considering equalities at each component $j_k$ for $1 \leq k \leq d$ with $C_2$ configured to perform the same move as a given move $(m_1, \ldots, m_d)$ of $C_1$ starting from either the critical start vertex or an intermediate vertex.

**Major Case 1:** We now analyze moves along $j_1$:

Dimension Case 1: Let $d = 2$. Focusing on $j_2$, we see

$$t + (2\mu + r + m_2)\lambda = t_1 + (2\mu_1 + r_1 + m_2)\lambda$$
$$\implies 2(\mu - \mu_1) + (r - r_1) = \frac{t_1 - t}{\lambda}.$$

Noting that $2(\mu - \mu_1) + (r - r_1) \in \mathbb{Z}$ and $|t - t_1| \leq \lambda - 1$ by construction, we see $t = t_1$. So our equation in $j_2$ implies

$$\mu - \mu_1 = \frac{r_1 - r}{2}.$$

Given $\mu - \mu_1 \in \mathbb{Z}$ and $|r - r_1| \leq 1$ by construction, we get $r = r_1$ and $\mu = \mu_1$. Going to our equation in $j_1$, we see

$$t + (2\nu + r + m_1)\lambda + (m_2 - m_1)x_1 = t_1 + (2\nu_1 + r_1 + m_1)\lambda$$
$$\implies 2(\nu - \nu_1) = \frac{(m_1 - m_2)x_1}{\lambda}.$$

Since $2(\nu - \nu_1) \in \mathbb{Z}$, $0 \leq x_1 \leq \lambda - 1$ and $m_2 - m_1 = \pm 1$ during moves along $j_1$ by construction, we must have $x_1 = 0$. Note that our $j_1$ move is arbitrary, which is why $m_2 - m_1$ is not fixed. With the above, our $j_1$ equation tells us $\nu = \nu_1$.

By our remark prior to the start of the proofs of Proposition 12, we get our results imply $\ell = \ell_1$. Hence, all identifying values $\ell, t$ agree, so $C_1 = C_{\ell, t} = C_{\ell_1, t_1} = C_2$ in this case.

Dimension Case 2: Let $d = 3$. Focusing on $j_3$, we have

$$s_1 + 2z\lambda + (-1)^r m_3 \lambda = s_1^* + 2z_1 \lambda + (-1)^{r_1} m_3 \lambda$$
$$\implies 2(z - z_1) + ((-1)^r - (-1)^{r_1})m_3 = \frac{s_1^* - s_1}{\lambda}.$$

Following from $2(z - z_1) + ((-1)^r - (-1)^{r_1})m_3 \in \mathbb{Z}$ and $|s_1 - s_1^*| \leq \lambda - 1$ by construction, we see that $s_1 = s_1^*$. This implies $p_2 = p_2^*$ by construction.

Then, in $j_2$ we see by our result $p_2 = p_2^*$,

$$t + (2\mu + r + m_2 + p_1)\lambda + p_2 = t_1 + (2\mu_1 + r_1 + m_2 + p_1^*)\lambda + p_2^*$$
$$\implies 2(\mu - \mu_1) + (r - r_1) + (p_1 - p_1^*) = \frac{t_1 - t}{\lambda}.$$

Noting that $2(\mu - \mu_1) + (r - r_1) + (p_1 - p_1^*) \in \mathbb{Z}$ and $|t - t_1| \leq \lambda - 1$ by construction, we have $t = t_1$. Going to $j_1$ with our result $t = t_1$ yields

$$t + (2\nu + r + m_1)\lambda + (m_2 - m_1)x_1 = t_1 + (2\nu_1 + r_1 + m_1)\lambda$$
$$\implies 2(\nu - \nu_1) + (r - r_1) = \frac{(m_1 - m_2)x_1}{\lambda}.$$

From the above, we get that $x_1 = 0$ as $2(\nu - \nu_1) + (r - r_1) \in \mathbb{Z}$, $0 \leq x_1 \leq \lambda - 1$ and $m_2 - m_1 = \pm 1$ during $j_1$ moves. Note that our $j_1$ move is arbitrary, which is why $m_2 - m_1$ is not fixed. So our new $j_1$ equation gives us

$$\nu - \nu_1 = \frac{r_1 - r}{2}.$$

Then, since $\nu - \nu_1 \in \mathbb{Z}$ and $|r - r_1| \leq 1$ by construction, it follows that $r = r_1$ and hence $\nu = \nu_1$. With the result $r = r_1$ applied to our equation in $j_3$, we see

$$z - z_1 = \frac{((-1)^{r_1} - (-1)^r)m_3}{2}$$
$$\implies z = z_1.$$

Hence, in $j_2$, our results imply

$$t + (2\mu + r + m_2 + p_1)\lambda + p_2 = t_1 + (2\mu_1 + r_1 + m_2 + p_1^*)\lambda + p_2^*$$
$$\implies \mu - \mu_1 = \frac{p_1^* - p_1}{2}.$$

Given $\mu - \mu_1 \in \mathbb{Z}$ and $|p_1 - p_1^*| \leq 1$ by construction, it follows that $p_1 = p_1^*$ and $\mu = \mu_1$. By our remark prior to the start of the proofs of Proposition 12, we get our results imply $\ell = \ell_1$. Hence, all identifying values $\ell, t, p_1, z, s_1$ agree, so $C_1 = C_{\ell, t, p_1, z, s_1} = C_{\ell_1, t_1, p_1^*, z_1, s_1^*} = C_2$ in this case.

<u>Dimension Case 3:</u> Let $d \geq 4$. Here we analyze $j_k$ for $4 \leq k \leq d$.

Looking at all $j_k$, we see

$$s_{k-2} + m_k\lambda = s_{k-2}^* + m_k\lambda$$
$$\implies s_{k-2} = s_{k-2}^*$$

for all $4 \leq k \leq d$. This then implies $p_{k-1} = p_{k-1}^*$ by construction. Now inspecting $j_3$, it follows that

$$s_1 + 2z\lambda + (-1)^r m_3\lambda = s_1^* + 2z_1\lambda + (-1)^{r_1} m_3\lambda$$
$$\implies 2(z - z_1) + ((-1)^r - (-1)^{r_1})m_3 = \frac{s_1^* - s_1}{\lambda}.$$

Consequently, $s_1 = s_1^*$ as $2(z - z_1) + ((-1)^r - (-1)^{r_1})m_3 \in \mathbb{Z}$ and $|s_1 - s_1^*| \leq \lambda - 1$ by construction. So the result above gives us that $p_2 = p_2^*$ by construction. Going to $j_2$, our results $p_2 = p_2^*$ and $p_{k-1} = p_{k-1}^*$ for all $4 \leq k \leq d$ give us

$$t + (2\mu + r + m_2 + p_1)\lambda + p_2 + \sum_{k=3}^{d-1} p_k = t_1 + (2\mu_1 + r_1 + m_2 + p_1^*)\lambda + p_2^* + \sum_{k=3}^{d-1} p_k^*$$

$$\implies 2(\mu - \mu_1) + (r - r_1) + (p_1 - p_1^*) = \frac{t_1 - t}{\lambda}.$$

Given $2(\mu - \mu_1) + (r - r_1) + (p_1 - p_1^*) \in \mathbb{Z}$ and $|t - t_1| \leq \lambda - 1$ by construction, it follows that $t = t_1$. Hence, since $t = t_1$, in $j_1$ we now see

$$t + (2\nu + r + m_1)\lambda + (m_2 - m_1)x_1 = t_1 + (2\nu_1 + r_1 + m_1)\lambda$$
$$\implies 2(\nu - \nu_1) + (r - r_1) = \frac{(m_1 - m_2)x_1}{\lambda}.$$

Following from $2(\nu - \nu_1) + (r - r_1) \in \mathbb{Z}$, $0 \leq x_1 \leq \lambda - 1$ and $m_2 - m_1 = \pm 1$ during $j_1$ moves, we have that $x_1 = 0$. So our $j_1$ equation implies

$$\nu - \nu_1 = \frac{r_1 - r}{2}.$$

Since $\nu - \nu_1 \in \mathbb{Z}$ and $|r - r_1| \leq 1$ by construction, it is the case that $r = r_1$ and $\nu = \nu_1$. From this, in $j_3$ we see

$$z - z_1 = \frac{((-1)^{r_1} - (-1)^r)m_3}{2}$$
$$\implies z = z_1.$$

Lastly, in $j_2$, the pertinent results from above imply

$$\mu - \mu_1 = \frac{p_1^* - p_1}{2}.$$

Given $\mu - \mu_1 \in \mathbb{Z}$ and $|p_1 - p_1^*| \leq 1$ by construction, it follows that $p_1 = p_1^*$ and so $\mu = \mu_1$.

From our results and the remark made prior to the start of the proofs of Proposition 12, it now follows that $\ell = \ell_1$. Given all identifying values $\ell, t, p_1, z, s_1, s_2, \ldots, s_{d-2}$ agree, it is the case that $C_1 = C_{\ell, t, p_1, z, s_1, \ldots, s_{d-2}} = C_{\ell_1, t_1, p_1^*, z_1, s_1^*, \ldots, s_{d-2}^*} = C_2$.

Thus, edges resulting from moves along $j_1$ as defined in the General Square Decomposition can only be shared between two cycles $C_1$ and $C_2$ if and only $C_1 = C_2$.

**Major Case 2:** Focusing on moves along $j_2$, we see in $j_1$

$$t + (2\nu + r + m_1)\lambda = t_1 + (2\nu_1 + r_1 + m_1)\lambda$$
$$\implies 2(\nu - \nu_1) + (r - r_1) = \frac{t_1 - t}{\lambda}.$$

From the above, we conclude $t = t_1$ as $2(\nu - \nu_1) + (r - r_1) \in \mathbb{Z}$ and $|t - t_1| \leq \lambda - 1$ by construction. So our equation in $j_1$ now implies

$$\nu - \nu_1 = \frac{r_1 - r}{2}.$$

Following from $\nu - \nu_1 \in \mathbb{Z}$ and $|r - r_1| \leq 1$ by construction, we get $r = r_1$ and so $\nu = \nu_1$. With this, we now proceed with our dimension cases:

<u>Dimension Case 1:</u> Let $d = 2$. Then, by our dimension assumption and our results in this major case prior to the dimension cases, we get in $j_2$

$$t + (2\mu + r + m_2)\lambda + \chi x_2 = t_1 + (2\mu_1 + r_1 + m_2)\lambda$$
$$\implies 2(\mu - \mu_1) = -\frac{\chi x_2}{\lambda}.$$

Then, since $2(\mu - \mu_1) \in \mathbb{Z}$, $\chi = \pm 1$ during $j_2$ moves and $0 \leq x_2 \leq \lambda - 1$ by construction, we have $x_2 = 0$ and hence $\mu = \mu_1$.

From our results and the remark made prior to the start of the proof of Proposition 12, it now follows that $\ell = \ell_1$. Since all identifying values $\ell, t$ agree, we get $C_1 = C_{\ell,t} = C_{\ell_1,t_1} = C_2$.

<u>Dimension Case 2:</u> Let $d = 3$. Focusing on $j_3$, we have

$$s_1 + 2z\lambda + (-1)^r m_3 \lambda = s_1^* + 2z_1\lambda + (-1)^{r_1} m_3 \lambda$$
$$\implies 2(z - z_1) + ((-1)^r - (-1)^{r_1})m_3 = \frac{s_1^* - s_1}{\lambda}.$$

Since $2(z - z_1) + ((-1)^r - (-1)^{r_1})m_3 \in \mathbb{Z}$ and $|s_1 - s_1^*| \leq \lambda - 1$ by construction, we get $s_1 = s_1^*$. This then implies $p_2 = p_2^*$ by construction. Now, our result $r = r_1$ from the start of this major case along with the above implies

$$z - z_1 = \frac{((-1)^{r_1} - (-1)^r)m_3}{2}$$
$$\implies z = z_1.$$

Going to $j_2$, we see by our results thus far

$$t + (2\mu + r + m_2 + p_1)\lambda + p_2 + \chi x_2 = t_1 + (2\mu_1 + r_1 + m_2 + p_1^*)\lambda + p_2^*$$
$$\implies 2(\mu - \mu_1) + (p_1 - p_1^*) = -\frac{\chi x_2}{\lambda}.$$

Now, observing that $2(\mu - \mu_1) + (p_1 - p_1^*) \in \mathbb{Z}$, $0 \leq x_2 \leq \lambda - 1$ and $\chi = \pm 1$ during $j_2$ moves by construction, we find $x_2 = 0$. So our equation in $j_2$ now implies

$$\mu - \mu_1 = \frac{p_1^* - p_1}{2}.$$

Then, since $\mu - \mu_1 \in \mathbb{Z}$ and $|p_1 - p_1^*| \leq 1$ by construction, it follows that $p_1 = p_1^*$ and so $\mu = \mu_1$.

From our results and the remark made prior to the start of the proof of Proposition 12, it now follows that $\ell = \ell_1$. Since all identifying values $\ell, t, p_1, z, s_1$ agree, we get $C_1 = C_{\ell,t,p_1,z,s_1} = C_{\ell_1,t_1,p_1^*,z_1,s_1^*} = C_2$.

<u>Dimension Case 3:</u> Let $d \geq 4$. Looking at $j_k$ for $4 \leq k \leq d$, it immediately follows that for all $4 \leq k \leq d$,

$$s_{k-2} + m_k \lambda = s_{k-2}^* + m_k \lambda$$
$$\implies s_{k-2} = s_{k-2}^*.$$

Then, by the above, we get $p_{k-1} = p^*_{k-1}$ for all $4 \leq k \leq d$ by construction.

Given $j_3$'s definition does not change when the dimension exceeds $d = 3$, we still have $s_1 = s^*_1$, $p_2 = p^*_2$, and $z = z_1$ by our argument in the previous dimension case. Note that the argument there relies strictly on the results obtained prior to the dimension cases in this major case.

Going to $j_2$ and applying the results $t = t_1$, $p_2 = p^*_2$ and $p_{k-1} = p^*_{k-1}$ for all $4 \leq k \leq d$, we get

$$t + (2\mu + r + m_2 + p_1)\lambda + p_2 + \chi x_2 + \sum_{k=3}^{d-1} p_k = t_1 + (2\mu_1 + r_1 + m_2 + p^*_1)\lambda + p^*_2 + \sum_{k=3}^{d-1} p^*_k$$

$$\implies 2(\mu - \mu_1) + (p_1 - p^*_1) = -\frac{\chi x_2}{\lambda}.$$

So we have $x_2 = 0$ as $2(\mu - \mu_1) + (p_1 - p^*_1) \in \mathbb{Z}$, $0 \leq x_2 \leq \lambda - 1$ and $\chi = \pm 1$ during $j_2$ moves by construction. Consequently, our equation in $j_2$ now implies

$$\mu - \mu_1 = \frac{p^*_1 - p_1}{2}.$$

Given $\mu - \mu_1 \in \mathbb{Z}$ and $|p_1 - p^*_1| \leq 1$ by construction, it is then the case that $p_1 = p^*_1$ and so $\mu = \mu_1$.

From our results and the remark made prior to the start of the proof of Proposition 12, it now follows that $\ell = \ell_1$. Since all identifying values $\ell, t, p_1, z, s_1, \ldots, s_{d-2}$ agree, we get $C_1 = C_{\ell, t, p_1, z, s_1, \ldots, s_{d-2}} = C_{\ell_1, t_1, p^*_1, z_1, s^*_1, \ldots, s^*_{d-2}} = C_2$.

Thus, edges resulting from moves along $j_2$ as defined in the General Square Decomposition can only be shared between two cycles $C_1$ and $C_2$ if and only if $C_1 = C_2$.

**Major Case 3:** Let $d \geq 3$. Here we analyze moves along $j_3$:

Looking at $j_1$, we have

$$t + (2\nu + r + m_1)\lambda = t_1 + (2\nu_1 + r_1 + m_1)\lambda$$

$$\implies 2(\nu - \nu_1) + (r - r_1) = \frac{t_1 - t}{\lambda}.$$

The above implies $t = t_1$ as $2(\nu - \nu_1) + (r - r_1) \in \mathbb{Z}$ and $|t - t_1| \leq \lambda - 1$ by construction. So our equation in $j_1$ now implies

$$\nu - \nu_1 = \frac{r_1 - r}{2}.$$

Since $\nu - \nu_1 \in \mathbb{Z}$ and $|r - r_1| \leq 1$ by construction, we get $r = r_1$ and so $\nu = \nu_1$. We now proceed with our dimension cases:

Dimension Case 1: Let $d = 3$. Looking at $j_2$, our results above imply

$$t + (2\mu + r + m_2 + p_1)\lambda + p_2 = t_1 + (2\mu_1 + r_1 + m_2 + p^*_1)\lambda + p^*_2$$

$$\implies 2(\mu - \mu_1) + (p_1 - p^*_1) = \frac{p^*_2 - p_2}{\lambda}.$$

Hence, $p_2 = p_2^*$ as $2(\mu - \mu_1) + (p_1 - p_1^*) \in \mathbb{Z}$ and $|p_2 - p_2^*| \leq \lambda - 1$ by construction. Then, since $0 \leq s_1, s_1^* \leq \lambda - 1$ by construction, we see $p_2$ and $p_2^*$ are injective over such $s_1$ and $s_1^*$, respectively by their definitions as

$$p_2 = p_2(s_1) = s_1 - \left\lfloor \frac{s_1}{\lambda} \right\rfloor \lambda \qquad \text{and} \qquad p_2^* = p_2^*(s_1^*) = s_1^* - \left\lfloor \frac{s_1^*}{\lambda} \right\rfloor \lambda.$$

Hence, $p_2 = p_2^*$ with $p_2 = p_2(s_1) = s_1$ and $p_2^* = p_2^*(s_1^*) = s_1^*$ implies $s_1 = s_1^*$. With this, in $j_3$ we have

$$s_1 + 2z\lambda + (-1)^r(m_3\lambda + (m_1 - m_3)x_3) = s_1^* + 2z_1\lambda + (-1)^{r_1} m_3 \lambda$$

$$\implies 2(z - z_1) + ((-1)^r - (-1)^{r_1})m_3 = \frac{(-1)^r(m_3 - m_1)x_3}{\lambda}.$$

Since $2(z - z_1) + ((-1)^r - (-1)^{r_1})m_3 \in \mathbb{Z}$, $m_1 - m_3 = \pm 1$ during $j_3$ moves and $0 \leq x_3 \leq \lambda - 1$ all by construction, we see it must be the case $x_3 = 0$. So our result $r = r_1$ from before the dimension cases in this major case along with the above yields

$$z - z_1 = \frac{((-1)^{r_1} - (-1)^r)m_3}{2}$$

$$\implies z = z_1.$$

Going to $j_2$ and applying the results $t = t_1$, $r = r_1$, and $p_2 = p_2^*$, we get

$$t + (2\mu + r + m_2 + p_1)\lambda + p_2 = t_1 + (2\mu_1 + r_1 + m_2 + p_1^*)\lambda + p_2^*$$

$$\implies \mu - \mu_1 = \frac{p_1^* - p_1}{2}.$$

Then, since $\mu - \mu_1 \in \mathbb{Z}$ and $|p_1 - p_1^*| \leq 1$ by construction, it follows that $p_1 = p_1^*$ and hence $\mu = \mu_1$.

From our results and the remark made prior to the start of the proof of Proposition 12, it now follows that $\ell = \ell_1$. Since all identifying values $\ell, t, p_1, z, s_1$ agree, we get $C_1 = C_{\ell, t, p_1, z, s_1} = C_{\ell_1, t_1, p_1^*, z_1, s_1^*} = C_2$.

<u>Dimension Case 2:</u> Let $d \geq 4$. Analyzing the components $j_k$ for $4 \leq k \leq d$, it immediately follows that

$$s_{k-2} + m_k \lambda = s_{k-2}^* + m_k \lambda$$

$$\implies s_{k-2} = s_{k-2}^*$$

for all $4 \leq k \leq d$. From here, we see $p_{k-1} = p_{k-1}^*$ for all $4 \leq k \leq d$ by construction.

Now, given $j_3$'s definition does not change due to the dimension exceeding three, we still have $s_1 = s_1^*$ and $z = z_1$ as shown in the dimension case $d = 3$ before this one. This is a consequence of the argument in the previous dimension case for $j_3$ relying strictly on the results obtained prior to the dimensions cases of this major case as the remaining pertinent results needed to begin from the same point as there are what we shown above. Note that $s_1 = s_1^*$ implies $p_2 = p_2^*$ by construction. Applying to $j_2$ the corresponding results from $j_3$ and our results $t = t_1$ and $r = r_1$ from before the dimension cases of this major case, we obtain

$$t + (2\mu + r + m_2 + p_1)\lambda + p_2 + \sum_{k=3}^{d-1} p_k = t_1 + (2\mu_1 + r_1 + m_2 + p_1^*)\lambda + p_2^* + \sum_{k=3}^{d-1} p_k^*$$

$$\implies \mu - \mu_1 = \frac{p_1^* - p_1}{2}.$$

Then, since $\mu - \mu_1 \in \mathbb{Z}$ and $|p_1 - p_1^*| \leq 1$ by construction, it is then the case that $p_1 = p_1^*$ and so $\mu = \mu_1$.

From our results and the remark made prior to the start of the proof of Proposition 12, it now follows that $\ell = \ell_1$. Since all identifying values $\ell, t, p_1, z, s_1, \ldots, s_{d-2}$ agree, we get $C_1 = C_{\ell, t, p_1, z, s_1, \ldots, s_{d-2}} = C_{\ell_1, t_1, p_1^*, z_1, s_1^*, \ldots, s_{d-2}^*} = C_2$.

Thus, edges resulting from moves along $j_3$ as defined in the General Square Decomposition can only be shared between two cycles $C_1$ and $C_2$ if and only if $C_1 = C_2$.

**Major Case 4:** Let $d \geq 4$. Focusing on moves along $j_k$ for $4 \leq k \leq d$, observe that for a move along $j_k$,

$$s_{k-2} + m_k \lambda + (m_{k-1} - m_k)x_k = s_{k-2}^* + m_k \lambda$$

$$\implies (s_{k-2} - s_{k-2}^*) + (m_{k-1} - m_k)x_k = 0.$$

Now, since $x_k$ needs to be such that the above equality holds for every $j_k$ move, which by construction has $(m_{k-1}, m_k) \in \{(0,1), (1,0)\}$, as our $j_k$ move is arbitrary, and we know $s_{k-2}$ and $s_{k-2}^*$ are fixed as they are identifying parameters of our cycles $C_1$ and $C_2$, we must have $x_k = 0$.

To see this, our observations about $s_{k-2}, s_{k-2}^*$, and $m_{k-1} - m_k$ tell us $m_{k-1} - m_k = \pm 1$ during $j_k$ moves and so we have the following system

$$\begin{cases} s_{k-2} - s_{k-2}^* = x_k \\ s_{k-2} - s_{k-2}^* = -x_k \end{cases}.$$

From this, we see that $x_k = 0$ is the only appropriate choice by the above and $0 \leq x_k \leq \lambda - 1$ by construction. So we get that $s_{k-2} = s_{k-2}^*$, which then implies $p_{k-1} = p_{k-1}^*$.

Looking at $j_1$, we get

$$t + (2\nu + r + m_1)\lambda = t_1 + (2\nu_1 + r_1 + m_1)\lambda$$

$$\implies 2(\nu - \nu_1) + (r - r_1) = \frac{t_1 - t}{\lambda}.$$

Given $2(\nu - \nu_1) + (r - r_1) \in \mathbb{Z}$ and $|t - t_1| \leq \lambda - 1$ by construction, we see $t = t_1$. Our equation in $j_1$ hence implies

$$\nu - \nu_1 = \frac{r_1 - r}{2}.$$

Given $\nu - \nu_1 \in \mathbb{Z}$ and $|r - r_1| \leq 1$ by construction, it is the case that $r = r_1$ and $\nu = \nu_1$.

Going to $j_3$, we have

$$s_1 + 2z\lambda + (-1)^r m_3 \lambda = s_1^* + 2z_1\lambda + (-1)^{r_1} m_3 \lambda$$
$$\implies 2(z - z_1) + ((-1)^r - (-1)^{r_1})m_3 = \frac{s_1^* - s_1}{\lambda}.$$

Since $2(z - z_1) + ((-1)^r - (-1)^{r_1})m_3 \in \mathbb{Z}$ and $|s_1 - s_1^*| \leq \lambda - 1$ by construction, it must be the case that $s_1 = s_1^*$, which implies $p_2 = p_2^*$ by construction. Then, by our results $s_1 = s_1^*$ and $r = r_1$, we get

$$s_1 + 2z\lambda + (-1)^r m_3 \lambda = s_1^* + 2z_1\lambda + (-1)^{r_1} m_3 \lambda$$
$$\implies z - z_1 = \frac{((-1)^{r_1} - (-1)^r)m_3}{2}$$
$$\implies z = z_1.$$

Before concluding, note that in the case $d > 4$, in the components $j_w$ for $4 \leq w \leq d$ with $w \neq k$, we have

$$s_{w-2} + m_w \lambda = s_{w-2}^* + m_w \lambda$$
$$\implies s_{w-2} = s_{w-2}^*$$

for all $4 \leq w \leq d$ with $w \neq k$. This then implies $p_{w-1} = p_{w-1}^*$ for all $4 \leq w \leq d$ with $w \neq k$.

Lastly, combining and applying all of our observations and results on $j_2$ yields

$$t + (2\mu + r + m_2 + p_1)\lambda + p_2 + \sum_{k=3}^{d-1} p_k = t_1 + (2\mu_1 + r_1 + m_2 + p_1^*)\lambda + p_2^* + \sum_{k=3}^{d-1} p_k^*$$
$$\implies \mu - \mu_1 = \frac{p_1^* - p_1}{2}.$$

Given $\mu - \mu_1 \in \mathbb{Z}$ and $|p_1 - p_1^*| \leq 1$ by construction, we have $p_1 = p_1^*$ and so $\mu = \mu_1$.

From our results and the remark made prior to the start of the proof of Proposition 12, it now follows that $\ell = \ell_1$. Since all identifying values $\ell, t, p_1, z, s_1, \ldots, s_{d-2}$ agree, we get $C_1 = C_{\ell,t,p_1,z,s_1,\ldots,s_{d-2}} = C_{\ell_1,t_1,p_1^*,z_1,s_1^*,\ldots,s_{d-2}^*} = C_2$.

Thus, edges resulting from moves along $j_k$ for $4 \leq k \leq d$ when $d \geq 4$ as defined in the General Square Decomposition can only be shared between two cycles $C_1$ and $C_2$ if and only if $C_1 = C_2$.

This concludes the proof of Proposition 12.

■

With this, we have shown that two cycles $C_1$ and $C_2$ defined by the General Square Decomposition share an edge if and only if $C_1 = C_2$, precisely Theorem 9 as desired.

## 6.2 Proof of all Cycles being of the Same Length:

Here, we show that all cycles defined by the edge set are of the same length. To do this, it suffices to show that $2d$ moves occur with each occurring exactly once and contributing $\lambda$ edges each time. Note that this is without distinguishing between partitioned and non-partitioned edges.

This will lead to the following result:

**Theorem 13** *Every cycle defined by the General Square Decomposition's edge set yields cycles of length $2d\lambda$, without distinguishing between partitioned and non-partitioned edges along each dimension.*

**Proof:** Let $d \in \mathbb{Z}^{\geq 2}$ and $1 \leq \lambda \leq \lambda^*$. Assume we have chosen the parameters $\ell, t, p_1, z, s_1, \ldots, s_{d-2}$, as applicable based on the dimension $d$ of our torus, for an arbitrary cycle defined by the General Square Decomposition's edge set. Also observe that the sliding parameters that allow for $\lambda$ edges along a given dimension per $(m_1, m_2, m_3, \ldots, m_d)$ move are the $0 \leq x_k \leq \lambda - 1$ for $1 \leq k \leq d$ with coefficients $(m_2 - m_1), \chi, (-1)^r(m_1 - m_3)$, and $(m_{w-1} - m_w)$ for $4 \leq w \leq d$ as applicable based on $d$. Throughout the proof, take note how only one $x_k$ has a non-zero coefficient at a time, meaning only moves along the $k$th dimension are being carried out for such $(m_1, m_2, m_3, \ldots, m_d)$ $d-$tuples. For moves along $j_3$, we will disregard the $(-1)^r$ component as that remains fixed given $r$ is the pairity of $\ell$, and $\ell$ is an identifying parameter of our cycle and hence fixed. We now case on the dimension of our torus with the torus as defined in Theorem 2:

**Dimension Case 1:** Let $d = 2$. Then, the definition of $M = 1$ as in Theorem 2 in terms of the truth value of $q$ tells us that appropriate $(m_1, m_2) \in \{0, 1\}^2$ are such that $(m_1 = m_2) \vee (m_1 \neq m_2)$. Hence, we see $(m_1, m_2) \in \{(0, 0), (0, 1), (1, 1), (1, 0)\}$.

<u>Move Case 1:</u> If our move is a $j_1$ move, observe $(m_1, m_2) \in \{(0, 1), (1, 0)\}$ as then $m_2 - m_1 = \pm 1$ and $\chi = 0$ as defined in Theorem 2. Following from the definition for $j_1$, $\lambda$ edges are contributed per $j_1$ move, so we get $2\lambda$ edges total here.

<u>Move Case 2:</u> If our move is a $j_2$ move, observe $(m_1, m_2) \in \{(0, 0), (1, 1)\}$ as then $\chi = \pm 1$ and $m_2 - m_1 = 0$. Following from the definition of $j_2$, $\lambda$ edges are contributed per $j_2$ move, so we get $2\lambda$ edges total here.

Note that all possible moves were used and each was used exactly once, so each occurs exactly once by the above.

Thus, our cycles are all of length $2\lambda + 2\lambda = 2 \cdot 2\lambda = d \cdot 2\lambda = 2d\lambda$ for this dimension case.

**Dimension Case 2:** Let $d = 3$. Then, the definition of $M = 1$ in terms of the truth value of $q$ tells us that $(m_1, m_2, m_3) \in \{0, 1\}^3$ are such that appropriate $(m_1 = m_2 = m_3) \vee (m_2 \neq m_1 = m_3) \vee (m_1 = m_2 \neq m_3)$. Hence, $(m_1, m_2, m_3) \in \{(0, 0, 0), (0, 1, 0), (1, 1, 0), (1, 1, 1), (1, 0, 1), (0, 0, 1)\}$.

<u>Move Case 1:</u> If our move is a $j_1$ move, observe $(m_1, m_2, m_3) \in \{(0, 1, 0), (1, 0, 1)\}$ as then $m_2 - m_1 = \pm 1, \chi = 0$, and $m_1 - m_3 = 0$. Following from the definition of $j_1$, $\lambda$ edges are contributed per $j_1$ move, so we get $2\lambda$ edges total here.

<u>Move Case 2:</u> If our move is a $j_2$ move, observe $(m_1, m_2, m_3) \in \{(0, 0, 0), (1, 1, 1)\}$ as then $\chi =$

$\pm 1, m_2 - m_1 = 0$, and $m_1 - m_3 = 0$. Following from the definition of $j_2$, $\lambda$ edges are contributed per $j_2$ move, so we get $2\lambda$ edges total here.

<u>Move Case 3:</u> If our move is a $j_3$ move, observe $(m_1, m_2, m_3) \in \{(1,1,0), (0,0,1)\}$ as then $m_1 - m_3 = \pm 1, m_2 - m_1 = 0$, and $\chi = 0$. Following from the definition of $j_3$, $\lambda$ edges are contributed per $j_3$ move, so we get $2\lambda$ total here.

Note that all possible moves were used and each was used exactly once, so each occurs exactly once by the above.

Thus, our cycles are all of length $2\lambda + 2\lambda + 2\lambda = 3 \cdot 2\lambda = d \cdot 2\lambda = 2d\lambda$ for this dimension case.

**Dimension Case 3:** Let $d \geq 4$. Then, the definition of $M = 1$ in terms of the truth value of $q$ tells us appropriate $(m_1, m_2, m_3, \ldots, m_d) \in \{0,1\}^d$ are such that one of the statements from $q(m_1, m_2, m_3, \ldots, m_d)$ below are true

$$(m_1 = m_2 = \cdots = m_d) \vee (m_2 \neq m_1 = m_3 = \cdots = m_d) \vee \bigvee_{k=2}^{d-1} (m_1 = m_2 = \cdots = m_k \neq m_{k+1} = \cdots = m_d).$$

Hence,

$(m_1, m_2, m_3, \ldots, m_d) \in \{(0,0,0,\ldots,0), (0,1,0,\ldots,0), (1,1,0,\ldots,0), (1,1,1,0,\ldots,0), (1,1,1,1,0,\ldots,0), \ldots,$
$(1,1,1,\ldots,1,0), (1,1,1,\ldots,1), (1,0,1,\ldots,1), (0,0,1,\ldots,1), (0,0,0,1,\ldots,1),$
$(0,0,0,0,1,\ldots,1), \ldots, (0,0,0,\ldots,0,1)\}$

Above we convey the general setting of appropriate moves admitted with the first 6 $d$-tuples enumerated above being the starting sequence of moves and the other half being the inverse sequence of moves to those moves to close up our cycle. By the logical statement $q$ from earlier, we see that there are $d$ possible statements each appropriate $(m_1, m_2, m_3, \ldots, m_d)$ can satisfy to be admitted as a move, and there are two moves $(m_1, m_2, m_3, \ldots, m_d), (1-m_1, 1-m_2, 1-m_3, \ldots, 1-m_d)$ for a given statement such that the truth value remains unaffected if either $d$-tuple is being applied. As a consequence, note that each statement corresponds to two $d$-tuples along a given dimension.

Given how each $m_j$ is to relate to another, choosing one $m_j$ dictates what we can choose for the others to satisfy one of the statements, and so we can think of each choice as turning on those necessary (1) and keeping the others off (0), and vice versa. Consequently, there are $2d$ moves with each $(m_1, m_2, m_3, \ldots, m_d)$ in the set above corresponding to one such move.

<u>Move Case 1:</u> If our move is a $j_1$ move, observe $(m_1, m_2, m_3, \ldots, m_d) \in \{(0,1,0,\ldots,0), (1,0,1,\ldots,1)\}$ as then $m_2 - m_1 = \pm 1, \chi = 0, m_1 - m_3 = 0$, and $m_{k-1} - m_k = 0$ for all $4 \leq k \leq d$. Following from the definition of $j_1$, $\lambda$ edges are contributed per $j_1$ move, so we get $2\lambda$ edges total here.

<u>Move Case 2:</u> If our move is a $j_2$ move, observe $(m_1, m_2, m_3, \ldots, m_d) \in \{(0,0,0,\ldots,0), (1,1,1,\ldots,1)\}$ as then $\chi = \pm 1, m_2 - m_1 = 0, m_1 - m_3 = 0$ and $m_{k-1} - m_k = 0$ for all $4 \leq k \leq d$. Following from the definition of $j_2$, $\lambda$ edges are contributed per $j_2$ move, so we get $2\lambda$ edges total here.

<u>Move Case 3:</u> If our move is a $j_3$ move, observe $(m_1, m_2, m_3, \ldots, m_d) \in \{(1,1,0,\ldots,0), (0,0,1,\ldots,1)\}$ as then $m_1 - m_3 = \pm 1, m_2 - m_1 = 0, \chi = 0$, and $m_{k-1} - m_k = 0$ for all $4 \leq k \leq d$. Following from the definition of $j_3$, $\lambda$ edges are contributed per $j_3$ move, so we get $2\lambda$ edges total here.

Move Case 4: If our move is a $j_k$ move for some $4 \leq k \leq d$, observe

$$(m_1, m_2, m_3, \ldots, m_{k-1}, m_k, m_{k+1}, \ldots, m_d) \in \{(1,1,1,\ldots,1,0,0,\ldots,0), (0,0,0,\ldots,0,1,1,\ldots,1)\}$$

as then $m_{k-1} - m_k = \pm 1, m_2 - m_1 = 0, \chi = 0, m_1 - m_3 = 0$, and $m_{w-1} - m_w = 0$ for all $4 \leq w \leq d$ with $w \neq k$ in the case $d > 4$. Following from the definition of $j_k$, $\lambda$ edges are contributed per $j_k$ move, so we get $2\lambda$ edges total here.

Note that all possible moves were used and each was used exactly once, so each occurs exactly once by the above.

Thus, since there are two moves along each dimension and there are $d$−dimensions in our torus, our cycles are all of length $\sum_{k=1}^{d} 2\lambda = d \cdot 2\lambda = 2d\lambda$ for this dimension case.

Given that in all dimension cases every move is used exactly once with each of the $d$ dimensions corresponding to two moves, each contributing $\lambda$ edges, we have that the General Square Decomposition defines cycles of length $2d\lambda$. Note that this is without distinguishing between partitioned and non-partitioned edges.

Now, letting $n \in \mathbb{Z}^+$, $a = 2^{i_1} + \cdots + 2^{i_d}$ with $i_1 > \cdots > i_d = 0$ and $\lambda = 2^{\alpha-1}$ with $\lambda^* = 2^{2n-1}$, we see that $Q_{2an}$ can be decomposed into cycles of length $a \cdot 2^{\alpha}$ for $1 \leq \alpha \leq 2n$ when distinguishing between partitioned and non-partitioned edges per Proposition 1 with $a_j = 2\lambda = 2^{\alpha}$ and $k_j = 2^{i_j}$ for $1 \leq j \leq d$.

With this, we have proven Theorem 13, as required.

∎

# 7 Proofs of Theorem 5 and Corollary 7

We will prove that the proposed General Lock-and-Key Decomposition edge set definition indeed defines edge-disjoint cycles. Following this, we will show that all cycles are of the same length and in particular of length $a \cdot 2^{\alpha}$ in the case of the cycle decompositions of $Q_{2an}$. Combining these results, we will have proven Theorem 5 and Corollary 7.

## 7.1 Edge-Disjoint Cycles Proofs:

Before we begin, we define a few terms that will be used in laying out the propositions to be proven. **Start vertex** refers to the vertex from which a given move $(m_1, \ldots, m_d)$ is to begin. **End vertex** refers to the vertex at which a move $(m_1, \ldots, m_d)$ leaves off once the referenced move has been carried out. **Major sets** refer to the four main sets from which all other edge sets are defined. We will number them $1 - 4$ from top to bottom as presented in Theorem 5 and refer to major sets as simply sets followed by the number corresponding to a given major set.

Our main result for this section will be **Theorem 14** *The General Lock-and-Key Decomposition's Edge set yields cycles that share an edge if and only if the two cycles are the same cycle.*

Now, let $C_1 = C_{\ell,\gamma,t,p_1,s_1,p_2,s_2,\ldots,p_{d-2},s_{d-2}}$ and $C_2 = C_{\ell_1,\gamma_1,t_1,p_1^*,s_1^*,p_2^*,s_2^*,\ldots,p_{d-2}^*,s_{d-2}^*}$ with the parameters ranging over the intervals established in Theorem 5. Note that any other parameters with no subscripts that are fixed for a given cycle, meaning they are identifiers of the cycle in some fashion, will have subscript 1 if they belong to $C_2$ and none if they belong to $C_1$. Further, if any of the parameters have a subscript prior to distinguishing between those belonging to $C_1$ and $C_2$, having no superscript will correspond to those parameters belonging to $C_1$ and a $*$ for a superscript will correspond to those parameters belonging to $C_2$.

Given the versatility of the edge set definition in that it performs a move $(m_1,\ldots,m_d)$ and those with the opposite orientation when required, we must prove the following:

**Proposition 15** *Two cycles $C_1$ and $C_2$ cannot share an edge as a consequence of $C_2$ performing a move with the opposite orientation to that of a given move $(m_1,\ldots,m_d)$ of $C_1$ starting from the end vertex of the edge belonging to $C_1$.*

**Proposition 16** *Two cycles $C_1$ and $C_2$ cannot share an edge as a consequence of $C_2$ performing a given move $(m_1^*,\ldots,m_d^*)$ with the same orientation as that which defines an edge belonging to $C_1$ starting from the same start vertex if the edge is present in two distinct major sets.*

**Proposition 17** *Two cycles $C_1$ and $C_2$ cannot share an edge as a consequence of $C_2$ performing a move with the same orientation as a given move $(m_1,\ldots,m_d)$ performed by $C_1$ from the same start vertex unless $C_2$ performs the same move from the same point $(x_1,x_2,\ldots,x_d)$ corresponding to the edge in $C_1$ and the edge comes from the same major set. Further, this is if and only if $C_1$ and $C_2$ are the same cycle i.e. $(C_1 = C_2)$.*

With this, we will have shown that the cycles resulting from the General Lock-and-Key Decomposition edge set are edge-disjoint, meaning the edges in a cycle's edge set belong uniquely to that cycle.

We will present a proof for each of the propositions. When we address the two sets that will be considered in a set comparison case throughout the proofs, the set number corresponding to the major set being considered for $C_1$ will always be less than or equal to that corresponding to $C_2$. Further, when we say that moves of $C_1$ along $j_k$ for $1 \leq k \leq d$ relate to those of $C_2$ in a given manner, we will mean that moves defined by the set assigned to $C_1$ in that set comparison case will relate to the moves defined by the set assigned to $C_2$ as established.

Note that the expressions defining the vertex pair in the edge set definition are $\pm 1$ during $j_1,\ldots,j_d$ moves respectively, and 0 otherwise by the construction $M$ in Theorem 5. $M$ only admits $d$-tuples $(m_1,\ldots,m_d) \in \{0,1\}^d$ with the properties that make the mappings in the vertex pair the case for a given move.

**Remark** When handling expressions of the form $\lambda_j(1-\lambda_k)$, where $\lambda_j = \lambda_k$ and $\lambda_j, \lambda_k \in \{0,1\}$ for $j,k \in \mathbb{Z}^+$, we see that given the parameters $\lambda_j$ and $\lambda_k$ are only defined at the roots of a polynomial of the same form, $\lambda_j(1-\lambda_k) = 0$. On a similar note, we also make use of the property $\lambda_j - \lambda_k = \text{sgn}(\lambda_j - \lambda_k)|\lambda_j - \lambda_k| = (-1)^{\lambda_k}|\lambda_j - \lambda_k|$ when $\lambda_j, \lambda_k \in \{0,1\}$ and hence the case that $|\lambda_j - \lambda_k| \leq 1$. Lastly, observe that when $\lambda_j \in \{0,1\}$ for $j \in \mathbb{Z}^+$, $\lambda_j^k = \lambda_j$ for every $k \in \mathbb{Z}^+$.

In the following proofs, we show that the General Lock-and-Key decomposition edge set definition holds more generally to decompose any torus defined by the $d-$fold Cartesian product of cycles at hand. When addressing moves with the opposite orientation, we may sometimes refer to them as inverses and it will be established if there is a unique inverse in a given case. Lastly, when referring to $j_k$ moves of a given cycle for $1 \leq k \leq d$, we will mean the $j_k$ moves defined by the set assigned to that cycle for that set comparison case.

We will be using all of the above without further mention in the following proofs.

### 7.1.1 Proof of Proposition 15:

Let $d \in \mathbb{Z}^{\geq 2}, \alpha_1 \leq \alpha^* \leq \alpha_d$, and for reference refer to Theorem 5 for specific bounds of each parameter, though we will bring some of them up in our arguments as necessary.

We would also like to make the reader aware that the form of the inverse move on $C_2$ to a move $(m_1, \ldots, m_d)$ on $C_1$ will vary based on which move and which two major sets are being considered. In some instances, there may be more than one move with the same orientation and so there would not be a unique inverse move in said case. The definition of the inverse move $(m_1^*, \ldots, m_d^*)$ to a given move $(m_1, \ldots, m_d)$ and all applicable equations that hold during every move will be introduced as required.

Note that we will be assuming for the sake of contradiction that two cycles share an edge with the move configurations established below, meaning all components of their vertices agree, and we then show that such scenarios do not occur by the definition of the General Lock-and-Key Decomposition. In doing so, we will have shown that these cycles cannot share edges by way of these configurations. We now proceed with the proof of Proposition 15 by considering equalities at each component $j_k$ for $1 \leq k \leq d$ with $C_2$ configured to perform a move $(m_1^*, \ldots, m_d^*)$ with the opposite orientation to that of a given move $(m_1, \ldots, m_d)$ of $C_1$ starting from the end vertex belonging to the edge in $C_1$:

**Dimension Case 1:** Let $d = 2$. Then, $C_1 = C_{\ell, \gamma}$ and $C_2 = C_{\ell_1, \gamma_1}$. Note that $t = 0 = t_1$ as $d = 2$ and so $\eta_3 = 0$. We now case on $\alpha^*$ for $\alpha_1 \leq \alpha^* \leq \alpha_2$:

**$\alpha^*-$Case 1:** Let $\alpha^* = \alpha_1$. In this case, only sets 3 and 4 are active as $A_1 = 0$ and $A_2 = 0$.

**Set Comparison Case 1:** We will consider whether edges are shared from moves with the opposite orientation via set 3:

Major Case 1: Focusing on moves along $j_1$, we see that an edge can be shared as a consequence of an inverse move if and only if $|\gamma - \gamma_1| = 1$ since $|\gamma - \gamma_1| \leq 1$ by construction and $\gamma \neq \gamma_1$. Then, following from $\alpha^* = \alpha_1$, we get $A_3 = 2$. Observing that during such $j_1$ moves $m_2 = 1 = m_2^*$, in $j_2$ we have

$$m_2 + \gamma + 4x_2 + (1 - m_2)m_1 = m_2^* + \gamma_1 + 4x_2^* + (1 - m_2^*)m_1^*$$
$$\implies x_2 - x_2^* = \frac{(-1)^\gamma}{4}.$$

Given that $x_2 - x_2^* \in \mathbb{Z}$ by construction and for every $\gamma \in \{0, 1\}$, $\frac{(-1)^\gamma}{4} \notin \mathbb{Z}$, we have a contradiction as the above equality is not possible.

<u>Major Case 2</u>: Focusing on moves along $j_2$, we see that no $j_2$ moves with the opposite orientation to a given $j_2$ move in $C_1$ are ever performed and so there is nothing to prove.

**Set Comparison Case 2:** We will consider whether edges are shared from moves with the opposite orientation via sets 3 and 4.

<u>Major Case 1</u>: Focusing on moves along $j_1$, we see that every $j_1$ move in $C_2$ is an inverse to every $j_1$ move in $C_1$ if and only if $\gamma = \gamma_1$. Looking at $j_2$ and applying the observations that $m_2 = 1 = m_2^*$ during all $j_1$ moves along with the above yields

$$m_2 + \gamma + 4x_2 + (1 - m_2)m_1 = 2 + m_1^* + \gamma_1 + m_2^*(1 - m_1^*) + 4x_2^*$$
$$\implies x_2 - x_2^* = \frac{1}{2}.$$

Given $x_2 - x_2^* \in \mathbb{Z}$ by construction and $\frac{1}{2} \notin \mathbb{Z}$, we have a contradiction as the above equality is not possible.

<u>Major Case 2</u>: Focusing on moves along $j_2$, we see that all $j_2$ moves defined by sets 3 and 4 in this case have the same orientation with respect to each other, and so there are no moves with the opposite orientation to consider.

**Set Comparison Case 3:** We will consider whether edges are shared from moves with the opposite orientation via set 4.

<u>Major Case 1</u>: Looking at moves along $j_1$, we see that $C_2$ has an inverse $j_1$ move to $C_1$ if and only if $|\gamma - \gamma_1| = 1$, meaning $\gamma \neq \gamma_1$, and $m_2 = 1 = m_2^*$ for all $j_1$ moves. Applying these observations to $j_2$, we get

$$2 + m_1 + \gamma + m_2(1 - m_1) + 4x_2 = 2 + m_1^* + \gamma_1 + m_2^*(1 - m_1^*) + 4x_2^*$$
$$\implies x_2 - x_2^* = \frac{(-1)^\gamma}{4}.$$

Since $x_2 - x_2^* \in \mathbb{Z}$ by construction and for every $\gamma \in \{0, 1\}$, $\frac{(-1)^\gamma}{4} \notin \mathbb{Z}$, we have a contradiction as the above equality is not possible.

<u>Major Case 2</u>: Focusing on moves along $j_2$, we see that there are no inverse moves to consider as all $j_2$ moves from $C_1$ have the same orientation with respect to those in $C_2$.

Thus, no edges are shared via moves with the opposite orientation for $\alpha^* = \alpha_1$ when $d = 2$.

**$\alpha^*$−Case 2:** Let $\alpha_1 < \alpha^* < \alpha_2$. Observe that here $A_1 = 1$ and $A_2 = 0$, and so sets $1, 2, 3$ and $4$ are all active. Lastly, note that $A_3 \geq 4$ is even.

**Set Comparison Case 1:** We will consider whether edges are shared from moves with the opposite orientation via set 1.

<u>Major Case 1</u>: Focusing on moves along $j_1$, we see that inverse $j_1$ moves between $C_1$ and $C_2$ occur if and only if $|\gamma - \gamma_1| = 1$. Looking at $j_2$ and applying the above observation along with the equations $m_1 + m_2 = 1$ and $m_1^* + m_2^* = 1$ that hold during all $j_1$ moves, we get

$$m_2 + \gamma + 4x_2 = m_2^* + \gamma_1 + 4x_2^*$$
$$\implies x_2 - x_2^* = \frac{(\gamma_1 - \gamma) + (m_2^* - m_2)}{4}.$$

From the above, we see $m_2^* - m_2 = -(\gamma_1 - \gamma)$ as $x_2 - x_2^* \in \mathbb{Z}$, $|\gamma_1 - \gamma| \leq 1$, and $|m_2^* - m_2| \leq 1$ by construction. It then follows that $m_1 - m_1^* = -(\gamma_1 - \gamma)$ and so in $j_1$, our results and assumptions give us

$$(\ell + \gamma)A_3 - \gamma + (-1)^\gamma(m_1 + 2x_1 + 1) = (\ell_1 + \gamma_1)A_3 - \gamma_1 + (-1)^{\gamma_1}(m_1^* + 2x_1^*)$$
$$\implies ((\ell - \ell_1) + (\gamma - \gamma_1))\frac{A_3}{2} - \frac{1 - (-1)^\gamma}{2}(\gamma - \gamma_1) + (-1)^\gamma(m_1^* + x_1 + x_1^*) = \frac{(-1)^{\gamma+1}}{2}.$$

Given $((\ell - \ell_1) + (\gamma - \gamma_1))\frac{A_3}{2} - \frac{1-(-1)^\gamma}{2}(\gamma - \gamma_1) + (-1)^\gamma(m_1^* + x_1 + x_1^*) \in \mathbb{Z}$ and for every $\gamma \in \{0, 1\}$, $\frac{(-1)^{\gamma+1}}{2} \notin \mathbb{Z}$ by construction, we have a contradiction as the above equality is not possible.

Major Case 2: Focusing on moves along $j_2$, we see that during such $j_2$ moves $(m_1^*, m_2^*) = (1-m_1, 1-m_2)$ and $m_1 = m_2$. Applying the above along with the equation $2m_2 + \chi = 1$ that holds during all $j_2$ moves in this case, in $j_2$ they give us

$$m_2 + \gamma + 4x_2 + \chi = m_2^* + \gamma_1 + 4x_2^*$$
$$\implies x_2 - x_2^* = \frac{\gamma_1 - \gamma}{4}.$$

Since $x_2 - x_2^* \in \mathbb{Z}$ and $|\gamma_1 - \gamma| \leq 1$ by construction, $\gamma = \gamma_1$. Going to $j_1$, we see

$$(\ell + \gamma)A_3 - \gamma + (-1)^\gamma(m_1 + 2x_1) = (\ell_1 + \gamma_1)A_3 - \gamma_1 + (-1)^{\gamma_1}(m_1^* + 2x_1^*)$$
$$\implies (\ell - \ell_1)\frac{A_3}{2} + (-1)^\gamma(m_1 + (x_1 - x_1^*)) = \frac{(-1)^\gamma}{2}.$$

Since $(\ell - \ell_1)\frac{A_3}{2} + (-1)^\gamma(m_1 + (x_1 - x_1^*)) \in \mathbb{Z}$ and for every $\gamma \in \{0, 1\}$, $\frac{(-1)^\gamma}{2} \notin \mathbb{Z}$, we have a contradiction as the above equality is not possible.

**Set Comparison Case 2:** We will consider whether edges are shared from moves with the opposite orientation via sets 1 and 2.

Major Case 1: Focusing on moves along $j_1$, we see that inverse $j_1$ moves can only take place in $C_2$ relative to $C_1$ if and only if $\gamma = \gamma_1$ since $|\gamma_1 - \gamma| \leq 1$ and $|\gamma_1 - \gamma| \neq 1$. Now, following from the equation $m_1^* = m_2^*$ that holds during $j_1$ moves, in $j_2$ we have

$$m_2 + \gamma + 4x_2 = 2 + m_1^* + \gamma_1 + 4x_2^*$$
$$\implies x_2 - x_2^* = \frac{2 + (m_2^* - m_2)}{4}.$$

Given $x_2 - x_2^* \in \mathbb{Z}$ and $|m_2^* - m_2| \leq 1$ by construction, for every $m_2, m_2^* \in \{0, 1\}$, it is the case $\frac{2+(m_2^*-m_2)}{4} \notin \mathbb{Z}$, giving us a contradiction as the above equality is not possible.

Major Case 2: Focusing on moves along $j_2$, we see that which $j_2$ moves of $C_2$ serve as moves with the opposite orientation to those of $C_1$ depends on $x_1^*$. Hence, we proceed by casing on $x_1^*$:

<u>Case 1</u>: Let $0 \leq x_1^* < \frac{A_3}{2} - 2$. Then, observing that during such $j_2$ moves the equations $m_1^* = 1 - m_2$ and $2m_2 + \chi = 1$ hold, in $j_2$ we get

$$m_2 + \gamma + 4x_2 + \chi = 2 + m_1^* + \gamma_1 + 4x_2^*$$
$$\implies x_2 - x_2^* = \frac{2 + (\gamma_1 - \gamma)}{4}.$$

Since $x_2 - x_2^* \in \mathbb{Z}$ by construction and for every $\gamma, \gamma_1 \in \{0, 1\}$, $\frac{2+(\gamma_1-\gamma)}{4} \notin \mathbb{Z}$ since $1 \leq 2+(\gamma_1-\gamma) \leq 3$ by construction, we have a contradiction as the above equality is not possible.

<u>Case 2</u>: Let $x_1^* = \frac{A_3}{2} - 2$. Here we see all of $C_2$'s $j_2$ moves serve as moves with the opposite orientation only for $C_1$'s $j_2$ move $(m_1, m_2) = (1, 1)$. Applying this to $j_2$ along with the equation $m_1^* + m_2^* = 1$, we see

$$m_2 + \gamma + 4x_2 - 1 = 2 + m_1^* + \gamma_1 + 4x_2^*$$
$$\implies x_2 - x_2^* = \frac{2 + m_1^* + (\gamma_1 - \gamma)}{4}.$$

Since $x_2 - x_2^* \in \mathbb{Z}$, $|\gamma_1 - \gamma| \leq 1$, and $m_1^* \in \{0, 1\}$ all by construction, we must have $m_1^* = 1$ and $\gamma_1 - \gamma = 1$, meaning $\gamma_1 = 1$ and $\gamma = 0$. Going to $j_1$, we have

$$(\ell + \gamma)A_3 - \gamma + (-1)^\gamma(m_1 + 2x_1) = (\ell_1 + 1 - \gamma_1)A_3 - \gamma_1 + (-1)^{\gamma_1+1}(2 + m_2^* + 2m_1^*(1 - m_2^*) + 2x_1^*)$$
$$\implies \ell - \ell_1 - 1 = -\frac{2(1 + x_1)}{A_3}.$$

Given $\ell - \ell_1 - 1 \in \mathbb{Z}$ and $2 \leq 2(1 + x_1) \leq A_3 - 2$ by construction, we get $-\frac{2(1+x_1)}{A_3} \notin \mathbb{Z}$. Hence, we have a contradiction as the above equality is not possible.

**Set Comparison Case 3:** We will consider whether edges are shared via sets 1 and 3 through moves with the opposite orientation.

<u>Major Case 1</u>: Focusing on moves along $j_1$, we see that inverse $j_1$ moves between $C_1$ and $C_2$ exist if and only if $|\gamma_1 - \gamma| = 1$. Observing that $m_1 + m_2 = 1$ and $m_2^* = 1$ during such moves, in $j_1$ we have

$$(\ell + \gamma)A_3 - \gamma + (-1)^\gamma(m_1 + 2x_1 + 1) = (\ell_1 + 1 - \gamma_1)A_3 - \gamma_1 + (-1)^{\gamma_1}(-2 + m_1^* + (1 - m_2^*)m_1^*)$$
$$\implies \ell - \ell_1 = \frac{(-1)^{\gamma+1}(1 + 2x_1 + (m_1^* - m_2))}{A_3}.$$

Given $\ell - \ell_1 \in \mathbb{Z}$ and $0 \leq 1 + 2x_1 + (m_1^* - m_2) \leq A_3 - 2$ by construction, we must have $x_1 = 0$ and $m_1^* - m_2 = -1$. This means $m_1^* = 0$ and $m_2 = 1$. Going to $j_2$ with these results, we get

$$m_2 + \gamma + 4x_2 = m_2^* + \gamma_1 + 4x_2^* + (1 - m_2^*)m_1^*$$
$$\implies x_2 - x_2^* = \frac{(-1)^\gamma}{4}.$$

Since $x_2 - x_2^* \in \mathbb{Z}$ and for all $\gamma \in \{0, 1\}$, $\frac{(-1)^\gamma}{4} \notin \mathbb{Z}$ by construction, we have a contradiction as the above equality is not possible.

<u>Major Case 2</u>: Focusing on moves along $j_2$, we observe that all $j_2$ moves belonging to $C_2$ here are moves with the opposite orientation for $C_1$'s $j_2$ move $(m_1, m_2) = (1, 1)$, and $m_2^* = 0$ for all such moves. Applying these observations to $j_2$, it follows that

$$m_2 + \gamma + 4x_2 - 1 = m_2^* + \gamma_1 + 4x_2^* + (1 - m_2^*)m_1^*$$
$$\implies x_2 - x_2^* = \frac{(\gamma_1 - \gamma) + m_1^*}{4}.$$

Given $x_2 - x_2^* \in \mathbb{Z}$, $|\gamma_1 - \gamma| \leq 1$, and $m_1^* \in \{0, 1\}$ all by construction, $m_1^* = -(\gamma_1 - \gamma)$ and so $\gamma_1 - \gamma = 0$ or $\gamma_1 - \gamma = -1$ as $\gamma_1 - \gamma = 1$ implies that for every $m_1^* \in \{0, 1\}$, $\frac{1+m_1^*}{4} \notin \mathbb{Z}$, giving us our contradiction as the above equality would not be possible in that case. We proceed by treating the remaining two cases together with $\gamma_1 - \gamma = -|\gamma_1 - \gamma|$ and $m_1^* = |\gamma_1 - \gamma|$.

Applying the above to $j_1$, we see

$$(\ell + \gamma)A_3 - \gamma + (-1)^\gamma(m_1 + 2x_1) = (\ell_1 + 1 - \gamma_1)A_3 - \gamma_1 + (-1)^{\gamma_1}(-2 + m_1^* + (1 - m_2^*)m_1^*)$$
$$\implies (\ell - \ell_1) + (\gamma + \gamma_1) - 1 = \frac{(-1)^{\gamma+1}(1 + 2x_1) + |\gamma_1 - \gamma| + 2(-1)^{\gamma_1+1}(1 - |\gamma_1 - \gamma|)}{A_3}.$$

In either case that $\gamma_1 - \gamma = -1$ or $\gamma_1 - \gamma = 0$, we see that $(\ell - \ell_1) + (\gamma + \gamma_1) - 1 \in \mathbb{Z}$ while the right hand side is not, giving us a contradiction as the above equality is not possible.

**Set Comparison Case 4:** We will consider whether edges are shared from moves with the opposite orientation via sets 1 and 4.

<u>Major Case 1</u>: Focusing on moves along $j_1$, we see that $C_2$ has edges defined by inverse $j_1$ moves to those defining $C_1$ if and only if $\gamma = \gamma_1$. Given that $m_2^* = 1$ during all $j_1$ moves in this case, in $j_1$ the above yields

$$(\ell + \gamma)A_3 - \gamma + (-1)^\gamma(m_1 + 2x_1 + 1) = (\ell_1 + 1 - \gamma_1)A_3 - \gamma_1 + (-1)^{\gamma_1+1}(m_1^* + m_1^*(1 - m_2^*))$$
$$\implies (\ell - \ell_1) + 2\gamma - 1 = \frac{(-1)^{\gamma+1}(m_1 + m_1^* + 2x_1 + 1)}{A_3}.$$

Given $(\ell - \ell_1) + 2\gamma - 1 \in \mathbb{Z}$ and $1 \leq m_1 + m_1^* + 2x_1 + 1 \leq A_3 - 1$ with $A_3 \geq 4$ all by construction, we see that $\frac{(-1)^{\gamma+1}(m_1+m_1^*+2x_1+1)}{A_3} \notin \mathbb{Z}$, a contradiction as the above equality is not possible.

<u>Major Case 2</u>: Focusing on moves along $j_2$, we see that $m_1 = m_2$, $2m_2 + \chi = 1$, $m_2^* = 0$, and the inverse $j_2$ moves have the form $(m_1^*, m_2^*) = (1 - m_1, 0)$. Applying these observations to $j_2$, we get

$$m_2 + \gamma + 4x_2 + \chi = 2 + m_1^* + \gamma_1 + m_2^*(1 - m_1^*) + 4x_2^*$$
$$\implies x_2 - x_2^* = \frac{2 + (\gamma_1 - \gamma)}{4}.$$

Given $x_2 - x_2^* \in \mathbb{Z}$ and $|\gamma_1 - \gamma| \leq 1$ by construction, we have that for all $\gamma, \gamma_1 \in \{0, 1\}$, $\frac{2+(\gamma_1-\gamma)}{4} \notin \mathbb{Z}$, a contradiction as the above equality is not possible.

**Set Comparison Case 5:** We will now consider whether edges are shared from moves with the opposite orientation via set 2.

Major Case 1: Focusing on moves along $j_1$, we see that $C_2$ is defined by inverse $j_1$ moves relative to those of $C_1$ if and only if $|\gamma_1 - \gamma| = 1$. Observing that $m_1 = m_2$ and $m_1^* = m_2^*$ during $j_1$ moves, in $j_1$ we get

$$(\ell + 1 - \gamma)A_3 - \gamma + (-1)^{\gamma+1}(2 + m_2 + 2m_1(1 - m_2) + 2x_1 + 1)$$
$$= (\ell_1 + 1 - \gamma_1)A_3 - \gamma_1 + (-1)^{\gamma_1+1}(2 + m_2^* + 2m_1^*(1 - m_2^*) + 2x_1^*)$$
$$\implies ((\ell - \ell_1) - (\gamma - \gamma_1))\frac{A_3}{2} + (-1)^{\gamma+1}(2 + (x_1 + x_1^*)) = \frac{(-1)^\gamma (m_2 + m_2^*)}{2}.$$

Given $((\ell - \ell_1) - (\gamma - \gamma_1))\frac{A_3}{2} + (-1)^{\gamma+1}(2 + (x_1 + x_1^*)) \in \mathbb{Z}$ and $0 \leq m_2 + m_2^* \leq 2$ by construction, we must have $m_2 = m_2^*$. Applying this to $j_2$, we have

$$2 + \gamma + m_1 + 4x_2 = 2 + \gamma_1 + m_1^* + 4x_2^*$$
$$\implies x_2 - x_2^* = \frac{(-1)^\gamma}{4}.$$

Since $x_2 - x_2^* \in \mathbb{Z}$ and for every $\gamma \in \{0, 1\}$, $\frac{(-1)^\gamma}{4} \notin \mathbb{Z}$, we have a contradiction as the above equality is not possible.

Major Case 2: Focusing on moves along $j_2$ observe that which $j_2$ moves of $C_2$ correspond as moves with the opposite orientation to those of $C_1$ depends on $x_1$ and $x_1^*$. Hence, we case on $x_1$ and then further case on $x_1^*$:

Case 1: Let $0 \leq x_1 < \frac{A_3}{2} - 2$. We now case on $x_1^*$:

Subcase 1: Let $0 \leq x_1^* < \frac{A_3}{2} - 2$. Then, during such moves it is the case $m_1 + m_2 = 1 = m_1^* + m_2^*$ and $(m_1^*, m_2^*) = (1 - m_1, m_1)$ during all $j_2$ moves. Applying these observations to $j_2$, we get

$$2 + m_1 + \gamma + 4x_2 + (m_2 - m_1) = 2 + m_1^* + \gamma_1 + 4x_2^*$$
$$\implies x_2 - x_2^* = \frac{\gamma_1 - \gamma}{4}.$$

Since $x_2 - x_2^* \in \mathbb{Z}$ and $|\gamma_1 - \gamma| \leq 1$ by construction, it must be the case $\gamma = \gamma_1$. Hence, in $j_1$ what we have thus far implies

$$(\ell + 1 - \gamma)A_3 - \gamma + (-1)^{\gamma+1}(2 + m_2 + 2m_1(1 - m_2) + 2x_1)$$
$$= (\ell_1 + 1 - \gamma_1)A_3 - \gamma_1 + (-1)^{\gamma_1+1}(2 + m_2^* + 2m_1^*(1 - m_2^*) + 2x_1^*)$$
$$\implies (\ell - \ell_1)\frac{A_3}{2} + (-1)^{\gamma+1}(m_1 + (x_1 - x_1^*)) = \frac{(-1)^{\gamma+1}}{2}.$$

Since $(\ell - \ell_1)\frac{A_3}{2} + (-1)^{\gamma+1}(m_1 + (x_1 - x_1^*)) \in \mathbb{Z}$ by construction and for all $\gamma \in \{0, 1\}$, $\frac{(-1)^{\gamma+1}}{2} \notin \mathbb{Z}$, we have a contradiction as the above equality is not possible.

Subcase 2: Let $x_1^* = \frac{A_3}{2} - 2$. Here all $j_2$ moves from $C_2$ serve as moves with the opposite orientation for $C_1$'s $j_2$ move $(m_1, m_2) = (1, 0)$, and during such moves it is the case $m_1^* + m_2^* = 1$. Applying this to $j_2$, it follows that

$$2 + m_1 + \gamma + 4x_2 + (m_2 - m_1) = 2 + m_1^* + \gamma_1 + 4x_2^*$$
$$\implies x_2 - x_2^* = \frac{(\gamma_1 - \gamma) + m_1^*}{4}.$$

Since $x_2 - x_2^* \in \mathbb{Z}$ and $-1 \leq (\gamma_1 - \gamma) + m_1^* \leq 2$ by construction, we see $\gamma_1 - \gamma = -1$ or $\gamma_1 - \gamma = 0$, and in both cases $m_1^* = |\gamma_1 - \gamma|$. Note that $\gamma_1 - \gamma = 1$ would give us a contradiction as for every $m_1^* \in \{0, 1\}$, $\frac{1+m_1^*}{4} \notin \mathbb{Z}$. Treating the remaining two cases of $\gamma_1 - \gamma$ together, in $j_1$ we see

$$(\ell + 1 - \gamma)A_3 - \gamma + (-1)^{\gamma+1}(2 + m_2 + 2m_1(1 - m_2) + 2x_1)$$
$$= (\ell_1 + 1 - \gamma_1)A_3 - \gamma_1 + (-1)^{\gamma_1+1}(2 + m_2^* + 2m_1^*(1 - m_2^*) + 2x_1^*)$$
$$\implies ((\ell - \ell_1) + (\gamma_1 - \gamma) - (-1)^{\gamma_1+1})\frac{A_3}{2} + (-1)^{\gamma+1}(2 + x_1) = \frac{|\gamma_1 - \gamma| + (-1)^{\gamma_1}(1 - |\gamma_1 - \gamma|)}{2}.$$

In either case $\gamma_1 - \gamma = -1$ or $\gamma_1 - \gamma = 0$, we have $((\ell - \ell_1) + (\gamma_1 - \gamma) - (-1)^{\gamma_1+1})\frac{A_3}{2} + (-1)^{\gamma+1}(2 + x_1) \in \mathbb{Z}$ and $\frac{|\gamma_1 - \gamma| + (-1)^{\gamma_1}(1 - |\gamma_1 - \gamma|)}{2} \notin \mathbb{Z}$, giving us a contradiction as the above equality is not possible.

<u>Case 2</u>: Let $x_1 = \frac{A_3}{2} - 2$. Given that the argument with $x_1 \neq x_1^*$ in case 1 is symmetrical, we need only consider when $x_1^* = \frac{A_3}{2} - 2$. In this case, all defined $j_2$ moves have the same orientation with respect to each other, and so there are no inverse $j_2$ moves to consider.

**Set Comparison Case 6:** We will now consider whether edges are shared via sets 2 and 3.

<u>Major Case 1</u>: Focusing on moves along $j_1$, we see that $C_2$'s $j_1$ moves only serve as moves with the opposite orientation to those of $C_1$ if and only if $\gamma = \gamma_1$. Applying this to $j_1$ along with $m_1 = m_2$ and $m_2^* = 1$, we have

$$(\ell + 1 - \gamma)A_3 - \gamma + (-1)^{\gamma+1}(2 + m_2 + 2m_1(1 - m_2) + 2x_1 + 1)$$
$$= (\ell_1 + 1 - \gamma_1)A_3 - \gamma_1 + (-1)^{\gamma_1}(-2 + m_1^* + (1 - m_2^*)m_1^*)$$
$$\implies \ell - \ell_1 = \frac{(-1)^{\gamma}(1 + 2x_1 + m_1 + m_1^*)}{A_3}.$$

Given $\ell - \ell_1 \in \mathbb{Z}$ and $1 \leq 1 + 2x_1 + m_1 + m_1^* \leq A_3 - 1$ by construction, we see that for all $\gamma \in \{0, 1\}$, $\frac{(-1)^{\gamma}(1+2x_1+m_1+m_1^*)}{A_3} \notin \mathbb{Z}$. Thus, we have a contradiction as the above equality is not possible.

<u>Major Case 2</u>: Focusing on moves along $j_2$, we see that the $j_2$ moves of $C_2$ that serve as moves with the opposite orientation to those of $C_1$ depend on $x_1$. Hence, we case on $x_1$:

<u>Case 1</u>: Let $0 \leq x_1 < \frac{A_3}{2} - 2$. Then, all $j_2$ moves of $C_2$ serve as moves with the opposite orientation only to $C_1$'s $j_2$ move $(m_1, m_2) = (1, 0)$ and during such moves $m_2^* = 0$. Applying this to $j_2$, we get

$$2 + m_1 + \gamma + 4x_2 + (m_2 - m_1) = m_2^* + \gamma_1 + 4x_2^* + (1 - m_2^*)m_1^*$$
$$\implies x_2 - x_2^* = \frac{-2 + (\gamma_1 - \gamma) + m_1^*}{4}.$$

Since $x_2 - x_2^* \in \mathbb{Z}$ and $-3 \leq -2 + (\gamma_1 - \gamma) + m_1^* \leq 0$, we must have $\gamma_1 - \gamma = 1 = m_1^*$. Applying our results to $j_1$, we have

$$(\ell + 1 - \gamma)A_3 - \gamma + (-1)^{\gamma+1}(2 + m_2 + 2m_1(1 - m_2) + 2x_1)$$
$$= (\ell_1 + 1 - \gamma_1)A_3 - \gamma_1 + (-1)^{\gamma_1}(-2 + m_1^* + (1 - m_2^*)m_1^*)$$
$$\implies ((\ell - \ell_1) + (\gamma_1 - \gamma))\frac{A_3}{2} + (-1)^{\gamma+1}(2 + x_1) = -\frac{1}{2}.$$

Given $((\ell - \ell_1) + (\gamma_1 - \gamma))\frac{A_3}{2} + (-1)^{\gamma+1}(2 + x_1) \in \mathbb{Z}$ and $-\frac{1}{2} \notin \mathbb{Z}$, we have a contradiction as the above equality is not possible.

<u>Subcase 2:</u> Let $x_1 = \frac{A_3}{2} - 2$. Then, we see that there are no inverse $j_2$ moves of $C_2$ to consider with respect to $C_1$ as all $j_2$ moves have the same orientation with respect to each other.

**Set Comparison Case 7:** We will consider whether edges are shared from moves with the opposite orientation via sets 2 and 4.

<u>Major Case 1:</u> Focusing on moves along $j_1$, we see that $C_2$ is defined by inverse $j_1$ moves with respect to $C_1$ if and only if $|\gamma_1 - \gamma| = 1$. Observing that $m_1 = m_2$ and $m_2^* = 1$ during all $j_1$ moves, in $j_1$ we see

$$(\ell + 1 - \gamma)A_3 - \gamma + (-1)^{\gamma+1}(2 + m_2 + 2m_1(1 - m_2) + 2x_1 + 1)$$
$$= (\ell_1 + 1 - \gamma_1)A_3 - \gamma_1 + (-1)^{\gamma_1+1}(m_1^* + m_1^*(1 - m_2^*))$$
$$\implies (\ell - \ell_1) + (-1)^{\gamma} = \frac{(-1)^{\gamma}(2 + 2x_1 + m_1 + m_1^*)}{A_3}.$$

Given $(\ell - \ell_1) + (-1)^{\gamma} \in \mathbb{Z}$ and $2 \leq 2 + 2x_1 + m_1 + m_1^* \leq A_3$, it must be the case $x_1 = \frac{A_3}{2} - 2$ and $m_1 = 1 = m_1^*$. Applying what we have so far to $j_1$, it follows that

$$2 + m_1 + \gamma + 4x_2 = 2 + m_1^* + \gamma_1 + m_2^*(1 - m_1^*) + 4x_2^*$$
$$\implies x_2 - x_2^* = \frac{(-1)^{\gamma}}{4}.$$

Since $x_2 - x_2^* \in \mathbb{Z}$ and for all $\gamma \in \{0, 1\}$, $\frac{(-1)^{\gamma}}{4} \notin \mathbb{Z}$ by construction, we have a contradiction as the above equality is not possible.

<u>Major Case 2:</u> Focusing on moves along $j_2$, observe that which $j_2$ moves from $C_2$ serve as moves with the opposite orientation to those of $C_1$ depends on $x_1$. So we case on $x_1$:

<u>Case 1:</u> Let $0 \leq x_1 < \frac{A_3}{2} - 2$. In this case, all $j_2$ moves satisfy $m_1 + m_2 = 1$ and $(m_1^*, m_2^*) = (1 - m_1, 0)$. Consequently, in $j_2$ we see

$$2 + m_1 + \gamma + 4x_2 + (m_2 - m_1) = 2 + m_1^* + \gamma_1 + m_2^*(1 - m_1^*) + 4x_2^*$$
$$\implies x_2 - x_2^* = \frac{\gamma_1 - \gamma}{4}$$

Given $x_2 - x_2^* \in \mathbb{Z}$ and $|\gamma_1 - \gamma| \leq 1$ by construction, we see that $\gamma = \gamma_1$. Applying this to $j_1$, we get

$$(\ell + 1 - \gamma)A_3 - \gamma + (-1)^{\gamma+1}(2 + m_2 + 2m_1(1 - m_2) + 2x_1)$$
$$= (\ell_1 + 1 - \gamma_1)A_3 - \gamma_1 + (-1)^{\gamma_1+1}(m_1^* + m_1^*(1 - m_2^*))$$

$$\implies \ell - \ell_1 = \frac{(-1)^\gamma(1 + 3m_1 + 2x_1)}{A_3}.$$

Following from $\ell - \ell_1 \in \mathbb{Z}$ and $1 \leq 1 + 3m_1 + 2x_1 < A_3$ by construction, we find that for every $0 \leq x_1 < \frac{A_3}{2} - 2$ and $\gamma, m_1 \in \{0, 1\}$, $\frac{(-1)^\gamma(1+3m_1+2x_1)}{A_3} \notin \mathbb{Z}$. Thus, we have a contradiction as the above equality is not possible.

<u>Case 2</u>: Let $x_1 = \frac{A_3}{2} - 2$. In this case, only $C_2$'s $j_2$ move $(m_1^*, m_2^*) = (1, 0)$ serves as an inverse to all of $C_1$'s $j_2$ moves, and during such moves $m_1 + m_2 = 1$. Applying this to $j_2$, we get

$$2 + m_1 + \gamma + 4x_2 + 1 = 2 + m_1^* + \gamma_1 + m_2^*(1 - m_1^*) + 4x_2^*$$

$$\implies x_2 - x_2^* = \frac{(\gamma_1 - \gamma) - m_1}{4}.$$

Since $x_2 - x_2^* \in \mathbb{Z}$, $|\gamma_1 - \gamma| \leq 1$, and $m_1 \in \{0, 1\}$, we see that either $\gamma_1 - \gamma = 1$ or $\gamma_1 - \gamma = 0$ with $m_1 = |\gamma_1 - \gamma|$ in both cases. Note that $\gamma_1 - \gamma = -1$ is not a possibility as that would require $m_1 = -1$, which is not possible as $m_1 \in \{0, 1\}$ by construction. With $\gamma_1 - \gamma = |\gamma_1 - \gamma|$ for the remaining two cases, our results and observations then give us in $j_1$ that

$$(\ell + 1 - \gamma)A_3 - \gamma + (-1)^{\gamma+1}(2 + m_2 + 2m_1(1 - m_2) + 2x_1)$$
$$= (\ell_1 + 1 - \gamma_1)A_3 - \gamma_1 + (-1)^{\gamma_1+1}(m_1^* + m_1^*(1 - m_2^*))$$

$$\implies (\ell - \ell_1) + (\gamma_1 - \gamma) + (-1)^{\gamma+1} = -\frac{|\gamma_1 - \gamma| + (-1)^\gamma(1 - |\gamma_1 - \gamma|) + 2(-1)^{\gamma_1}}{A_3}.$$

In either case that $\gamma_1 - \gamma = 1$ or $\gamma_1 - \gamma = 0$, we see that $(\ell - \ell_1) + (\gamma_1 - \gamma) + (-1)^{\gamma+1} \in \mathbb{Z}$ by construction while $-\frac{|\gamma_1-\gamma|+(-1)^\gamma(1-|\gamma_1-\gamma|)+2(-1)^{\gamma_1}}{A_3} \notin \mathbb{Z}$, giving us our contradiction as the above equality is thus not possible.

**Set Comparison Case 8:** We will consider whether edges are shared from moves with the opposite orientation via set 3.

<u>Major Case 1</u>: Focusing on moves along $j_1$, observe that $C_2$ has inverse $j_1$ moves with respect to those of $C_1$ if and only if $|\gamma_1 - \gamma| = 1$. Following from $m_2 = 1 = m_2^*$ during $j_1$ moves, in $j_2$ we get

$$m_2 + \gamma + 4x_2 + (1 - m_2)m_1 = m_2^* + \gamma_1 + 4x_2^* + (1 - m_2^*)m_1^*$$

$$\implies x_2 - x_2^* = \frac{(-1)^\gamma}{4}.$$

Given $x_2 - x_2^* \in \mathbb{Z}$ by construction and for all $\gamma \in \{0, 1\}$, $\frac{(-1)^\gamma}{4} \notin \mathbb{Z}$, we have a contradiction as the above equality is not possible.

<u>Major Case 2</u>: Focusing on moves along $j_2$, we see that there are no moves with the opposing orientation in $C_2$ with respect to those of $C_1$, and so there is nothing to prove.

**Set Comparison Case 9:** We will consider whether edges are shared from moves with the opposite orientation via sets 3 and 4.

<u>Major Case 1</u>: Focusing on moves along $j_1$, we see that $C_2$ has inverse $j_1$ moves with respect to those of $C_1$ if and only if $\gamma_1 = \gamma$. Then, given $m_2 = 1 = m_2^*$ during $j_1$ moves, in $j_2$ we find

$$m_2 + \gamma + 4x_2 + (1 - m_2)m_1 = 2 + m_1^* + \gamma_1 + m_2^*(1 - m_1^*) + 4x_2^*$$
$$\implies x_2 - x_2^* = \frac{1}{2}.$$

Since $x_2 - x_2^* \in \mathbb{Z}$ by construction and $\frac{1}{2} \notin \mathbb{Z}$, we have a contradiction as the above equality is not possible.

<u>Major Case 2</u>: Focusing on moves along $j_2$, we see that $C_2$'s $j_2$ move $(m_1^*, m_2^*) = (1, 0)$ is the only inverse with respect to all of $C_1$'s $j_2$ moves. Hence, the above along with the observation that $m_2 = 0$ during $j_2$ moves yields in $j_2$

$$m_2 + \gamma + 4x_2 + (1 - m_2)m_1 + 1 = 2 + m_1^* + \gamma_1 + (1 - m_1^*)m_2^* + 4x_2^*$$
$$\implies x_2 - x_2^* = \frac{2 + (\gamma_1 - \gamma) - m_1}{4}.$$

Since $x_2 - x_2^* \in \mathbb{Z}$ and $0 \leq 2 + (\gamma_1 - \gamma) - m_1 \leq 3$ by construction, it follows that $\gamma_1 - \gamma = -1$, meaning $\gamma = 1$ and $\gamma_1 = 0$, and $m_1 = 1$. Applying this to $j_1$, we have

$$(\ell + 1 - \gamma)A_3 - \gamma + (-1)^\gamma(-2 + m_1 + (1 - m_2)m_1) = (\ell_1 + 1 - \gamma_1)A_3 - \gamma_1 + (-1)^{\gamma_1 + 1}(m_1^* + m_1^*(1 - m_2^*))$$
$$\implies (\ell - \ell_1) - 1 = -\frac{1}{A_3}.$$

Given $(\ell - \ell_1) - 1 \in \mathbb{Z}$ and $-\frac{1}{A_3} \notin \mathbb{Z}$ since $A_3 \geq 4$, we have a contradiction as the above equality is not possible.

**Set Comparison Case 10:** We will consider whether edges are shared from moves with the opposite orientation via set 4.

<u>Major Case 1</u>: Focusing on moves along $j_1$, we see that $C_2$ has inverse $j_1$ moves with respect to $C_1$ if and only if $|\gamma_1 - \gamma| = 1$. Observing that $m_2 = 1 = m_2^*$ during $j_1$ moves, in $j_2$ we have

$$2 + m_1 + \gamma + m_2(1 - m_1) + 4x_2 = 2 + m_1^* + \gamma_1 + m_2^*(1 - m_1^*) + 4x_2^*$$
$$\implies x_2 - x_2^* = \frac{(-1)^\gamma}{4}.$$

Given $x_2 - x_2^* \in \mathbb{Z}$ by construction and for all $\gamma \in \{0, 1\}$, $\frac{(-1)^\gamma}{4} \notin \mathbb{Z}$, we have a contradiction as the equality above is not possible.

<u>Major Case 2</u>: Focusing on moves along $j_2$, we see that during $j_2$ moves $m_2 = 0$ and $(m_1^*, m_2^*) = (1 - m_1, 0)$. Consequently, noting that $(-1)^{m_1} = 1 - 2m_1$ since $m_1 \in \{0, 1\}$, in $j_2$ we get

$$2 + m_1 + \gamma + m_2(1 - m_1) + 4x_2 + (-1)^{m_1} = 2 + m_1^* + \gamma_1 + m_2^*(1 - m_1^*) + 4x_2^*$$
$$\implies x_2 - x_2^* = \frac{\gamma_1 - \gamma}{4}.$$

Given $x_2 - x_2^* \in \mathbb{Z}$ and $|\gamma_1 - \gamma| \leq 1$, we see that it must be the case that $\gamma = \gamma_1$. Applying this to $j_1$, we obtain

$$(\ell+1-\gamma)A_3 - \gamma + (-1)^{\gamma+1}(m_1 + m_1(1-m_2)) = (\ell_1+1-\gamma_1)A_3 - \gamma_1 + (-1)^{\gamma_1+1}(m_1^* + m_1^*(1-m_2^*))$$

$$\implies \ell - \ell_1 = \frac{2(-1)^{\gamma+m_1+1}}{A_3}.$$

Since $\ell - \ell_1 \in \mathbb{Z}$ by construction and for every $\gamma, m_1 \in \{0,1\}$, $\frac{2(-1)^{\gamma+m_1+1}}{A_3} \notin \mathbb{Z}$, we have a contradiction as the above equality is not possible.

Thus, no edges are shared from moves with the opposite orientation for $\alpha_1 < \alpha^* < \alpha_2$ when $d = 2$.

$\boldsymbol{\alpha^* - }$**Case 3:** Let $\alpha^* = \alpha_2$. In this case, sets 1 and 3 are active as $A_1 = 1$ and $A_2 = 1$.

**Set Comparison Case 1:** We will consider whether edges are shared from moves with the opposite orientation via set 1.

Major Case 1: Focusing on moves along $j_1$, we observe that there are no inverse $j_1$ moves to consider as all $j_1$ moves are defined with the same orientation with respect to each other.

Major Case 2: Focusing on moves along $j_2$, we see that during all $j_2$ moves $m_1 = m_2$ and $(m_1^*, m_2^*) = (1-m_1, 1-m_2)$. Consequently, in $j_1$ we have

$$\ell A_3 + (m_1 + 2x_1) = \ell_1 A_3 + (m_1^* + 2x_1^*)$$

$$\implies (\ell - \ell_1)\frac{A_3}{2} + (m_1 + (x_1 - x_1^*)) = \frac{1}{2}.$$

Given $(\ell - \ell_1)\frac{A_3}{2} + (m_1 + (x_1 - x_1^*)) \in \mathbb{Z}$ by construction and $\frac{1}{2} \notin \mathbb{Z}$, we have a contradiction as the above equality is not possible.

**Set Comparison Case 2:** We will consider whether edges are shared from moves with the opposite orientation via sets 1 and 3.

Major Case 1: Focusing on moves along $j_1$, we see all $j_1$ moves defined have the same orientation with respect to each other and hence there are no inverse $j_1$ moves in $C_2$ to consider with respect to $C_1$.

Major Case 2: Focusing on moves along $j_2$, observe that all $j_2$ moves from $C_2$ serve as moves with the opposite orientation only to $C_1$'s $j_2$ move $(m_1, m_2) = (1, 1)$. Combining this with the observation that $m_1^* = m_2^*$ during $j_2$ moves, we get in $j_1$

$$\ell A_3 + (m_1 + 2x_1) = (\ell_1 + 1)A_3 + (-2 + m_1^*)$$

$$\implies \ell - \ell_1 - 1 = -\frac{2x_1 + 3 - m_1^*}{A_3}.$$

Given $\ell - \ell_1 - 1 \in \mathbb{Z}$ and $2 \leq 2x_1 + 3 - m_1^* \leq A_3 - 1$ by construction, it follows that for every $0 \leq x_1 \leq \frac{A_3}{2} - 2$ and $m_1^* \in \{0, 1\}$, $-\frac{2x_1+3-m_1^*}{A_3} \notin \mathbb{Z}$, giving us our contradiction as the above equality is not possible.

**Set Comparison Case 3:** We will consider whether edges are shared from moves with the opposite orientation via set 3.

Major Case 1: Focusing on moves along $j_1$, we see that there are no inverse $j_1$ moves to consider as all $j_1$ moves are defined with the same orientation with respect to each other.

Major Case 2: Focusing on moves along $j_2$, observe that there are no inverse $j_2$ moves to consider as all $j_2$ moves have the same orientation with respect to each other.

Thus, no edges are shared from moves with the opposite orientation for $\alpha^* = \alpha_2$ when $d = 2$, allowing us to conclude that no edges are shared from moves with the opposite orientation for $\alpha_1 \leq \alpha^* \leq \alpha_2$ when $d = 2$.

**Dimension Case 2:** Let $d \geq 3$. Then, $C_1 = C_{\ell,\gamma,t,p_1,s_1,\ldots,p_{d-2},s_{d-2}}$ and $C_2 = C_{\ell_1,\gamma_1,t_1,p_1^*,s_1^*,\ldots,p_{d-2}^*,s_{d-2}^*}$. We now case on $\alpha^*$ for $\alpha_1 \leq \alpha^* \leq \alpha_d$.

**$\alpha^*$−Case 1:** Let $\alpha^* = \alpha_1$. In this case, only sets 3 and 4 are active as $A_1 = 0$ and $A_2 = 0$. Note that $0 \leq t \leq 1$, $p_1 = 0 = p_1^*$, $A_4 = 1$, and $x_3 = 0 = x_3^*$ as a consequence of our assumption on $d$ and $\alpha^*$.

**Set Comparison Case 1:** We will consider whether edges are shared from moves with the opposite orientation via set 3.

Major Case 1: Focusing on moves along $j_1$, we see that $C_2$ has inverse $j_1$ moves with respect to $C_1$ if and only if $|\gamma_1 - \gamma| = 1$. Observing that $m_2 = 1 = m_2^*$, $m_1 = m_d$ and $m_1^* = m_d^*$ during $j_1$ moves, in $j_2$ we have

$$m_2 + \gamma + 4x_2 + (1 - m_2)m_1 + 2(1 - m_1)m_d = m_2^* + \gamma_1 + 4x_2^* + (1 - m_2^*)m_1^* + 2(1 - m_1^*)m_d^*$$

$$\implies x_2 - x_2^* = \frac{(-1)^\gamma}{4}.$$

Given $x_2 - x_2^* \in \mathbb{Z}$ by construction and for all $\gamma \in \{0, 1\}$, $\frac{(-1)^\gamma}{4} \notin \mathbb{Z}$, we have a contradiction as the above equality is not possible.

Major Case 2: Focusing on moves along $j_2$, we see that there are no inverse $j_2$ moves to consider as all $j_2$ moves have the same orientation.

Major Case 3: Focusing on moves along $j_3$, we see that during $j_3$ moves $m_1 = m_2$, $m_d = 1 - m_1$, and

$$(m_1^*, \ldots, m_d^*) = ((1 - |\gamma_1 - \gamma|)(1 - m_1) + |\gamma_1 - \gamma|m_1, \ldots, (1 - |\gamma_1 - \gamma|)(1 - m_d) + |\gamma_1 - \gamma|m_d).$$

Noting that $2m_2 - 1 = (-1)^{m_2+1}$, in $j_2$ we get

$$m_2 + \gamma + 4x_2 + (1 - m_2)m_1 + 2(1 - m_1)m_d = m_2^* + \gamma_1 + 4x_2^* + (1 - m_2^*)m_1^* + 2(1 - m_1^*)m_d^*$$

$$\implies x_2 - x_2^* = \frac{(-1)^\gamma|\gamma_1 - \gamma| + (-1)^{m_2+1}(1 - |\gamma_1 - \gamma|)}{4}.$$

Given $x_2 - x_2^* \in \mathbb{Z}$ by construction and for every $m_2, \gamma, \gamma_1 \in \{0, 1\}$, $\frac{(-1)^\gamma|\gamma_1-\gamma|+(-1)^{m_2+1}(1-|\gamma_1-\gamma|)}{4} \notin \mathbb{Z}$, we have a contradiction as the above equality is not possible.

<u>Major Case 4</u>: Focusing on $j_k$ moves for $4 \leq k \leq d$ if $d \geq 4$, we see that during such moves $m_3^* = 1 - m_3$, giving us in $j_3$

$$p_1 + 2A_4 s_1 + (-1)^\gamma m_3 + 2x_3 = p_1^* + 2A_4 s_1^* + (-1)^{\gamma_1} m_3^* + 2x_3^*$$

$$\implies (s_1 - s_1^*) + \frac{(-1)^{\gamma_1} + (-1)^\gamma}{2} m_3 = \frac{(-1)^{\gamma_1}}{2}.$$

Since $(s_1 - s_1^*) + \frac{(-1)^{\gamma_1}+(-1)^\gamma}{2} m_3 \in \mathbb{Z}$ by construction and for every $\gamma_1 \in \{0, 1\}$, $\frac{(-1)^{\gamma_1}}{2} \notin \mathbb{Z}$, we have a contradiction as the above equality is not possible.

**Set Comparison Case 2:** We will consider whether edges are shared from moves with the opposite orientation via sets 3 and 4.

<u>Major Case 1</u>: Focusing on moves along $j_1$, we see that $C_2$ has inverse $j_1$ moves with respect to $C_1$ if and only if $\gamma = \gamma_1$. Observing that during $j_1$ moves $m_1 = m_d$, $m_2 = 1 = m_2^*$ and $m_1^* = m_d^*$, in $j_2$ we have

$$m_2 + \gamma + 4x_2 + (1 - m_2)m_1 + 2(1 - m_1)m_d = 2 + m_1^* + \gamma_1 + (1 - m_1^*)m_2^* + 4x_2^* + 2(1 - m_1^*)m_d^*$$

$$\implies x_2 - x_2^* = \frac{1}{2}.$$

Since $x_2 - x_2^* \in \mathbb{Z}$ by construction and $\frac{1}{2} \notin \mathbb{Z}$, we have a contradiction as the above equality is not possible.

<u>Major Case 2</u>: Focusing on moves along $j_2$, observe that there are no inverse $j_2$ moves to consider as all $j_2$ moves are defined to have the same orientation.

<u>Major Case 3</u>: Focusing on moves along $j_3$, it follows that during $j_3$ moves $m_1 = m_2$, $m_d = 1 - m_1$, and

$$(m_1^*, \ldots, m_d^*) = ((1 - |\gamma_1 - \gamma|)(1 - m_1) + |\gamma_1 - \gamma|m_1, \ldots, (1 - |\gamma_1 - \gamma|)(1 - m_d) + |\gamma_1 - \gamma|m_d).$$

Applying this to $j_2$ along with $(-1)^{m_2+1} = 2m_2 - 1$, we get

$$m_2 + \gamma + 4x_2 + (1 - m_2)m_1 + 2(1 - m_1)m_d = 2 + m_1^* + \gamma_1 + (1 - m_1^*)m_2^* + 4x_2^* + 2(1 - m_1^*)m_d^*$$

$$\implies x_2 - x_2^* = \frac{2 + (-1)^\gamma|\gamma_1 - \gamma| + (-1)^{m_2+1}(1 - |\gamma_1 - \gamma|)}{4}.$$

Given $x_2 - x_2^* \in \mathbb{Z}$ by construction, we have a contradiction as for every $m_2, \gamma, \gamma_1 \in \{0, 1\}$, $\frac{2+(-1)^\gamma|\gamma_1-\gamma|+(-1)^{m_2+1}(1-|\gamma_1-\gamma|)}{4} \notin \mathbb{Z}$, hence making the above equality not possible.

**Major Case 4:** Focusing on moves along $j_k$ for $4 \leq k \leq d$ when $d \geq 4$, observe that $m_3^* = 1 - m_3$ during all such moves, giving us in $j_3$

$$p_1 + 2A_4 s_1 + (-1)^\gamma m_3 + 2x_3 = 2s_1^* + (-1)^{\gamma_1} m_3^*$$
$$\implies (s_1 - s_1^*) + \frac{(-1)^{\gamma_1} + (-1)^\gamma}{2} m_3 = \frac{(-1)^{\gamma_1}}{2}.$$

Following from $(s_1 - s_1^*) + \frac{(-1)^{\gamma_1} + (-1)^\gamma}{2} m_3 \in \mathbb{Z}$ by construction and for every $\gamma_1 \in \{0,1\}$, $\frac{(-1)^{\gamma_1}}{2} \notin \mathbb{Z}$, we have a contradiction as the above equality is not possible.

**Set Comparison 3:** We will consider whether edges are shared from moves with the opposite orientation via set 4.

**Major Case 1:** Focusing on moves along $j_1$, we see that $C_2$ has inverse $j_1$ moves to those of $C_1$ if and only if $|\gamma_1 - \gamma| = 1$. Observing that during $j_1$ moves $m_1 = m_d$, $m_1^* = m_d^*$ and $m_2 = 1 = m_2^*$, in $j_2$ these observations yield

$$2 + m_1 + \gamma + m_2(1 - m_1) + 4x_2 + 2(1 - m_1)m_d = 2 + m_1^* + \gamma_1 + m_2^*(1 - m_1^*) + 4x_2^* + 2(1 - m_1^*)m_d^*$$
$$\implies x_2 - x_2^* = \frac{(-1)^\gamma}{4}.$$

Given $x_2 - x_2^* \in \mathbb{Z}$ by construction and for all $\gamma \in \{0,1\}$, $\frac{(-1)^\gamma}{4} \notin \mathbb{Z}$, we have a contradiction as the above equality is not possible.

**Major Case 2:** Focusing on moves along $j_2$, we see that there are no $j_2$ moves to consider as every $j_2$ move defined in $C_2$ has the same orientation with respect to the other $j_2$ moves defined in $C_1$.

**Major Case 3:** Focusing on moves along $j_3$, we see that during $j_3$ moves $m_1 = m_2$, $m_d = 1 - m_1$, and

$$(m_1^*, \ldots, m_d^*) = ((1 - |\gamma_1 - \gamma|)(1 - m_1) + |\gamma_1 - \gamma|m_1, \ldots, (1 - |\gamma_1 - \gamma|)(1 - m_d) + |\gamma_1 - \gamma|m_d).$$

Applying the above to $j_2$ along with $(-1)^{m_2+1} = 2m_2 - 1$, we have

$$2 + m_1 + \gamma + m_2(1 - m_1) + 4x_2 + 2(1 - m_1)m_d = 2 + m_1^* + \gamma_1 + m_2^*(1 - m_1^*) + 4x_2^* + 2(1 - m_1^*)m_d^*$$
$$\implies x_2 - x_2^* = \frac{(-1)^\gamma |\gamma_1 - \gamma| + (-1)^{m_2+1}(1 - |\gamma_1 - \gamma|)}{4}.$$

Since $x_2 - x_2^* \in \mathbb{Z}$ by construction and for every $m_2, \gamma, \gamma_1 \in \{0,1\}$, $\frac{(-1)^\gamma |\gamma_1 - \gamma| + (-1)^{m_2+1}(1 - |\gamma_1 - \gamma|)}{4} \notin \mathbb{Z}$, we have a contradiction as the above equality is not possible.

**Major Case 4:** Focusing on moves along $j_k$ for $4 \leq k \leq d$ when $d \geq 4$, observe that $m_3^* = 1 - m_3$ during all such moves, giving us in $j_3$

$$2s_1 + (-1)^\gamma m_3 = 2s_1^* + (-1)^{\gamma_1} m_3^*$$

$$\implies (s_1 - s_1^*) + \frac{(-1)^{\gamma_1} + (-1)^{\gamma}}{2} m_3 = \frac{(-1)^{\gamma_1}}{2}.$$

Following from $(s_1 - s_1^*) + \frac{(-1)^{\gamma_1}+(-1)^{\gamma}}{2} m_3 \in \mathbb{Z}$ by construction and for every $\gamma_1 \in \{0,1\}$, $\frac{(-1)^{\gamma_1}}{2} \notin \mathbb{Z}$, we have a contradiction as the above equality is not possible.

Thus, no edges are shared from moves with the opposite orientation for $\alpha^* = \alpha_1$ when $d \geq 3$.

**$\alpha^*$−Case 2:** Let $\alpha_1 < \alpha^* < \alpha_2$. Then, $A_1 = 1$ and $A_2 = 0$, meaning sets $1, 2, 3,$ and $4$ are active. Note that here $t = 0 = t_1$, $p_1 = 0 = p_1^*$, $A_4 = 1$, and $x_3 = 0 = x_3^*$ due to our assumption on $\alpha^*$.

**Set Comparison Case 1:** We will consider whether edges are shared from moves with the opposite orientation via set 1.

Major Case 1: Focusing on moves along $j_1$, observe that all $j_1$ moves are defined with the same orientation, and so they have the same orientation with respect to each other. Hence, there are no inverse $j_1$ moves to consider.

Major Case 2: Focusing on moves along $j_2$, observe that during all $j_2$ moves $m_1 = m_d$ and $(m_1^*, m_d^*) = (1 - m_1, 1 - m_d)$. Applying this to $j_1$, we find

$$\ell A_3 + 1 + (m_1 + 2x_1 + 2m_d(1 - m_1)) = \ell_1 A_3 + 1 + (m_1^* + 2x_1^* + 2m_d^*(1 - m_1^*))$$
$$\implies (\ell - \ell_1)\frac{A_3}{2} + (x_1 - x_1^*) + m_1 = \frac{1}{2}.$$

Given $(\ell - \ell_1)\frac{A_3}{2} + (x_1 - x_1^*) + m_1 \in \mathbb{Z}$ by construction and $\frac{1}{2} \notin \mathbb{Z}$, we have a contradiction as the above equality is not possible.

Major Case 3: Focusing on moves along $j_k$ for $3 \leq k \leq d$, observe that during all such moves $m_d = 1 - m_1$, and $(m_1^*, m_d^*) = (1 - m_1, 1 - m_d)$. Applying this to $j_1$, we see

$$\ell A_3 + 1 + (m_1 + 2x_1 + 2m_d(1 - m_1)) = \ell_1 A_3 + 1 + (m_1^* + 2x_1^* + 2m_d^*(1 - m_1^*))$$
$$\implies (\ell - \ell_1)\frac{A_3}{2} + m_1 + (2m_d - 1) + (x_1 - x_1^*) = \frac{1}{2}.$$

Since $(\ell-\ell_1)\frac{A_3}{2} + m_1 + (2m_d - 1) + (x_1 - x_1^*) \in \mathbb{Z}$ by construction and $\frac{1}{2} \notin \mathbb{Z}$, we have a contradiction as the above equality is not possible.

**Set Comparison Case 2:** We will consider whether edges are shared from moves with the opposite orientation via sets 1 and 2.

Major Case 1: Focusing on moves along $j_1$, observe that during such moves $m_1 = m_3 = m_d$, $m_2 = 1 - m_1$, and $m_1^* = m_2^* = m_d^*$. Noting that $r_1^* = 0$ since $m_1^* = m_3^*$ when $d = 3$ and $m_{d-1}^* = m_d^*$ when $d \geq 4$ during these $j_1$ moves, in $j_1$ we get

$$\ell A_3 + 1 + (m_1 + 2x_1 + 2m_d(1 - m_1) + 1)$$
$$= (\ell_1 + 1)A_3 - (2 + m_2^* + 2m_1^*(1 - m_2^*) + 2x_1^* + 2m_d^*(1 - m_1^*) - r_1^*(A_3 - 1))$$
$$\implies (\ell - \ell_1 - 1)\frac{A_3}{2} + (x_1 + x_1^*) + 2 = -\frac{m_1 + m_2^*}{2}.$$

Since $(\ell - \ell_1 - 1)\frac{A_3}{2} + (x_1 + x_1^*) + 2 \in \mathbb{Z}$ and $0 \leq m_1 + m_2^* \leq 2$, we have $m_1 = m_2^*$. Applying this to $j_3$, we see

$$p_1 + s_1 A_4 + 1 - m_3 + 2x_3 = s_1^* + m_1^* + r_1^*$$
$$\implies s_1 - s_1^* = 2m_1 - 1.$$

Lastly, applying this to $j_2$, we find

$$m_2 + \gamma + (s_1 + 1) + 4x_2 = 2 + \gamma_1 + (1 + s_1^* + m_d^* - r_1^*) + 4x_2^*$$
$$\implies x_2 - x_2^* = \frac{2 + (\gamma_1 - \gamma)}{4}.$$

Given $x_2 - x_2^* \in \mathbb{Z}$ and $1 \leq 2 + (\gamma_1 - \gamma) \leq 3$ by construction and for every $\gamma, \gamma_1 \in \{0,1\}$, $\frac{2+(\gamma_1-\gamma)}{4} \notin \mathbb{Z}$, we have a contradiction as the above equality is not possible.

<u>Major Case 2</u>: Focusing on moves along $j_2$, we see that during these moves $m_1 = m_2 = m_d$, $m_1^* = m_2^*$, and $m_d^* = 1 - m_1^*$. Applying this to $j_1$, we have

$$\ell A_3 + 1 + (m_1 + 2x_1 + 2m_d(1 - m_1))$$
$$= (\ell_1 + 1)A_3 - (2 + m_2^* + 2m_1^*(1 - m_2^*) + 2x_1^* + 2m_d^*(1 - m_1^*) - r_1^*(A_3 - 1))$$
$$\implies (\ell - \ell_1 - 1 - r_1^*)A_3 + 2(x_1 + x_1^*) + (4 + m_1 + m_d^* + r_1^*) = 0.$$

Given that which $j_2$ moves of $C_2$ will serve as moves with the opposite orientation to those of $C_1$ depends on $x_1^*$, we case on $x_1^*$:

<u>Case 1</u>: Let $0 \leq x_1^* < \frac{A_3}{2} - 2$. Then, we also get in particular $(m_1^*, m_2^*, m_d^*) = (m_1, m_2, 1 - m_1)$, meaning $r_1^* = 0$ since $x_1^* \neq \frac{A_3}{2} - 2$, and so our equation in $j_1$ gives us

$$(\ell - \ell_1 - 1)\frac{A_3}{2} + (x_1 + x_1^*) = -\frac{5}{2}.$$

Since $(\ell - \ell_1 - 1)\frac{A_3}{2} + (x_1 + x_1^*) \in \mathbb{Z}$ by construction and $-\frac{5}{2} \notin \mathbb{Z}$, we have a contradiction as the above equality is not possible.

<u>Case 2</u>: Let $x_1^* = \frac{A_3}{2} - 2$. Then, all $j_2$ moves have the same orientation and so all serve as moves with the opposite orientation only to $C_1$'s $j_1$ move with $(m_1, m_2, m_d) = (1, 1, 1)$. Hence, in $j_1$ we get

$$\ell - \ell_1 - r_1^* = -\frac{2x_1 + 1 + m_d^* + r_1^*}{A_3}.$$

Following from $\ell - \ell_1 - r_1^* \in \mathbb{Z}$ and $1 \leq 2x_1 + 1 + m_d^* + r_1^* \leq A_3 - 1$ by construction, it follows that for every $0 \leq x_1 \leq \frac{A_3}{2} - 2$ and $m_d^*, r_1^* \in \{0,1\}$, $-\frac{2x_1+1+m_d^*+r_1^*}{A_3} \notin \mathbb{Z}$, giving us a contradiction as the above equality is not possible.

<u>Major Case 3</u>: Focusing on moves along $j_3$, observe that during such moves $m_1 = m_2$, $m_d = 1 - m_1$, and $(m_1^*, m_2^*, m_d^*) = (1 - m_1, m_2, m_d)$. Note that the above imply $m_1^* = m_d^*$, $m_2^* = 1 - m_1^*$, and $r_1^* = 0$ since $m_1^* = m_3^*$ when $d = 3$ and $m_{d-1}^* = m_d^*$ when $d \geq 4$ during such moves. Consequently, in $j_1$ we find

$$\ell A_3 + 1 + (m_1 + 2x_1 + 2m_d(1 - m_1))$$
$$= (\ell_1 + 1)A_3 - (2 + m_2^* + 2m_1^*(1 - m_2^*) + 2x_1^* + 2m_d^*(1 - m_1^*) - r_1^*(A_3 - 1))$$
$$\implies (\ell - \ell_1 - 1)\frac{A_3}{2} + m_d + (x_1 + x_1^*) = -\frac{5}{2}.$$

Since $(\ell - \ell_1 - 1)\frac{A_3}{2} + m_d + (x_1 + x_1^*) \in \mathbb{Z}$ by construction and $-\frac{5}{2} \notin \mathbb{Z}$, we have a contradiction as the above equality is not possible.

<u>Major Case 4</u>: Observe that during all $j_k$ moves for $4 \leq k \leq d$ when $d \geq 4$ it is the case $m_1 = m_2 = m_3$, $m_d = 1 - m_1$, and $(m_1^*, m_2^*, m_d^*) = (m_2, m_3, m_d)$. Hence, $r_1^* = 0$ since $m_{d-1}^* = m_d^*$, and so in $j_1$ we obtain

$$\ell A_3 + 1 + (m_1 + 2x_1 + 2m_d(1 - m_1))$$
$$= (\ell_1 + 1)A_3 - (2 + m_2^* + 2m_1^*(1 - m_2^*) + 2x_1^* + 2m_d^*(1 - m_1^*) - r_1^*(A_3 - 1))$$
$$\implies (\ell - \ell_1 - 1)\frac{A_3}{2} + m_d + (x_1 + x_1^*) = -\frac{5}{2}.$$

Given $(\ell - \ell_1 - 1)\frac{A_3}{2} + m_d + (x_1 + x_1^*) \in \mathbb{Z}$ by construction and $-\frac{5}{2} \notin \mathbb{Z}$, we have a contradiction as the above equality is not possible.

**Set Comparison Case 3:** We will consider whether edges are shared from moves with the opposite orientation via sets 1 and 3.

<u>Major Case 1</u>: Focusing on moves along $j_1$, observe that all $j_1$ moves defined have the same orientation with respect to each other. Hence, there are no inverse $j_1$ moves to consider.

<u>Major Case 2</u>: Focusing on moves along $j_2$, observe that all of $C_2$'s $j_2$ moves serve as moves with the opposite orientation only to $C_1$'s $j_2$ move with $(m_1, m_d) = (1, 1)$. Noting that it is also the case $m_1^* = m_d^*$, in $j_1$ we get

$$\ell A_3 + 1 + (m_1 + 2x_1 + 2m_d(1 - m_1)) = (\ell_1 + 1)A_3 + (-2 + m_1^* + 2m_d^*(1 - m_1^*) - (A_3 - 2)(1 - m_d^*))$$
$$\implies \ell - \ell_1 - m_d^* = -\frac{2x_1 + 2 + m_1^*}{A_3}.$$

Given $\ell - \ell_1 - m_d^* \in \mathbb{Z}$ and $2 \leq 2x_1 + 2 + m_1^* \leq A_3 - 1$ by construction, we see that for every $0 \leq x_1 \leq \frac{A_3}{2} - 2$ and $m_1^* \in \{0, 1\}$, $-\frac{2x_1 + 2 + m_1^*}{A_3} \notin \mathbb{Z}$, giving us a contradiction as the above equality is not possible.

<u>Major Case 3</u>: Focusing on moves along $j_k$ for $3 \leq k \leq d$, observe that during all such moves $m_d = 1 - m_1$ and $(m_1^*, m_d^*) = (m_1, m_d)$, and so in $j_1$ we have

$$\ell A_3 + 1 + (m_1 + 2x_1 + 2m_d(1 - m_1)) = (\ell_1 + 1)A_3 + (-2 + m_1^* + 2m_d^*(1 - m_1^*) - (A_3 - 2)(1 - m_d^*))$$
$$\implies (\ell - \ell_1 - m_d^*)\frac{A_3}{2} + m_d^* + x_1 = -\frac{1}{2}.$$

Given $(\ell - \ell_1 - m_d^*)\frac{A_3}{2} + m_d^* + x_1 \in \mathbb{Z}$ by construction and $-\frac{1}{2} \notin \mathbb{Z}$, we have a contradiction as the above equality is not possible.

**Set Comparison Case 4:** We will consider whether edges are shared from moves with the opposite orientation via sets 1 and 4.

<u>Major Case 1</u>: Focusing on moves along $j_1$, observe that during these moves $m_1 = m_d$, $m_2 = 1 - m_1$, $m_1^* = m_d^*$, and $m_2^* = 1$. Applying this to $j_1$, we get

$$\ell A_3 + 1 + (m_1 + 2x_1 + 2m_d(1-m_1) + 1) = (\ell_1 + 1)A_3 - (m_1^* + m_1^*(1-m_2^*) + 2m_d^*(1-m_1^*))$$
$$\implies \ell - \ell_1 - 1 = -\frac{2x_1 + 2 + m_1 + m_1^*}{A_3}.$$

Since $\ell - \ell_1 - 1 \in \mathbb{Z}$ and $2 \leq 2x_1 + 2 + m_1 + m_1^* \leq A_3$, it must be the case $x_1 = \frac{A_3}{2} - 2$ and $m_1 + m_1^* = 2$, meaning $m_1 = 1 = m_1^*$. Applying this to $j_3$ along with the observation that $m_3 = m_1$ and $m_3^* = m_1^*$, we obtain

$$p_1 + s_1 A_4 + 1 - m_3 + 2x_3 = s_1^* + m_1^* m_3^*$$
$$\implies s_1 - s_1^* = 1.$$

Hence, in $j_2$ we now have

$$m_2 + \gamma + (s_1 + 1) + 4x_2 = 2 + m_1^* + (s_1^* + (1-m_1^*)m_d^* + m_d^*) + \gamma_1 + m_2^*(1-m_1^*) + 4x_2^*$$
$$\implies x_2 - x_2^* = \frac{2 + (\gamma_1 - \gamma)}{4}.$$

Given $x_2 - x_2^* \in \mathbb{Z}$ and $1 \leq 2 + (\gamma_1 - \gamma) \leq 3$ by construction, it follows that for every $\gamma, \gamma_1 \in \{0, 1\}$, $\frac{2+(\gamma_1-\gamma)}{4} \notin \mathbb{Z}$, giving us a contradiction as the above equality is not possible.

<u>Major Case 2</u>: Focusing on moves along $j_2$, observe that during such moves $m_1 = m_d$ and $(m_1^*, m_2^*, m_d^*) = (0, 0, 1 - m_d)$, giving us in $j_1$

$$\ell A_3 + 1 + (m_1 + 2x_1 + 2m_d(1-m_1)) = (\ell_1 + 1)A_3 - (m_1^* + m_1^*(1-m_2^*) + 2m_d^*(1-m_1^*))$$
$$\implies \ell - \ell_1 - 1 = -\frac{3 + 2x_1 - m_1}{A_3}.$$

Given $\ell - \ell_1 - 1 \in \mathbb{Z}$ and $2 \leq 3 + 2x_1 - m_1 \leq A_3 - 1$ by construction, then for every $0 \leq x_1 \leq \frac{A_3}{2} - 2$ and $m_1 \in \{0, 1\}$, $-\frac{3+2x_1-m_1}{A_3} \notin \mathbb{Z}$, yielding a contradiction as the above equality is not possible.

<u>Major Case 3</u>: Focusing on moves along $j_3$, observe that during these moves $m_1 = m_2$, $m_d = 1 - m_1$, and $(m_1^*, m_2^*, m_d^*) = (1, m_2, m_d)$. Applying these observations to $j_1$, we get

$$\ell A_3 + 1 + (m_1 + 2x_1 + 2m_d(1-m_1)) = (\ell_1 + 1)A_3 - (m_1^* + m_1^*(1-m_2^*) + 2m_d^*(1-m_1^*))$$
$$\implies (\ell - \ell_1 - 1)\frac{A_3}{2} + m_d + x_1 = -\frac{3}{2}.$$

Since $(\ell - \ell_1 - 1)\frac{A_3}{2} + m_d + x_1 \in \mathbb{Z}$ by construction and $-\frac{3}{2} \notin \mathbb{Z}$, we have a contradiction as the above equality is not possible.

<u>Major Case 4</u>: Focusing on moves along $j_k$ for $4 \leq k \leq d$ when $d \geq 4$, observe that during all such moves $m_1 = m_2$, $m_d = 1 - m_1$, and $(m_1^*, m_2^*, m_d^*) = (m_1, m_2, m_d)$. Thus, in $j_1$ we get

$$\ell A_3 + 1 + (m_1 + 2x_1 + 2m_d(1 - m_1)) = (\ell_1 + 1)A_3 - (m_1^* + m_1^*(1 - m_2^*) + 2m_d^*(1 - m_1^*))$$

$$\implies (\ell - \ell_1 - 1)\frac{A_3}{2} + m_d + x_1 = -\frac{3}{2}.$$

Following from $(\ell - \ell_1 - 1)\frac{A_3}{2} + m_d + x_1 \in \mathbb{Z}$ by construction and $-\frac{3}{2} \notin \mathbb{Z}$, we have a contradiction as the above equality is not possible.

**Set Comparison Case 5:** We will consider whether edges are shared from moves with the opposite orientation via set 2.

<u>Major Case 1</u>: Focusing on moves along $j_1$, observe that all $j_1$ moves defined have the same orientation and so there are no inverse $j_1$ moves to consider.

<u>Major Case 2</u>: Focusing on moves along $j_2$, we see that which $j_2$ moves of $C_2$ serve as moves with the opposite orientation to those of $C_1$ depends on $x_1$ and $x_1^*$. Consequently, we case on $x_1$ and case further on $x_1^*$ in each case:

<u>Case 1</u>: Let $0 \leq x_1 < \frac{A_3}{2} - 2$. To determine which $j_2$ moves in $C_2$ serve as moves with the opposite orientation to those of $C_1$, we further case on $x_1^*$:

<u>Subcase 1</u>: Let $0 \leq x_1^* < \frac{A_3}{2} - 2$. Then, it is the case that $m_1 = m_2$, $m_d = 1 - m_1 = m_3$, and $(m_1^*, m_2^*, m_d^*) = (1 - m_1, 1 - m_2, 1 - m_d)$. Noting that $r_1 = 0 = r_1^*$ since $x_1 \neq \frac{A_3}{2} - 2 \neq x_1^*$, in $j_1$ we find

$$(\ell + 1)A_3 - (2 + m_2 + 2m_1(1 - m_2) + 2x_1 + 2m_d(1 - m_1) - r_1(A_3 - 1))$$
$$= (\ell_1 + 1)A_3 - (2 + m_2^* + 2m_1^*(1 - m_2^*) + 2x_1^* + 2m_d^*(1 - m_1^*) - r_1^*(A_3 - 1))$$

$$\implies (\ell - \ell_1)\frac{A_3}{2} - m_2 + (x_1^* - x_1) + (m_d^* - m_d) = -\frac{1}{2}.$$

Since $(\ell - \ell_1)\frac{A_3}{2} - m_2 + (x_1^* - x_1) + (m_d^* - m_d) \in \mathbb{Z}$ by construction and $-\frac{1}{2} \notin \mathbb{Z}$, we have a contradiction as the above equality is not possible.

<u>Subcase 2</u>: Let $x_1^* = \frac{A_3}{2} - 2$. Now observe that $(m_1, m_2, m_d) = (0, 0, 1)$, $m_1^* = \cdots = m_{d-1}^*$, $m_3^* = \eta_3(1 - \eta_4)(1 - m_1^*) + \eta_4 m_1^*$, and $m_d^* = 1 - m_1^*$. Hence, $r_1 = 0$ since $x_1 \neq \frac{A_3}{2} - 2$ and so using the definition of $r_1^*$ for $d \geq 3$, in $j_1$ we have

$$(\ell + 1)A_3 - (2 + m_2 + 2m_1(1 - m_2) + 2x_1 + 2m_d(1 - m_1) - r_1(A_3 - 1))$$
$$= (\ell_1 + 1)A_3 - (2 + m_2^* + 2m_1^*(1 - m_2^*) + 2x_1^* + 2m_d^*(1 - m_1^*) - r_1^*(A_3 - 1))$$

$$\implies (\ell - \ell_1 - r_1^* + 1)\frac{A_3}{2} - (2 + x_1) = \frac{(-1)^{m_1^*+1}}{2}.$$

Since $(\ell - \ell_1 - r_1^* + 1)\frac{A_3}{2} - (2 + x_1) \in \mathbb{Z}$ by construction and for every $m_1^* \in \{0, 1\}$, $\frac{(-1)^{m_1^*+1}}{2} \notin \mathbb{Z}$, we have a contradiction as the above equality is not possible.

<u>Case 2:</u> Let $x_1 = \frac{A_3}{2} - 2$. Given that the argument with $x_1 \neq x_1^*$ in case 1 is symmetrical, we need only consider when $x_1^* = \frac{A_3}{2} - 2$. In this case, all defined $j_2$ moves have the same orientation with respect to each other, and so there are no inverse $j_2$ moves to consider.

<u>Major Case 3:</u> Focusing on moves along $j_3$, observe that during such moves $m_1 = m_d$, $m_2 = 1 - m_1$, and $(m_1^*, m_2^*, m_d^*) = (1 - m_1, 1 - m_2, 1 - m_d)$. Note that $m_1 = m_3$ and $m_1^* = m_3^*$ when $d = 3$, and $m_{d-1} = m_d$ and $m_{d-1}^* = m_d^*$ when $d \geq 4$. It now follows that $r_1 = 0 = r_1^*$ and so in $j_1$ we see

$$(\ell + 1)A_3 - (2 + m_2 + 2m_1(1 - m_2) + 2x_1 + 2m_d(1 - m_1) - r_1(A_3 - 1))$$
$$= (\ell_1 + 1)A_3 - (2 + m_2^* + 2m_1^*(1 - m_2^*) + 2x_1^* + 2m_d^*(1 - m_1^*) - r_1^*(A_3 - 1))$$
$$\implies (\ell - \ell_1)\frac{A_3}{2} - m_2 + (1 - 2m_1) + (x_1^* - x_1) = -\frac{1}{2}.$$

Given $(\ell - \ell_1)\frac{A_3}{2} - m_2 + (1 - 2m_1) + (x_1^* - x_1) \in \mathbb{Z}$ by construction and $-\frac{1}{2} \notin \mathbb{Z}$, we have a contradiction as the above equality is not possible.

<u>Major Case 4:</u> Focusing on moves along $j_k$ for $4 \leq k \leq d$ when $d \geq 4$, we see that during all such moves $m_1 = m_2$, $m_d = 1 - m_1$, and $(m_1^*, m_2^*, m_d^*) = (1 - m_1, 1 - m_2, 1 - m_d)$. Note that $m_{d-1} = m_d$ and $m_{d-1}^* = m_d^*$, meaning $r_1 = 0 = r_1^*$. Applying these observations to $j_1$, we find

$$(\ell + 1)A_3 - (2 + m_2 + 2m_1(1 - m_2) + 2x_1 + 2m_d(1 - m_1) - r_1(A_3 - 1))$$
$$= (\ell_1 + 1)A_3 - (2 + m_2^* + 2m_1^*(1 - m_2^*) + 2x_1^* + 2m_d^*(1 - m_1^*) - r_1^*(A_3 - 1))$$
$$\implies (\ell - \ell_1)\frac{A_3}{2} - m_2 + (x_1^* - x_1) + (m_d^* - m_d) = -\frac{1}{2}.$$

Since $(\ell - \ell_1)\frac{A_3}{2} - m_2 + (x_1^* - x_1) + (m_d^* - m_d) \in \mathbb{Z}$ by construction and $-\frac{1}{2} \notin \mathbb{Z}$, we have a contradiction as the above equality is not possible.

**Set Comparison Case 6:** We will consider whether edges are shared from moves with the opposite orientation via sets 2 and 3.

<u>Major Case 1:</u> Focusing on moves along $j_1$, observe that during these moves $m_1 = \cdots = m_d$, $m_1^* = m_3^* = m_d^*$, and $m_2^* = 1 - m_1^*$. In either case that $d = 3$ or $d \geq 4$, it follows that $r_1 = 0$ as $m_1 = m_3$ when $d = 3$ and $m_{d-1} = m_d$ when $d \geq 4$. Applying all of the above to $j_1$, we get

$$(\ell + 1)A_3 - (2 + m_2 + 2m_1(1 - m_2) + 2x_1 + 2m_d(1 - m_1) - r_1(A_3 - 1) + 1)$$
$$= (\ell_1 + 1)A_3 + (-2 + m_1^* + 2m_d^*(1 - m_1^*) - (A_3 - 2)(1 - m_d^*))$$
$$\implies \ell - \ell_1 + m_2^* = \frac{2x_1 + 2 + m_2 + m_2^*}{A_3}.$$

Since $\ell - \ell_1 + m_2^* \in \mathbb{Z}$ and $2 \leq 2x_1 + 2 + m_2 + m_2^* \leq A_3$ by construction, it must be the case $m_2 + m_2^* = 2$, meaning $m_2 = 1 = m_2^*$. Applying this to $j_3$, we find

$$s_1 + m_1 + r_1 = p_1^* + s_1^*A_4 + m_3^* + 2x_3^*$$
$$\implies s_1 - s_1^* = -1.$$

Lastly, in $j_2$ we see

$$2 + \gamma + (1 + s_1 + m_d - r_1) + 4x_2 = m_2^* + \gamma_1 + 4x_2^* + 2m_1^*(1 - m_2^*) + (s_1^* + 2m_d^*(1 - m_1^*))$$
$$\implies x_2 - x_2^* = \frac{-2 + (\gamma_1 - \gamma)}{4}.$$

Given $x_2 - x_2^* \in \mathbb{Z}$ and $-3 \leq -2 + (\gamma_1 - \gamma) \leq -1$ by construction, it follows that $\frac{-2+(\gamma_1-\gamma)}{4} \notin \mathbb{Z}$. Consequently, we have a contradiction as the above equality is not possible.

<u>Major Case 2</u>: Focusing on moves along $j_2$, we see that $C_2$'s $j_2$ moves serve as moves with the opposite orientation to those of $C_1$ if and only if $0 \leq x_1 < \frac{A_3}{2} - 2$. We exclude $x_1 = \frac{A_3}{2} - 2$ as then $C_1$'s $j_2$ moves would be defined with the same orientation as those of $C_2$, meaning there would be no moves with the opposite orientation to consider.

So suppose $0 \leq x_1 < \frac{A_3}{2} - 2$. Then, observe that all of $C_2$'s $j_2$ moves serve as moves with the opposite orientation only to $C_1$'s $j_2$ move with $(m_1, m_2, m_d) = (0, 0, 1)$ and during such moves $m_1^* = m_2^* = m_d^*$. Using the above along with the observation $r_1 = 0$ since $x_1 \neq \frac{A_3}{2} - 2$, in $j_1$ we get

$$(\ell + 1)A_3 - (2 + m_2 + 2m_1(1 - m_2) + 2x_1 + 2m_d(1 - m_1) - r_1(A_3 - 1))$$
$$= (\ell_1 + 1)A_3 + (-2 + m_1^* + 2m_d^*(1 - m_1^*) - (A_3 - 2)(1 - m_d^*))$$
$$\implies (\ell - \ell_1) + (1 - m_d^*) = \frac{2x_1 + 4 - m_1^*}{A_3}.$$

Since $(\ell - \ell_1) + (1 - m_d^*) \in \mathbb{Z}$ and $3 \leq 2x_1 + 4 - m_1^* < A_3$ by construction and our case assumption, we find $\frac{2x_1+4-m_1^*}{A_3} \notin \mathbb{Z}$. Hence, we have a contradiction as the above equality is not possible.

<u>Major Case 3</u>: Focusing on moves along $j_3$, we see that it is then the case $m_1 = m_d$, $m_2 = 1 - m_1$, and $(m_1^*, m_2^*, m_d^*) = (m_1, 1 - m_2, 1 - m_d)$. Noting that $m_1 = m_3$ when $d = 3$ and $m_{d-1} = m_d$ when $d \geq 4$, we find $r_1 = 0$. Applying this to $j_1$, we have

$$(\ell + 1)A_3 - (2 + m_2 + 2m_1(1 - m_2) + 2x_1 + 2m_d(1 - m_1) - r_1(A_3 - 1))$$
$$= (\ell_1 + 1)A_3 + (-2 + m_1^* + 2m_d^*(1 - m_1^*) - (A_3 - 2)(1 - m_d^*))$$
$$\implies (\ell - \ell_1 + m_1^*)\frac{A_3}{2} - (m_1 + x_1) = \frac{3}{2}.$$

Given $(\ell - \ell_1 + m_1^*)\frac{A_3}{2} - (m_1 + x_1) \in \mathbb{Z}$ by construction and $\frac{3}{2} \notin \mathbb{Z}$, we have a contradiction as the above equality is not possible.

<u>Major Case 4</u>: Focusing on moves along $j_k$ for $4 \leq k \leq d$ when $d \geq 4$, we see that during all such moves $m_1 = m_2$, $m_d = 1 - m_1$, and $(m_1^*, m_2^*, m_d^*) = (1 - m_1, 1 - m_2, 1 - m_d)$. Note that $r_1 = 0$ as $m_{d-1} = m_d$. Thus, in $j_1$ we have

$$(\ell + 1)A_3 - (2 + m_2 + 2m_1(1 - m_2) + 2x_1 + 2m_d(1 - m_1) - r_1(A_3 - 1))$$
$$= (\ell_1 + 1)A_3 + (-2 + m_1^* + 2m_d^*(1 - m_1^*) - (A_3 - 2)(1 - m_d^*))$$
$$\implies (\ell - \ell_1 + m_1^*)\frac{A_3}{2} - (m_d + x_1) = \frac{3}{2}.$$

Following from $(\ell - \ell_1 + m_1^*)\frac{A_3}{2} - (m_d + x_1) \in \mathbb{Z}$ by construction and $\frac{3}{2} \notin \mathbb{Z}$, we see we have a contradiction as the above equality is not possible.

**Set Comparison Case 7:** We will consider whether edges are shared from moves with the opposite orientation via sets 2 and 4.

<u>Major Case 1:</u> Focusing on moves along $j_1$, we see that all $j_1$ moves have the same orientation with respect to each other, and so there are no inverse $j_1$ moves to consider.

<u>Major Case 2:</u> Focusing on moves along $j_2$, observe that which of $C_2$'s $j_2$ moves serve as moves with the opposite orientation to those of $C_1$ depend on $x_1$. Hence, we proceed by casing on $x_1$:

<u>Case 1:</u> Let $0 \leq x_1 < \frac{A_3}{2} - 2$. Observing that in this case $m_1 = m_2$, $m_d = 1 - m_1$, $(m_1^*, m_2^*, m_d^*) = (0, 0, 1 - m_d)$, and $r_1 = 0$ since $x_1 \neq \frac{A_3}{2} - 2$, in $j_1$ we find

$$(\ell + 1)A_3 - (2 + m_2 + 2m_1(1 - m_2) + 2x_1 + 2m_d(1 - m_1) - r_1(A_3 - 1))$$
$$= (\ell_1 + 1)A_3 - (m_1^* + m_1^*(1 - m_2^*) + 2m_d^*(1 - m_1^*))$$
$$\implies \ell - \ell_1 = \frac{4 + 2x_1 - 3m_1}{A_3}.$$

Since $\ell - \ell_1 \in \mathbb{Z}$ and $1 \leq 4 + 2x_1 - 3m_1 < A_3$ by construction and our case assumption, we see that for every $0 \leq x_1 < \frac{A_3}{2} - 2$ and $m_1 \in \{0, 1\}$, $\frac{4+2x_1-3m_1}{A_3} \notin \mathbb{Z}$, giving us a contradiction as the above equality is not possible.

<u>Case 2:</u> Let $x_1 = \frac{A_3}{2} - 2$. Then, we see that $m_1 = m_2$, $m_d = 1 - m_1$, and in particular only $C_2$'s $j_2$ move with $(m_1^*, m_2^*, m_d^*) = (0, 0, 1)$ serves as an inverse to all of $C_1$'s $j_2$ moves. Hence, in $j_1$ we get

$$(\ell + 1)A_3 - (2 + m_2 + 2m_1(1 - m_2) + 2x_1 + 2m_d(1 - m_1) - r_1(A_3 - 1))$$
$$= (\ell_1 + 1)A_3 - (m_1^* + m_1^*(1 - m_2^*) + 2m_d^*(1 - m_1^*))$$
$$\implies \ell - \ell_1 + r_1 - 1 = \frac{m_d + r_1 - 3}{A_3}.$$

Given $\ell - \ell_1 + r_1 - 1 \in \mathbb{Z}$ and $-3 \leq m_d + r_1 - 3 \leq -1$ with $A_3 \geq 4$ by construction and for every $m_d, r_1 \in \{0, 1\}$, $\frac{m_d + r_1 - 3}{A_3} \notin \mathbb{Z}$, we have a contradiction as the above equality is not possible.

<u>Major Case 3:</u> Focusing on moves along $j_3$, observe that during these moves $m_1 = m_d$, $m_2 = 1 - m_1$, and $(m_1^*, m_2^*, m_d^*) = (1, 1 - m_2, 1 - m_d)$. Note that $r_1 = 0$ as $m_1 = m_3$ when $d = 3$ and $m_{d-1} = m_d$ when $d \geq 4$. Applying these observations to $j_1$, we obtain

$$(\ell + 1)A_3 - (2 + m_2 + 2m_1(1 - m_2) + 2x_1 + 2m_d(1 - m_1) - r_1(A_3 - 1))$$
$$= (\ell_1 + 1)A_3 - (m_1^* + m_1^*(1 - m_2^*) + 2m_d^*(1 - m_1^*))$$
$$\implies (\ell - \ell_1)\frac{A_3}{2} - (m_1 + x_1) = \frac{1}{2}.$$

Given $(\ell - \ell_1)\frac{A_3}{2} - (m_1 + x_1) \in \mathbb{Z}$ by construction and $\frac{1}{2} \notin \mathbb{Z}$, we have a contradiction as the above equality is not possible.

<u>Major Case 4:</u> Observing that $m_1 = m_2$, $m_d = 1 - m_1$, $(m_1^*, m_2^*, m_d^*) = (1 - m_1, 1 - m_2, 1 - m_d)$, and $r_1 = 0$ as $m_{d-1} = m_d$ during all $j_k$ moves for $4 \leq k \leq d$ when $d \geq 4$, in $j_1$ we see

$$(\ell + 1)A_3 - (2 + m_2 + 2m_1(1 - m_2) + 2x_1 + 2m_d(1 - m_1) - r_1(A_3 - 1))$$

$$= (\ell_1 + 1)A_3 - (m_1^* + m_1^*(1 - m_2^*) + 2m_d^*(1 - m_1^*))$$
$$\implies (\ell - \ell_1)\frac{A_3}{2} - (m_2 + (m_d - m_d^*) + x_1) = \frac{1}{2}.$$

Following from $(\ell - \ell_1)\frac{A_3}{2} - (m_2 + (m_d - m_d^*) + x_1) \in \mathbb{Z}$ by construction and $\frac{1}{2} \notin \mathbb{Z}$, we have a contradiction as the above equality is not possible.

**Set Comparison Case 8:** We will consider whether edges are shared from moves with the opposite orientation via set 3.

Major Case 1: Focusing on moves along $j_1$, observe that all $j_1$ moves are defined to have the same orientation with respect to each other and so there are no inverse $j_1$ moves to consider.

Major Case 2: Focusing on moves along $j_2$, note that all $j_2$ moves are defined to have the same orientation with respect to each other and so there are no inverse $j_2$ moves to consider.

Major Case 3: Focusing on moves along $j_k$ for $3 \leq k \leq d$, we find that during all such moves $m_1 = m_2$, $m_d = 1 - m_1$, and $(m_1^*, m_2^*, m_d^*) = (1 - m_1, 1 - m_2, 1 - m_d)$. Applying this to $j_1$, we see

$$(\ell + 1)A_3 + (-2 + m_1 + 2m_d(1 - m_1) - (A_3 - 2)(1 - m_d))$$
$$= (\ell_1 + 1)A_3 + (-2 + m_1^* + 2m_d^*(1 - m_1^*) - (A_3 - 2)(1 - m_d^*))$$
$$\implies (\ell - \ell_1 + 2m_d - 1)\frac{A_3}{2} + m_1 = \frac{1}{2}.$$

Given $(\ell - \ell_1 + 2m_d - 1)\frac{A_3}{2} + m_1 \in \mathbb{Z}$ by construction and $\frac{1}{2} \notin \mathbb{Z}$, we have a contradiction as the above equality is not possible.

**Set Comparison Case 9:** We will consider whether edges are shared from moves with the opposite orientation via sets 3 and 4.

Major Case 1: Focusing on moves along $j_1$, we see that during such moves $m_1 = m_3 = m_d$, $m_2 = 1 - m_1$, $m_1^* = m_3^* = m_d^*$, and $m_2^* = 1$. Applying this to $j_1$, we have

$$(\ell + 1)A_3 + (-2 + m_1 + 2m_d(1 - m_1) - (A_3 - 2)(1 - m_d) + 1)$$
$$= (\ell_1 + 1)A_3 - (m_1^* + m_1^*(1 - m_2^*) + 2m_d^*(1 - m_1^*))$$
$$\implies \ell - \ell_1 + m_d - 1 = -\frac{1 - m_1 + m_1^*}{A_3}.$$

Since $\ell - \ell_1 + m_d - 1 \in \mathbb{Z}$ and $0 \leq 1 - m_1 + m_1^* \leq 2$ by construction, it follows that $m_1^* - m_1 = -1$, meaning $m_1 = 1$ and $m_1^* = 0$. Applying this to $j_3$, we get

$$p_1 + s_1 A_4 + m_3 + 2x_3 = s_1^* + m_1^* m_3^*$$
$$\implies s_1 - s_1^* = -1.$$

Consequently, in $j_2$ we obtain

$$m_2 + \gamma + 4x_2 + 2m_1(1 - m_2) + (s_1 + 2m_d(1 - m_1))$$
$$= 2 + m_1^* + (s_1^* + m_d^*(1 - m_1^*) + m_d^*) + \gamma_1 + m_2^*(1 - m_1^*) + 4x_2^*$$

$$\implies x_2 - x_2^* = \frac{2 + (\gamma_1 - \gamma)}{4}.$$

Given $x_2 - x_2^* \in \mathbb{Z}$ and $1 \leq 2 + (\gamma_1 - \gamma) \leq 3$ by construction, it follows that for every $\gamma, \gamma_1 \in \{0, 1\}$, $\frac{2+(\gamma_1-\gamma)}{4} \notin \mathbb{Z}$, giving us a contradiction as the above equality is not possible.

<u>Major Case 2</u>: Focusing on moves along $j_2$, we see that only $C_2$'s $j_2$ move with $(m_1^*, m_2^*, m_d^*) = (0, 0, 1)$ serves as an inverse $j_2$ move to those of $C_1$, and during such moves $m_1 = m_2 = m_d$. Applying this to $j_1$, we find

$$(\ell + 1)A_3 + (-2 + m_1 + 2m_d(1 - m_1) - (A_3 - 2)(1 - m_d))$$
$$= (\ell_1 + 1)A_3 - (m_1^* + m_1^*(1 - m_2^*) + 2m_d^*(1 - m_1^*))$$
$$\implies \ell - \ell_1 + m_d - 1 = \frac{m_1 - 2}{A_3}.$$

Following from $\ell - \ell_1 + m_d - 1 \in \mathbb{Z}$ and $-2 \leq m_1 - 2 \leq -1$ by construction, it follows that for every $m_1 \in \{0, 1\}$, $\frac{m_1-2}{A_3} \notin \mathbb{Z}$, giving us a contradiction as the above equality is not possible.

<u>Major Case 3</u>: Focusing on moves along $j_3$, observe that during such moves $m_1 = m_2$, $m_d = 1 - m_1$, and $(m_1^*, m_2^*, m_d^*) = (1, 1 - m_2, 1 - m_d)$. Applying this to $j_1$, we see

$$(\ell + 1)A_3 + (-2 + m_1 + 2m_d(1 - m_1) - (A_3 - 2)(1 - m_d))$$
$$= (\ell_1 + 1)A_3 - (m_1^* + m_1^*(1 - m_2^*) + 2m_d^*(1 - m_1^*))$$
$$\implies (\ell - \ell_1 + m_d - 1)\frac{A_3}{2} + m_1 = -\frac{1}{2}.$$

Since $(\ell - \ell_1 + m_d - 1)\frac{A_3}{2} + m_1 \in \mathbb{Z}$ by construction and $-\frac{1}{2} \notin \mathbb{Z}$, we have a contradiction as the above equality is not possible.

<u>Major Case 4</u>: Focusing on moves along $j_k$ for $4 \leq k \leq d$ when $d \geq 4$, we see that during such moves $m_1 = m_2$, $m_d = 1 - m_1$, and $(m_1^*, m_2^*, m_d^*) = (1 - m_1, 1 - m_2, 1 - m_d)$. Applying this to $j_1$, we get

$$(\ell + 1)A_3 + (-2 + m_1 + 2m_d(1 - m_1) - (A_3 - 2)(1 - m_d))$$
$$= (\ell_1 + 1)A_3 - (m_1^* + m_1^*(1 - m_2^*) + 2m_d^*(1 - m_1^*))$$
$$\implies (\ell - \ell_1 + m_d - 1)\frac{A_3}{2} + m_d^* = -\frac{1}{2}.$$

Given $(\ell - \ell_1 + m_d - 1)\frac{A_3}{2} + m_d^* \in \mathbb{Z}$ by construction and $-\frac{1}{2} \notin \mathbb{Z}$, we obtain a contradiction as the above equality is not possible.

**Set Comparison Case 10:** We will consider whether edges are shared from moves with the opposite orientation via set 4.

<u>Major Case 1</u>: Focusing on moves along $j_1$, note that there are no inverse $j_1$ moves to consider as all $j_1$ moves are defined with the same orientation with respect to each other.

<u>Major Case 2</u>: Focusing on moves along $j_2$, observe that during such moves $m_1 = m_2 = 0$ and $(m_1^*, m_2^*, m_d^*) = (0, 0, 1 - m_d)$. Noting that $2m_d - 1 = (-1)^{m_d+1}$, in $j_1$ we find

$$(\ell+1)A_3 - (m_1 + m_1(1-m_2) + 2m_d(1-m_1)) = (\ell_1+1)A_3 - (m_1^* + m_1^*(1-m_2^*) + 2m_d^*(1-m_1^*))$$
$$\implies \ell - \ell_1 = \frac{2(-1)^{m_d+1}}{A_3}.$$

Since $\ell - \ell_1 \in \mathbb{Z}$ by construction and for every $m_d \in \{0,1\}$, $\frac{2(-1)^{m_d+1}}{A_3} \notin \mathbb{Z}$, we have a contradiction as the above equality is not possible.

Major Case 3: Focusing on moves along $j_3$, observe that during such moves $m_1 = 1$, $m_2 = 1 - m_d$, and $(m_1^*, m_2^*, m_d^*) = (1, 1-m_2, 1-m_d)$. Thus, in $j_1$ we have

$$(\ell+1)A_3 - (m_1 + m_1(1-m_2) + 2m_d(1-m_1)) = (\ell_1+1)A_3 - (m_1^* + m_1^*(1-m_2^*) + 2m_d^*(1-m_1^*))$$
$$\implies (\ell - \ell_1)\frac{A_3}{2} + m_2 = \frac{1}{2}.$$

Following from $(\ell-\ell_1)\frac{A_3}{2} + m_2 \in \mathbb{Z}$ by construction and $\frac{1}{2} \notin \mathbb{Z}$, we see we have a contradiction as the above equality is not possible.

Major Case 4: Focusing on moves along $j_k$ for $4 \leq k \leq d$ when $d \geq 4$, we see that during all such moves $m_1 = m_2$, $m_d = 1 - m_1$, and $(m_1^*, m_2^*, m_d^*) = (1-m_1, 1-m_2, 1-m_d)$. Applying this to $j_1$, we get

$$(\ell+1)A_3 - (m_1 + m_1(1-m_2) + 2m_d(1-m_1)) = (\ell_1+1)A_3 - (m_1^* + m_1^*(1-m_2^*) + 2m_d^*(1-m_1^*))$$
$$\implies (\ell - \ell_1)\frac{A_3}{2} - (m_1 + (m_d - m_d^*)) = -\frac{1}{2}.$$

Given $(\ell-\ell_1)\frac{A_3}{2} - (m_1 + (m_d - m_d^*)) \in \mathbb{Z}$ by construction and $-\frac{1}{2} \notin \mathbb{Z}$, we have a contradiction as the above equality is not possible.

Thus, no edges are shared from moves with the opposite orientation for $\alpha_1 < \alpha^* < \alpha_2$ when $d \geq 3$.

$\boldsymbol{\alpha^*}$**–Case 3:** Let $\alpha^* = \alpha_2$. Then, $A_1 = 1$ and $A_2 = 1$, meaning sets 1 and 3 are active. Note that $A_3 \geq 4$ is even, $t = 0 = t_1$, $p_1 = 0 = p_1^*$, and $A_4 = 1$.

**Set Comparison Case 1:** We will consider whether edges are shared from moves with the opposite orientation via set 1.

Major Case 1: Focusing on moves along $j_1$, observe that all $j_1$ moves are defined with the same orientation with respect to each other and so there are no inverse $j_1$ moves to consider.

Major Case 2: Focusing on moves along $j_2$, observe that during such moves $m_1 = m_d$ and $(m_1^*, m_d^*) = (1-m_1, 1-m_d)$. Hence, in $j_1$ we have

$$\ell A_3 + (m_1 + 2x_1 + 2m_d(1-m_1)) = \ell_1 A_3 + (m_1^* + 2x_1^* + 2m_d^*(1-m_1^*))$$
$$\implies (\ell - \ell_1)\frac{A_3}{2} + m_1 + (x_1 - x_1^*) = \frac{1}{2}$$

Given $(\ell - \ell_1)\frac{A_3}{2} + m_1 + (x_1 - x_1^*) \in \mathbb{Z}$ and $\frac{1}{2} \notin \mathbb{Z}$, we have a contradiction as the above equality is not possible.

**Major Case 3:** Focusing on moves along $j_k$ for $3 \leq k \leq d$, note that it is the case $m_d = 1 - m_1$ and $(m_1^*, m_d^*) = (1 - m_1, 1 - m_d)$ during all such moves. Applying this to $j_1$, we get

$$\ell A_3 + (m_1 + 2x_1 + 2m_d(1 - m_1)) = \ell_1 A_3 + (m_1^* + 2x_1^* + 2m_d^*(1 - m_1^*))$$
$$\implies (\ell - \ell_1)\frac{A_3}{2} + m_1 + (x_1 - x_1^*) + (2m_d - 1) = \frac{1}{2}.$$

Since $(\ell - \ell_1)\frac{A_3}{2} + m_1 + (x_1 - x_1^*) + (2m_d - 1) \in \mathbb{Z}$ and $\frac{1}{2} \notin \mathbb{Z}$, we have a contradiction as the above equality is not possible.

**Set Comparison Case 2:** We will consider whether edges are shared from moves with the opposite orientation via sets 1 and 3.

**Major Case 1:** Focusing on moves along $j_1$, we see that there are no inverse $j_1$ moves to consider as all $j_1$ moves are defined with the same orientation with respect to each other.

**Major Case 2:** Focusing on moves along $j_2$, observe that all of $C_2$'s $j_2$ moves serve as moves with the opposite orientation only to $C_1$'s $j_2$ move with $(m_1, m_d) = (1, 1)$, and during these moves $m_1^* = m_d^*$. Applying this to $j_1$, we get

$$\ell A_3 + (m_1 + 2x_1 + 2m_d(1 - m_1)) = (\ell_1 + 1)A_3 + (-2 + m_1^* + 2m_d^*(1 - m_1^*))$$
$$\implies \ell - \ell_1 - 1 = -\frac{3 + 2x_1 - m_1^*}{A_3}.$$

Given $\ell - \ell_1 - 1 \in \mathbb{Z}$ and $2 \leq 3 + 2x_1 - m_1^* \leq A_3 - 1$ by construction, it follows that for every $0 \leq x_1 \leq \frac{A_3}{2} - 2$ and $m_1^* \in \{0, 1\}$, $-\frac{3 + 2x_1 - m_1^*}{A_3} \notin \mathbb{Z}$. This gives us a contradiction as the above equality is not possible.

**Major Case 3:** Focusing on moves along $j_k$ for $3 \leq k \leq d$, observe that during all such moves it is the case $m_d = 1 - m_1$ and $(m_1^*, m_d^*) = (1 - m_1, 1 - m_d)$. Applying this to $j_1$, we have

$$\ell A_3 + (m_1 + 2x_1 + 2m_d(1 - m_1)) = (\ell_1 + 1)A_3 + (-2 + m_1^* + 2m_d^*(1 - m_1^*))$$
$$\implies (\ell - \ell_1 - 1)\frac{A_3}{2} + (x_1 + m_d) = -\frac{1}{2}.$$

Since $(\ell - \ell_1 - 1)\frac{A_3}{2} + (x_1 + m_d) \in \mathbb{Z}$ by construction and $-\frac{1}{2} \notin \mathbb{Z}$, we have a contradiction as the above equality is not possible.

**Set Comparison Case 3:** We will consider whether edges are shared from moves with the opposite orientation via set 3.

**Major Case 1:** Focusing on moves along $j_1$, observe that there are no inverse $j_1$ moves to consider as all $j_1$ moves are defined with the same orientation with respect to each other.

<u>Major Case 2</u>: Focusing on moves along $j_2$, we see that there are no inverse $j_2$ moves to consider as all $j_2$ moves are defined with the same orientation with respect to each other.

<u>Major Case 3</u>: Focusing on moves along $j_k$ for $3 \leq k \leq d$, observe that during all such moves $m_d = 1 - m_1$ and $(m_1^*, m_d^*) = (1 - m_1, 1 - m_d)$. Applying this to $j_1$, we find

$$(\ell+1)A_3 + (-2 + m_1 + 2m_d(1-m_1)) = (\ell_1 + 1)A_3 + (-2 + m_1^* + 2m_d^*(1-m_1^*))$$
$$\implies (\ell - \ell_1)\frac{A_3}{2} + m_d = \frac{1}{2}.$$

Since $(\ell - \ell_1)\frac{A_3}{2} + m_d \in \mathbb{Z}$ by construction and $\frac{1}{2} \notin \mathbb{Z}$, we have a contradiction as the above equality is not possible.

Consequently, no edges are shared from moves with the opposite orientation for $\alpha^* = \alpha_2$ when $d \geq 3$.

**$\alpha^*$−Case 4:** Let $\alpha_2 < \alpha^* \leq \alpha_d$. Then, $A_1 = 1$, $A_2 = 1$, and $\eta_{\alpha_2} = 1$, meaning only sets 1 and 3 are active. Note that $t = 0 = t_1$, and for $4 \leq k \leq d+1$, $A_k = 1$ for $\alpha_2 < \alpha^* \leq \alpha_{k-2}$ and $A_k \geq 2$ for $\alpha_{k-2} < \alpha^* \leq \alpha_d$.

Before we proceed, we wish to present the three general states corresponding to a move along $j_k$ for $3 \leq k \leq d$ as defined by set 3:

1. **Stair-casing Moves:** For moves along $j_k$ for $3 \leq k \leq d$, we see that stair-casing moves take place when all of the following are true: $\alpha_{k-1} < \alpha^* \leq \alpha_d$, $R_1 = 0$, and $R_{k-1} = 1$. Note that $R_1 = 0$ and $R_{k-1} = 1$ together imply that $x_z = (1-r_2)(A_{z+1} - 1)$ for all $3 \leq z \leq k-1$ and that there exists at least one $k \leq w \leq d$ such that $x_w \neq (1-r_2)(A_{w+1} - 1)$. Since $x_w \neq (1-r_2)(A_{w+1} - 1)$ is to hold for all $r_2 \in \{0, 1\}$, letting $k \leq w \leq d$ be the smallest integer, by the Well-Ordering Principle, such that the above is satisfied, it must be the case $A_{w+1} \geq 2$, which then implies $\alpha_{w-1} < \alpha^* \leq \alpha_d$.

2. **Normal Moves:** For moves along $j_k$ for $4 \leq k \leq d$, observe that normal moves take place when either $\alpha_2 < \alpha^* \leq \alpha_{k-1}$ or it is the case that all of the following are true: $\alpha_{k-1} < \alpha^* \leq \alpha_d$, $R_{k-1} = 0$ and $R_1 = 0$. It follows that $R_{k-1} = 0$ and $R_1 = 0$ together imply that there exists at least one $3 \leq w \leq k-1$ such that $x_w \neq (1-r_2)(A_{w+1} - 1)$. Note that we exclude $j_3$ moves as $R_2 = 1$ during all $j_3$ moves, even when $R_1 = 0$.

3. **Column Transitions:** For moves along $j_k$ for $3 \leq k \leq d$, it follows that column transitions take place when it is the case that both $\alpha_{k-1} < \alpha^* \leq \alpha_d$ and $R_1 = 1$. From $R_1 = 1$, we see that $x_z = (1-r_2)(A_{z+1} - 1)$ for all $3 \leq z \leq d$.

**Set Comparison Case 1:** Considering whether edges are shared from moves with the opposite orientation via set 1, we see that the same arguments from set comparison case 1 of $\alpha^*$−case 3 with $d \geq 3$ hold here as all arguments there make use only of the $j_1$ component that remains the same for all $\alpha_2 \leq \alpha^* \leq \alpha_d$. By same, we mean that no additional terms are introduced from the definition of $j_1$ as a consequence of increasing $\alpha^*$.

**Set Comparison Case 2:** We will consider whether edges are shared from moves with the opposite orientation via sets 1 and 3.

**Major Case 1:** Focusing on moves along $j_1$, observe that all $j_1$ moves defined have the same orientation with respect to each other and so there are no inverse $j_1$ moves to consider.

**Major Case 2:** Focusing on moves along $j_2$, we see that which $j_2$ moves in $C_2$ serve as moves with the opposite orientation to those of $C_1$ depend on $r_2^*$ and $x_k^*$ for $3 \leq k \leq d$. To case on $r_2^*$ and $x_k^*$ for $3 \leq k \leq d$, we proceed by casing on $R_1^*$:

**Case 1:** Let $R_1^* = 0$. In this case, we see that during all $j_2$ moves $m_1 = m_d$ and $(m_1^*, m_d^*) = (1 - m_1, 1 - m_d)$. Applying this to $j_1$, we get

$$\ell A_3 + (m_1 + 2x_1 + 2m_d(1 - m_1)) = (\ell_1 + 1)A_3 + (-2 + m_1^* + 2m_d^*(1 - m_1^*))$$

$$\implies (\ell - \ell_1 - 1)\frac{A_3}{2} + m_1 + (1 + x_1) = \frac{1}{2}.$$

Given $(\ell - \ell_1 - 1)\frac{A_3}{2} + m_1 + (1 + x_1) \in \mathbb{Z}$ by construction and $\frac{1}{2} \notin \mathbb{Z}$, we have a contradiction as the above equality is not possible.

**Case 2:** Let $R_1^* = 1$. Then, we see that all of $C_2$'s $j_2$ moves serve as moves with the opposite orientation only to that of $C_1$ with $(m_1, m_d) = (1, 1)$, and during such moves $m_1^* = m_d^*$. Applying this to $j_1$, we find

$$\ell A_3 + (m_1 + 2x_1 + 2m_d(1 - m_1)) = (\ell_1 + 1)A_3 + (-2 + m_1^* + 2m_d^*(1 - m_1^*))$$

$$\implies \ell - \ell_1 - 1 = -\frac{3 + 2x_1 - m_1^*}{A_3}.$$

Since $\ell - \ell_1 - 1 \in \mathbb{Z}$ and $2 \leq 3 + 2x_1 - m_1^* \leq A_3 - 1$ by construction, it follows that for every $0 \leq x_1 \leq \frac{A_3}{2} - 2$ and $m_1^* \in \{0, 1\}$, $-\frac{3 + 2x_1 - m_1^*}{A_3} \notin \mathbb{Z}$, giving us a contradiction as the above equality is not possible.

**Major Case 3:** Focusing on moves along $j_k$ for $3 \leq k \leq d$, we will consider the three general states that define $C_2$ in this case as each state has a different group of $j_k$ moves serving as moves with the opposite orientation to the $j_k$ moves of $C_1$:

**Case 1:** Considering stair-casing moves along $j_k$ for $3 \leq k \leq d$, we see that $\alpha_{k-1} < \alpha^* \leq \alpha_d$ and for all $3 \leq z \leq k - 1$, $x_z^* = (1 - r_2^*)(A_{z+1} - 1)$. Then, $(m_1, m_d) = (r_2^*, 1 - r_2^*)$ and $m_d^* = 1 - m_1^*$, giving us in $j_1$

$$\ell A_3 + (m_1 + 2x_1 + 2m_d(1 - m_1)) = (\ell_1 + 1)A_3 + (-2 + m_1^* + 2m_d^*(1 - m_1^*))$$

$$\implies \ell - \ell_1 - 1 = -\frac{3 - r_2^* + 2x_1 - m_d^*}{A_3}.$$

Given $\ell - \ell_1 - 1 \in \mathbb{Z}$ and $1 \leq 3 - r_2^* + 2x_1 - m_d^* \leq A_3 - 1$ by construction, it follows that for every $0 \leq x_1 \leq \frac{A_3}{2} - 2$ and $r_2^*, m_d^* \in \{0, 1\}$, $-\frac{3 - r_2^* + 2x_1 - m_d^*}{A_3} \notin \mathbb{Z}$. Thus, we have a contradiction as the above equality is not possible.

**Case 2:** Considering normal moves along $j_k$ for $4 \leq k \leq d$, we know that either $\alpha_2 < \alpha^* \leq \alpha_{k-1}$ or $\alpha_{k-1} < \alpha^* \leq \alpha_d$, $R_1^* = 0$, and $R_{k-1}^* = 1$. Observing that during all such moves $m_d = 1 - m_1$ and $(m_1^*, m_d^*) = (1 - m_1, 1 - m_d)$, in $j_1$ we get

$$\ell A_3 + (m_1 + 2x_1 + 2m_d(1 - m_1)) = (\ell_1 + 1)A_3 + (-2 + m_1^* + 2m_d^*(1 - m_1^*))$$
$$\implies (\ell - \ell_1 - 1)\frac{A_3}{2} + m_1 + x_1 + (m_d - m_d^*) = -\frac{1}{2}.$$

Following from $(\ell - \ell_1 - 1)\frac{A_3}{2} + m_1 + x_1 + (m_d - m_d^*) \in \mathbb{Z}$ by construction and $-\frac{1}{2} \notin \mathbb{Z}$, we see that we have a contradiction as the above equality is not possible.

<u>Case 3:</u> Focusing on column transitions along $j_k$ for $3 \leq k \leq d$, observe that it must be the case $\alpha_{k-1} < \alpha^* \leq \alpha_d$ and $R_1^* = 1$. Noting that during all these moves $m_d = 1 - m_1$ and $(m_1^*, m_d^*) = ((1 - r_2^*)(1 - m_1) + r_2^* m_1, (1 - r_2^*)(1 - m_d) + r_2^* m_d)$, in $j_1$ we have

$$\ell A_3 + (m_1 + 2x_1 + 2m_d(1 - m_1)) = (\ell_1 + 1)A_3 + (-2 + m_1^* + 2m_d^*(1 - m_1^*))$$
$$\implies \ell - \ell_1 - 1 = -\frac{2 + 2x_1 + (1 - r_2^*)(2m_d - 1)}{A_3}.$$

Since $\ell - \ell_1 - 1 \in \mathbb{Z}$ and $1 \leq 2 + 2x_1 + (1 - r_2^*)(2m_d - 1) \leq A_3 - 1$ by construction, we see that for every $0 \leq x_1 \leq \frac{A_3}{2} - 2$ and $r_2^*, m_d \in \{0, 1\}$, $-\frac{2+2x_1+(1-r_2^*)(2m_d-1)}{A_3} \notin \mathbb{Z}$. Consequently, we have a contradiction as the above equality is not possible.

**Set Comparison Case 3:** We will consider whether edges are shared from moves with the opposite orientation via set 3.

<u>Major Case 1:</u> Focusing on moves along $j_1$, observe that all $j_1$ moves have the same orientation with respect to each other and hence there are no inverse $j_1$ moves to consider.

<u>Major Case 2:</u> Focusing on moves along $j_2$, we see that which $j_2$ moves in $C_2$ serve as moves with the opposite orientation to those of $C_1$ depend on $R_1$ and $R_1^*$. Hence, we proceed by casing on $R_1$ and $R_1^*$:

<u>Case 1:</u> Let $R_1 = 0$. To determine which of $C_2$'s $j_2$ moves serves as an inverse to those of $C_1$, we further case on $R_1^*$:

<u>Subcase 1:</u> Let $R_1^* = 0$. Then, $m_1 = m_d$ and $(m_1^*, m_d^*) = (1 - m_1, 1 - m_d)$, giving us in $j_1$

$$(\ell + 1)A_3 + (-2 + m_1 + 2m_d(1 - m_1)) = (\ell_1 + 1)A_3 + (-2 + m_1^* + 2m_d^*(1 - m_1^*))$$
$$\implies (\ell - \ell_1)\frac{A_3}{2} + m_1 = \frac{1}{2}.$$

Given $(\ell - \ell_1)\frac{A_3}{2} + m_1 \in \mathbb{Z}$ by construction and $\frac{1}{2} \notin \mathbb{Z}$, we have a contradiction as the above equality is not possible.

<u>Subcase 2:</u> Let $R_1^* = 1$. Then, $(m_1, m_2, \ldots, m_d) = (1, 1, \ldots, 1)$, $m_1^* = \cdots = m_d^*$, and for all $3 \leq z \leq d$, $x_z^* = (1 - r_2^*)(A_{z+1} - 1)$. Applying this to $j_1$, we get

$$(\ell + 1)A_3 + (-2 + m_1 + 2m_d(1 - m_1)) = (\ell_1 + 1)A_3 + (-2 + m_1^* + 2m_d^*(1 - m_1^*))$$
$$\implies \ell - \ell_1 = \frac{m_1^* - 1}{A_3}.$$

Since $\ell - \ell_1 \in \mathbb{Z}$ and $-1 \leq m_1^* - 1 \leq 0$ by construction, it must be the case $m_1^* = 1$ and so $(m_1, \ldots, m_d) = (m_1^*, \ldots, m_d^*)$. Note that since $R_1 = 0$, $C_1$ is stair-casing along $j_2$, meaning $\alpha_{w-1} < \alpha^* \leq \alpha_d$, where $3 \leq w \leq d$ is the smallest integer such that $x_w \neq (1 - r_2)(A_{w+1} - 1)$.

To complete our treatment of this subcase, we will require the following result we will prove using induction:

We will prove that for every $3 \leq w \leq d$ such that $w$ is the smallest integer with $x_w \neq (1-r_2)(A_{w+1}-1)$, it is the case that for all $3 \leq z \leq w-1$, $x_z = (1-r_2)(A_{z+1}-1)$. We proceed by induction on $w$:

<u>Base Case 1</u>: Let $w = 3$. Then, the statement is vacuously true as there does not exist an integer $z$ such that $3 \leq z \leq 2$.

<u>Base Case 2</u>: Let $w = 4$. Then, $z = 3$ and so our assumptions for this subcase give us in $j_3$ that

$$p_1 + 2A_4 s_1 + (1 - 2r_2\Psi_0)m_3 + 2x_3 + 2\eta_{\alpha_3}(1 - R_1)((1 - r_2)X_0 - r_2\Psi_0)(1 - m_3)m_d$$
$$= p_1^* + 2A_4 s_1^* + (1 - 2r_2^*\Psi_0^*)m_3^* + 2x_3^* + 2\eta_{\alpha_3}(1 - R_1^*)((1 - r_2^*)X_0^* - r_2^*\Psi_0^*)(1 - m_3^*)m_d^*$$
$$\implies A_4(s_1 - s_1^*) + (r_2^* - r_2) + (x_3 - x_3^*) = \frac{p_1^* - p_1}{2}.$$

Since $A_4(s_1 - s_1^*) + (r_2^* - r_2) + (x_3 - x_3^*) \in \mathbb{Z}$ and $|p_1^* - p_1| \leq 1$ by construction, it must be the case $p_1 = p_1^*$. Recalling that $x_3^* = (1 - r_2^*)(A_4 - 1)$ since $R_1^* = 1$, we are left with

$$(s_1 - s_1^*) - (1 - r_2^*) = -\frac{x_3 + (1 - r_2)}{A_4}.$$

Following from $(s_1 - s_1^*) - (1 - r_2^*) \in \mathbb{Z}$ and $1 - r_2 \leq x_3 + (1 - r_2) \leq A_4 - r_2$, we see that we must have $x_3 = (1 - r_2)(A_4 - 1)$.

<u>Induction Step</u>: Let $3 \leq w \leq d$ and $d \geq 3$ be such that $3 \leq w + 1 \leq d$ so that the following makes sense to consider. Otherwise, we will have shown that $x_z = (1-r_2)(A_{z+1}-1)$ for all $3 \leq z \leq d$ when $w = d$, which is not possible as $R_1 = 0$ by assumption. Further, assume that for all $3 \leq z \leq y - 1$ with $3 \leq y \leq w$, $x_z = (1 - r_2)(A_{z+1} - 1)$. We will show that our statement holds for $3 \leq w+1 \leq d$ by showing that $x_w = (1 - r_2)(A_{w+1} - 1)$ since our induction hypothesis grants us the above for all $3 \leq z \leq w - 1$. Note that our assumptions on $x_z$ and $x_{z^*}^*$, with $x_{z^*}^* = (1 - r_2^*)(A_{z^*+1} - 1)$ for all $3 \leq z^* \leq d$ since $R_1^* = 1$, imply $(1-r_2)X_{w-3} = (1-r_2)$, $r_2\Psi_{w-3} = r_2$, $(1-r_2^*)X_{w-3}^* = (1-r_2^*)$, and $r_2^*\Psi_{w-3}^* = r_2^*$. The corresponding assumptions and results we have so far, along with the observations $\eta_{w+1} = 1$ since $w + 1 \leq d$ and $\eta_{\alpha_w} = \eta_{\alpha_{w-1}} = 1$ since $\alpha^* > \alpha_w > \alpha_{w-1}$, yield in $j_w$

$$p_{w-2} + 2A_{w+1}s_{w-2} + (1 - 2r_2\Psi_{w-3})m_w + 2x_w + 2\eta_{w+1}\eta_{\alpha_w}(1 - R_1)((1 - r_2)X_{w-3} - r_2\Psi_{w-3})(1 - m_w)m_d$$
$$= p_{w-2}^* + 2A_{w+1}s_{w-2}^* + (1-2r_2^*\Psi_{w-3}^*)m_w^* + 2x_w^* + 2\eta_{w+1}\eta_{\alpha_w}(1-R_1^*)((1-r_2^*)X_{w-3}^* - r_2^*\Psi_{w-3}^*)(1-m_w^*)m_d^*$$
$$\implies A_{w+1}(s_{w-2} - s_{w-2}^*) + (r_2^* - r_2) + (x_w - x_w^*) = \frac{p_{w-2}^* - p_{w-2}}{2}.$$

Since $A_{w+1}(s_{w-2} - s_{w-2}^*) + (r_2^* - r_2) + (x_w - x_w^*) \in \mathbb{Z}$ and $|p_{w-2}^* - p_{w-2}| \leq 1$, it must be the case $p_{w-2} = p_{w-2}^*$. Hence, observing that $x_w^* = (1 - r_2^*)(A_{w+1} - 1)$ since $R_1^* = 1$, we are left with

$$(s_{w-2} - s_{w-2}^*) - (1 - r_2^*) = -\frac{x_w + (1 - r_2)}{A_{w+1}}.$$

Given $(s_{w-2} - s^*_{w-2}) - (1 - r^*_2) \in \mathbb{Z}$ and $1 - r_2 \leq x_w + (1 - r_2) \leq A_{w+1} - r_2$ by construction, it follows that $x_w = (1 - r_2)(A_{w+1} - 1)$.

By the Principle of Strong Mathematical Induction, it follows that for every $3 \leq w \leq d$, it is the case that $x_z = (1 - r_2)(A_{z+1} - 1)$ for every $3 \leq z \leq w - 1$.

From this, we know that $x_z = (1 - r_2)(A_{z+1} - 1)$ for all $3 \leq z \leq w - 1$. Given the definitions of $j_d$ and $j_k$ for $3 \leq k \leq d - 1$ differ slightly, we case on whether $w = d$:

Sub-subcase 1: Let $w = d$. Then, $x_d \neq (1 - r_2)(A_{d+1} - 1)$, and for all $3 \leq k \leq d - 1$, $x_k = (1 - r_2)(A_{k+1} - 1)$. Note that $\eta_{\alpha_{d-1}} = 1$ since $\alpha_{d-1} < \alpha^* \leq \alpha_d$ as $A_{d+1} \geq 2$ for the assumption $x_d \neq (1 - r_2)(A_{d+1} - 1)$ to hold for all $r_2 \in \{0, 1\}$. Using what we have so far along with the observations $r_2 \Psi_{d-3} = r_2$ and $r^*_2 \Psi^*_{d-3} = r^*_2$ granted by our assumptions $x_z = (1 - r_2)(A_{z+1} - 1)$ and $x^*_z = (1 - r^*_2)(A_{z+1} - 1)$ for all $3 \leq z \leq w - 1$, in $j_d$ we obtain

$$p_{d-2} + 2A_{d+1}s_{d-2} + (1 - 2r_2\Psi_{d-3})m_d + 2x_d = p^*_{d-2} + 2A_{d+1}s^*_{d-2} + (1 - 2r^*_2\Psi^*_{d-3})m^*_d + 2x^*_d$$

$$\implies A_{d+1}(s_{d-2} - s^*_{d-2}) + (r^*_2 - r_2) + (x_d - x^*_d) = \frac{p^*_{d-2} - p_{d-2}}{2}.$$

Since $A_{d+1}(s_{d-2} - s^*_{d-2}) + (r^*_2 - r_2) + (x_d - x^*_d) \in \mathbb{Z}$ and $|p^*_{d-2} - p_{d-2}| \leq 1$ by construction, it must be the case $p_{d-2} = p^*_{d-2}$. Hence, we now have

$$(s_{d-2} - s^*_{d-2}) - (1 - r^*_2) = -\frac{x_d + (1 - r_2)}{A_{d+1}}.$$

Given $(s_{d-2} - s^*_{d-2}) - (1 - r^*_2) \in \mathbb{Z}$ by construction and $x_d + (1 - r_2) \neq (1 - r_2)A_{d+1}$ by assumption, we see that $-\frac{x_d + (1 - r_2)}{A_{d+1}} \notin \mathbb{Z}$. Consequently, we have a contradiction as the above equality is not possible. It now follows that this scenario is impossible as the smallest $w$ is also the largest it could be, in this case, and so there does not exist such a $w$.

Sub-subcase 2: Let $3 \leq w < d$. Then, $x_z = (1 - r_2)(A_{z+1} - 1)$ for all $3 \leq z \leq w - 1$ and $x_w \neq (1 - r_2)(A_{w+1} - 1)$, meaning $\alpha_{w-1} < \alpha^* \leq \alpha_d$ and $\eta_{\alpha_{w-1}} = 1$ by the same reasoning as before. Since $R^*_1 = 1$, we know $x^*_k = (1 - r^*_2)(A_{k+1} - 1)$ for all $3 \leq k \leq d$. Now, observe that our assumptions on $x_z$ and $x^*_z$ for $3 \leq z \leq w - 1$ imply $(1 - r_2)X_{w-3} = (1 - r_2)$, $r_2\Psi_{w-3} = r_2$, $(1 - r^*_2)X^*_{w-3} = (1 - r^*_2)$, and $r^*_2\Psi^*_{w-3} = r^*_2$. Recalling that $R_1 = 0$ and $R^*_1 = 1$, all of the above, along with our previous results for this subcase, in $j_w$ yield

$$p_{w-2} + 2A_{w+1}s_{w-2} + (1 - 2r_2\Psi_{w-3})m_w + 2x_w + 2\eta_{w+1}\eta_{\alpha_w}(1 - R_1)((1 - r_2)X_{w-3} - r_2\Psi_{w-3})(1 - m_w)m_d$$
$$= p^*_{w-2} + 2A_{w+1}s^*_{w-2} + (1 - 2r^*_2\Psi^*_{w-3})m^*_w + 2x^*_w + 2\eta_{w+1}\eta_{\alpha_w}(1 - R^*_1)((1 - r^*_2)X^*_{w-3} - r^*_2\Psi^*_{w-3})(1 - m^*_w)m^*_d$$

$$\implies A_{w+1}(s_{w-2} - s^*_{w-2}) + (r^*_2 - r_2) + (x_w - x^*_w) = \frac{p^*_{w-2} - p_{w-2}}{2}.$$

Since $A_{w+1}(s_{w-2} - s^*_{w-2}) + (r^*_2 - r_2) + (x_w - x^*_w) \in \mathbb{Z}$ and $|p^*_{w-2} - p_{w-2}| \leq 1$ by construction, we must have $p_{w-2} = p^*_{w-2}$. We now get

$$(s_{w-2} - s^*_{w-2}) - (1 - r^*_2) = -\frac{x_w + (1 - r_2)}{A_{w+1}}.$$

Following from $(s_{w-2} - s^*_{w-2}) - (1 - r^*_2) \in \mathbb{Z}$ by construction and $x_w + (1 - r_2) \neq (1 - r_2)\mathrm{A}_{w+1}$ by assumption, we have $-\frac{x_w + (1 - r_2)}{\mathrm{A}_{w+1}} \notin \mathbb{Z}$, giving us a contradiction as the above equality is not possible.

Further, since $w$ was the arbitrary smallest and all components $j_k$ for $3 \leq k \leq d-1$ are defined as above, the above holds for all such $w$. So there is no such $3 \leq w \leq d-1$. Now, if $w = d$, then by sub-subcase 1 we conclude $w \neq d$ as we get a contradiction. Having exhausted all possible $3 \leq w \leq d$, we conclude that there does not exist a smallest integer $3 \leq w \leq d$ such that $x_w \neq (1 - r_2)(\mathrm{A}_{w+1} - 1)$. This means that the scenario we have considered never occurs.

<u>Case 2:</u> Let $R_1 = 1$. Since the argument made in subcase 2 of case 1 is symmetrical, the case with both $R_1 = 1$ and $R^*_1 = 0$ is subject to the same conclusion as in subcase 2 of case 1.

Letting $R^*_1 = 1$, we see that all $j_2$ moves in $C_1$ and $C_2$ have the same orientation with respect to each other and so there are no inverse $j_2$ moves to consider.

<u>Major Case 3:</u> Focusing on moves along $j_k$ for $3 \leq k \leq d$, we proceed by casing on the three general states defined by set 3 along $j_k$:

<u>Case 1:</u> Let $R_1 = 0$ and $R_{k-1} = 1$ for $3 \leq k \leq d$. This means that $C_1$ is stair-casing, so $x_z = (1 - r_2)(\mathrm{A}_{z+1} - 1)$ for all $3 \leq z \leq k-1$ and there exists at least one integer $k \leq y \leq d$ such that $x_y \neq (1 - r_2)(\mathrm{A}_{y+1} - 1)$. Note that the latter deduction implies there exists a smallest integer $k \leq w \leq d$ such that $x_w \neq (1-r_2)(\mathrm{A}_{w+1}-1)$. However, since this is to hold for all $r_2 \in \{0, 1\}$, this is possible if and only if $\mathrm{A}_{w+1} \geq 2$, which by its definition is if and only if $\alpha_{w-1} < \alpha^* \leq \alpha_d$. Hence, for this case, assume $\alpha_{w-1} < \alpha^* \leq \alpha_d$. In what follows, we case on $R^*_1$ and $R^*_{k-1}$ as appropriate:

<u>Subcase 1:</u> Let $R^*_1 = 0$ and $R^*_{k-1} = 1$ for $3 \leq k \leq d$. Then, $C_2$ is stair-casing, meaning $x^*_z = (1 - r^*_2)(\mathrm{A}_{z+1} - 1)$ for all $3 \leq z \leq k-1$ and there exists at least one integer $k \leq y^* \leq d$ such that $x^*_{y^*} \neq (1-r^*_2)(\mathrm{A}_{y^*+1}-1)$. If $r_2 = r^*_2$, then $C_2$ has no inverse $j_k$ moves with respect to $C_1$'s $j_k$ moves.

So suppose $r_2 \neq r^*_2$. From this, we see that $m_1 = m_2$, $m_d = 1 - m_1$, $m^*_1 = m^*_2$, and $m^*_d = 1 - m^*_1$ during all such $j_k$ moves. Applying this to $j_1$, we get

$$(\ell + 1)\mathrm{A}_3 + (-2 + m_1 + 2m_d(1 - m_1)) = (\ell_1 + 1)\mathrm{A}_3 + (-2 + m^*_1 + 2m^*_d(1 - m^*_1))$$
$$\implies \ell - \ell_1 = \frac{m^*_d - m_d}{\mathrm{A}_3}.$$

Since $\ell - \ell_1 \in \mathbb{Z}$ and $|m^*_d - m_d| \leq 1$ with $\mathrm{A}_3 \geq 4$ by construction, we must have $m_d = m^*_d$, meaning $(m_1, m_2, m_d) = (m^*_1, m^*_2, m^*_d)$. Consequently, in $j_2$ we see

$$m_2 + \gamma + 2x_2 + 2m_1(1 - m_2) + 2((1 - m_1)m_d - (1 - R_1)((1 - m_2)m_1 + m_d(1 - m_1)))$$
$$= m^*_2 + \gamma_1 + 2x^*_2 + 2m^*_1(1 - m^*_2) + 2((1 - m^*_1)m^*_d - (1 - R^*_1)((1 - m^*_2)m^*_1 + m^*_d(1 - m^*_1)))$$
$$\implies x_2 - x^*_2 = \frac{\gamma_1 - \gamma}{2}.$$

Given $x_2 - x^*_2 \in \mathbb{Z}$ and $|\gamma_1 - \gamma| \leq 1$ by construction, it must be the case $\gamma = \gamma_1$ and so $x_2 = x^*_2$. This last deduction implies $r_2 = r^*_2$, a contradiction to our assumption that $r_2 \neq r^*_2$.

<u>Subcase 2:</u> Since $\alpha_{w-1} < \alpha^* \leq \alpha_d$, let $R_1^* = 0$ and $R_{k-1}^* = 0$ for all $4 \leq k \leq d$. This means that $C_2$ is being defined by normal moves, so there exists at least one $3 \leq w^* \leq k-1$ such that $x_{w^*}^* \neq (1-r_2^*)(A_{w^*+1} - 1)$ during $j_k$ moves for all $4 \leq k \leq d$. Note that $j_3$ is not subject to normal moves as it is always the case $R_2^* = 1$ when $R_1^* = 0$.

Observing that during all such moves $m_1 = m_2$, $m_d = 1 - m_1$, and $(m_1^*, \ldots, m_{k-1}^*, m_d^*) = (r_2, \ldots, r_2, 1 - r_2)$, in $j_1$ we have

$$(\ell + 1)A_3 + (-2 + m_1 + 2m_d(1 - m_1)) = (\ell_1 + 1)A_3 + (-2 + m_1^* + 2m_d^*(1 - m_1^*))$$

$$\implies \ell - \ell_1 = \frac{m_d^* - m_d}{A_3}.$$

Since $\ell - \ell_1 \in \mathbb{Z}$ and $|m_d^* - m_d| \leq 1$ by construction, we must have $m_d = m_d^*$ and so $(m_1, \ldots, m_{k-1}, m_d) = (m_1^*, \ldots, m_{k-1}^*, m_d^*)$. Consequently, in $j_2$ we find

$$m_2 + \gamma + 2x_2 + 2m_1(1 - m_2) + 2((1 - m_1)m_d - (1 - R_1)((1 - m_2)m_1 + m_d(1 - m_1)))$$
$$= m_2^* + \gamma_1 + 2x_2^* + 2m_1^*(1 - m_2^*) + 2((1 - m_1^*)m_d^* - (1 - R_1^*)((1 - m_2^*)m_1^* + m_d^*(1 - m_1^*)))$$

$$\implies x_2 - x_2^* = \frac{\gamma_1 - \gamma}{2}.$$

Given $x_2 - x_2^* \in \mathbb{Z}$ and $|\gamma_1 - \gamma| \leq 1$ by construction, it must be the case $\gamma = \gamma_1$ and so $x_2 = x_2^*$. From this last deduction, it follows that $r_2 = r_2^*$.

To complete our treatment of this subcase, we will require the following result we will prove by way of induction:

We aim to show that for every smallest $3 \leq w^* \leq k-1$ such that $x_{w^*}^* \neq (1 - r_2^*)(A_{w^*+1} - 1)$ during a given $j_k$ move for $4 \leq k \leq d$, it is the case that for every $3 \leq z \leq w^* - 1$, $x_z^* = (1 - r_2^*)(A_{z+1} - 1)$.

<u>Base Case 1:</u> Let $w^* = 3$. Then, the statement is vacuously true as there does not exist an integer $z$ such that $3 \leq z \leq 2$.

<u>Base Case 2:</u> Let $w^* = 4$. Then, $z = 3$. Using our results and assumptions from this subcase along with the observation $\eta_{\alpha_3} = 1$ since $\alpha^* > \alpha_{w-1} > \alpha_{w^*-1}$, in $j_3$ we obtain

$$p_1 + 2A_4 s_1 + (1 - 2r_2 \Psi_0)m_3 + 2x_3 + 2\eta_{\alpha_3}(1 - R_1)((1 - r_2)X_0 - r_2\Psi_0)(1 - m_3)m_d$$
$$= p_1^* + 2A_4 s_1^* + (1 - 2r_2^* \Psi_0^*)m_3^* + 2x_3^* + 2\eta_{\alpha_3}(1 - R_1^*)((1 - r_2^*)X_0^* - r_2^*\Psi_0^*)(1 - m_3^*)m_d^*$$

$$\implies A_4(s_1 - s_1^*) + (x_3 - x_3^*) = \frac{p_1^* - p_1}{2}.$$

Following from $A_4(s_1 - s_1^*) + (x_3 - x_3^*) \in \mathbb{Z}$ and $|p_1^* - p_1| \leq 1$ by construction, it must be the case $p_1 = p_1^*$. Note that $x_3 = (1 - r_2)(A_4 - 1)$ since $R_{k-1} = 1$. This leaves us with

$$(s_1 - s_1^*) + (1 - r_2) = \frac{x_3^* + (1 - r_2^*)}{A_4}.$$

Since $(s_1 - s_1^*) + (1 - r_2) \in \mathbb{Z}$ and $1 - r_2^* \leq x_3^* + (1 - r_2^*) \leq A_4 - r_2^*$ with $A_4 \geq 2$, it follows that $x_3^* = (1 - r_2^*)(A_4 - 1)$.

Induction Step: Assume that $3 \leq w^* \leq k-1$ with $4 \leq k \leq d$ and $d \geq 3$ such that $3 \leq w+1 \leq k-1$ so that the following makes sense to consider. Otherwise, we will have shown that $x_z^* = (1-r_2^*)(A_{z+1}-1)$ for all $3 \leq z \leq k-1$ when $w^* = k-1$, which is not possible as $R_{k-1}^* = 0$ by assumption. Further assume that for all $3 \leq z \leq y^* - 1$ with $3 \leq y^* \leq w^*$, $x_z^* = (1 - r_2^*)(A_{z+1} - 1)$. We will show that our statement holds for $3 \leq w^* + 1 \leq k - 1$ by showing that $x_{w^*}^* = (1 - r_2^*)(A_{w^*+1} - 1)$ since our induction hypothesis already grants us the above for all $3 \leq z \leq w^* - 1$. Note that $r_2 = r_2^*$, and that our assumptions on $x_z$ and $x_z^*$ for $3 \leq z \leq w^* - 1$ imply $(1 - r_2)X_{w^*-3} = (1 - r_2)$, $r_2\Psi_{w^*-3} = r_2$, $(1 - r_2^*)X_{w^*-3}^* = (1 - r_2^*)$, and $r_2^*\Psi_{w^*-3}^* = r_2^*$. Applying the pertinent assumptions and results we have so far along with the observations $\eta_{w^*+1} = 1$ since $w^* + 1 \leq k - 1 < d$ and $\eta_{\alpha_{w^*}} = 1 = \eta_{\alpha_{w^*-1}}$ since $\alpha^* > \alpha_{w-1} > \alpha_{w^*-1}$, in $j_{w^*}$ we have

$p_{w^*-2} + 2A_{w^*+1}s_{w^*-2} + (1 - 2r_2\Psi_{w^*-3})m_{w^*} + 2x_{w^*} + 2\eta_{w^*+1}\eta_{\alpha_{w^*}}(1 - R_1)((1 - r_2)X_{w^*-3} - r_2\Psi_{w^*-3})(1 - m_{w^*})m_d$

$= p_{w^*-2}^* + 2A_{w^*+1}s_{w^*-2}^* + (1 - 2r_2^*\Psi_{w^*-3}^*)m_{w^*}^* + 2x_{w^*}^* + 2\eta_{w^*+1}\eta_{\alpha_{w^*}}(1 - R_1^*)((1 - r_2^*)X_{w^*-3}^* - r_2^*\Psi_{w^*-3}^*)(1 - m_{w^*}^*)m_d^*$

$$\implies A_{w^*+1}(s_{w^*-2} - s_{w^*-2}^*) + (x_{w^*} - x_{w^*}^*) = \frac{p_{w^*-2}^* - p_{w^*-2}}{2}.$$

Since $A_{w^*+1}(s_{w^*-2} - s_{w^*-2}^*) + (x_{w^*} - x_{w^*}^*) \in \mathbb{Z}$ and $|p_{w^*-2}^* - p_{w^*-2}| \leq 1$, it must be the case $p_{w^*-2} = p_{w^*-2}^*$. Note that $x_{w^*} = (1 - r_2)(A_{w^*+1} - 1)$ since $R_{k-1} = 1$. Hence, this leaves us with

$$(s_{w^*-2} - s_{w^*-2}^*) + (1 - r_2) = \frac{x_{w^*}^* + (1 - r_2^*)}{A_{w^*+1}}.$$

Given $(s_{w^*-2} - s_{w^*-2}^*) + (1 - r_2^*) \in \mathbb{Z}$ and $1 - r_2 \leq x_{w^*}^* + (1 - r_2^*) \leq A_{w^*+1} - r_2^*$ by construction, it follows that $x_{w^*}^* = (1 - r_2^*)(A_{w^*+1} - 1)$.

By the Principle of Strong Mathematical Induction, it follows that for every $3 \leq w^* \leq k - 1$, it is the case that $x_z^* = (1 - r_2^*)(A_{z+1} - 1)$ for every $3 \leq z \leq w^* - 1$.

Now, we see that during all $j_k$ moves for $4 \leq k \leq d$, it is the case that the smallest $3 \leq w^* \leq k - 1$ such that $x_{w^*}^* \neq (1 - r_2^*)(A_{w^*+1} - 1)$ satisfies $x_z^* = (1 - r_2)(A_{z+1} - 1)$ for every $3 \leq z \leq w^* - 1$. By the same argument made in the induction step, in $j_{w^*}$ we get

$p_{w^*-2} + 2A_{w^*+1}s_{w^*-2} + (1 - 2r_2\Psi_{w^*-3})m_{w^*} + 2x_{w^*} + 2\eta_{w^*+1}\eta_{\alpha_{w^*}}(1 - R_1)((1 - r_2)X_{w^*-3} - r_2\Psi_{w^*-3})(1 - m_{w^*})m_d$

$= p_{w^*-2}^* + 2A_{w^*+1}s_{w^*-2}^* + (1 - 2r_2^*\Psi_{w^*-3}^*)m_{w^*}^* + 2x_{w^*}^* + 2\eta_{w^*+1}\eta_{\alpha_{w^*}}(1 - R_1^*)((1 - r_2^*)X_{w^*-3}^* - r_2^*\Psi_{w^*-3}^*)(1 - m_{w^*}^*)m_d^*$

$$\implies A_{w^*+1}(s_{w^*-2} - s_{w^*-2}^*) + (x_{w^*} - x_{w^*}^*) = \frac{p_{w^*-2}^* - p_{w^*-2}}{2}.$$

Since $A_{w^*+1}(s_{w^*-2} - s_{w^*-2}^*) + (x_{w^*} - x_{w^*}^*) \in \mathbb{Z}$ and $|p_{w^*-2}^* - p_{w^*-2}| \leq 1$, it must be the case $p_{w^*-2} = p_{w^*-2}^*$. Note that $x_{w^*} = (1 - r_2)(A_{w^*+1} - 1)$ since $R_{k-1} = 1$. Hence, this leaves us with

$$(s_{w^*-2} - s_{w^*-2}^*) + (1 - r_2) = \frac{x_{w^*}^* + (1 - r_2^*)}{A_{w^*+1}}.$$

Given $(s_{w^*-2} - s_{w^*-2}^*) + (1 - r_2) \in \mathbb{Z}$ and $x_{w^*}^* + (1 - r_2^*) \neq (1 - r_2^*)A_{w^*+1}$ by construction and assumption for this subcase, it follows that $\frac{x_{w^*}^* + (1-r_2^*)}{A_{w^*+1}} \notin \mathbb{Z}$. Consequently, we have a contradiction as the above equality is not possible. Further, since $3 \leq w^* \leq k - 1$ is the arbitrary smallest integer assumed to satisfy $x_{w^*}^* \neq (1 - r_2^*)(A_{w^*+1} - 1)$, it follows that this holds for all such $w^*$, meaning this scenario does not ever occur as there does not exist such a $w^*$ under these assumptions.

Subcase 3: Let $R_1^* = 1$. Note that here we treat $j_k$ moves for all $3 \leq k \leq d$. Then, $C_2$ is being defined by column transitions, and so $x_z^* = (1 - r_2^*)(A_{z+1} - 1)$ for all $3 \leq z \leq d$. From our case assumption, recall that $x_z = (1 - r_2)(A_{z+1} - 1)$ for all $3 \leq z \leq k - 1$ and that there exists at least one $k \leq y \leq d$ such that $x_y \neq (1 - r_2)(A_{y+1} - 1)$. In particular, there exists a smallest integer $k \leq w \leq d$, by the Well-Ordering Principle, such that $x_w \neq (1 - r_2)(A_{w+1} - 1)$, leaving $x_z = (1 - r_2)(A_{z+1} - 1)$ for every $k \leq z \leq w - 1$. Observing that $m_1 = m_2$, $m_d = 1 - m_1$, $m_1^* = m_2^* = (1 - r_2)r_2^* + r_2(1 - r_2^*)$, and $m_d^* = 1 - m_1^* = (1 - r_2)(1 - r_2^*) + r_2 r_2^*$, in $j_1$ we get

$$(\ell + 1)A_3 + (-2 + m_1 + 2m_d(1 - m_1)) = (\ell_1 + 1)A_3 + (-2 + m_1^* + 2m_d^*(1 - m_1^*))$$

$$\implies (\ell - \ell_1)\frac{A_3}{2} + (m_d - m_d^*) = \frac{m_1^* - m_1}{2}.$$

Since $(\ell - \ell_1)\frac{A_3}{2} + (m_d^* - m_d) \in \mathbb{Z}$ and $|m_d^* - m_d| \leq 1$ by construction, we see $m_1 = m_1^*$ and hence $(m_1, m_2, m_d) = (m_1^*, m_2^*, m_d^*)$. Applying this to $j_2$, we find

$$m_2 + \gamma + 2x_2 + 2m_1(1 - m_2) + 2((1 - m_1)m_d - (1 - R_1)((1 - m_2)m_1 + m_d(1 - m_1)))$$
$$= m_2^* + \gamma_1 + 2x_2^* + 2m_1^*(1 - m_2^*) + 2((1 - m_1^*)m_d^* - (1 - R_1^*)((1 - m_2^*)m_1^* + m_d^*(1 - m_1^*)))$$

$$\implies (x_2 - x_2^*) - m_d^* = \frac{\gamma_1 - \gamma}{2}.$$

Since $(x_2 - x_2^*) - m_d^* \in \mathbb{Z}$ and $|\gamma_1 - \gamma| \leq 1$ by construction, we must have $\gamma = \gamma_1$, allowing us to conclude $x_2 - x_2^* = m_d^*$. The last conclusion gives us

$$(x_2 - x_2^*) - m_d^* = 0$$
$$\implies (r_2 - r_2^*) - (1 - r_2)(1 - r_2^*) - r_2 r_2^* \equiv 0 \pmod{2}$$
$$\implies 2r_2 - 2r_2 r_2^* - 1 \equiv 0 \pmod{2}$$
$$\implies -1 \equiv 0 \pmod{2}.$$

Thus, we have a contradiction as the above is not possible.

Case 2: Focusing on moves along $j_k$ for $4 \leq k \leq d$ since $j_3$ is not subject to normal moves, either $\alpha_2 < \alpha^* \leq \alpha_{k-1}$ or let $\alpha_{k-1} < \alpha^* \leq \alpha_d$ with $R_1 = 0$ and $R_{k-1} = 0$. Then, $C_1$ is being defined by normal moves and so, in the case $\alpha_{k-1} < \alpha^* \leq \alpha_d$, with $R_1 = 0$ and $R_{k-1} = 0$, there exists at least one integer $3 \leq y \leq k - 1$ such that $x_y \neq (1 - r_2)(A_{y+1} - 1)$. Let $3 \leq w \leq k - 1$ be the smallest such integer and observe that $x_z = (1 - r_2)(A_{z+1} - 1)$ for all $3 \leq z \leq w - 1$.

Subcase 1: Assume $\alpha_2 < \alpha^* \leq \alpha_{k-1}$ if this is already being assumed and otherwise assume all of the following are the case: $\alpha_{k-1} < \alpha^* \leq \alpha_d$, $R_{k-1}^* = 0$ and $R_1^* = 0$. Now, observing that during all such moves $m_d = 1 - m_1$ and $(m_1^*, m_d^*) = (1 - m_1, 1 - m_d)$, in $j_1$ we obtain

$$(\ell + 1)A_3 + (-2 + m_1 + 2m_d(1 - m_1)) = (\ell_1 + 1)A_3 + (-2 + m_1^* + 2m_d^*(1 - m_1^*))$$

$$\implies (\ell - \ell_1)\frac{A_3}{2} + m_1 + (m_d - m_d^*) = \frac{1}{2}.$$

Given $(\ell - \ell_1)\frac{A_3}{2} + m_1 + (m_d - m_d^*) \in \mathbb{Z}$ by construction and $\frac{1}{2} \notin \mathbb{Z}$, we have a contradiction as the above equality is not possible.

Subcase 2: Let $R_1^* = 1$, meaning $C_2$ is being defined by column transitions. As a consequence, we must have $\alpha_{k-1} < \alpha^* \leq \alpha_d$ and $x_z^* = (1 - r_2^*)(A_{z+1} - 1)$ for all $3 \leq z \leq d$. Further, we have $R_1 = 0$ and $R_{k-1} = 0$. Observing that during all such moves $m_1 = \cdots = m_{k-1}$, $m_k = \cdots = m_d = 1 - m_1$, $m_1^* = \cdots = m_{k-1}^* = (1-r_2^*)(1-m_1)+r_2^*m_1$ and $m_k^* = \cdots = m_d^* = 1-m_1^* = (1-r_2^*)(1-m_d)+r_2^*m_d$, in $j_1$ we see

$$(\ell+1)A_3 + (-2 + m_1 + 2m_d(1-m_1)) = (\ell_1+1)A_3 + (-2 + m_1^* + 2m_d^*(1-m_1^*))$$
$$\implies (\ell - \ell_1)\frac{A_3}{2} + (m_d - m_d^*) = \frac{m_1^* - m_1}{2}.$$

Since $(\ell-\ell_1)\frac{A_3}{2} + (m_d - m_d^*) \in \mathbb{Z}$ and $|m_1^* - m_1| \le 1$ by construction, it follows that $m_1 = m_1^*$. This last deduction implies

$$m_1 - m_1^* = 0$$
$$\implies (1 - r_2^*)(2m_1 - 1) = 0.$$

Given $r_2^*, m_1 \in \{0,1\}$, we see that it must be the case $r_2^* = 1$. So we get $(m_1^*, \ldots, m_d^*) = (m_1, \ldots, m_d)$. Applying this to $j_2$, we find

$$m_2 + \gamma + 2x_2 + 2m_1(1-m_2) + 2((1-m_1)m_d - (1-R_1)((1-m_2)m_1 + m_d(1-m_1)))$$
$$= m_2^* + \gamma_1 + 2x_2^* + 2m_1^*(1-m_2^*) + 2((1-m_1^*)m_d^* - (1-R_1^*)((1-m_2^*)m_1^* + m_d^*(1-m_1^*)))$$
$$\implies (x_2 - x_2^*) - m_d^* = \frac{\gamma_1 - \gamma}{2}.$$

Since $(x_2 - x_2^*) - m_d^* \in \mathbb{Z}$ and $|\gamma_1 - \gamma| \le 1$ by construction, we must have $\gamma = \gamma_1$, allowing us to conclude $x_2 - x_2^* = m_d^*$. Hence, it follows that

$$(x_2 - x_2^*) - m_d^* = 0$$
$$\implies (r_2 - r_2^*) - m_d^* \equiv 0 \pmod 2$$
$$\implies r_2 - m_d^* \equiv 1 \pmod 2.$$

Then, since $|r_2 - m_d| \le 1$ by construction, we must have $|r_2 - m_d^*| = 1$. From this, it immediately follows that $r_2 + m_d^* = 1$. We will invoke this last result where applicable without further mention.

To complete our treatment of this subcase, we will require the following result we will prove using induction:

We will prove that for every $3 \le w \le k-1$ such that $w$ is the smallest integer with $x_w \ne (1-r_2)(A_{w+1}-1)$, it is the case that for all $3 \le z \le w-1$, $x_z = (1-r_2)(A_{z+1}-1)$. We proceed by induction on $w$:

<u>Base Case 1:</u> Let $w = 3$. Then, the statement is vacuously true as there does not exist an integer $z$ such that $3 \le z \le 2$.

<u>Base Case 2:</u> Let $w = 4$. Then, $z = 3$ and so our assumptions for this subcase give us in $j_3$ that

$$p_1 + 2A_4s_1 + (1 - 2r_2\Psi_0)m_3 + 2x_3 + 2\eta_{\alpha_3}(1-R_1)((1-r_2)X_0 - r_2\Psi_0)(1-m_3)m_d$$
$$= p_1^* + 2A_4s_1^* + (1 - 2r_2^*\Psi_0^*)m_3^* + 2x_3^* + 2\eta_{\alpha_3}(1-R_1^*)((1-r_2^*)X_0^* - r_2^*\Psi_0^*)(1-m_3^*)m_d^*$$
$$\implies A_4(s_1 - s_1^*) + (1-r_2)m_3 + (x_3 - x_3^*) + (1-2r_2)(1-m_3)m_d = \frac{p_1^* - p_1}{2}.$$

Since $A_4(s_1 - s_1^*) + (1 - r_2)m_3 + (x_3 - x_3^*) + (1 - 2r_2)(1 - m_3)m_d \in \mathbb{Z}$ and $|p_1^* - p_1| \leq 1$ by construction, it must be the case $p_1 = p_1^*$. Noting that $x_3^* = 0$ since $r_2^* = 1$, the above leaves us with

$$s_1 - s_1^* = -\frac{x_3 + (1 - r_2)}{A_4}.$$

Following from $s_1 - s_1^* \in \mathbb{Z}$ and $1 - r_2 \leq x_3 + (1 - r_2) \leq A_4 - r_2$, we see that we must have $x_3 = (1 - r_2)(A_4 - 1)$.

<u>Induction Step:</u> Let $3 \leq w \leq k - 1$ with $4 \leq k \leq d$ and $d \geq 3$ such that $3 \leq w + 1 \leq d$ so that the following makes sense to consider. Otherwise, we will have shown that $x_z = (1 - r_2)(A_{z+1} - 1)$ for all $3 \leq z \leq k - 1$ when $w = k - 1$, which is not possible as $R_{k-1} = 0$ by assumption. Further assume that for all $3 \leq z \leq y - 1$ with $3 \leq y \leq w$, $x_z = (1 - r_2)(A_{z+1} - 1)$. We will show that our statement holds for $3 \leq w + 1 \leq k - 1$ by showing that $x_w = (1 - r_2)(A_{w+1} - 1)$ since our induction hypothesis grants us the above for all $3 \leq z \leq w - 1$. Note that our assumptions on $x_z$ and $x_z^*$ imply $(1 - r_2)X_{w-3} = (1 - r_2)$, $r_2\Psi_{w-3} = r_2$, $(1 - r_2^*)X_{w-3}^* = (1 - r_2^*)$, and $r_2^*\Psi_{w-3}^* = r_2^*$. The corresponding assumptions and results we have so far, along with the observations $\eta_{w+1} = 1$ since $w + 1 \leq k - 1 < d$ and $\eta_{\alpha_w} = \eta_{\alpha_{w-1}} = 1$ since $\alpha^* > \alpha_{k-1} \geq \alpha_w$, yield in $j_w$

$$p_{w-2} + 2A_{w+1}s_{w-2} + (1 - 2r_2\Psi_{w-3})m_w + 2x_w + 2\eta_{w+1}\eta_{\alpha_w}(1 - R_1)((1 - r_2)X_{w-3} - r_2\Psi_{w-3})(1 - m_w)m_d$$
$$= p_{w-2}^* + 2A_{w+1}s_{w-2}^* + (1 - 2r_2^*\Psi_{w-3}^*)m_w^* + 2x_w^* + 2\eta_{w+1}\eta_{\alpha_w}(1 - R_1^*)((1 - r_2^*)X_{w-3}^* - r_2^*\Psi_{w-3}^*)(1 - m_w^*)m_d^*$$
$$\implies A_{w+1}(s_{w-2} - s_{w-2}^*) + (1 - r_2)m_w + (x_w - x_w^*) + (1 - 2r_2)(1 - m_w)m_d = \frac{p_{w-2}^* - p_{w-2}}{2}.$$

Since $A_{w+1}(s_{w-2} - s_{w-2}^*) + (1 - r_2)m_w + (x_w - x_w^*) + (1 - 2r_2)(1 - m_w)m_d \in \mathbb{Z}$ and $|p_{w-2}^* - p_{w-2}| \leq 1$, it must be the case $p_{w-2} = p_{w-2}^*$. Hence, observing that $x_w^* = (1 - r_2^*)(A_{w+1} - 1) = 0$ since $r_2^* = 1$, we are left with

$$s_{w-2} - s_{w-2}^* = -\frac{x_w + (1 - r_2)}{A_{w+1}}.$$

Given $s_{w-2} - s_{w-2}^* \in \mathbb{Z}$ and $1 - r_2 \leq x_w + (1 - r_2) \leq A_{w+1} - r_2$ by construction, it follows that $x_w = (1 - r_2)(A_{w+1} - 1)$.

By the Principle of Strong Mathematical Induction, it follows that for every $3 \leq w \leq k - 1$, it is the case that $x_z = (1 - r_2)(A_{z+1} - 1)$ for every $3 \leq z \leq w - 1$.

Now, during every $j_k$ move for $4 \leq k \leq d$, we see that for the smallest integer $3 \leq w \leq k - 1$ such that $x_w \neq (1 - r_2)(A_{w+1} - 1)$ granted by our case assumption, our induction argument above gives us that for all $3 \leq z \leq w - 1$, $x_z = (1 - r_2)(A_{z+1} - 1)$. Applying the same argument made in the induction step, in $j_w$ we obtain

$$p_{w-2} + 2A_{w+1}s_{w-2} + (1 - 2r_2\Psi_{w-3})m_w + 2x_w + 2\eta_{w+1}\eta_{\alpha_w}(1 - R_1)((1 - r_2)X_{w-3} - r_2\Psi_{w-3})(1 - m_w)m_d$$
$$= p_{w-2}^* + 2A_{w+1}s_{w-2}^* + (1 - 2r_2^*\Psi_{w-3}^*)m_w^* + 2x_w^* + 2\eta_{w+1}\eta_{\alpha_w}(1 - R_1^*)((1 - r_2^*)X_{w-3}^* - r_2^*\Psi_{w-3}^*)(1 - m_w^*)m_d^*$$
$$\implies A_{w+1}(s_{w-2} - s_{w-2}^*) + (1 - r_2)m_w + (x_w - x_w^*) + (1 - 2r_2)(1 - m_w)m_d = \frac{p_{w-2}^* - p_{w-2}}{2}.$$

Since $A_{w+1}(s_{w-2} - s_{w-2}^*) + (1 - r_2)m_w + (x_w - x_w^*) + (1 - 2r_2)(1 - m_w)m_d \in \mathbb{Z}$ and $|p_{w-2}^* - p_{w-2}| \leq 1$, it must be the case $p_{w-2} = p_{w-2}^*$. Hence, observing that $x_w^* = (1 - r_2^*)(A_{w+1} - 1) = 0$ since $r_2^* = 1$, we are left with

$$s_{w-2} - s^*_{w-2} = -\frac{x_w + (1-r_2)}{A_{w+1}}.$$

Given $s_{w-2} - s^*_{w-2} \in \mathbb{Z}$ and $x_w + (1-r_2) \neq (1-r_2)A_{w+1}$ by construction and assumption for this subcase, it follows that $-\frac{x_w+(1-r_2)}{A_{w+1}} \notin \mathbb{Z}$. Consequently, we have a contradiction as the above equality is not possible. Further, since $3 \leq w \leq k-1$ is the arbitrary smallest integer assumed to satisfy $x_w \neq (1-r_2)(A_{w+1} - 1)$, it follows that this holds for all such $w$, meaning this scenario does not ever occur as there does not exist such a $w$ under these assumptions.

<u>Case 3:</u> Let $R_1 = 1$. Then, $C_1$ is being defined by column transitions along $j_k$ for $3 \leq k \leq d$, and so we must have $\alpha_{k-1} < \alpha^* \leq \alpha_d$. Further, we know that $x_z = (1-r_2)(A_{z+1} - 1)$ for all $3 \leq z \leq d$. In this case, we will only need to consider when $C_2$ is being defined by column transitions.

Hence, let $R^*_1 = 1$ so that $x^*_z = (1-r^*_2)(A_{z+1} - 1)$ for all $3 \leq z \leq d$, and observe that during all such moves $m_1 = m_2$, $m_d = 1 - m_1$, $m^*_1 = m^*_2 = (1 - |r^*_2 - r_2|)(1 - m_1) + |r^*_2 - r_2|m_1$ and $m^*_d = 1 - m^*_1 = (1 - |r^*_2 - r_2|)(1 - m_d) + |r^*_2 - r_2|m_d$. Applying this to $j_1$, we see

$$(\ell+1)A_3 + (-2 + m_1 + 2m_d(1-m_1)) = (\ell_1 + 1)A_3 + (-2 + m^*_1 + 2m^*_d(1-m^*_1))$$
$$\implies (\ell - \ell_1)\frac{A_3}{2} + (m_d - m^*_d) = \frac{m^*_1 - m_1}{2}.$$

Since $(\ell - \ell_1)\frac{A_3}{2} + (m_d - m^*_d) \in \mathbb{Z}$ and $|m^*_1 - m_1| \leq 1$ by construction, it follows that $m_1 = m^*_1$ and so $(m_1, m_2, m_d) = (m^*_1, m^*_2, m^*_d)$, as we will see next. The former deduction above implies

$$m_1 - m^*_1 = 0$$
$$\implies (1 - |r^*_2 - r_2|)(2m_1 - 1) = 0.$$

Given $|r^*_2 - r_2| \leq 1$ and $m_1 \in \{0, 1\}$, it must be the case $|r^*_2 - r_2| = 1$, meaning $r_2 \neq r^*_2$. Going to $j_2$, we find

$$m_2 + \gamma + 2x_2 + 2m_1(1 - m_2) + 2((1-m_1)m_d - (1-R_1)((1-m_2)m_1 + m_d(1-m_1)))$$
$$= m^*_2 + \gamma_1 + 2x^*_2 + 2m^*_1(1 - m^*_2) + 2((1-m^*_1)m^*_d - (1-R^*_1)((1-m^*_2)m^*_1 + m^*_d(1-m^*_1)))$$
$$\implies x_2 - x^*_2 = \frac{\gamma_1 - \gamma}{2}.$$

Following from $x_2 - x^*_2 \in \mathbb{Z}$ and $|\gamma_1 - \gamma| \leq 1$ by construction, we must have $\gamma = \gamma_1$, meaning $x_2 = x^*_2$. From this, it is immediate that $r_2 = r^*_2$, a contradiction to our deduction that $r_2 \neq r^*_2$. Hence, our initial equality is not possible.

Thus, for all $\alpha_1 \leq \alpha^* \leq \alpha_d$ with $d \geq 2$, it follows that $C_1$ and $C_2$ do not share edges from moves with the opposite orientation, proving Proposition 15 as desired.

∎

### 7.1.2 Proof of Proposition 16:

Let $d \in \mathbb{Z}^{\geq 2}$, $\alpha_1 \leq \alpha^* \leq \alpha_d$, and for reference refer to Theorem 5 for specific bounds of each parameter, though we will bring some of them up in our arguments as necessary.

We would also like to make the reader aware that the form of the moves with same orientation in $C_2$ to a given move $(m_1, \ldots, m_d)$ on $C_1$ will vary based on which move and which two major sets are being considered. In some instances, there may be more than one move with the same orientation and so there would not be a unique inverse move in said case. The definition of the inverse move $(m_1^*, \ldots, m_d^*)$ to a given move $(m_1, \ldots, m_d)$ and all applicable equations that hold during every move will be introduced as required.

Note that we will be assuming for the sake of contradiction that two cycles share an edge with the move configurations established below, meaning all components of their vertices agree, and we then show that such scenarios do not occur by the definition of the General Lock-and-Key Decomposition. In particular, we will be assuming that the same edge is present in two distinct major sets at a time and show that this not possible. In doing so, we will have shown that these cycles cannot share edges by way of these configurations. We now proceed with the proof of Proposition 16 by considering equalities at each component $j_k$ for $1 \leq k \leq d$ with $C_2$ configured to perform a move $(m_1^*, \ldots, m_d^*)$ with the same orientation as that of a given move $(m_1, \ldots, m_d)$ of $C_1$ starting from the same vertex belonging to the edge in $C_1$.

**Dimension Case 1:** Let $d = 2$. Then, $C_1 = C_{\ell,\gamma}$ and $C_2 = C_{\ell_1,\gamma_1}$. Note that $t = 0 = t_1$ since $d = 2$. We now case on $\alpha^*$ for $\alpha_1 \leq \alpha^* \leq \alpha_2$:

$\boldsymbol{\alpha^*}$**—Case 1:** Let $\alpha^* = \alpha_1$. In this case, only sets 3 and 4 are active as $A_1 = 0$ and $A_2 = 0$.

**Set Comparison Case 1:** We will consider whether edges are shared from moves with the same orientation via sets 3 and 4:

Major Case 1: Focusing on moves along $j_1$, we see that $C_1$ and $C_2$ have $j_1$ moves with the same orientation if and only if $\gamma \neq \gamma_1$, meaning $|\gamma_1 - \gamma| = 1$. Noting that $m_2 = 1 = m_2^*$ during these moves, in $j_2$ we get

$$m_2 + \gamma + 4x_2 + (1 - m_2)m_1 = 2 + m_1^* + \gamma_1 + m_2^*(1 - m_1^*) + 4x_2^*$$
$$\implies x_2 - x_2^* = \frac{2 + (-1)^\gamma}{4}.$$

Since $x_2 - x_2^* \in \mathbb{Z}$ and for every $\gamma \in \{0, 1\}$, $\frac{2+(-1)^\gamma}{4} \notin \mathbb{Z}$ by construction, we have a contradiction as the above equality is not possible.

Major Case 2: Focusing on moves along $j_2$, we observe that during all such moves $m_2 = 0 = m_2^*$. Applying this to $j_2$, we find

$$m_2 + \gamma + 4x_2 + (1 - m_2)m_1 = 2 + m_1^* + \gamma_1 + m_2^*(1 - m_1^*) + 4x_2^*$$
$$\implies x_2 - x_2^* = \frac{2 + (\gamma_1 - \gamma) + (m_1^* - m_1)}{4}.$$

Since $x_2 - x_2^* \in \mathbb{Z}$ and $0 \leq 2 + (\gamma_1 - \gamma) + (m_1^* - m_1) \leq 4$ by construction, we must have $(\gamma_1 - \gamma) = (m_1^* - m_1)$ with $|\gamma_1 - \gamma| = |m_1^* - m_1| = 1$. Hence, in $j_1$ we have

$$(\ell+1-\gamma)A_3-\gamma+(-1)^\gamma(-2+m_1+(1-m_2)m_1) = (\ell_1+1-\gamma_1)A_3-\gamma_1+(-1)^{\gamma_1+1}(m_1^*+m_1^*(1-m_2^*))$$

$$\implies ((\ell-\ell_1)+(\gamma_1-\gamma))\frac{A_3}{2}+((-1)^\gamma(m_1-1)+(-1)^{\gamma_1}m_1^*) = \frac{(-1)^{\gamma+1}}{2}.$$

Given $((\ell-\ell_1)+(\gamma_1-\gamma))\frac{A_3}{2}+((-1)^\gamma(m_1-1)+(-1)^{\gamma_1}m_1^*) \in \mathbb{Z}$ and for every $\gamma \in \{0,1\}$, $\frac{(-1)^{\gamma+1}}{2} \notin \mathbb{Z}$ by construction, we have a contradiction as the above equality is not possible.

Thus, for $\alpha^* = \alpha_1$, no edges are shared from moves with the same orientation via two distinct sets.

**$\alpha^*$−Case 2:** Let $\alpha_1 < \alpha^* < \alpha_2$. Then, $A_1 = 1$ and $A_2 = 0$, so sets $1, 2, 3$ and $4$ are all active. Lastly, note that $A_3 \geq 4$ is even.

**Set Comparison Case 1:** We will consider whether edges are shared from moves with the same orientation via sets 1 and 2.

Major Case 1: Focusing on moves along $j_1$, we see that $C_1$ and $C_2$ have $j_1$ moves with the same orientation if and only if $\gamma \neq \gamma_1$, meaning $|\gamma_1-\gamma| = 1$. Applying this to $j_2$ along with the observations $m_1 + m_2 = 1$ and $m_1^* = m_2^*$, we see

$$m_2 + \gamma + 4x_2 = 2 + m_1^* + \gamma_1 + 4x_2^*$$

$$\implies x_2 - x_2^* = \frac{2+(\gamma_1-\gamma)+(m_2^*-m_2)}{4}.$$

Since $x_2 - x_2^* \in \mathbb{Z}$ and $0 \leq 2+(\gamma_1-\gamma)+(m_2^*-m_2) \leq 4$ by construction, we must have $(\gamma_1-\gamma) = (m_2^*-m_2)$ with $|\gamma_1-\gamma| = |m_2^*-m_2| = 1$. Consequently, in $j_1$ we obtain

$$(\ell+\gamma)A_3-\gamma+(-1)^\gamma(m_1+2x_1) = (\ell_1+1-\gamma_1)A_3-\gamma_1+(-1)^{\gamma_1+1}(2+m_2^*+2m_1^*(1-m_2^*)+2x_1^*)$$

$$\implies ((\ell-\ell_1)+(\gamma+\gamma_1)-1)\frac{A_3}{2}+(-1)^\gamma(x_1-x_1^*) = \frac{(-1)^{m_2}}{2}.$$

Given $((\ell-\ell_1)+(\gamma+\gamma_1)-1)\frac{A_3}{2}+(-1)^\gamma(x_1-x_1^*) \in \mathbb{Z}$ and for all $m_2 \in \{0,1\}$, $\frac{(-1)^{m_2}}{2} \notin \mathbb{Z}$ by construction, we have a contradiction as the equality above is not possible.

Major Case 2: Focusing on moves along $j_2$, we see that which $j_2$ moves in $C_2$ have the same orientation as those in $C_1$ depends on $x_1^*$. Hence, we case on $x_1^*$:

Case 1: Let $0 \leq x_1^* < \frac{A_3}{2} - 2$. Then, $m_1 = m_2$ and $(m_1^*, m_2^*) = (m_1, 1-m_1)$, and so in $j_2$ we get

$$m_2 + \gamma + 4x_2 = 2 + m_1^* + \gamma_1 + 4x_2^*$$

$$\implies x_2 - x_2^* = \frac{2+(\gamma_1-\gamma)}{4}.$$

Given $x_2 - x_2^* \in \mathbb{Z}$ and $1 \leq 2+(\gamma_1-\gamma) \leq 3$ by construction, it follows that for every $\gamma, \gamma_1 \in \{0,1\}$, $\frac{2+(\gamma_1-\gamma)}{4} \notin \mathbb{Z}$, giving us a contradiction as the above equality is not possible.

<u>Case 2:</u> Let $x_1^* = \frac{A_3}{2} - 2$. Here, all of $C_2$'s $j_2$ moves have the same orientation as $C_1$'s $j_2$ move $(m_1, m_2) = (0,0)$, and during all such moves $m_1^* + m_2^* = 1$. So in $j_2$, we have

$$m_2 + \gamma + 4x_2 = 2 + m_1^* + \gamma_1 + 4x_2^*$$
$$\implies x_2 - x_2^* = \frac{2 + m_1^* + (\gamma_1 - \gamma)}{4}.$$

Since $x_2 - x_2^* \in \mathbb{Z}$ and $1 \leq 2 + m_1^* + (\gamma_1 - \gamma) \leq 4$ by construction, we must have $(\gamma_1 - \gamma) = m_1^* = 1$. Applying this to $j_1$, we get

$$(\ell + \gamma)A_3 - \gamma + (-1)^\gamma(m_1 + 2x_1) = (\ell_1 + 1 - \gamma_1)A_3 - \gamma_1 + (-1)^{\gamma_1+1}(2 + m_2^* + 2m_1^*(1 - m_2^*) + 2x_1^*)$$
$$\implies (\ell - \ell_1)\frac{A_3}{2} + (x_1 - x_1^*) - 1 = \frac{1}{2}.$$

Since $(\ell - \ell_1)\frac{A_3}{2} + (x_1 - x_1^*) - 1 \in \mathbb{Z}$ by construction and $\frac{1}{2} \notin \mathbb{Z}$, we have a contradiction as the above equality is not possible.

**Set Comparison Case 2:** We will consider whether edges are shared from moves with the same orientation via sets 1 and 3.

<u>Major Case 1:</u> Focusing on moves along $j_1$, we see that $C_1$ and $C_2$ have $j_1$ moves with the same orientation if and only if $\gamma = \gamma_1$. Observing that during all such moves $m_1 + m_2 = 1$ and $m_2^* = 1$, in $j_1$ we get

$$(\ell + \gamma)A_3 - \gamma + (-1)^\gamma(m_1 + 2x_1) = (\ell_1 + 1 - \gamma_1)A_3 - \gamma_1 + (-1)^{\gamma_1}(-2 + m_1^* + (1 - m_2^*)m_1^*)$$
$$\implies (\ell - \ell_1) + 2\gamma - 1 = \frac{(-1)^{\gamma_1+1}(2 + (m_1 - m_1^*) + 2x_1)}{A_3}.$$

Given $(\ell - \ell_1) + 2\gamma - 1 \in \mathbb{Z}$ and $1 \leq 2 + (m_1 - m_1^*) + 2x_1 \leq A_3 - 1$ by construction, it follows that for every $0 \leq x_1 \leq \frac{A_3}{2} - 2$ and $m_1, m_1^* \in \{0, 1\}$, $\frac{(-1)^{\gamma_1+1}(2+(m_1-m_1^*)+2x_1)}{A_3} \notin \mathbb{Z}$. Thus, we have a contradiction as the above equality is not possible.

<u>Major Case 2:</u> Focusing on moves along $j_2$, observe that all of $C_2$'s $j_2$ moves have the same orientation as $C_1$'s $j_2$ move $(m_1, m_2) = (0,0)$ only, and during all such moves $m_2^* = 0$. Applying this to $j_2$, we see

$$m_2 + \gamma + 4x_2 = m_2^* + \gamma_1 + 4x_2^* + (1 - m_2^*)m_1^*$$
$$\implies x_2 - x_2^* = \frac{m_1^* + (\gamma_1 - \gamma)}{4}.$$

Given $x_2 - x_2^* \in \mathbb{Z}$ and $-1 \leq m_1^* + (\gamma_1 - \gamma) \leq 2$ by construction, we must have $m_1^* = -(\gamma_1 - \gamma)$ with either $\gamma_1 - \gamma = -1$ or $\gamma_1 - \gamma = 0$. Note that if $\gamma_1 - \gamma = 1$, then for every $m_1^* \in \{0, 1\}$, $\frac{m_1^* + (\gamma_1 - \gamma)}{4} \notin \mathbb{Z}$. We will treat the remaining two cases together with $\gamma_1 - \gamma = -|\gamma_1 - \gamma|$ and $m_1^* = |\gamma_1 - \gamma|$. Now, in $j_1$ we find

$$(\ell + \gamma)A_3 - \gamma + (-1)^\gamma(m_1 + 2x_1) = (\ell_1 + 1 - \gamma_1)A_3 - \gamma_1 + (-1)^{\gamma_1}(-2 + m_1^* + (1 - m_2^*)m_1^*)$$

$$\implies (\ell - \ell_1) + (\gamma + \gamma_1) - 1 = \frac{|\gamma_1 - \gamma| + 2(-1)^{\gamma_1+1}(1 - |\gamma_1 - \gamma|) + 2(-1)^{\gamma+1}x_1}{A_3}.$$

In either case that $\gamma_1 - \gamma \in \{-1, 0\}$, it follows that $(\ell - \ell_1) + (\gamma + \gamma_1) - 1 \in \mathbb{Z}$ and for every $0 \leq x_1 \leq \frac{A_3}{2} - 2$, $\frac{|\gamma_1-\gamma|+2(-1)^{\gamma_1+1}(1-|\gamma_1-\gamma|)+2(-1)^{\gamma+1}x_1}{A_3} \notin \mathbb{Z}$ by construction. Hence, we have a contradiction as the above equality is not possible.

**Set Comparison Case 3:** We will consider whether edges are shared from the moves with the same orientation via sets 1 and 4.

Major Case 1: Focusing on moves along $j_1$, observe that $C_2$ has $j_1$ moves with the same orientation as those of $C_1$ if and only if $|\gamma_1 - \gamma| = 1$. Following from $m_1 + m_2 = 1$ and $m_2^* = 1$ during all such moves, in $j_1$ we get

$$(\ell + \gamma)A_3 - \gamma + (-1)^\gamma(m_1 + 2x_1) = (\ell_1 + 1 - \gamma_1)A_3 - \gamma_1 + (-1)^{\gamma_1+1}(m_1^* + m_1^*(1 - m_2^*))$$
$$\implies (\ell - \ell_1) + (\gamma + \gamma_1) - 1 = \frac{(-1)^{\gamma+1}(1 + 2x_1 + (m_1 - m_1^*))}{A_3}.$$

Since $(\ell - \ell_1) + (\gamma + \gamma_1) - 1 \in \mathbb{Z}$ and $0 \leq 1 + 2x_1 + (m_1 - m_1^*) \leq A_3 - 2$ by construction, it must be the case $x_1 = 0$ and $m_1 - m_1^* = -1$, meaning $m_1 = 0$ and $m_1^* = 1$. Applying this to $j_2$, we have

$$m_2 + \gamma + 4x_2 = 2 + m_1^* + \gamma_1 + m_2^*(1 - m_1^*) + 4x_2^*$$
$$\implies x_2 - x_2^* = \frac{2 + (-1)^\gamma}{4}.$$

Following from $x_2 - x_2^* \in \mathbb{Z}$ and for every $\gamma \in \{0, 1\}$, $\frac{2+(-1)^\gamma}{4} \notin \mathbb{Z}$, we have a contradiction as the above equality is not possible.

Major Case 2: Focusing on moves along $j_2$, we see that during all such moves $m_1 = m_2$ and $(m_1^*, m_2^*) = (m_1, 0)$, giving us in $j_2$

$$m_2 + \gamma + 4x_2 = 2 + m_1^* + \gamma_1 + m_2^*(1 - m_1^*) + 4x_2^*$$
$$\implies x_2 - x_2^* = \frac{2 + (\gamma_1 - \gamma)}{4}.$$

Given $x_2 - x_2^* \in \mathbb{Z}$ and $1 \leq 2 + (\gamma_1 - \gamma) \leq 3$ by construction, it follows that for every $\gamma, \gamma_1 \in \{0, 1\}$, $\frac{2+(\gamma_1-\gamma)}{4} \notin \mathbb{Z}$, giving us a contradiction as the above equality is not possible.

**Set Comparison Case 4:** We will consider whether edges are shared from moves with the same orientation via sets 2 and 3.

Major Case 1: Focusing on moves along $j_1$, observe that $C_2$'s $j_1$ moves have the same orientation as those of $C_1$ if and only if $|\gamma_1 - \gamma| = 1$. Since $m_1 = m_2$ and $m_2^* = 1$ during all such $j_1$ moves, in $j_1$ we see

$$(\ell + 1 - \gamma)A_3 - \gamma + (-1)^{\gamma+1}(2 + m_2 + 2m_1(1 - m_2) + 2x_1)$$
$$= (\ell_1 + 1 - \gamma_1)A_3 - \gamma_1 + (-1)^{\gamma_1}(-2 + m_1^* + (1 - m_2^*)m_1^*)$$
$$\implies (\ell - \ell_1) + (\gamma_1 - \gamma) = \frac{(-1)^\gamma(3 + (m_1 - m_1^*) + 2x_1)}{A_3}.$$

Given $(\ell - \ell_1) + (\gamma_1 - \gamma) \in \mathbb{Z}$ and $2 \leq 3 + (m_1 - m_1^*) + 2x_1 \leq A_3$ by construction, we must have $x_1 = \frac{A_3}{2} - 2$ and $m_1 - m_1^* = 1$, meaning $m_1 = 1$ and $m_1^* = 0$. Applying this to $j_2$, we obtain

$$2 + m_1 + \gamma + 4x_2 = m_2^* + \gamma_1 + 4x_2^* + (1 - m_2^*)m_1^*$$
$$\implies x_2 - x_2^* = \frac{-2 + (-1)^\gamma}{4}.$$

Since $x_2 - x_2^* \in \mathbb{Z}$ and $-3 \leq -2 + (-1)^\gamma \leq -1$ by construction, it follows that for every $\gamma \in \{0, 1\}$, $\frac{-2+(-1)^\gamma}{4} \notin \mathbb{Z}$, giving us a contradiction as the above equality is not possible.

Major Case 2: Focusing on moves along $j_2$, we see that which $j_2$ moves in $C_1$ have the same orientation as those in $C_2$ depends on $x_1$. Hence, we proceed by casing on $x_1$:

Case 1: Let $0 \leq x_1 < \frac{A_3}{2} - 2$. Then, we observe that $m_1 + m_2 = 1$ with $(m_1, m_2) = (0, 1)$, and $m_2^* = 0$, giving us in $j_2$

$$2 + m_1 + \gamma + 4x_2 = m_2^* + \gamma_1 + 4x_2^* + (1 - m_2^*)m_1^*$$
$$\implies x_2 - x_2^* = \frac{-2 + (\gamma_1 - \gamma) + m_1^*}{4}.$$

Following from $x_2 - x_2^* \in \mathbb{Z}$ and $-3 \leq -2 + (\gamma_1 - \gamma) + m_1^* \leq 0$ by construction, it must be the case $\gamma_1 - \gamma = m_1^* = 1$, meaning $\gamma_1 = 1$ and $\gamma = 0$. Consequently, in $j_1$ we find

$$(\ell + 1 - \gamma)A_3 - \gamma + (-1)^{\gamma+1}(2 + m_2 + 2m_1(1 - m_2) + 2x_1)$$
$$= (\ell_1 + 1 - \gamma_1)A_3 - \gamma_1 + (-1)^{\gamma_1}(-2 + m_1^* + (1 - m_2^*)m_1^*)$$
$$\implies (\ell - \ell_1) + 1 = \frac{2(1 + x_1)}{A_3}.$$

Since $(\ell - \ell_1) + 1 \in \mathbb{Z}$ and $2 \leq 2(1 + x_1) < A_3 - 2$ by construction and our case assumption, we see that $\frac{2(1+x_1)}{A_3} \notin \mathbb{Z}$. So we have a contradiction as the above equality is not possible.

Case 2: Let $x_1 = \frac{A_3}{2} - 2$. Then, observing that during all $j_2$ moves $m_1 + m_2 = 1$ and $m_2^* = 0$, in $j_1$ we obtain

$$(\ell + 1 - \gamma)A_3 - \gamma + (-1)^{\gamma+1}(2 + m_2 + 2m_1(1 - m_2) + 2x_1)$$
$$= (\ell_1 + 1 - \gamma_1)A_3 - \gamma_1 + (-1)^{\gamma_1}(-2 + m_1^* + (1 - m_2^*)m_1^*)$$
$$\implies (\ell - \ell_1) + (\gamma_1 - \gamma) + (-1)^{\gamma+1} = \frac{(-1)^\gamma(m_1 - 1 - |\gamma_1 - \gamma|) + 2(-1)^{\gamma_1}(m_1^* - 1)}{A_3}.$$

Observing that $(\ell - \ell_1) + (\gamma_1 - \gamma) + (-1)^{\gamma+1} \in \mathbb{Z}$ by construction, the above implies $m_1 = m_1^*$ in either case that $|\gamma_1 - \gamma| \in \{0, 1\}$. Applying this to $j_2$, we see

$$2 + m_1 + \gamma + 4x_2 = m_2^* + \gamma_1 + 4x_2^* + (1 - m_2^*)m_1^*$$
$$\implies x_2 - x_2^* = \frac{-2 + (\gamma_1 - \gamma)}{4}.$$

Since $x_2 - x_2^* \in \mathbb{Z}$ and $-3 \leq 2 + (\gamma_1 - \gamma) \leq -1$ by construction, it follows that $\frac{-2+(\gamma_1-\gamma)}{4} \notin \mathbb{Z}$, giving us a contradiction as the above equality is not possible.

**Set Comparison Case 5:** We will consider whether edges are shared from moves with the same orientation via sets 2 and 4.

Major Case 1: Focusing on moves along $j_1$, we see that $C_2$ has $j_1$ moves with the same orientation as those of $C_1$ if and only if $\gamma = \gamma_1$. Observing that $m_1 = m_2$ and $m_2^* = 1$ during all such moves, in $j_1$ we have

$$(\ell + 1 - \gamma)A_3 - \gamma + (-1)^{\gamma+1}(2 + m_2 + 2m_1(1 - m_2) + 2x_1)$$
$$= (\ell_1 + 1 - \gamma_1)A_3 - \gamma_1 + (-1)^{\gamma_1+1}(m_1^* + m_1^*(1 - m_2^*))$$
$$\implies \ell - \ell_1 = \frac{(-1)^{\gamma}(2 + 2x_1 + (m_1 - m_1^*))}{A_3}.$$

Given $\ell - \ell_1 \in \mathbb{Z}$ and $1 \leq 2 + 2x_1 + (m_1 - m_1^*) \leq A_3 - 1$ by construction, it must be the case that for all $0 \leq x_1 \leq \frac{A_3}{2} - 2$ and $m_1, m_1^* \in \{0, 1\}$, $\frac{(-1)^{\gamma}(2+2x_1+(m_1-m_1^*))}{A_3} \notin \mathbb{Z}$, giving us a contradiction as the above equality is not possible.

Major Case 2: Focusing on moves along $j_2$, observe that which $j_2$ moves in $C_1$ have the same orientation as those in $C_2$ depends on $x_1$. Hence, we case on $x_1$:

Case 1: Let $0 \leq x_1 < \frac{A_3}{2} - 2$. In this case, note that $m_1 + m_2 = 1$ and $(m_1^*, m_2^*) = (m_1, 0)$. Consequently, in $j_2$ we find

$$2 + m_1 + \gamma + 4x_2 = 2 + m_1^* + \gamma_1 + m_2^*(1 - m_1^*) + 4x_2^*$$
$$\implies x_2 - x_2^* = \frac{\gamma_1 - \gamma}{4}.$$

Given $x_2 - x_2^* \in \mathbb{Z}$ and $|\gamma_1 - \gamma| \leq 1$ by construction, it follows that $\gamma = \gamma_1$. Applying this to $j_1$, we see

$$(\ell + 1 - \gamma)A_3 - \gamma + (-1)^{\gamma+1}(2 + m_2 + 2m_1(1 - m_2) + 2x_1)$$
$$= (\ell_1 + 1 - \gamma_1)A_3 - \gamma_1 + (-1)^{\gamma_1+1}(m_1^* + m_1^*(1 - m_2^*))$$
$$\implies \ell - \ell_1 = \frac{(-1)^{\gamma}(2 + m_2 + 2x_1)}{A_3}.$$

Since $\ell - \ell_1 \in \mathbb{Z}$ and $2 \leq 2 + m_2 + 2x_1 < A_3 - 1$ by construction, it follows that for every $\gamma, m_2 \in \{0, 1\}$ and $0 \leq x_1 < \frac{A_3}{2} - 2$, $\frac{(-1)^{\gamma}(2+m_2+2x_1)}{A_3} \notin \mathbb{Z}$. So we have a contradiction as the above equality is not possible.

Case 2: Let $x_1 = \frac{A_3}{2} - 2$. Then, during all $j_2$ moves it is the case $m_1 + m_2 = 1$ and $(m_1^*, m_2^*) = (0, 0)$, and so in $j_2$ we have

$$2 + m_1 + \gamma + 4x_2 = 2 + m_1^* + \gamma_1 + m_2^*(1 - m_1^*) + 4x_2^*$$
$$\implies x_2 - x_2^* = \frac{(\gamma_1 - \gamma) - m_1}{4}.$$

Following from $x_2 - x_2^* \in \mathbb{Z}$ and $-2 \leq (\gamma_1 - \gamma) - m_1 \leq 1$ by construction, it must be the case $m_1 = (\gamma_1 - \gamma)$ with $\gamma_1 - \gamma \in \{0, 1\}$. Note that $\gamma_1 - \gamma = -1$ would give us a contradiction as for every $m_1 \in \{0, 1\}$, we would have $-\frac{1+m_1}{4} \notin \mathbb{Z}$. Applying our previous result to $j_1$ for the remaining two cases with $\gamma_1 - \gamma \in \{0, 1\}$, we obtain

$$(\ell + 1 - \gamma)A_3 - \gamma + (-1)^{\gamma+1}(2 + m_2 + 2m_1(1 - m_2) + 2x_1)$$
$$= (\ell_1 + 1 - \gamma_1)A_3 - \gamma_1 + (-1)^{\gamma_1+1}(m_1^* + m_1^*(1 - m_2^*))$$
$$\implies (\ell - \ell_1) + (\gamma_1 - \gamma) + (-1)^{\gamma+1} = \frac{(-1)^{\gamma+1}(1 - (\gamma_1 - \gamma)) - (\gamma_1 - \gamma)}{A_3}.$$

Since $(\ell - \ell_1) + (\gamma_1 - \gamma) + (-1)^{\gamma+1} \in \mathbb{Z}$ and for every $\gamma_1 - \gamma \in \{0, 1\}$, $\frac{(-1)^{\gamma+1}(1-(\gamma_1-\gamma))-(\gamma_1-\gamma)}{A_3} \notin \mathbb{Z}$, we have a contradiction as the above equality is not possible.

**Set Comparison Case 6:** We will consider whether edges are shared from moves with the same orientation via sets 3 and 4.

Major Case 1: Focusing on moves along $j_1$, we see that $C_2$ has $j_1$ moves with the same orientation as those in $C_1$ if and only if $|\gamma_1 - \gamma| = 1$. Then, given $m_2 = 1 = m_2^*$ during all such moves, in $j_1$ we get

$$(\ell+1-\gamma)A_3 - \gamma + (-1)^{\gamma}(-2 + m_1 + (1-m_2)m_1) = (\ell_1 + 1 - \gamma_1)A_3 - \gamma_1 + (-1)^{\gamma_1+1}(m_1^* + (1-m_2^*)m_1^*)$$
$$\implies (\ell - \ell_1) + (-1)^{\gamma} = \frac{(-1)^{\gamma}(1 + (m_1^* - m_1))}{A_3}.$$

Since $(\ell - \ell_1) + (-1)^{\gamma} \in \mathbb{Z}$ and $0 \leq 1 + (m_1^* - m_1) \leq 2$ by construction, we must have $m_1^* - m_1 = -1$, meaning $m_1^* = 0$ and $m_1 = 1$. Applying this to $j_2$, we find

$$m_2 + \gamma + 4x_2 + (1-m_2)m_1 = 2 + m_1^* + \gamma_1 + m_2^*(1 - m_1^*) + 4x_2^*$$
$$\implies x_2 - x_2^* = \frac{2 + (-1)^{\gamma}}{4}.$$

Given $x_2 - x_2^* \in \mathbb{Z}$ and for every $\gamma \in \{0, 1\}$, $\frac{2+(-1)^{\gamma}}{4} \notin \mathbb{Z}$ by construction, we have a contradiction as the above equality is not possible.

Major Case 2: Focusing on moves along $j_2$, we see that only $C_2$'s $j_2$ move $(m_1^*, m_2^*) = (0, 0)$ has the same orientation as those in $C_1$ and during such moves $m_2 = 0$. Applying this to $j_2$, we get

$$m_2 + \gamma + 4x_2 + (1-m_2)m_1 = 2 + m_1^* + \gamma_1 + (1-m_1^*)m_2^* + 4x_2^*$$
$$\implies x_2 - x_2^* = \frac{2 + (\gamma_1 - \gamma) - m_1}{4}.$$

Following from $x_2 - x_2^* \in \mathbb{Z}$ and $0 \leq 2 + (\gamma_1 - \gamma) - m_1 \leq 3$ by construction, it must be the case $m_1 = 1$ and $\gamma_1 - \gamma = -1$, meaning $\gamma = 1$ and $\gamma_1 = 0$. Applying this to $j_1$, it follows that

$$(\ell+1-\gamma)A_3 - \gamma + (-1)^{\gamma}(-2 + m_1 + (1-m_2)m_1) = (\ell_1 + 1 - \gamma_1)A_3 - \gamma_1 + (-1)^{\gamma_1+1}(m_1^* + m_1^*(1-m_2^*))$$
$$\implies (\ell - \ell_1) - 1 = \frac{1}{A_3}.$$

Given $(\ell - \ell_1) - 1 \in \mathbb{Z}$ and $\frac{1}{A_3} \notin \mathbb{Z}$ since $A_3 \geq 4$, we have a contradiction as the above equality is not possible.

Thus, no edges are shared from moves with the same orientation via two distinct sets for $\alpha_1 < \alpha^* < \alpha_2$ when $d = 2$.

**$\alpha^*$−Case 3:** Let $\alpha^* = \alpha_2$. Then, $A_1 = 1$ and $A_2 = 1$, meaning only sets 1 and 3 are active. Note that $A_3 \geq 4$ is even.

**Set Comparison Case 1:** We will consider whether edges are shared from moves with the same orientation via sets 1 and 3.

Major Case 1: Focusing on moves along $j_1$, in $j_1$ we get

$$\ell A_3 + (m_1 + 2x_1) = (\ell_1 + 1)A_3 + (-2 + m_1^*)$$
$$\implies \ell - \ell_1 - 1 = -\frac{2 + 2x_1 + (m_1 - m_1^*)}{A_3}.$$

Since $\ell - \ell_1 - 1 \in \mathbb{Z}$ and $1 \leq 2 + 2x_1 + (m_1 - m_1^*) \leq A_3 - 1$ by construction, it follows that for every $0 \leq x_1 \leq \frac{A_3}{2} - 2$ and $m_1, m_1^* \in \{0, 1\}$, $-\frac{2+2x_1+(m_1-m_1^*)}{A_3} \notin \mathbb{Z}$, giving us a contradiction as the above equality is not possible.

Major Case 2: Focusing on moves along $j_2$, we see that all of $C_2$'s $j_2$ moves have the same orientation as $C_1$'s $j_2$ move $(m_1, m_2) = (0, 0)$. Applying this to $j_1$, we find

$$\ell A_3 + (m_1 + 2x_1) = (\ell_1 + 1)A_3 + (-2 + m_1^*)$$
$$\implies \ell - \ell_1 - 1 = -\frac{2 + 2x_1 - m_1^*}{A_3}.$$

Since $\ell - \ell_1 - 1 \in \mathbb{Z}$ and $1 \leq 2 + 2x_1 - m_1^* \leq A_3 - 2$ by construction, it follows that for every $0 \leq x_1 \leq \frac{A_3}{2} - 2$ and $m_1^* \in \{0, 1\}$, $-\frac{2+2x_1-m_1^*}{A_3} \notin \mathbb{Z}$, giving us a contradiction as the above equality is not possible.

Hence, no edges are shared from moves with the same orientation via two distinct sets for $\alpha^* = \alpha_2$ when $d = 2$.

**Dimension Case 2:** Let $d \geq 3$. Then, $C_1 = C_{\ell,\gamma,t,p_1,s_1,\ldots,p_{d-2},s_{d-2}}$ and $C_2 = C_{\ell_1,\gamma_1,t_1,p_1^*,s_1^*,\ldots,p_{d-2}^*,s_{d-2}^*}$. We now case on $\alpha_1 \leq \alpha^* \leq \alpha_d$.

**$\alpha^*$−Case 1:** Let $\alpha^* = \alpha_1$. In this case, only sets 3 and 4 are active as $A_1 = 0$ and $A_2 = 0$. In particular, note that $p_1 = 0 = p_1^*$, $A_4 = 1$, and $x_3 = 0 = x_3^*$.

Major Case 1: Focusing on moves along $j_1$, we see that $C_2$ has $j_1$ moves with the same orientation as those in $C_1$ if and only if $|\gamma_1 - \gamma| = 1$. Observing that during such moves $m_1 = m_d$, $m_2 = 1 = m_2^*$ and $m_1^* = m_d^*$, in $j_2$ we find

$$m_2 + \gamma + 4x_2 + (1 - m_2)m_1 + 2(1 - m_1)m_d = 2 + m_1^* + \gamma_1 + (1 - m_1^*)m_2^* + 4x_2^* + 2(1 - m_1^*)m_d^*$$
$$\implies x_2 - x_2^* = \frac{2 + (-1)^\gamma}{4}.$$

Since $x_2 - x_2^* \in \mathbb{Z}$ by construction and for every $\gamma \in \{0,1\}$, $\frac{2+(-1)^\gamma}{4} \notin \mathbb{Z}$, we have a contradiction as the above equality is not possible.

<u>Major Case 2</u>: Focusing on moves along $j_2$, observe that during all such moves $m_1 = m_3 = m_d$, $m_1^* = m_3^* = m_d^*$, and $m_2 = 0 = m_2^*$. These observations in $j_2$ give us

$$m_2 + \gamma + 4x_2 + (1-m_2)m_1 + 2(1-m_1)m_d = 2 + m_1^* + \gamma_1 + (1-m_1^*)m_2^* + 4x_2^* + 2(1-m_1^*)m_d^*$$

$$\implies x_2 - x_2^* = \frac{2 + (\gamma_1 - \gamma) + (m_1^* - m_1)}{4}.$$

Since $x_2 - x_2^* \in \mathbb{Z}$ and $0 \leq 2 + (\gamma_1 - \gamma) + (m_1^* - m_1) \leq 4$ by construction, we must have in particular $(\gamma_1 - \gamma) = (m_1^* - m_1)$ with $|\gamma_1 - \gamma| = |m_1^* - m_1| = 1$. Now, observe that in $j_3$ we get

$$p_1 + 2A_4 s_1 + (-1)^\gamma m_3 + 2x_3 = 2s_1^* + (-1)^{\gamma_1} m_3^*$$

$$\implies s_1 - s_1^* = \frac{(-1)^{\gamma+1}(m_3 + m_3^*)}{2}.$$

Since $s_1 - s_1^* \in \mathbb{Z}$ and $0 \leq m_3 + m_3^* \leq 2$ with $m_3, m_3^* \in \{0,1\}$ by construction, we see $m_3 = m_3^*$. This then implies $m_1 = m_3 = m_3^* = m_1^*$, but this contradicts our deduction $|m_1^* - m_1| = 1$ as this implies $m_1 \neq m_1^*$. Hence, we have a contradiction as the above is not possible.

<u>Major Case 3</u>: Focusing on moves along $j_3$, we see that during $j_3$ moves $m_1 = m_2$, $m_3 = m_d = 1 - m_1$, $m_1^* = m_2^*$, and $m_3^* = m_d^* = 1 - m_1^*$. Applying this to $j_2$, we have

$$m_2 + \gamma + 4x_2 + (1-m_2)m_1 + 2(1-m_1)m_d = 2 + m_1^* + \gamma_1 + (1-m_1^*)m_2^* + 4x_2^* + 2(1-m_1^*)m_d^*$$

$$\implies x_2 - x_2^* = \frac{2 + (\gamma_1 - \gamma) + (m_d^* - m_d)}{4}.$$

Since $x_2 - x_2^* \in \mathbb{Z}$ and $0 \leq 2 + (\gamma_1 - \gamma) + (m_d^* - m_d) \leq 4$ by construction, it must be the case $(\gamma_1 - \gamma) = (m_d^* - m_d)$ with $|\gamma_1 - \gamma| = |m_d^* - m_d| = 1$. Applying this to $j_3$, we find

$$p_1 + 2A_4 s_1 + (-1)^\gamma m_3 + 2x_3 = 2s_1^* + (-1)^{\gamma_1} m_3^*$$

$$\implies s_1 - s_1^* = \frac{(-1)^{\gamma+1}(m_3 + m_3^*)}{2}.$$

Following from $s_1 - s_1^* \in \mathbb{Z}$ and $0 \leq m_3 + m_3^* \leq 2$ with $m_3, m_3^* \in \{0,1\}$ by construction, it must be the case $m_3 = m_3^*$. This implies $m_d = m_3 = m_3^* = m_d^*$, contradicting our assumption $|m_d^* - m_d| = 1$ that implies $m_d \neq m_d^*$. Consequently, the above is not possible.

<u>Major Case 4</u>: Focusing on moves along $j_k$ for $4 \leq k \leq d$ when $d \geq 4$, we se that during all such moves $m_1 = m_2$, $m_d = 1 - m_1$ and $(m_1^*, m_2^*, m_d^*) = (m_1, m_2, m_d)$, giving us in $j_2$

$$m_2 + \gamma + 4x_2 + (1-m_2)m_1 + 2(1-m_1)m_d = 2 + m_1^* + \gamma_1 + (1-m_1^*)m_2^* + 4x_2^* + 2(1-m_1^*)m_d^*$$

$$\implies x_2 - x_2^* = \frac{2 + (\gamma_1 - \gamma)}{4}.$$

Following from $x_2 - x_2^* \in \mathbb{Z}$ and $|\gamma_1 - \gamma| \leq 1$ by construction, it follows that for every $\gamma, \gamma_1 \in \{0, 1\}$, $\frac{2+(\gamma_1-\gamma)}{4} \notin \mathbb{Z}$, giving us a contradiction as the above equality is not possible.

Thus, no edges are shared from moves with the same orientation via two distinct sets for $\alpha^* = \alpha_1$ when $d \geq 3$.

**$\alpha^*$−Case 2:** Let $\alpha_1 < \alpha^* < \alpha_2$. Then, $A_1 = 1$ and $A_2 = 0$, meaning sets $1, 2, 3$, and $4$ are all active. Note that $A_3 \geq 4$ is even, $p_1 = 0 = p_1^*$, $t = 0 = t_1$, $A_4 = 1$, and $x_3 = 0 = x_3^*$.

**Set Comparison Case 1:** We will consider whether edges are shared from moves with the same orientation via sets 1 and 2.

Major Case 1: Focusing on moves along $j_1$, we see that $C_2$ has no $j_1$ moves with the same orientation as those of $C_1$ to consider.

Major Case 2: Focusing on moves along $j_2$, observe that which $j_2$ moves in $C_2$ have the same orientation as those in $C_1$ depends on $x_1^*$. So we proceed by casing on $x_1^*$:

Case 1: Let $0 \leq x_1^* < \frac{A_3}{2} - 2$. Then, it is the case that $m_1 = m_2 = m_3 = m_d$, $(m_1^*, m_2^*, m_d^*) = (1 - m_1, 1 - m_2, m_1)$, and $r_1^* = 0$ since $x_1^* \neq \frac{A_3}{2} - 2$. So our equation in $j_3$ implies

$$p_1 + A_4 s_1 + (1 - m_3) + 2x_3 = s_1^* + m_1^* + r_1^*$$
$$\implies s_1 = s_1^*.$$

Applying this to $j_2$, we find

$$m_2 + \gamma + (s_1 + 1) + 4x_2 = 2 + \gamma_1 + (1 + s_1^* + m_d^* - r_1^*) + 4x_2^*$$
$$\implies x_2 - x_2^* = \frac{2 + (\gamma_1 - \gamma)}{4}.$$

Following from $x_2 - x_2^* \in \mathbb{Z}$ and $1 \leq 2 + (\gamma_1 - \gamma) \leq 3$ by construction, we get the for every $\gamma, \gamma_1 \in \{0, 1\}$, $\frac{2+(\gamma_1-\gamma)}{4} \notin \mathbb{Z}$, giving us a contradiction as the above equality is not possible.

Case 2: Let $x_1^* = \frac{A_3}{2} - 2$. Then, observe that all of $C_2$'s $j_2$ moves have the same orientation only as $C_1$'s $j_2$ move with $(m_1, m_2, m_3, m_d) = (0, 0, 0, 0)$, and during such moves $m_1^* = m_2^*$ and $m_d^* = 1 - m_1^*$. Hence, in $j_3$ we see

$$p_1 + A_4 s_1 + (1 - m_3) + 2x_3 = s_1^* + m_1^* + r_1^*$$
$$\implies s_1 - s_1^* = m_1^* + r_1^* - 1.$$

Now, invoking the above, our equation in $j_2$ implies

$$m_2 + \gamma + (s_1 + 1) + 4x_2 = 2 + \gamma_1 + (1 + s_1^* + m_d^* - r_1^*) + 4x_2^*$$
$$\implies 2(x_2 - x_2^*) + (m_1^* + r_1^*) - 2 = \frac{\gamma_1 - \gamma}{2}.$$

Since $2(x_2 - x_2^*) + (m_1^* + r_1^*) - 2 \in \mathbb{Z}$ and $|\gamma_1 - \gamma| \leq 1$ by construction, we must have $\gamma = \gamma_1$. This leaves us with

$$(x_2 - x_2^*) - 1 = -\frac{m_1^* + r_1^*}{2}.$$

Given $(x_2 - x_2^*) - 1 \in \mathbb{Z}$ and $0 \leq m_1^* + r_1^* \leq 2$ by construction, it follows that $m_1^* = r_1^*$. This implies $m_1^* = m_3^*(1 - m_1^*)$ when $d = 3$. Noting that $m_d^* = m_3^* = 1 - m_1^*$ in this case, our result $m_1^* = r_1^*$ implies $m_1^* = \frac{1}{2}$, but this is not possible as $m_1^* \in \{0, 1\}$ by construction. If $d \geq 4$, observe that $m_1^* = m_d^*(1 - m_{d-1}^*)$, $m_1^* = m_{d-1}^*$, $m_d^* = 1 - m_1^*$, and so we see $m_1^* = \frac{1}{2}$. However this is not possible as $m_1^* \in \{0, 1\}$ by construction.

Thus, for all $d \geq 3$, we see that the above yields a contradiction.

<u>Major Case 3</u>: Focusing on moves along $j_3$, observe that during all such moves $m_1 = m_2$, $m_3 = m_d = 1 - m_1$, $(m_1^*, m_2^*, m_d^*) = (m_1, 1 - m_2, 1 - m_d)$, and $r_1^* = 0$ since $m_1^* = m_3^*$ when $d = 3$ and $m_{d-1}^* = m_d^*$ when $d \geq 4$. Applying this to $j_3$, we find

$$p_1 + A_4 s_1 + (1 - m_3) + 2x_3 = s_1^* + m_1^* + r_1^*$$
$$\implies s_1 = s_1^*.$$

Hence, in $j_2$ we see

$$m_2 + \gamma + (s_1 + 1) + 4x_2 = 2 + \gamma_1 + (1 + s_1^* + m_d^* - r_1^*) + 4x_2^*$$
$$\implies x_2 - x_2^* = \frac{2 + (\gamma_1 - \gamma)}{4}.$$

Given $x_2 - x_2^* \in \mathbb{Z}$ and $1 \leq 2 + (\gamma_1 - \gamma) \leq 3$ by construction, it follows that for every $\gamma, \gamma_1 \in \{0, 1\}$, $\frac{2 + (\gamma_1 - \gamma)}{4} \notin \mathbb{Z}$, giving us a contradiction as the above equality is not possible.

<u>Major Case 4</u>: Focusing on moves along $j_k$ for $4 \leq k \leq d$ when $d \geq 4$, we see that during all such moves it is the case $m_1 = m_2 = m_3$, $m_d = 1 - m_1$, and $(m_1^*, m_2^*, m_d^*) = (1 - m_2, 1 - m_3, 1 - m_d)$. Consequently, noting that $r_1^* = 0$ since $m_{d-1}^* = m_d^*$, in $j_3$ we get

$$p_1 + A_4 s_1 + (1 - m_3) + 2x_3 = s_1^* + m_1^* + r_1^*$$
$$\implies s_1 = s_1^*.$$

Applying this to $j_2$, it follows that

$$m_2 + \gamma + (s_1 + 1) + 4x_2 = 2 + \gamma_1 + (1 + s_1^* + m_d^* - r_1^*) + 4x_2^*$$
$$\implies x_2 - x_2^* = \frac{2 + (\gamma_1 - \gamma)}{4}.$$

Since $x_2 - x_2^* \in \mathbb{Z}$ and $1 \leq 2 + (\gamma_1 - \gamma) \leq 3$ by construction, it follows that for every $\gamma, \gamma_1 \in \{0, 1\}$, $\frac{2 + (\gamma_1 - \gamma)}{4} \notin \mathbb{Z}$, giving us a contradiction as the equality above is not possible.

**Set Comparison Case 2:** We will consider whether edges are shared from moves with the same orientation via sets 1 and 3.

<u>Major Case 1</u>: Focusing on moves along $j_1$, observe that all $j_1$ moves have the same orientation and that during these moves $m_1 = m_d$, $m_2 = 1 - m_1$, $m_1^* = m_d^*$, and $m_2^* = 1 - m_1^*$. Applying these observations to $j_1$, we have

$$\ell A_3 + 1 + (m_1 + 2x_1 + 2m_d(1-m_1)) = (\ell_1 + 1)A_3 + (-2 + m_1^* + 2m_d^*(1-m_1^*) - (A_3 - 2)(1-m_d^*))$$

$$\implies (\ell - \ell_1) + (m_2^* - 1) = -\frac{1 + 2x_1 + (m_1 + m_1^*)}{A_3}.$$

Given $(\ell - \ell_1) + (m_2^* - 1) \in \mathbb{Z}$ and $1 \leq 1 + 2x_1 + (m_1 + m_1^*) \leq A_3 - 1$ by construction, it is the case that for all $0 \leq x_1 \leq \frac{A_3}{2} - 2$ and $m_1, m_1^* \in \{0, 1\}$, $-\frac{1+2x_1+(m_1+m_1^*)}{A_3} \notin \mathbb{Z}$, giving us a contradiction as the above equality is not possible.

Major Case 2: Focusing on moves along $j_2$, we see that all of $C_2$'s $j_2$ moves have the same orientation only as $C_1$'s $j_2$ move with $(m_1, m_2, m_d) = (0, 0, 0)$, and during all such moves $m_1^* = m_2^* = m_d^*$. Hence, in $j_1$ we obtain

$$\ell A_3 + 1 + (m_1 + 2x_1 + 2m_d(1-m_1)) = (\ell_1 + 1)A_3 + (-2 + m_1^* + 2m_d^*(1-m_1^*) - (A_3 - 2)(1-m_d^*))$$

$$\implies \ell - \ell_1 - m_1^* = -\frac{1 + 2x_1 + m_1^*}{A_3}.$$

Following from $\ell - \ell_1 - m_1^* \in \mathbb{Z}$ and $1 \leq 1 + 2x_1 + m_1^* \leq A_3 - 2$ by construction, we see that for every $0 \leq x_1 \leq \frac{A_3}{2} - 2$ and $m_1^* \in \{0, 1\}$, $-\frac{1+2x_1+m_1^*}{A_3} \notin \mathbb{Z}$, giving us a contradiction as the above equality is not possible.

Major Case 3: Focusing on moves along $j_k$ for $3 \leq k \leq d$, observe that during all such moves $m_1 = m_2, m_d = 1 - m_1$, and $(m_1^*, m_2^*, m_d^*) = (1 - m_1, 1 - m_2, 1 - m_d)$. Consequently, in $j_1$ we obtain

$$\ell A_3 + 1 + (m_1 + 2x_1 + 2m_d(1-m_1)) = (\ell_1 + 1)A_3 + (-2 + m_1^* + 2m_d^*(1-m_1^*) - (A_3 - 2)(1-m_d^*))$$

$$\implies (\ell - \ell_1) + (m_1^* - 1) = -\frac{2(1 + x_1)}{A_3}.$$

Since $(\ell - \ell_1) + (m_1^* - 1) \in \mathbb{Z}$ and $2 \leq 2(1 + x_1) \leq A_3 - 2$ by construction, it follows that for all $0 \leq x_1 \leq \frac{A_3}{2} - 2$, $-\frac{2(1+x_1)}{A_3} \notin \mathbb{Z}$, giving us a contradiction as the above equality is not possible.

**Set Comparison Case 3:** We will consider whether edges are shared from moves with the same orientation via sets 1 and 4.

Major Case 1: Focusing on moves along $j_1$, we see that $C_2$ is not defined by any $j_1$ moves with the same orientation as those in $C_1$, and so there are no $j_1$ moves in $C_2$ to consider.

Major Case 2: Focusing on moves along $j_2$, observe that during such moves $m_1 = m_2 = m_3 = m_d$ and $(m_1^*, m_2^*, m_d^*) = (0, 0, m_d)$, giving us in $j_3$

$$p_1 + A_4 s_1 + (1 - m_3) + 2x_3 = s_1^* + m_1^* m_3^*$$

$$\implies s_1 - s_1^* = m_3 - 1.$$

Following from the above, in $j_2$ we see

$$m_2 + \gamma + (s_1 + 1) + 4x_2 = 2 + m_1^* + (s_1^* + (1 - m_1^*)m_d^* + m_d^*) + \gamma_1 + m_2^*(1 - m_1^*) + 4x_2^*$$

$$\implies x_2 - x_2^* = \frac{2 + (\gamma_1 - \gamma)}{4}.$$

Given $x_2 - x_2^* \in \mathbb{Z}$ and $1 \leq 2 + (\gamma_1 - \gamma) \leq 3$ by construction, it is the case that for all $\gamma, \gamma_1 \in \{0, 1\}$, $\frac{2+(\gamma_1-\gamma)}{4} \notin \mathbb{Z}$. Hence, we have contradiction as the above equality is not possible.

<u>Major Case 3</u>: Focusing on moves along $j_3$, we see that $m_1 = m_2$, $m_3 = m_d = 1 - m_1$, and $(m_1^*, m_2^*, m_d^*) = (1, 1 - m_2, 1 - m_d)$ with $m_3^* = m_d^*$. Consequently, in $j_3$ we get

$$p_1 + A_4 s_1 + (1 - m_3) + 2x_3 = s_1^* + m_1^* m_3^*$$

$$\implies s_1 = s_1^*.$$

Applying this to $j_2$, we have

$$m_2 + \gamma + (s_1 + 1) + 4x_2 = 2 + m_1^* + (s_1^* + (1 - m_1^*)m_d^* + m_d^*) + \gamma_1 + m_2^*(1 - m_1^*) + 4x_2^*$$

$$\implies x_2 - x_2^* = \frac{2 + (\gamma_1 - \gamma)}{4}.$$

Since $x_2 - x_2^* \in \mathbb{Z}$ and $1 \leq 2 + (\gamma_1 - \gamma) \leq 3$ by construction, it must be the case that for every $\gamma, \gamma_1 \in \{0, 1\}$, $\frac{2+(\gamma_1-\gamma)}{4} \notin \mathbb{Z}$, giving us a contradiction as the above equality is not possible.

<u>Major Case 4</u>: Focusing on moves along $j_k$ for all $4 \leq k \leq d$ when $d \geq 4$, we see that during all of these moves $m_1 = m_2 = m_3$, $m_d = 1 - m_1$ and $(m_1^*, m_2^*, m_3^*, m_d^*) = (1 - m_1, 1 - m_2, m_d + (1 - m_d)(1 - m_4), 1 - m_d)$. Applying these observations to $j_3$, we find

$$p_1 + A_4 s_1 + (1 - m_3) + 2x_3 = s_1^* + m_1^* m_3^*$$

$$\implies s_1 = s_1^*.$$

Applying this to $j_2$, we get

$$m_2 + \gamma + (s_1 + 1) + 4x_2 = 2 + m_1^* + (s_1^* + (1 - m_1^*)m_d^* + m_d^*) + \gamma_1 + m_2^*(1 - m_1^*) + 4x_2^*$$

$$\implies x_2 - x_2^* = \frac{2 + (\gamma_1 - \gamma)}{4}.$$

Since $x_2 - x_2^* \in \mathbb{Z}$ and $1 \leq 2 + (\gamma_1 - \gamma) \leq 3$ by construction, it must be the case that for all $\gamma, \gamma_1 \in \{0, 1\}$, $\frac{2+(\gamma_1-\gamma)}{4} \notin \mathbb{Z}$, giving us a contradiction as the above equality is not possible.

**Set Comparison Case 4:** We will consider whether edges are shared from moves with the same orientation via sets 2 and 3.

<u>Major Case 1</u>: Focusing on moves along $j_1$, observe that there are no $j_1$ moves to consider as there are no $j_1$ moves in $C_2$ with the same orientation as those in $C_1$.

<u>Major Case 2</u>: Focusing on moves along $j_2$, we see that which $j_2$ moves in $C_1$ have the same orientation as those in $C_2$ depends on $x_1$. Hence, we proceed by casing on $x_1$:

<u>Case 1</u>: Let $0 \leq x_1 < \frac{A_3}{2} - 2$. Then, $(m_1, m_2, m_d) = (1, 1, 0)$ during the only $j_1$ move in $C_1$ with the same orientation as the $j_1$ moves in $C_2$, and $m_1^* = m_2^* = m_d^*$ during all such moves. Noting that $r_1 = 0$ since $x_1 \neq \frac{A_3}{2} - 2$, the above in $j_1$ yields

$$(\ell + 1)A_3 - (2 + m_2 + 2m_1(1 - m_2) + 2x_1 + 2m_d(1 - m_1) - r_1(A_3 - 1))$$
$$= (\ell_1 + 1)A_3 + (-2 + m_1^* + 2m_d^*(1 - m_1^*) - (A_3 - 2)(1 - m_d^*))$$
$$\implies (\ell - \ell_1) + (1 - m_d^*) = \frac{3 + 2x_1 - m_1^*}{A_3}.$$

Given $(\ell - \ell_1) + (1 - m_d^*) \in \mathbb{Z}$ and $2 \leq 3 + 2x_1 - m_1^* < A_3 - 1$ by construction and our case assumption, we get that for every $0 \leq x_1 < \frac{A_3}{2} - 2$ and $m_1^* \in \{0, 1\}$, $\frac{3 + 2x_1 - m_1^*}{A_3} \notin \mathbb{Z}$. Thus, we have a contradiction as the equality above is not possible.

<u>Case 2:</u> Let $x_1 = \frac{A_3}{2} - 2$. Then, all of $C_1$'s $j_2$ moves have the same orientation as those in $C_2$, and $m_1 = \cdots = m_{d-1}$, $m_d = 1 - m_1$ and $m_1^* = \cdots = m_d^*$ during all such moves. Applying this to $j_1$, we obtain

$$(\ell + 1)A_3 - (2 + m_2 + 2m_1(1 - m_2) + 2x_1 + 2m_d(1 - m_1) - r_1(A_3 - 1))$$
$$= (\ell_1 + 1)A_3 + (-2 + m_1^* + 2m_d^*(1 - m_1^*) - (A_3 - 2)(1 - m_d^*))$$
$$\implies (\ell - \ell_1) + r_1 - m_d^* = \frac{r_1 + m_d - (1 + m_1^*)}{A_3}.$$

Since $(\ell - \ell_1) + r_1 - m_d^* \in \mathbb{Z}$ and $-2 \leq r_1 + m_d - (1 + m_1^*) \leq 2$ with $m_1^*, m_d, r_1 \in \{0, 1\}$, it must be the case $r_1 = 1 + m_1^* - m_d$. Applying this to $j_3$, we get

$$s_1 + m_1 + r_1 = p_1^* + A_4 s_1^* + m_3^* + 2x_3^*$$
$$\implies s_1 - s_1^* = -2m_1.$$

Lastly, in $j_2$ we have

$$2 + \gamma + (1 + s_1 + m_d - r_1) + 4x_2 = m_2^* + \gamma_1 + 2(1 - m_2^*)m_1^* + 4x_2^* + (s_1^* + 2(1 - m_1^*)m_d^*)$$
$$\implies 2(x_2 - x_2^*) - 2m_1 - m_1^* + 2 = \frac{\gamma_1 - \gamma}{2}.$$

Since $2(x_2 - x_2^*) - 2m_1 - m_1^* + 2 \in \mathbb{Z}$ and $|\gamma_1 - \gamma| \leq 1$ by construction, it follows that $\gamma = \gamma_1$. This leaves us with

$$(x_2 - x_2^*) - m_1 + 1 = \frac{m_1^*}{2}.$$

Given $(x_2 - x_2^*) - m_1 + 1 \in \mathbb{Z}$ and $m_1^* \in \{0, 1\}$ by construction, it must be the case $m_1^* = 0$, meaning $r_1 = 1 - m_d$. If $d = 3$, then $m_3(1 - m_1) = 1 - m_3$ implies $m_3 = \frac{1}{2}$, but this is not possible as $m_3 \in \{0, 1\}$ by construction. If $d \geq 4$, then $m_d(1 - m_{d-1}) = 1 - m_d$ implies $m_d = \frac{1}{2}$, but this is not possible either since $m_d \in \{0, 1\}$ by construction.

Thus, for all $d \geq 3$, we have a contradiction as the above is not possible.

<u>Major Case 3:</u> Focusing on moves along $j_3$, we see that during these moves $m_1 = m_d$, $m_2 = 1 - m_1$, $(m_1^*, m_2^*, m_d^*) = (1, m_2, m_d)$, and $m_3^* = m_d^*$. Noting that $r_1 = 0$ since $m_3 = m_1$ when $d = 3$ and $m_{d-1} = m_d$ when $d \geq 4$, in $j_1$ we get

$$(\ell + 1)A_3 - (2 + m_2 + 2m_1(1 - m_2) + 2x_1 + 2m_d(1 - m_1) - r_1(A_3 - 1))$$

$$= (\ell_1 + 1)A_3 + (-2 + m_1^* + 2m_d^*(1 - m_1^*) - (A_3 - 2)(1 - m_d^*))$$
$$\implies (\ell - \ell_1) + (1 - m_d^*) = \frac{4 + 2x_1 - m_1}{A_3}.$$

Since $(\ell - \ell_1) + (1 - m_d^*) \in \mathbb{Z}$ and $3 \leq 4 + 2x_1 - m_1 \leq A_3$ by construction, it follows that $x_1 = \frac{A_3}{2} - 2$ and $m_1 = 0$. In $j_3$, our assumptions imply

$$s_1 + m_1 + r_1 = p_1^* + A_4 s_1^* + m_3^* + 2x_3^*$$
$$\implies s_1 = s_1^*.$$

Hence, in $j_2$ we find

$$2 + \gamma + (1 + s_1 + m_d - r_1) + 4x_2 = m_2^* + \gamma_1 + 2(1 - m_2^*)m_1^* + 4x_2^* + (s_1^* + 2(1 - m_1^*)m_d^*)$$
$$\implies x_2 - x_2^* = \frac{-2 + (\gamma_1 - \gamma)}{4}.$$

Given $x_2 - x_2^* \in \mathbb{Z}$ and $-3 \leq -2 + (\gamma_1 - \gamma) \leq -1$ by construction, we get that $\frac{-2+(\gamma_1-\gamma)}{4} \notin \mathbb{Z}$, giving us a contradiction as the above equality is not possible.

<u>Major Case 4</u>: Focusing on moves along $j_k$ for $4 \leq k \leq d$ when $d \geq 4$, we see that during all such moves $m_1 = m_2$, $m_{d-1} = m_d = 1 - m_1$, $(m_1^*, m_2^*, m_d^*) = (m_1, m_2, m_d)$, and $m_3^* = m_1$. Noting that $r_1 = 0$ since $d \geq 4$ and $m_d = m_{d-1}$, in $j_3$ we find

$$s_1 + m_1 + r_1 = p_1^* + A_4 s_1^* + m_3^* + 2x_3^*$$
$$\implies s_1 = s_1^*.$$

Applying this to $j_2$, it follows that

$$2 + \gamma + (1 + s_1 + m_d - r_1) + 4x_2 = m_2^* + \gamma_1 + 2(1 - m_2^*)m_1^* + 4x_2^* + (s_1^* + 2(1 - m_1^*)m_d^*)$$
$$\implies x_2 - x_2^* = \frac{-2 + (\gamma_1 - \gamma)}{4}.$$

Following from $x_2 - x_2^* \in \mathbb{Z}$ and $-3 \leq -2 + (\gamma_1 - \gamma) \leq -1$ by construction, we see that for every $\gamma, \gamma_1 \in \{0, 1\}$, $\frac{-2+(\gamma_1-\gamma)}{4} \notin \mathbb{Z}$, giving us a contradiction as the above equality is not possible.

**Set Comparison Case 5:** We will consider whether edges are shared from moves with the same orientation via sets 2 and 4.

<u>Major Case 1</u>: Focusing on moves along $j_1$, we see that all of $C_2$'s $j_1$ moves have the same orientation as those in $C_1$, and during these moves $m_1 = m_2 = m_d$, $m_1^* = m_d^*$, and $m_2^* = 1$. Note that $r_1 = 0$ as $m_1 = m_3$ when $d = 3$ and $m_{d-1} = m_d$ when $d \geq 4$. Applying these observations to $j_1$, we get

$$(\ell + 1)A_3 - (2 + m_2 + 2m_1(1 - m_2) + 2x_1 + 2m_d(1 - m_1) - r_1(A_3 - 1))$$
$$= (\ell_1 + 1)A_3 - (m_1^* + m_1^*(1 - m_2^*) + 2(1 - m_1^*)m_d^*)$$
$$\implies \ell - \ell_1 = \frac{2 + 2x_1 + (m_1 - m_1^*)}{A_3}.$$

Given $\ell - \ell_1 \in \mathbb{Z}$ and $1 \leq 2 + 2x_1 + (m_1 - m_1^*) \leq A_3 - 1$ by construction, it must then be the case that for every $0 \leq x_1 \leq \frac{A_3}{2} - 2$ and $m_1, m_1^* \in \{0, 1\}$, $\frac{2+2x_1+(m_1-m_1^*)}{A_3} \notin \mathbb{Z}$, giving us a contradiction as the above equality is not possible.

Major Case 2: Focusing on moves along $j_2$, observe that which of $C_1$'s $j_2$ moves have the same orientation as those in $C_2$ depends on $x_1$. Consequently, we proceed by casing on $x_1$:

Case 1: Let $0 \leq x_1 < \frac{A_3}{2} - 2$. Then, observing that during these moves $m_1 = m_2$, $m_d = 1 - m_1$, $(m_1^*, m_2^*, m_d^*) = (0, 0, m_d)$ and $r_1 = 0$ since $x_1 \neq \frac{A_3}{2} - 2$, in $j_1$ we find

$$(\ell + 1)A_3 - (2 + m_2 + 2m_1(1 - m_2) + 2x_1 + 2m_d(1 - m_1) - r_1(A_3 - 1))$$
$$= (\ell_1 + 1)A_3 - (m_1^* + m_1^*(1 - m_2^*) + 2(1 - m_1^*)m_d^*)$$
$$\implies \ell - \ell_1 = \frac{2 + 2x_1 + m_2}{A_3}.$$

Since $\ell - \ell_1 \in \mathbb{Z}$ and $2 \leq 2 + 2x_1 + m_2 < A_3 - 1$ by construction and our case assumption, it follows that for every $0 \leq x_1 < \frac{A_3}{2} - 2$ and $m_2 \in \{0, 1\}$, $\frac{2+2x_1+m_2}{A_3} \notin \mathbb{Z}$, giving us a contradiction as the above equality is not possible.

Case 2: Let $x_1 = \frac{A_3}{2} - 2$. Then, observe that during all of these moves $m_1 = \cdots = m_{d-1}$, $m_d = 1 - m_1$, and $(m_1^*, m_2^*, m_d^*) = (0, 0, 0)$. So in $j_1$ we obtain

$$(\ell + 1)A_3 - (2 + m_2 + 2m_1(1 - m_2) + 2x_1 + 2m_d(1 - m_1) - r_1(A_3 - 1))$$
$$= (\ell_1 + 1)A_3 - (m_1^* + m_1^*(1 - m_2^*) + 2(1 - m_1^*)m_d^*)$$
$$\implies (\ell - \ell_1) + r_1 - 1 = \frac{m_d + r_1 - 1}{A_3}.$$

Following from $(\ell - \ell_1) + r_1 - 1 \in \mathbb{Z}$ and $-1 \leq m_d + r_1 - 1 \leq 1$, we must have have $m_d + r_1 = 1$. When $d = 3$, our result $1 - m_3 = m_3(1 - m_1)$ implies $m_3 = \frac{1}{2}$, which is not possible as $m_3 \in \{0, 1\}$ by construction. When $d \geq 4$, our result $1 - m_d = m_d(1 - m_{d-1})$ implies $m_d = \frac{1}{2}$, which is not possible as $m_d \in \{0, 1\}$ by construction.

Thus, for all $d \geq 3$, the above yields a contradiction.

Major Case 3: Focusing on moves along $j_3$, observe that during these moves $m_1 = m_3 = m_d$, $m_2 = 1 - m_1$, $m_{d-1} = m_d$ when $d \geq 4$, and $(m_1^*, m_2^*, m_d^*) = (1, m_2, m_d)$. Since $m_1 = m_3$ when $d = 3$ and $m_{d-1} = m_d$ when $d \geq 4$, $r_1 = 0$ and so in $j_1$ we have

$$(\ell + 1)A_3 - (2 + m_2 + 2m_1(1 - m_2) + 2x_1 + 2m_d(1 - m_1) - r_1(A_3 - 1))$$
$$= (\ell_1 + 1)A_3 - (m_1^* + m_1^*(1 - m_2^*) + 2(1 - m_1^*)m_d^*)$$
$$\implies \ell - \ell_1 = \frac{2(1 + x_1)}{A_3}.$$

Given $\ell - \ell_1 \in \mathbb{Z}$ and $2 \leq 2(1 + x_1) \leq A_3 - 2$ by construction, it must then be the case that for every $0 \leq x_1 \leq \frac{A_3}{2} - 2$, $\frac{2(1+x_1)}{A_3} \notin \mathbb{Z}$, giving us a contradiction as the equality above is not possible.

Major Case 4: Focusing on moves along $j_k$ for $4 \leq k \leq d$ when $d \geq 4$, we see that during all such moves $m_1 = m_2$, $m_{d-1} = m_d = 1 - m_1$, and $(m_1^*, m_2^*, m_d^*) = (m_1, m_2, m_d)$. Observing that $r_1 = 0$ since $m_{d-1} = m_d$ when $d \geq 4$, in $j_1$ we get

$$(\ell+1)A_3 - (2 + m_2 + 2m_1(1-m_2) + 2x_1 + 2m_d(1-m_1) - r_1(A_3-1))$$
$$= (\ell_1+1)A_3 - (m_1^* + m_1^*(1-m_2^*) + 2(1-m_1^*)m_d^*)$$
$$\implies \ell - \ell_1 = \frac{2(1+x_1)}{A_3}.$$

Since $\ell - \ell_1 \in \mathbb{Z}$ and $2 \leq 2(1+x_1) \leq A_3 - 2$ by construction, we see that for every $0 \leq x_1 \leq \frac{A_3}{2} - 2$, $\frac{2(1+x_1)}{A_3} \notin \mathbb{Z}$, giving us a contradiction as the above equality is not possible.

**Set Comparison Case 6:** We will consider whether edges are shared from moves with the same orientation via sets 3 and 4.

Major Case 1: Focusing on moves along $j_1$, we see that $C_2$ has no $j_1$ moves with the same orientation as those in $C_1$, and so there are no $j_1$ moves to consider.

Major Case 2: Focusing on moves along $j_2$, observe that only $C_2$'s $j_2$ move with $(m_1^*, m_2^*, m_d^*) = (0,0,0)$ has the same orientation as the $j_2$ moves in $C_1$, and during these moves $m_1 = m_2 = m_3 = m_d$. Applying this to $j_3$, we find

$$p_1 + A_4 s_1 + m_3 + 2x_3 = s_1^* + m_1^* m_3^*$$
$$\implies s_1 - s_1^* = -m_3.$$

Lastly, in $j_2$ we have

$$m_2 + \gamma + 4x_2 + 2(1-m_2)m_1 + (s_1 + 2(1-m_1)m_d)$$
$$= 2 + m_1^* + (s_1^* + (1-m_1^*)m_d^* + m_d^*) + \gamma_1 + m_2^*(1-m_1^*) + 4x_2^*$$
$$\implies x_2 - x_2^* = \frac{2 + (\gamma_1 - \gamma)}{4}.$$

Following from $x_2 - x_2^* \in \mathbb{Z}$ and $1 \leq 2+(\gamma_1 - \gamma) \leq 3$ by construction, we see that for all $\gamma, \gamma_1 \in \{0,1\}$, $\frac{2+(\gamma_1-\gamma)}{4} \notin \mathbb{Z}$, giving us a contradiction as the above equality is not possible.

Major Case 3: Focusing on moves along $j_3$, observe that during such moves $m_1 = m_2$, $m_3 = \cdots = m_d = 1 - m_1$, and $(m_1^*, m_2^*, m_d^*) = (1, m_2, m_d)$ with $m_3^* = m_d^*$. Applying these observations to $j_3$, we obtain

$$p_1 + A_4 s_1 + m_3 + 2x_3 = s_1^* + m_1^* m_3^*$$
$$\implies s_1 = s_1^*.$$

Consequently, in $j_2$ we get

$$m_2 + \gamma + 4x_2 + 2(1-m_2)m_1 + (s_1 + 2(1-m_1)m_d)$$
$$= 2 + m_1^* + (s_1^* + (1-m_1^*)m_d^* + m_d^*) + \gamma_1 + m_2^*(1-m_1^*) + 4x_2^*$$
$$\implies x_2 - x_2^* = \frac{2 + (\gamma_1 - \gamma)}{4}.$$

Since $x_2 - x_2^* \in \mathbb{Z}$ and $1 \leq 2 + (\gamma_1 - \gamma) \leq 3$ by construction, we see that for all $\gamma, \gamma_1 \in \{0, 1\}$, $\frac{2+(\gamma_1-\gamma)}{4} \notin \mathbb{Z}$, giving us a contradiction as the above equality is not possible.

**Major Case 4**: Focusing on moves along $j_k$ for $4 \leq k \leq d$ when $d \geq 4$, we see that during all such moves $m_1 = m_2 = m_3$, $m_d = 1 - m_1$, and $(m_1^*, m_2^*, m_3^*, m_d^*) = (m_1, m_2, (1 - m_d) + m_d m_4, m_d)$. Hence, in $j_3$ we obtain

$$p_1 + A_4 s_1 + m_3 + 2x_3 = s_1^* + m_1^* m_3^*$$
$$\implies s_1 = s_1^*.$$

Consequently, in $j_2$ we get

$$m_2 + \gamma + 4x_2 + 2(1-m_2)m_1 + (s_1 + 2(1-m_1)m_d)$$
$$= 2 + m_1^* + (s_1^* + (1-m_1^*)m_d^* + m_d^*) + \gamma_1 + m_2^*(1-m_1^*) + 4x_2^*$$
$$\implies x_2 - x_2^* = \frac{2 + (\gamma_1 - \gamma)}{4}.$$

Given $x_2 - x_2^* \in \mathbb{Z}$ and $1 \leq 2 + (\gamma_1 - \gamma) \leq 3$ by construction, we find that for all $\gamma, \gamma_1 \in \{0, 1\}$, $\frac{2+(\gamma_1-\gamma)}{4} \notin \mathbb{Z}$. Thus, we have a contradiction as the above equality is not possible.

Thus, no edges are shared from moves with the same orientation via two distinct sets for $\alpha_1 < \alpha^* < \alpha_2$ when $d \geq 3$.

$\boldsymbol{\alpha^*}$**−Case 3:** Let $\alpha^* = \alpha_2$. Then, $A_1 = 1$ and $A_2 = 1$, meaning only sets 1 and 3 are active. Note that $A_3 \geq 4$ is even.

**Set Comparison Case 1:** We will consider whether edges are shared from moves with the same orientation via sets 1 and 3.

**Major Case 1**: Focusing on moves along $j_1$, we see that during all such moves $m_1 = m_d$ and $m_1^* = m_d^*$, giving us in $j_1$

$$\ell A_3 + (m_1 + 2x_1 + 2m_d(1-m_1)) = (\ell_1 + 1)A_3 + (-2 + m_1^* + 2m_d^*(1-m_1^*))$$
$$\implies \ell - \ell_1 - 1 = -\frac{2 + 2x_1 + (m_1 - m_1^*)}{A_3}.$$

Following from $\ell - \ell_1 - 1 \in \mathbb{Z}$ and $1 \leq 2 + 2x_1 + (m_1 - m_1^*) \leq A_3 - 1$ by construction, we see that for every $0 \leq x_1 \leq \frac{A_3}{2} - 2$ and $m_1, m_1^* \in \{0, 1\}$, $-\frac{2+2x_1+(m_1-m_1^*)}{A_3} \notin \mathbb{Z}$, giving us a contradiction as the equality above is not possible.

**Major Case 2**: Focusing on moves along $j_2$, observe that only $C_1$'s $j_2$ move with $(m_1, m_d) = (0, 0)$ has the same orientation as the $j_2$ moves in $C_2$, and during these moves $m_1^* = m_d^*$. Consequently, these observations in $j_1$ imply

$$\ell A_3 + (m_1 + 2x_1 + 2m_d(1-m_1)) = (\ell_1 + 1)A_3 + (-2 + m_1^* + 2m_d^*(1-m_1^*))$$
$$\implies \ell - \ell_1 - 1 = -\frac{2 + 2x_1 - m_1^*}{A_3}.$$

Given $\ell - \ell_1 - 1 \in \mathbb{Z}$ and $1 \leq 2 + 2x_1 - m_1^* \leq A_3 - 2$ by construction, it follows that for every $0 \leq x_1 \leq \frac{A_3}{2} - 2$ and $m_1^* \in \{0, 1\}$, $-\frac{2+2x_1-m_1^*}{A_3} \notin \mathbb{Z}$, giving us a contradiction as the above equality is not possible.

Major Case 3: Focusing on moves along $j_k$ for $3 \leq k \leq d$, we see that during all such moves $m_d = 1 - m_1$ and $(m_1^*, m_d^*) = (m_1, m_d)$, giving us in $j_1$ that

$$\ell A_3 + (m_1 + 2x_1 + 2m_d(1 - m_1)) = (\ell_1 + 1)A_3 + (-2 + m_1^* + 2m_d^*(1 - m_1^*))$$

$$\implies \ell - \ell_1 - 1 = -\frac{2(1 + x_1)}{A_3}.$$

Since $\ell - \ell_1 - 1 \in \mathbb{Z}$ and $2 \leq 2(1 + x_1) \leq A_3 - 2$ by construction, it follows that for every $0 \leq x_1 \leq \frac{A_3}{2} - 2$, $-\frac{2(1+x_1)}{A_3} \notin \mathbb{Z}$, giving us a contradiction as the above equality is not possible.

Thus, no edges are shared from moves with the same orientation from two distinct sets when $\alpha^* = \alpha_2$ when $d \geq 3$.

$\boldsymbol{\alpha^*-}$**Case 4:** Let $\alpha_2 < \alpha^* \leq \alpha_d$. Then, $A_1 = 1$ and $A_2 = 1$, meaning only sets 1 and 3 are active as in the previous $\alpha^*-$case. Note that $\eta_{\alpha_2} = 1$.

**Set Comparison Case 1:** We will consider whether edges are shared from moves with the same orientation via sets 1 and 3.

Major Case 1: Focusing on moves along $j_1$, we see that during all such moves $m_1 = m_d$ and $(m_1^*, m_d^*) = (m_1, m_d)$, giving us in $j_1$

$$\ell A_3 + (m_1 + 2x_1 + 2m_d(1 - m_1)) = (\ell_1 + 1)A_3 + (-2 + m_1^* + 2m_d^*(1 - m_1^*))$$

$$\implies \ell - \ell_1 - 1 = -\frac{2(1 + x_1)}{A_3}.$$

Since $\ell - \ell_1 - 1 \in \mathbb{Z}$ and $2 \leq 2(1 + x_1) \leq A_3 - 2$ by construction, it follows that for every $0 \leq x_1 \leq \frac{A_3}{2} - 2$, $-\frac{2(1+x_1)}{A_3} \notin \mathbb{Z}$, giving us a contradiction as the above equality is not possible.

Major Case 2: Focusing on moves along $j_2$, observe that there are two general states to consider with $j_2$ moves in $C_2$ differing in orientation from those in $C_1$ depending on the state. Hence, we case on whether $C_2$ is being defined by either stair-casing/normal moves or column transitions:

Case 1: Let $R_1^* = 0$. Then, observing that during these stair-casing/normal moves along $j_2$ it is the case $m_1 = m_d$ and $(m_1^*, m_d^*) = (m_1, m_d)$, in $j_1$ we obtain

$$\ell A_3 + (m_1 + 2x_1 + 2m_d(1 - m_1)) = (\ell_1 + 1)A_3 + (-2 + m_1^* + 2m_d^*(1 - m_1^*))$$

$$\implies \ell - \ell_1 - 1 = -\frac{2(1 + x_1)}{A_3}.$$

Since $\ell - \ell_1 - 1 \in \mathbb{Z}$ and $2 \leq 2(1+x_1) \leq A_3 - 2$ by construction, it follows that for every $0 \leq x_1 \leq \frac{A_3}{2} - 2$ and $m_1^* \in \{0,1\}$, $-\frac{2(1+x_1)}{A_3} \notin \mathbb{Z}$, giving us a contradiction as the above equality is not possible.

<u>Case 2</u>: Let $R_1^* = 1$. Then, during the $j_2$ column transitions, only $C_1$'s $j_2$ move with $(m_1, m_d) = (0, 0)$ has the same orientation as the $j_2$ moves in $C_2$. Noting that during $C_1$'s $j_2$ moves $m_1^* = m_d^*$, our observations for this case in $j_1$ yield

$$\ell A_3 + (m_1 + 2x_1 + 2m_d(1 - m_1)) = (\ell_1 + 1)A_3 + (-2 + m_1^* + 2m_d^*(1 - m_1^*))$$
$$\implies \ell - \ell_1 - 1 = -\frac{2 + 2x_1 - m_1^*}{A_3}.$$

Following from $\ell - \ell_1 - 1 \in \mathbb{Z}$ and $1 \leq 2 + 2x_1 - m_1^* \leq A_3 - 2$ by construction, it is then the case that for all $0 \leq x_1 \leq \frac{A_3}{2} - 2$ and $m_1^* \in \{0, 1\}$, $-\frac{2+2x_1-m_1^*}{A_3} \notin \mathbb{Z}$, giving us a contradiction as the above equality is not possible.

<u>Major Case 3</u>: Focusing on moves along $j_k$ for $3 \leq k \leq d$, we proceed by casing on whether either stair-casing moves, normal moves, or column transitions are being carried out in defining $C_2$ as applicable:

<u>Case 1</u>: Let $R_1^* = 0$ and $R_{k-1}^* = 1$. Then, $C_2$ is stair-casing along $j_k$ for $3 \leq k \leq d$, and so we see that $m_d = 1 - m_1$ and $m_d^* = 1 - m_1^*$ during all such moves. Hence, in $j_1$ we see

$$\ell A_3 + (m_1 + 2x_1 + 2m_d(1 - m_1)) = (\ell_1 + 1)A_3 + (-2 + m_1^* + 2m_d^*(1 - m_1^*))$$
$$\implies \ell - \ell_1 - 1 = -\frac{2 + 2x_1 + (m_d - m_d^*)}{A_3}.$$

Given $\ell - \ell_1 - 1 \in \mathbb{Z}$ and $1 \leq 2 + 2x_1 + (m_d - m_d^*) \leq A_3 - 1$ by construction, we see that for all $0 \leq x_1 \leq \frac{A_3}{2} - 2$ and $m_d, m_d^* \in \{0, 1\}$, $-\frac{2+2x_1+(m_d-m_d^*)}{A_3} \notin \mathbb{Z}$, giving us a contradiction as the above equality is not possible.

<u>Case 2</u>: Let either $\alpha_2 < \alpha^* \leq \alpha_{k-1}$ or $\alpha_{k-1} < \alpha^* \leq \alpha_d$, $R_1^* = 0$ and $R_{k-1}^* = 0$. Then, $C_2$ is being defined by normal moves along $j_k$ for $4 \leq k \leq d$. Observing that during all such moves $m_d = 1 - m_1$ and $(m_1^*, m_d^*) = (m_1, m_d)$, in $j_1$ we see

$$\ell A_3 + (m_1 + 2x_1 + 2m_d(1 - m_1)) = (\ell_1 + 1)A_3 + (-2 + m_1^* + 2m_d^*(1 - m_1^*))$$
$$\implies \ell - \ell_1 - 1 = -\frac{2(1 + x_1)}{A_3}.$$

Since $\ell - \ell_1 - 1 \in \mathbb{Z}$ and $2 \leq 2(1 + x_1) \leq A_3 - 2$ by construction, it must be the case that for every $0 \leq x_1 \leq \frac{A_3}{2} - 2$, $-\frac{2(1+x_1)}{A_3} \notin \mathbb{Z}$, giving us a contradiction as the above equality is not possible.

<u>Case 3</u>: Let $R_1^* = 1$. Then, $C_2$ is being defined by column transitions along $j_k$ for $3 \leq k \leq d$, and so we must have $\alpha_{k-1} < \alpha^* \leq \alpha_d$. Observing that during all of these moves $m_d = 1 - m_1$ and $m_d^* = 1 - m_1^*$, in $j_1$ we have

$$\ell A_3 + (m_1 + 2x_1 + 2m_d(1 - m_1)) = (\ell_1 + 1)A_3 + (-2 + m_1^* + 2m_d^*(1 - m_1^*))$$
$$\implies \ell - \ell_1 - 1 = -\frac{2 + 2x_1 + (m_d - m_d^*)}{A_3}.$$

Since $\ell - \ell_1 - 1 \in \mathbb{Z}$ and $1 \leq 2 + 2x_1 + (m_d - m_d^*) \leq A_3 - 1$ by construction, we see that for all $0 \leq x_1 \leq \frac{A_3}{2} - 2$ and $m_d, m_d^* \in \{0, 1\}$, $-\frac{2+2x_1+(m_d-m_d^*)}{A_3} \notin \mathbb{Z}$, giving us a contradiction as the above equality is not possible.

Hence, no edges are shared from moves with the same orientation via two distinct sets for $\alpha_2 < \alpha^* \leq \alpha_d$ when $d \geq 3$.

Thus, for all $d \in \mathbb{Z}^{\geq 2}$ and $\alpha_1 \leq \alpha^* \leq \alpha_d$, no edges are shared from moves with the same orientation coming from two distinct sets via the General Lock-and-Key Decomposition's Edge Set Definition, proving Proposition 16 as required.

∎

### 7.1.3 Proof of Proposition 17:

Let $d \in \mathbb{Z}^{\geq 2}$, $\alpha_1 \leq \alpha^* \leq \alpha_d$, and for reference refer to Theorem 5 for specific bounds of each parameter, though we will bring some of them up in our arguments as necessary.

We would also like to make the reader aware that the form of the moves with same orientation in $C_2$ to a given move $(m_1, \ldots, m_d)$ in $C_1$ will vary based on which move and which major set is being considered. In some instances, there may be more than one move with the same orientation and so there would not be a unique move with the same orientation in said case. The definition of the move $(m_1^*, \ldots, m_d^*)$ with the same orientation as a given move $(m_1, \ldots, m_d)$ and all applicable equations that hold during every move will be introduced as required.

Note that we will be assuming that two cycles share an edge with the move configurations established below, meaning all components of their vertices agree, and we then show that such scenarios do not occur unless the conditions outlined below hold by the definition of the General Lock-and-Key Decomposition. In particular, we will be assuming that the same edge is defined by the same major set for both $C_1$ and $C_2$, and we show that this is possible if and only if $(m_1, \ldots, m_d) = (m_1^*, \ldots, m_d^*)$ and the move is to be performed from the same point $(x_1, x_2, \ldots, x_d)$. In doing so, we will have shown that these cycles share edges by way of these configurations if and only if $C_1$ and $C_2$ are the same cycle. Indeed, the following proof alone gives us the forward direction while Propositions 15 and 16 together with this proof gives us the backwards direction. We now proceed with the proof of Proposition 17 by considering equalities at each component $j_k$ for $1 \leq k \leq d$ with $C_2$ configured to perform a move $(m_1^*, \ldots, m_d^*)$ with the same orientation as that of a given move $(m_1, \ldots, m_d)$ of $C_1$ starting from the same vertex belonging to the edge in $C_1$.

**Dimension Case 1:** Let $d = 2$. Then, $C_1 = C_{\ell,\gamma}$ and $C_2 = C_{\ell_1,\gamma_1}$. We now case on $\alpha^*$ for $\alpha_1 \leq \alpha^* \leq \alpha_2$:

$\boldsymbol{\alpha^*}$**−Case 1:** Let $\alpha^* = \alpha_1$. In this case, only sets 3 and 4 are active as $A_1 = 0$ and $A_2 = 0$.

**Set Comparison Case 1:** We will consider when edges are shared from moves with the same orientation via set 3.

Major Case 1: Focusing on moves along $j_1$, we see that $C_1$ and $C_2$ have $j_1$ moves with the same orientation if and only if $\gamma = \gamma_1$. Observing that $m_2 = 1 = m_2^*$ during all $j_1$ moves, in $j_2$ we have

$$m_2 + \gamma + 4x_2 + (1 - m_2)m_1 = m_2^* + \gamma_1 + 4x_2^* + (1 - m_2^*)m_1^*$$
$$\implies x_2 = x_2^*.$$

Going to $j_1$, our assumptions give us

$$(\ell+1-\gamma)A_3-\gamma+(-1)^\gamma(-2+m_1+(1-m_2)m_1) = (\ell_1+1-\gamma_1)A_3-\gamma_1+(-1)^{\gamma_1}(-2+m_1^*+(1-m_2^*)m_1^*)$$

$$\implies \ell-\ell_1 = \frac{(-1)^\gamma(m_1^*-m_1)}{A_3}.$$

Since $\ell - \ell_1 \in \mathbb{Z}$ and $|m_1^* - m_1| \leq 1$ with $A_3 \geq 2$ by construction, it must be the case $m_1 = m_1^*$, meaning $(m_1, m_2) = (m_1^*, m_2^*)$, and so $\ell = \ell_1$. Hence, $C_1 = C_{\ell,\gamma} = C_{\ell_1,\gamma_1} = C_2$.

<u>Major Case 2</u>: Focusing on moves along $j_2$, observe that during all $j_2$ moves $m_2 = 0 = m_2^*$. Applying this to $j_2$, we get

$$m_2 + \gamma + 4x_2 + (1-m_2)m_1 = m_2^* + \gamma_1 + 4x_2^* + (1-m_2^*)m_1^*$$

$$\implies x_2 - x_2^* = \frac{(m_1^* - m_1) + (\gamma_1 - \gamma)}{4}.$$

Following from $x_2 - x_2^* \in \mathbb{Z}$ and $-2 \leq (m_1^* - m_1) + (\gamma_1 - \gamma) \leq 2$ by construction, we must have $(m_1^* - m_1) = -(\gamma_1 - \gamma)$ and so $x_2 = x_2^*$. Consequently, in $j_1$ we see

$$(\ell+1-\gamma)A_3-\gamma+(-1)^\gamma(-2+m_1+(1-m_2)m_1) = (\ell_1+1-\gamma_1)A_3-\gamma_1+(-1)^{\gamma_1}(-2+m_1^*+(1-m_2^*)m_1^*)$$

$$\implies ((\ell-\ell_1)+(\gamma_1-\gamma))\frac{A_3}{2} + (-1)^\gamma(m_1-1) + (-1)^{\gamma_1+1}(m_1^*-1) = \frac{\gamma-\gamma_1}{2}.$$

Since $((\ell-\ell_1)+(\gamma_1-\gamma))\frac{A_3}{2}+(-1)^\gamma(m_1-1)+(-1)^{\gamma_1+1}(m_1^*-1) \in \mathbb{Z}$ and $|\gamma_1-\gamma| \leq 1$ by construction, it follows that $\gamma = \gamma_1$, which by our previous result in $j_2$ implies $m_1 = m_1^*$ and so $(m_1, m_2) = (m_1^*, m_2^*)$. Thus, our last equality above implies $\ell = \ell_1$ and so $C_1 = C_{\ell,\gamma} = C_{\ell_1,\gamma_1} = C_2$.

**Set Comparison Case 2:** We will consider when edges are shared from moves with the same orientation via set 4.

<u>Major Case 1</u>: Focusing on moves along $j_1$, we see that $C_1$ and $C_2$ have $j_1$ moves with the same orientation if and only if $\gamma = \gamma_1$. Noting that $m_2 = 1 = m_2^*$ during all $j_1$ moves, in $j_2$ we find

$$2 + m_1 + \gamma + m_2(1-m_1) + 4x_2 = 2 + m_1^* + \gamma_1 + m_2^*(1-m_1^*) + 4x_2^*$$

$$\implies x_2 = x_2^*.$$

Now, in $j_1$ our observations yield

$$(\ell+1-\gamma)A_3-\gamma+(-1)^{\gamma+1}(m_1+m_1(1-m_2)) = (\ell_1+1-\gamma_1)A_3-\gamma_1+(-1)^{\gamma_1+1}(m_1^*+m_1^*(1-m_2^*))$$

$$\implies \ell-\ell_1 = \frac{(-1)^{\gamma+1}(m_1^*-m_1)}{A_3}.$$

Following from $\ell - \ell_1 \in \mathbb{Z}$ and $|m_1^* - m_1| \leq 1$ with $A_3 \geq 2$ by construction, we get that $m_1 = m_1^*$, meaning $(m_1, m_2) = (m_1^*, m_2^*)$ and so $\ell = \ell_1$. Consequently, $C_1 = C_{\ell,\gamma} = C_{\ell_1,\gamma_1} = C_2$.

Major Case 2: Focusing on moves along $j_2$, observe that during all $j_2$ moves $m_2 = 0 = m_2^*$, giving us in $j_2$

$$2 + m_1 + \gamma + m_2(1-m_1) + 4x_2 = 2 + m_1^* + \gamma_1 + m_2^*(1-m_1^*) + 4x_2^*$$
$$\implies x_2 - x_2^* = \frac{(m_1^* - m_1) + (\gamma_1 - \gamma)}{4}.$$

Given $x_2 - x_2^* \in \mathbb{Z}$ and $-2 \leq (m_1^* - m_1) + (\gamma_1 - \gamma) \leq 2$ by construction, we must have $(m_1^* - m_1) = -(\gamma_1 - \gamma)$ and so $x_2 = x_2^*$. Thus, in $j_1$ our results and observations imply

$$(\ell + 1 - \gamma)A_3 - \gamma + (-1)^{\gamma+1}(m_1 + m_1(1-m_2)) = (\ell_1 + 1 - \gamma_1)A_3 - \gamma_1 + (-1)^{\gamma_1+1}(m_1^* + m_1^*(1-m_2^*))$$
$$\implies ((\ell - \ell_1) + (\gamma_1 - \gamma))\frac{A_3}{2} + (-1)^{\gamma+1}m_1 - (-1)^{\gamma_1+1}m_1^* = \frac{\gamma - \gamma_1}{2}.$$

Since $((\ell - \ell_1) + (\gamma_1 - \gamma))\frac{A_3}{2} + (-1)^{\gamma+1}m_1 - (-1)^{\gamma_1+1}m_1^* \in \mathbb{Z}$ and $|\gamma_1 - \gamma| \leq 1$ by construction, it is the case that $\gamma = \gamma_1$, which by our previous main result implies $m_1 = m_1^*$. So $(m_1, m_2) = (m_1^*, m_2^*)$ and $\ell = \ell_1$. Thus, $C_1 = C_{\ell,\gamma} = C_{\ell_1,\gamma_1} = C_2$.

Hence, edges are shared from moves with the same orientation for $\alpha^* = \alpha_1$ when $d = 2$ if and only if $C_1 = C_2$.

**$\alpha^*$−Case 2:** Let $\alpha_1 < \alpha^* < \alpha_2$. Then, sets 1, 2, 3, and 4 are active as $A_1 = 1$ and $A_2 = 0$. Note that $A_3 \geq 4$ is even.

**Set Comparison Case 1:** We will consider when edges are shared from moves with the same orientation via set 1.

Major Case 1: Focusing on moves along $j_1$, we see that $C_1$ and $C_2$ have $j_1$ moves with the same orientation if and only if $\gamma = \gamma_1$. Observing that $m_1 + m_2 = 1$ and $m_1^* + m_2^* = 1$ during these moves, in $j_2$ we get

$$m_2 + \gamma + 4x_2 = m_2^* + \gamma_1 + 4x_2^*$$
$$\implies x_2 - x_2^* = \frac{m_2^* - m_2}{4}.$$

Since $x_2 - x_2^* \in \mathbb{Z}$ and $|m_2^* - m_2| \leq 1$ by construction, it follows that $m_2 = m_2^*$ and so $x_2 = x_2^*$. Noting that by our previous main result $m_1 + m_2 = 1 = m_1^* + m_2^*$ implies $m_1 = m_1^*$, meaning $(m_1, m_2) = (m_1^*, m_2^*)$, in $j_1$ we have

$$(\ell + \gamma)A_3 - \gamma + (-1)^\gamma(m_1 + 2x_1) = (\ell_1 + \gamma_1)A_3 - \gamma_1 + (-1)^{\gamma_1}(m_1^* + 2x_1^*)$$
$$\implies \ell - \ell_1 = \frac{2(-1)^\gamma(x_1^* - x_1)}{A_3}.$$

Since $\ell - \ell_1 \in \mathbb{Z}$ and $2|x_1^* - x_1| \leq A_3 - 4$ by construction, it follows that $x_1 = x_1^*$ and so $\ell = \ell_1$. Thus, $C_1 = C_{\ell,\gamma} = C_{\ell_1,\gamma_1} = C_2$.

Major Case 2: Focusing on moves along $j_2$, observe that during all $j_2$ moves $m_1 = m_2$ and $(m_1^*, m_2^*) = (m_1, m_2)$, giving us in $j_2$

$$m_2 + \gamma + 4x_2 = m_2^* + \gamma_1 + 4x_2^*$$
$$\implies x_2 - x_2^* = \frac{\gamma_1 - \gamma}{4}.$$

Since $x_2 - x_2^* \in \mathbb{Z}$ and $|\gamma_1 - \gamma| \leq 1$ by construction, we must have $\gamma = \gamma_1$, giving us $x_2 = x_2^*$. Applying this to $j_1$, we see

$$(\ell + \gamma)A_3 - \gamma + (-1)^\gamma(m_1 + 2x_1) = (\ell_1 + \gamma_1)A_3 - \gamma_1 + (-1)^{\gamma_1}(m_1^* + 2x_1^*)$$
$$\implies \ell - \ell_1 = \frac{2(-1)^\gamma(x_1^* - x_1)}{A_3}$$

Given $\ell - \ell_1 \in \mathbb{Z}$ and $2|x_1^* - x_1| \leq A_3 - 4$ by construction, it must be the case $x_1 = x_1^*$ and so $\ell = \ell_1$. Hence, $C_1 = C_{\ell,\gamma} = C_{\ell_1,\gamma_1} = C_2$.

**Set Comparison Case 2:** We will consider when edges are shared from moves with the same orientation via set 2.

Major Case 1: Focusing on moves along $j_1$, observe that $C_1$ and $C_2$ have $j_1$ moves with the same orientation if and only if $\gamma = \gamma_1$. Since $m_1 = m_2$ and $m_1^* = m_2^*$ during all $j_1$ moves, in $j_2$ we get

$$2 + \gamma + m_1 + 4x_2 = 2 + \gamma_1 + m_1^* + 4x_2^*$$
$$\implies x_2 - x_2^* = \frac{m_1^* - m_1}{4}.$$

Following from $x_2 - x_2^* \in \mathbb{Z}$ and $|m_1^* - m_1| \leq 1$ by construction, it follows that $m_1 = m_1^*$, meaning $(m_1, m_2) = (m_1^*, m_2^*)$, and $x_2 = x_2^*$. Applying this to $j_1$, we see

$$(\ell + 1 - \gamma)A_3 - \gamma + (-1)^{\gamma+1}(2 + m_2 + 2m_1(1 - m_2) + 2x_1)$$
$$= (\ell_1 + 1 - \gamma_1)A_3 - \gamma_1 + (-1)^{\gamma_1+1}(2 + m_2^* + 2m_1^*(1 - m_2^*) + 2x_1^*)$$
$$\implies \ell - \ell_1 = \frac{2(-1)^{\gamma+1}(x_1^* - x_1)}{A_3}.$$

Since $\ell - \ell_1 \in \mathbb{Z}$ and $2|x_1^* - x_1| \leq A_3 - 4$ by construction, we must have $x_1 = x_1^*$ and so $\ell = \ell_1$. Hence, $C_1 = C_{\ell,\gamma} = C_{\ell_1,\gamma_1} = C_2$.

Major Case 2: Focusing on moves along $j_2$, observe that which moves in $C_1$ have the same orientation as those in $C_2$ depends on $x_1$ and $x_1^*$. Hence, we case on $x_1$ and further case on $x_1^*$:

Case 1: Let $0 \leq x_1 < \frac{A_3}{2} - 2$. We now case on $x_1^*$:

Subcase 1: Let $0 \leq x_1^* < \frac{A_3}{2} - 2$. Observing that $m_1 + m_2 = 1$, $m_1^* + m_2^* = 1$ and $(m_1^*, m_2^*) = (m_1, m_2)$, in $j_2$ we have

$$2 + m_1 + \gamma + 4x_2 = 2 + m_1^* + \gamma_1 + 4x_2^*$$
$$\implies x_2 - x_2^* = \frac{\gamma_1 - \gamma}{4}.$$

Given $x_2 - x_2^* \in \mathbb{Z}$ and $|\gamma_1 - \gamma| \leq 1$ by construction, it follows that $\gamma = \gamma_1$ and so $x_2 = x_2^*$. Applying this to $j_1$, we get

$$(\ell + 1 - \gamma)A_3 - \gamma + (-1)^{\gamma+1}(2 + m_2 + 2m_1(1 - m_2) + 2x_1)$$
$$= (\ell_1 + 1 - \gamma_1)A_3 - \gamma_1 + (-1)^{\gamma_1+1}(2 + m_2^* + 2m_1^*(1 - m_2^*) + 2x_1^*)$$
$$\implies \ell - \ell_1 = \frac{2(-1)^{\gamma+1}(x_1^* - x_1)}{A_3}.$$

Following from $\ell - \ell_1 \in \mathbb{Z}$ and $2|x_1^* - x_1| < A_3 - 4$ by construction, it must be the case that $x_1 = x_1^*$ and so $\ell = \ell_1$. Consequently, $C_1 = C_{\ell,\gamma} = C_{\ell_1,\gamma_1} = C_2$.

<u>Subcase 2:</u> Let $x_1^* = \frac{A_3}{2} - 2$. Then, we see that only $C_1$'s $j_2$ move $(m_1, m_2) = (0, 1)$ has the same orientation as all of $C_2$'s $j_2$ moves, and during these moves $m_1^* + m_2^* = 1$. So in $j_2$ we find

$$2 + m_1 + \gamma + 4x_2 = 2 + m_1^* + \gamma_1 + 4x_2^*$$
$$\implies x_2 - x_2^* = \frac{m_1^* + (\gamma_1 - \gamma)}{4}.$$

Given $x_2 - x_2^* \in \mathbb{Z}$ and $-1 \leq m_1^* + (\gamma_1 - \gamma) \leq 2$ by construction, it must be the case $m_1^* = |\gamma_1 - \gamma|$ with $\gamma_1 - \gamma \in \{-1, 0\}$. Note that if $\gamma_1 - \gamma = 1$, then for every $m_1^* \in \{0, 1\}$, $\frac{m_1^*+1}{4} \notin \mathbb{Z}$, giving us a contradiction as the above equality would not be possible in said case. Treating the remaining two cases for $\gamma_1 - \gamma \in \{-1, 0\}$ with $\gamma_1 - \gamma = -|\gamma_1 - \gamma|$, in $j_1$ we obtain

$$(\ell + 1 - \gamma)A_3 - \gamma + (-1)^{\gamma+1}(2 + m_2 + 2m_1(1 - m_2) + 2x_1)$$
$$= (\ell_1 + 1 - \gamma_1)A_3 - \gamma_1 + (-1)^{\gamma_1+1}(2 + m_2^* + 2m_1^*(1 - m_2^*) + 2x_1^*)$$
$$\implies (\ell - \ell_1) + (\gamma_1 - \gamma) + (-1)^{\gamma_1} = \frac{(-1)^{\gamma}(2x_1 + 4 - 2|\gamma_1 - \gamma|)}{A_3}.$$

In either case that $\gamma_1 - \gamma \in \{-1, 0\}$, we see $(\ell - \ell_1) + (\gamma_1 - \gamma) + (-1)^{\gamma_1} \in \mathbb{Z}$ and for every $0 \leq x_1 < \frac{A_3}{2} - 2$, $\frac{(-1)^{\gamma}(2x_1+4-2|\gamma_1-\gamma|)}{A_3} \notin \mathbb{Z}$, giving us a contradiction as the above equality is not possible.

<u>Case 2:</u> Let $x_1 = \frac{A_3}{2} - 2$. Since the subcase with $x_1$ as assumed and $0 \leq x_1^* < \frac{A_3}{2} - 2$ is symmetric to that of subcase 2 of case 1, the same argument addresses the subcase outlined above. So let $x_1^* = \frac{A_3}{2} - 2$. Then, observing that during such moves $m_1 + m_2 = 1$ and $m_1^* + m_2^* = 1$, in $j_2$ we see

$$2 + m_1 + \gamma + 4x_2 = 2 + m_1^* + \gamma_1 + 4x_2^*$$
$$\implies x_2 - x_2^* = \frac{(m_1^* - m_1) + (\gamma_1 - \gamma)}{4}.$$

Following from $x_2 - x_2^* \in \mathbb{Z}$ and $-2 \leq (m_1^* - m_1) + (\gamma_1 - \gamma) \leq 2$ by construction, it follows that $(m_1^* - m_1) = -(\gamma_1 - \gamma)$ and so $x_2 = x_2^*$. Consequently, in $j_1$ we find

$$(\ell + 1 - \gamma)A_3 - \gamma + (-1)^{\gamma+1}(2 + m_2 + 2m_1(1 - m_2) + 2x_1)$$
$$= (\ell_1 + 1 - \gamma_1)A_3 - \gamma_1 + (-1)^{\gamma_1+1}(2 + m_2^* + 2m_1^*(1 - m_2^*) + 2x_1^*)$$
$$\implies (\ell - \ell_1) + (\gamma_1 - \gamma) + (-1)^{\gamma+1} + (-1)^{\gamma_1} = \frac{-2((m_1 - 1)\gamma - \gamma_1(m_1^* - 1))}{A_3}.$$

Given $(\ell - \ell_1) + (\gamma_1 - \gamma) + (-1)^{\gamma+1} + (-1)^{\gamma_1} \in \mathbb{Z}$ and $-2 \leq -2((m_1 - 1)\gamma - \gamma_1(m_1^* - 1)) \leq 2$ with $A_3 \geq 4$ by construction, we must have $\gamma(m_1 - 1) = \gamma_1(m_1^* - 1)$. From this, observe that

$$m_1 - m_1^* = \gamma_1 - \gamma = \gamma_1 m_1^* - \gamma m_1.$$

<u>Subcase 1:</u> Let $\gamma = 0$. Then, our equality above implies $\gamma_1(1 - m_1^*) = 0$, which is if and only if $\gamma_1 = 0$ or $m_1^* = 1$. If $\gamma_1 = 0$, then $\gamma = \gamma_1$ and our previous equality implies $m_1 = m_1^*$. By our initial observations, this implies $m_2 = m_2^*$, meaning $(m_1, m_2) = (m_1^*, m_2^*)$. If $m_1^* = 1$, then our previous equality implies $\gamma_1 = m_1 - 1$. Since $m_1, \gamma_1 \in \{0, 1\}$, it must be the case $m_1 = 1$ and $\gamma_1 = 0$, so $\gamma = \gamma_1$ and $m_1 = m_1^*$, meaning $(m_1, m_2) = (m_1^*, m_2^*)$.

<u>Subcase 2:</u> Let $\gamma = 1$. Then, our equality above implies $m_1 - 1 = \gamma_1(m_1^* - 1)$. Since $\gamma_1 - \gamma = m_1 - m_1^*$, we see

$$m_1 - 1 = (m_1 - m_1^* + 1)(m_1^* - 1)$$
$$\implies m_1 = m_1^*.$$

Given $m_1 = m_1^*$, we have $\gamma_1 - 1 = m_1 - m_1^*$ implies $\gamma_1 = 1$. Hence, $\gamma = \gamma_1$ and $m_1 = m_1^*$, giving us $(m_1, m_2) = (m_1^*, m_2^*)$.

Since $(m_1, m_2) = (m_1^*, m_2^*)$ and $\gamma = \gamma_1$ in all cases, we see that our last equality in $j_1$ implies $\ell = \ell_1$. Thus, $C_1 = C_{\ell, \gamma} = C_{\ell_1, \gamma_1} = C_2$.

**Set Comparison Case 3:** We will consider when edges are shared from moves with the same orientation via set 3.

<u>Major Case 1:</u> Focusing on moves along $j_1$, note that $C_1$ and $C_2$ have $j_1$ moves with the same orientation if and only if $\gamma = \gamma_1$. Since $m_2 = 1 = m_2^*$ during all $j_1$ moves, in $j_2$ we get

$$m_2 + \gamma + 4x_2 + (1 - m_2)m_1 = m_2^* + \gamma_1 + 4x_2^* + (1 - m_2^*)m_1^*$$
$$\implies x_2 = x_2^*.$$

Going to $j_1$, it follows that

$$(\ell+1-\gamma)A_3 - \gamma + (-1)^\gamma(-2 + m_1 + m_1(1-m_2)) = (\ell_1+1-\gamma_1)A_3 - \gamma_1 + (-1)^{\gamma_1}(-2 + m_1^* + m_1^*(1-m_2^*))$$
$$\implies \ell - \ell_1 = \frac{(-1)^\gamma(m_1^* - m_1)}{A_3}.$$

Since $\ell - \ell_1 \in \mathbb{Z}$ and $|m_1^* - m_1| \leq 1$ by construction, we see $m_1 = m_1^*$, meaning $(m_1, m_2) = (m_1^*, m_2^*)$, and so $\ell = \ell_1$. Thus, $C_1 = C_{\ell, \gamma} = C_{\ell_1, \gamma_1} = C_2$.

<u>Major Case 2:</u> Focusing on moves along $j_2$, it is the case that $m_2 = 0 = m_2^*$ during these moves, giving us in $j_2$

$$m_2 + \gamma + 4x_2 + (1 - m_2)m_1 = m_2^* + \gamma_1 + 4x_2^* + (1 - m_2^*)m_1^*$$
$$\implies x_2 - x_2^* = \frac{(m_1^* - m_1) + (\gamma_1 - \gamma)}{4}.$$

Since $x_2 - x_2^* \in \mathbb{Z}$ and $-2 \leq (m_1^* - m_1) + (\gamma_1 - \gamma) \leq 2$ by construction, it must be the case $(m_1^* - m_1) = -(\gamma_1 - \gamma)$ and so $x_2 = x_2^*$. In $j_1$, we now have

$$(\ell + 1 - \gamma)A_3 - \gamma + (-1)^\gamma(-2 + m_1 + m_1(1 - m_2))$$
$$= (\ell_1 + 1 - \gamma_1)A_3 - \gamma_1 + (-1)^{\gamma_1}(-2 + m_1^* + m_1^*(1 - m_2^*))$$
$$\implies ((\ell - \ell_1) + (\gamma_1 - \gamma))\frac{A_3}{2} + (-1)^\gamma(m_1 - 1) + (-1)^{\gamma_1+1}(m_1^* - 1) = \frac{\gamma - \gamma_1}{2}.$$

Given $((\ell-\ell_1)+(\gamma_1-\gamma))\frac{A_3}{2}+(-1)^\gamma(m_1-1)+(-1)^{\gamma_1+1}(m_1^*-1) \in \mathbb{Z}$ and $|\gamma_1-\gamma| \leq 1$ by construction, we must have $\gamma = \gamma_1$, which by our previous result implies $m_1 = m_1^*$, meaning $(m_1, m_2) = (m_1^*, m_2^*)$. So our last equality in $j_1$ now implies $\ell = \ell_1$, giving us $C_1 = C_{\ell,\gamma} = C_{\ell_1,\gamma_1} = C_2$.

**Set Comparison Case 4:** We will consider when edges are shared from moves with the same orientation via set 4.

<u>Major Case 1</u>: Focusing on moves along $j_1$, we see that $C_1$ and $C_2$ have $j_1$ moves with the same orientation if and only if $\gamma = \gamma_1$. Noting that $m_2 = 1 = m_2^*$ during these moves, in $j_2$ we obtain

$$2 + m_1 + \gamma + m_2(1 - m_1) + 4x_2 = 2 + m_1^* + \gamma_1 + m_2^*(1 - m_1^*) + 4x_2^*$$
$$\implies x_2 = x_2^*.$$

Going to $j_1$, our assumptions give us

$$(\ell + 1 - \gamma)A_3 - \gamma + (-1)^{\gamma+1}(m_1 + m_1(1 - m_2))$$
$$= (\ell_1 + 1 - \gamma_1)A_3 - \gamma_1 + (-1)^{\gamma_1+1}(m_1^* + m_1^*(1 - m_2^*))$$
$$\implies \ell - \ell_1 = \frac{(-1)^{\gamma+1}(m_1^* - m_1)}{A_3}.$$

Since $\ell - \ell_1 \in \mathbb{Z}$ and $|m_1^* - m_1| \leq 1$ with $A_3 \geq 4$ by construction and our $\alpha^*$−case assumption, it follows that $m_1 = m_1^*$ and so $\ell = \ell_1$. Thus, $(m_1, m_2) = (m_1^*, m_2^*)$ and $C_1 = C_{\ell,\gamma} = C_{\ell_1,\gamma_1} = C_2$.

<u>Major Case 2</u>: Focusing on moves along $j_2$, it is the case that $m_2 = 0$ and $(m_1^*, m_2^*) = (m_1, m_2)$. Hence, in $j_2$ this yields

$$2 + m_1 + \gamma + m_2(1 - m_1) + 4x_2 = 2 + m_1^* + \gamma_1 + m_2^*(1 - m_1^*) + 4x_2^*$$
$$\implies x_2 - x_2^* = \frac{\gamma_1 - \gamma}{4}.$$

Following from $x_2 - x_2^* \in \mathbb{Z}$ and $|\gamma_1 - \gamma| \leq 1$ by construction, it must be the case $\gamma = \gamma_1$ and so $x_2 = x_2^*$. Applying this to $j_1$, we get

$$(\ell + 1 - \gamma)A_3 - \gamma + (-1)^{\gamma+1}(m_1 + m_1(1 - m_2))$$
$$= (\ell_1 + 1 - \gamma_1)A_3 - \gamma_1 + (-1)^{\gamma_1+1}(m_1^* + m_1^*(1 - m_2^*))$$
$$\implies \ell = \ell_1.$$

Thus, $C_1 = C_{\ell,\gamma} = C_{\ell_1,\gamma_1} = C_2$.

Consequently, edges are shared from moves with the same orientation from the same set for $\alpha_1 < \alpha^* < \alpha_2$ when $d = 2$ if and only if $C_1 = C_2$.

**α\*−Case 3:** Let $\alpha^* = \alpha_2$. Then, $A_1 = 1$ and $A_2 = 1$, meaning only sets 1 and 3 are active. Note that $A_3 \geq 4$ is even.

**Set Comparison Case 1:** We will consider when edges are shared from moves with the same orientation via set 1.

Major Case 1: Focusing on moves along $j_1$, we see that during all $j_1$ moves $m_1 + m_2 = 1$ and $m_1^* + m_2^* = 1$, giving us in $j_1$

$$\ell A_3 + (m_1 + 2x_1) = \ell_1 A_3 + (m_1^* + 2x_1^*)$$
$$\implies (\ell - \ell_1)\frac{A_3}{2} + (x_1 - x_1^*) = \frac{m_1^* - m_1}{2}.$$

Given $(\ell - \ell_1)\frac{A_3}{2} + (x_1 - x_1^*) \in \mathbb{Z}$ and $|m_1^* - m_1| \leq 1$ by construction, it must be that $m_1 = m_1^*$, meaning $(m_1, m_2) = (m_1^*, m_2^*)$. This leaves us with

$$\ell - \ell_1 = \frac{2(x_1^* - x_1)}{A_3}.$$

Following from $\ell - \ell_1 \in \mathbb{Z}$ and $2|x_1^* - x_1| \leq A_3 - 4$, we get $x_1 = x_1^*$ and so $\ell = \ell_1$. Lastly, our deduction $m_2 = m_2^*$ in $j_2$ yields

$$m_2 + \gamma + 2x_2 = m_2^* + \gamma_1 + 2x_2^*$$
$$\implies x_2 - x_2^* = \frac{\gamma_1 - \gamma}{2}.$$

Since $x_2 - x_2^* \in \mathbb{Z}$ and $|\gamma_1 - \gamma| \leq 1$ by construction, it follows that $\gamma = \gamma_1$ and so $x_2 = x_2^*$. Thus, $C_1 = C_{\ell,\gamma} = C_{\ell_1,\gamma_1} = C_2$.

Major Case 2: Focusing on moves along $j_2$, observe that during all $j_2$ moves $m_1 = m_2$ and $(m_1^*, m_2^*) = (m_1, m_2)$. Applying this to $j_1$, it immediately follows that

$$\ell A_3 + (m_1 + 2x_1) = \ell_1 A_3 + (m_1^* + 2x_1^*)$$
$$\implies \ell - \ell_1 = \frac{2(x_1^* - x_1)}{A_3}.$$

Since $\ell - \ell_1 \in \mathbb{Z}$ and $2|x_1^* - x_1| \leq A_3 - 4$ by construction, we get $x_1 = x_1^*$ and so $\ell = \ell_1$. Our assumptions in $j_2$ imply

$$m_2 + \gamma + 2x_2 = m_2^* + \gamma_1 + 2x_2^*$$
$$\implies x_2 - x_2^* = \frac{\gamma_1 - \gamma}{2}.$$

Given $x_2 - x_2^* \in \mathbb{Z}$ and $|\gamma_1 - \gamma| \leq 1$ by construction, it must be the case $\gamma = \gamma_1$ and so $x_2 = x_2^*$. Thus, $C_1 = C_{\ell,\gamma} = C_{\ell_1,\gamma_1} = C_2$.

**Set Comparison Case 2:** We will consider when edges are shared from moves with the same orientation via set 3.

**Major Case 1:** Focusing on moves along $j_1$, we see that during all $j_1$ moves $m_1 + m_2 = 1$ and $m_1^* + m_2^* = 1$. Consequently, in $j_1$ we have

$$(\ell + 1)A_3 + (-2 + m_1) = (\ell_1 + 1)A_3 + (-2 + m_1^*)$$
$$\implies \ell - \ell_1 = \frac{m_1^* - m_1}{A_3}.$$

Since $\ell - \ell_1 \in \mathbb{Z}$ and $|m_1^* - m_1| \leq 1$ with $A_3 \geq 4$, it must be the case $m_1 = m_1^*$, meaning $(m_1, m_2) = (m_1^*, m_2^*)$, and so $\ell = \ell_1$. Applying the deduction $m_2 = m_2^*$ along with our result $m_1 = m_1^*$, in $j_2$ we see

$$m_2 + \gamma + 2x_2 + 2(1 - m_2)m_1 = m_2^* + \gamma_1 + 2x_2^* + 2(1 - m_2^*)m_1^*$$
$$\implies x_2 - x_2^* = \frac{\gamma_1 - \gamma}{2}.$$

Following from $x_2 - x_2^* \in \mathbb{Z}$ and $|\gamma_1 - \gamma| \leq 1$ by construction, it must be the case $\gamma = \gamma_1$ and so $x_2 = x_2^*$. Thus, $C_1 = C_{\ell,\gamma} = C_{\ell_1,\gamma_1} = C_2$.

**Major Case 2:** Focusing on moves along $j_1$, we see that during all $j_1$ moves $m_1 = m_2$ and $m_1^* = m_2^*$. Consequently, in $j_1$ we have

$$(\ell + 1)A_3 + (-2 + m_1) = (\ell_1 + 1)A_3 + (-2 + m_1^*)$$
$$\implies \ell - \ell_1 = \frac{m_1^* - m_1}{A_3}.$$

Since $\ell - \ell_1 \in \mathbb{Z}$ and $|m_1^* - m_1| \leq 1$ with $A_3 \geq 4$, it must be the case $m_1 = m_1^*$, meaning $(m_1, m_2) = (m_1^*, m_2^*)$, and so $\ell = \ell_1$. Applying the deduction $m_2 = m_2^*$ along with our result $m_1 = m_1^*$, in $j_2$ we find

$$m_2 + \gamma + 2x_2 + 2(1 - m_2)m_1 = m_2^* + \gamma_1 + 2x_2^* + 2(1 - m_2^*)m_1^*$$
$$\implies x_2 - x_2^* = \frac{\gamma_1 - \gamma}{2}.$$

Following from $x_2 - x_2^* \in \mathbb{Z}$ and $|\gamma_1 - \gamma| \leq 1$ by construction, it must be the case $\gamma = \gamma_1$ and so $x_2 = x_2^*$. Thus, $C_1 = C_{\ell,\gamma} = C_{\ell_1,\gamma_1} = C_2$.

Hence, edges are shared from moves with the same orientation from the same set for $\alpha^* = \alpha_2$ when $d = 2$ if and only if $C_1 = C_2$.

**Dimension Case 2:** Let $d \geq 3$. Then, $C_1 = C_{\ell,\gamma,t,p_1,s_1,\ldots,p_{d-2},s_{d-2}}$ and $C_2 = C_{\ell_1,\gamma_1,t_1,p_1^*,s_1^*,\ldots,p_{d-2}^*,s_{d-2}^*}$.

**$\alpha^*$−Case 1:** Let $\alpha^* = \alpha_1$. Then, $A_1 = 0$ and $A_2 = 0$, meaning sets 3 and 4 are active. Note that $p_{k-2} = 0 = p_{k-2}^*$, $A_{k+1} = 1$, and $x_k = 0 = x_k^*$ for all $3 \leq k \leq d$.

**Set Comparison Case 1:** We will consider when edges are shared from moves with the same orientation via set 3.

**Major Case 1:** Focusing on moves along $j_1$, we see that $C_1$ and $C_2$ have $j_1$ moves with the same orientation if and only if $\gamma = \gamma_1$. Observing that during these moves $m_2 = 1 = m_2^*$, $m_1 = m_3 = \cdots = m_d$ and $m_1^* = m_3^* = \cdots = m_d^*$, in $j_2$ we see

$$m_2 + \gamma + 4x_2 + (1-m_2)m_1 + 2(1-m_1)m_d = m_2^* + \gamma_1 + 4x_2^* + (1-m_2^*)m_1^* + 2(1-m_1^*)m_d^*$$
$$\implies x_2 = x_2^*.$$

Now, our observations in $j_3$ imply

$$p_1 + 2A_4 s_1 + (-1)^\gamma m_3 + 2x_3 = p_1^* + 2A_4 s_1^* + (-1)^{\gamma_1} m_3^* + 2x_3^*$$
$$\implies s_1 - s_1^* = \frac{(-1)^\gamma (m_3^* - m_3)}{2}.$$

Given $s_1 - s_1^* \in \mathbb{Z}$ and $|m_3^* - m_3| \leq 1$ by construction, it must be the case $m_3 = m_3^*$, so $s_1 = s_1^*$ and $(m_1, \ldots, m_d) = (m_1^*, \ldots, m_d^*)$.

<u>Major Case 2</u>: Focusing on moves along $j_2$, observe that during these moves $m_2 = 0 = m_2^*$, $m_1 = m_3 = \cdots = m_d$, and $m_1^* = m_3^* = \cdots = m_d^*$. Applying this to $j_2$, we obtain

$$m_2 + \gamma + 4x_2 + (1-m_2)m_1 + 2m_d(1-m_1) = m_2^* + \gamma_1 + 4x_2^* + (1-m_2^*)m_1^* + 2m_d^*(1-m_1^*)$$
$$\implies x_2 - x_2^* = \frac{(\gamma_1 - \gamma) + (m_1^* - m_1)}{4}.$$

Since $x_2 - x_2^* \in \mathbb{Z}$ and $-2 \leq (\gamma_1 - \gamma) + (m_1^* - m_1) \leq 2$ by construction, it must be the case $(m_1^* - m_1) = -(\gamma_1 - \gamma)$ and $x_2 = x_2^*$. Applying this to $j_3$, we find

$$p_1 + 2A_4 s_1 + (-1)^\gamma m_3 + 2x_3 = p_1^* + 2A_4 s_1^* + (-1)^{\gamma_1} m_3^* + 2x_3^*$$
$$\implies (s_1 - s_1^*) - \frac{(-1)^{\gamma_1} - (-1)^\gamma}{2} m_3 = \frac{(-1)^{\gamma_1}(m_3^* - m_3)}{2}.$$

Following from $(s_1 - s_1^*) - \frac{(-1)^{\gamma_1} - (-1)^\gamma}{2} m_3 \in \mathbb{Z}$ and $|m_3^* - m_3| \leq 1$ by construction, it must be the case $m_3 = m_3^*$, meaning $(m_1, \ldots, m_d) = (m_1^*, \ldots, m_d^*)$, and so $\gamma = \gamma_1$. Hence, $s_1 = s_1^*$.

<u>Major Case 3</u>: Focusing on moves along $j_3$, we see that during such moves $m_1 = m_2$, $m_3 = \cdots = m_d = 1 - m_1$, and $(m_1^*, \ldots, m_d^*) = ((1-|\gamma_1 - \gamma|)m_1 + |\gamma_1 - \gamma|(1-m_1), \ldots, (1-|\gamma_1 - \gamma|)m_d + |\gamma_1 - \gamma|(1-m_d))$. Consequently, in $j_3$ we have

$$p_1 + 2A_4 s_1 + (-1)^\gamma m_3 + 2x_3 = p_1^* + 2A_4 s_1^* + (-1)^{\gamma_1} m_3^* + 2x_3^*$$
$$\implies (s_1 - s_1^*) - \frac{(-1)^{\gamma_1} - (-1)^\gamma}{2} m_3 = \frac{(-1)^{\gamma_1}(m_3^* - m_3)}{2}.$$

Since $(s_1 - s_1^*) - \frac{(-1)^{\gamma_1} - (-1)^\gamma}{2} m_3 \in \mathbb{Z}$ and $|m_3^* - m_3| \leq 1$ by construction, it must be the case $m_3 = m_3^*$, meaning $(m_1, \ldots, m_d) = (m_1^*, \ldots, m_d^*)$. Our deduction $m_3 = m_3^*$ implies

$$m_3^* - m_3 = 0$$
$$\implies (-1)^{m_3} |\gamma_1 - \gamma| = 0.$$

From the above, it must be the case $\gamma = \gamma_1$. So our last equality in $j_3$ implies $s_1 = s_1^*$. Now, in $j_2$ our results give us

$$m_2 + \gamma + 4x_2 + (1 - m_2)m_1 + 2m_d(1 - m_1) = m_2^* + \gamma_1 + 4x_2^* + (1 - m_2^*)m_1^* + 2m_d^*(1 - m_1^*)$$

$$\implies x_2 = x_2^*.$$

**Major Case 4:** Focusing on moves along $j_k$ for $4 \leq k \leq d$ when $d \geq 4$, we see that during all such moves $m_1 = \cdots = m_{k-1}$, $m_k = \cdots = m_d = 1 - m_1$ and $(m_1^*, \ldots, m_d^*) = (m_1, \ldots, m_d)$, giving us in $j_2$

$$m_2 + \gamma + 4x_2 + (1 - m_2)m_1 + 2m_d(1 - m_1) = m_2^* + \gamma_1 + 4x_2^* + (1 - m_2^*)m_1^* + 2m_d^*(1 - m_1^*)$$

$$\implies x_2 - x_2^* = \frac{\gamma_1 - \gamma}{4}.$$

Following from $x_2 - x_2^* \in \mathbb{Z}$ and $|\gamma_1 - \gamma| \leq 1$ by construction, it must be the case $\gamma = \gamma_1$ and so $x_2 = x_2^*$. Knowing that $m_3 = m_3^*$ and $\gamma = \gamma_1$, in $j_3$ we find

$$p_1 + 2A_4 s_1 + (-1)^\gamma m_3 + 2x_3 = p_1^* + 2A_4 s_1^* + (-1)^{\gamma_1} m_3^* + 2x_3^*$$

$$\implies s_1 = s_1^*.$$

To conclude our treatment of all the major cases, we will only work with the following assumptions that we deduced in all major cases: $(m_1, \ldots, m_d) = (m_1^*, \ldots, m_d^*)$, $\gamma = \gamma_1$, $x_2 = x_2^*$, and $s_1 = s_1^*$.

Applying our assumptions to $j_1$, we obtain

$$(\ell + 1 - \gamma)A_3 - \gamma + t + (-1)^\gamma(-2 + m_1 + (1 - m_2)m_1 + 2m_d(1 - m_1))$$
$$= (\ell_1 + 1 - \gamma_1)A_3 - \gamma_1 + t_1 + (-1)^{\gamma_1}(-2 + m_1^* + (1 - m_2^*)m_1^* + 2m_d^*(1 - m_1^*))$$

$$\implies \ell - \ell_1 = \frac{t_1 - t}{A_3}.$$

Since $\ell - \ell_1 \in \mathbb{Z}$ and $|t_1 - t| \leq 1$ with $A_3 \geq 2$ by construction, it must be the case $t = t_1$ and so $\ell = \ell_1$. Now, for all remaining components $j_k$ for $4 \leq k \leq d$ when $d \geq 4$, observe

$$p_{k-2} + A_{k+1} s_{k-2} + m_k + 2x_k = p_{k-2}^* + A_{k+1} s_{k-2}^* + m_k^* + 2x_k^*$$

$$\implies s_{k-2} = s_{k-2}^*.$$

Thus, for all $4 \leq k \leq d$, $s_{k-2} = s_{k-2}^*$. So $C_1 = C_{\ell, \gamma, t, p_1, s_1, \ldots, p_{d-2}, s_{d-2}} = C_{\ell_1, \gamma_1, t_1, p_1^*, s_1^*, \ldots, p_{d-2}^*, s_{d-2}^*} = C_2$.

**Set Comparison Case 2:** We will consider when edges are shared from moves with the same orientation via set 4.

**Major Case 1:** Focusing on moves along $j_1$, we see that $C_1$ and $C_2$ have $j_1$ moves with the same orientation if and only if $\gamma = \gamma_1$. Noting that during these moves $m_1 = m_3 = \cdots = m_d$, $m_2 = 1 = m_2^*$ and $m_1^* = m_3^* = \cdots = m_d^*$, in $j_2$ we see

$$2 + m_1 + \gamma + m_2(1 - m_1) + 4x_2 + 2(1 - m_1)m_d = 2 + m_1^* + \gamma_1 + m_2^*(1 - m_1^*) + 4x_2^* + 2(1 - m_1^*)m_d^*$$

$$\implies x_2 = x_2^*.$$

Now, in $j_3$ our assumptions yield

$$2s_1 + (-1)^\gamma m_3 = 2s_1^* + (-1)^{\gamma_1} m_3^*$$

$$\implies s_1 - s_1^* = \frac{(-1)^\gamma (m_3^* - m_3)}{2}.$$

Following from $s_1 - s_1^* \in \mathbb{Z}$ and $|m_3^* - m_3| \leq 1$ by construction, it must be the case $m_3 = m_3^*$ and so $s_1 = s_1^*$. Hence, $(m_1, \ldots, m_d) = (m_1^*, \ldots, m_d^*)$.

<u>Major Case 2</u>: Focusing on moves along $j_2$, observe that during all $j_2$ moves $m_1 = m_3 = \cdots = m_d$, $m_2 = 0 = m_2^*$ and $m_1^* = m_3^* = \cdots = m_d^*$, giving us in $j_2$

$$2 + m_1 + \gamma + m_2(1 - m_1) + 4x_2 + 2(1 - m_1)m_d = 2 + m_1^* + \gamma_1 + m_2^*(1 - m_1^*) + 4x_2^* + 2(1 - m_1^*)m_d^*$$

$$\implies x_2 - x_2^* = \frac{(m_1^* - m_1) + (\gamma_1 - \gamma)}{4}.$$

Since $x_2 - x_2^* \in \mathbb{Z}$ and $-2 \leq (m_1^* - m_1) + (\gamma_1 - \gamma) \leq 2$ by construction, we must have $(m_1^* - m_1) = -(\gamma_1 - \gamma)$, and so $x_2 = x_2^*$. Now, in $j_3$ we get

$$2s_1 + (-1)^\gamma m_3 = 2s_1^* + (-1)^{\gamma_1} m_3^*$$

$$\implies (s_1 - s_1^*) + \frac{(-1)^\gamma - (-1)^{\gamma_1}}{2} m_3 = \frac{(-1)^{\gamma_1}(m_3^* - m_3)}{2}.$$

Following from $(s_1 - s_1^*) + \frac{(-1)^\gamma - (-1)^{\gamma_1}}{2} m_3 \in \mathbb{Z}$ and $|m_3^* - m_3| \leq 1$ by construction, it must be the case $m_3 = m_3^*$, which by our result in $j_2$ implies $\gamma = \gamma_1$. So our last equality in $j_3$ implies $s_1 = s_1^*$.

<u>Major Case 3</u>: Focusing on moves along $j_3$, we see that during these moves $m_1 = m_2$, $m_3 = \cdots = m_d = 1 - m_1$, and $(m_1^*, \ldots, m_d^*) = ((1 - |\gamma_1 - \gamma|)m_1 + |\gamma_1 - \gamma|(1 - m_1), \ldots, (1 - |\gamma_1 - \gamma|)m_d + |\gamma_1 - \gamma|(1 - m_d))$. Hence, in $j_3$ we see

$$2s_1 + (-1)^\gamma m_3 = 2s_1^* + (-1)^{\gamma_1} m_3^*$$

$$\implies (s_1 - s_1^*) + \frac{(-1)^\gamma - (-1)^{\gamma_1}}{2} m_3 = \frac{(-1)^{\gamma_1}(m_3^* - m_3)}{2}.$$

Following from $(s_1 - s_1^*) + \frac{(-1)^\gamma - (-1)^{\gamma_1}}{2} m_3 \in \mathbb{Z}$ and $|m_3^* - m_3| \leq 1$ by construction, it must be the case $m_3 = m_3^*$. This implies

$$m_3^* - m_3 = 0$$

$$\implies (-1)^{m_3} |\gamma_1 - \gamma| = 0.$$

Since $\gamma = \gamma_1$ by the equality above, our last equality in $j_3$ implies $s_1 = s_1^*$ and $(m_1, \ldots, m_d) = (m_1^*, \ldots, m_d^*)$. Hence, in $j_2$ we have

$$2 + m_1 + \gamma + m_2(1 - m_1) + 4x_2 + 2(1 - m_1)m_d = 2 + m_1^* + \gamma_1 + m_2^*(1 - m_1^*) + 4x_2^* + 2(1 - m_1^*)m_d^*$$

$$\implies x_2 = x_2^*.$$

<u>Major Case 4</u>: Focusing on moves along $j_k$ for $4 \leq k \leq d$ when $d \geq 4$, we see that $m_1 = \cdots = m_{k-1}$, $m_k = \cdots = m_d = 1 - m_1$, and $(m_1, \ldots, m_d) = (m_1^*, \ldots, m_d^*)$. These observations in $j_2$ give us

$$2 + m_1 + \gamma + m_2(1 - m_1) + 4x_2 + 2(1 - m_1)m_d = 2 + m_1^* + \gamma_1 + m_2^*(1 - m_1^*) + 4x_2^* + 2(1 - m_1^*)m_d^*$$

$$\implies x_2 - x_2^* = \frac{\gamma_1 - \gamma}{4}.$$

Since $x_2 - x_2^* \in \mathbb{Z}$ and $|\gamma_1 - \gamma| \leq 1$ by construction, it follows that $\gamma = \gamma_1$ and so $x_2 = x_2^*$. Applying this to $j_3$, we get

$$2s_1 + (-1)^\gamma m_3 = 2s_1^* + (-1)^{\gamma_1} m_3^*$$

$$\implies s_1 = s_1^*.$$

To conclude our treatment of all major cases, we will only work with the following assumptions that we deduced in all major cases: $(m_1, \ldots, m_d) = (m_1^*, \ldots, m_d^*)$, $\gamma = \gamma_1$, $x_2 = x_2^*$, and $s_1 = s_1^*$.

Applying our assumptions to $j_1$, we have

$$(\ell + 1 - \gamma)A_3 - \gamma + t + (-1)^{\gamma+1}(m_1 + m_1(1 - m_2) + 2(1 - m_1)m_d)$$
$$= (\ell_1 + 1 - \gamma_1)A_3 - \gamma_1 + t_1 + (-1)^{\gamma_1+1}(m_1^* + m_1^*(1 - m_2^*) + 2(1 - m_1^*)m_d^*)$$

$$\implies \ell - \ell_1 = \frac{t_1 - t}{A_3}.$$

Given $\ell - \ell_1 \in \mathbb{Z}$ and $|t_1 - t| \leq 1$ with $A_3 \geq 2$ by construction, we must have $t = t_1$ and so $\ell = \ell_1$. Inspecting the remaining components $j_k$ for $4 \leq k \leq d$ when $d \geq 4$, observe that $m_k = m_k^*$ by our initial deductions and so

$$s_{k-2} + m_k = s_{k-2}^* + m_k^*$$

$$\implies s_{k-2} = s_{k-2}^*.$$

Thus, for all $4 \leq k \leq d$, $s_{k-2} = s_{k-2}^*$. So $C_1 = C_{\ell, \gamma, t, p_1, s_1, \ldots, p_{d-2}, s_{d-2}} = C_{\ell_1, \gamma_1, t_1, p_1^*, s_1^*, \ldots, p_{d-2}^*, s_{d-2}^*} = C_2$.

Consequently, edges are shared from moves with the same orientation from the same set for $\alpha^* = \alpha_1$ when $d \geq 3$ if and only if $C_1 = C_2$.

$\boldsymbol{\alpha^* -}$**Case 2**: Let $\alpha_1 < \alpha^* < \alpha_2$. Then, $A_1 = 1$ and $A_2 = 0$, meaning sets 1, 2, 3, and 4 are all active. Note that $A_3 \geq 4$ is even, $t = 0 = t_1$, and for all $3 \leq k \leq d$, $p_{k-2} = 0 = p_{k-2}^*$, $A_{k+1} = 1$ and $x_k = 0 = x_k^*$.

**Set Comparison Case 1:** We will consider when edges are shared from moves with the same orientation via set 1.

<u>Major Case 1</u>: Focusing on moves along $j_1$, we see that during all $j_1$ moves $m_1 = m_3 = \cdots = m_d$, $m_2 = 1 - m_1$, $m_1^* = m_3^* = \cdots = m_d^*$, and $m_2^* = 1 - m_1^*$. Consequently, in $j_1$ we get

$$\ell A_3 + 1 + (m_1 + 2x_1 + 2m_d(1 - m_1)) = \ell_1 A_3 + 1 + (m_1^* + 2x_1^* + 2m_d^*(1 - m_1^*))$$

$$\implies (\ell - \ell_1)\frac{A_3}{2} + (x_1 - x_1^*) = \frac{m_1^* - m_1}{2}.$$

Since $(\ell - \ell_1)\frac{A_3}{2} + (x_1 - x_1^*) \in \mathbb{Z}$ and $|m_1^* - m_1| \leq 1$ by construction, it must be the case $m_1 = m_1^*$, meaning $(m_1, \ldots, m_d) = (m_1^*, \ldots, m_d^*)$. So we are left with

$$\ell - \ell_1 = \frac{2(x_1^* - x_1)}{A_3}.$$

Since $\ell - \ell_1 \in \mathbb{Z}$ and $2|x_1^* - x_1| \leq A_3 - 4$, it must be the case $x_1 = x_1^*$ and so $\ell = \ell_1$.

<u>Major Case 2:</u> Focusing on moves along $j_z$ for all $2 \leq z \leq d$, we see that during these moves $(m_1^*, \ldots, m_d^*) = (m_1, \ldots, m_d)$. Applying these observations to $j_1$, we obtain

$$\ell A_3 + 1 + (m_1 + 2x_1 + 2m_d(1 - m_1)) = \ell_1 A_3 + 1 + (m_1^* + 2x_1^* + 2m_d^*(1 - m_1^*))$$
$$\implies \ell - \ell_1 = \frac{2(x_1^* - x_1)}{A_3}.$$

Following from $\ell - \ell_1 \in \mathbb{Z}$ and $2|x_1 - x_1^*| \leq A_3 - 4$ by construction, we must have $x_1 = x_1^*$ and so $\ell = \ell_1$.

To conclude our treatment of all the major cases, we will work with the following assumptions that we were able to deduce in all major cases: $(m_1, \ldots, m_d) = (m_1^*, \ldots, m_d^*)$, $x_1 = x_1^*$, and $\ell = \ell_1$.

Our assumptions in $j_k$ for all $3 \leq k \leq d$ give us

$$p_{k-2} + A_{k+1}s_{k-2} + (1 - m_k) + 2x_k = p_{k-2}^* + A_{k+1}s_{k-2}^* + (1 - m_k^*) + 2x_k^*$$
$$\implies s_{k-2} = s_{k-2}^*.$$

So for all $3 \leq k \leq d$, $s_{k-2} = s_{k-2}^*$. Observing in particular that $s_1 = s_1^*$, in $j_2$ we find

$$m_2 + \gamma + (s_1 + 1) + 4x_2 = m_2^* + \gamma_1 + (s_1^* + 1) + 4x_2^*$$
$$\implies x_2 - x_2^* = \frac{\gamma_1 - \gamma}{4}.$$

Since $x_2 - x_2^* \in \mathbb{Z}$ and $|\gamma_1 - \gamma| \leq 1$ by construction, it must be the case $\gamma = \gamma_1$ and so $x_2 = x_2^*$. Thus, $(m_1, \ldots, m_d) = (m_1^*, \ldots, m_d^*)$ at the same start vertex, and $C_1 = C_{\ell, \gamma, t, p_1, s_1, \ldots, p_{d-2}, s_{d-2}} = C_{\ell_1, \gamma_1, t_1, p_1^*, s_1^*, \ldots, p_{d-2}^*, s_{d-2}^*} = C_2$.

**Set Comparison Case 2:** We will consider when edges are shared from moves with the same orientation via set 2.

<u>Major Case 1:</u> Focusing on moves along $j_1$, we see that during these moves $m_1 = \cdots = m_d$ and $m_1^* = \cdots = m_d^*$. Consequently, noting that $r_1 = 0 = r_1^*$ as $m_1 = m_3$ and $m_1^* = m_3^*$ when $d = 3$ and $m_{d-1} = m_d$ and $m_{d-1}^* = m_d^*$ when $d \geq 4$, in $j_1$ we obtain

$$(\ell + 1)A_3 - (2 + m_2 + 2m_1(1 - m_2) + 2x_1 + 2m_d(1 - m_1) - r_1(A_3 - 1))$$
$$= (\ell_1 + 1)A_3 - (2 + m_2^* + 2m_1^*(1 - m_2^*) + 2x_1^* + 2m_d^*(1 - m_1^*) - r_1^*(A_3 - 1))$$
$$\implies (\ell - \ell_1)\frac{A_3}{2} + (x_1^* - x_1) = \frac{m_2 - m_2^*}{2}.$$

Given $(\ell - \ell_1)\frac{A_3}{2} + (x_1^* - x_1) \in \mathbb{Z}$ and $|m_2^* - m_2| \leq 1$ by construction, we must have $m_2 = m_2^*$, meaning $(m_1, \ldots, m_d) = (m_1^*, \ldots, m_d^*)$. This leaves us with

$$\ell - \ell_1 = \frac{2(x_1 - x_1^*)}{A_3}.$$

Following from $\ell - \ell_1 \in \mathbb{Z}$ and $2|x_1 - x_1^*| \leq A_3 - 4$ by construction, we see $x_1 = x_1^*$ and so $\ell = \ell_1$.

<u>Major Case 2</u>: Focusing on moves along $j_2$, we see that which $j_2$ moves in $C_1$ and $C_2$ have the same orientation depend on $x_1$ and $x_1^*$. Consequently, we case on $x_1$ and further case on $x_1^*$:

<u>Case 1</u>: Let $0 \leq x_1 < \frac{A_3}{2} - 2$. We now case on $x_1^*$:

<u>Subcase 1</u>: Let $0 \leq x_1^* < \frac{A_3}{2} - 2$. Observing that $(m_1, \ldots, m_d) = (m_1^*, \ldots, m_d^*)$ and $r_1 = 0 = r_1^*$ since $x_1 \neq \frac{A_3}{2} - 2 \neq x_1^*$, in $j_1$ we get

$$(\ell + 1)A_3 - (2 + m_2 + 2m_1(1 - m_2) + 2x_1 + 2m_d(1 - m_1) - r_1(A_3 - 1))$$
$$= (\ell_1 + 1)A_3 - (2 + m_2^* + 2m_1^*(1 - m_2^*) + 2x_1^* + 2m_d^*(1 - m_1^*) - r_1^*(A_3 - 1))$$
$$\implies \ell - \ell_1 = \frac{2(x_1 - x_1^*)}{A_3}.$$

Following from $\ell - \ell_1 \in \mathbb{Z}$ and $2|x_1 - x_1^*| < A_3 - 4$ by construction, it must be the case $x_1 = x_1^*$ and so $\ell = \ell_1$.

<u>Subcase 2</u>: Let $x_1^* = \frac{A_3}{2} - 2$. Observe that $(m_1, \ldots, m_{d-1}, m_d) = (1, \ldots, 1, 0)$ during $C_1$'s only $j_2$ move with the same orientation as those of $C_2$, and $m_1^* = \cdots = m_{d-1}^*$ and $m_d^* = 1 - m_1^*$. Now, since $r_1 = 0$ as $x_1 \neq \frac{A_3}{2} - 2$, in $j_1$ we see

$$(\ell + 1)A_3 - (2 + m_2 + 2m_1(1 - m_2) + 2x_1 + 2m_d(1 - m_1) - r_1(A_3 - 1))$$
$$= (\ell_1 + 1)A_3 - (2 + m_2^* + 2m_1^*(1 - m_2^*) + 2x_1^* + 2m_d^*(1 - m_1^*) - r_1^*(A_3 - 1))$$
$$\implies (\ell - \ell_1) + (1 - r_1^*) = \frac{3 + 2x_1 + (m_1^* - r_1^*)}{A_3}.$$

Given $(\ell - \ell_1) + (1 - r_1^*) \in \mathbb{Z}$ and $2 \leq 3 + 2x_1 + (m_1^* - r_1^*) < A_3$ by construction and our case assumption, it follows that for every $0 \leq x_1 < \frac{A_3}{2} - 2$ and $m_1^*, r_1^* \in \{0, 1\}$, $\frac{3+2x_1+(m_1^*-r_1^*)}{A_3} \notin \mathbb{Z}$. Hence, we have a contradiction as the above equality is not possible.

<u>Case 2</u>: Let $x_1 = \frac{A_3}{2} - 2$. Then, since the case with $0 \leq x_1^* < \frac{A_3}{2} - 2$ is symmetric to subcase 2 of the previous case, it is already addressed by the argument there. So let $x_1^* = \frac{A_3}{2} - 2$. Now observe that $m_1 = \cdots = m_{d-1}$, $m_d = 1 - m_1$, $m_1^* = \cdots = m_{d-1}^*$, and $m_d^* = 1 - m_1^*$. Applying this to $j_1$, we get

$$(\ell + 1)A_3 - (2 + m_2 + 2m_1(1 - m_2) + 2x_1 + 2m_d(1 - m_1) - r_1(A_3 - 1))$$
$$= (\ell_1 + 1)A_3 - (2 + m_2^* + 2m_1^*(1 - m_2^*) + 2x_1^* + 2m_d^*(1 - m_1^*) - r_1^*(A_3 - 1))$$
$$\implies (\ell - \ell_1) + (r_1 - r_1^*) = \frac{(r_1 - r_1^*) + (m_d - m_d^*)}{A_3}.$$

Since $(\ell - \ell_1) + (r_1 - r_1^*) \in \mathbb{Z}$ and $-2 \leq (r_1 - r_1^*) + (m_d - m_d^*) \leq 2$ with $A_3 \geq 4$ by construction, it must be the case $(r_1 - r_1^*) = -(m_d - m_d^*)$. This implies

$$(r_1 - r_1^*) + (m_d - m_d^*) = 0$$
$$\implies (r_1 + m_d) - (r_1^* + m_d^*) = 0.$$

If $d = 3$, then $r_1 = m_3(1 - m_1)$ and $r_1^* = m_3^*(1 - m_1^*)$. Recalling that $m_3 = 1 - m_1$ and $m_3^* = 1 - m_1^*$ in this case, we get
$$m_3(2 - m_1) - m_3^*(2 - m_1^*) = 0$$
$$\implies m_3 = m_3^*.$$

If $d \geq 4$, then $r_1 = m_d(1 - m_{d-1})$ and $r_1^* = m_d^*(1 - m_{d-1}^*)$. Noting that $m_d = 1 - m_{d-1}$ and $m_d^* = 1 - m_{d-1}^*$ in this case, we obtain
$$m_d(2 - m_{d-1}) - m_d^*(2 - m_{d-1}^*) = 0$$
$$\implies m_d = m_d^*.$$

For all $d \geq 3$, we have $(m_1, \ldots, m_d) = (m_1^*, \ldots, m_d^*)$. So $r_1 = r_1^*$ and hence our last equality implies $\ell = \ell_1$.

Major Case 3: Focusing on moves along $j_k$ for $3 \leq k \leq d$, observe that during all such moves $(m_1^*, \ldots, m_d^*) = (m_1, \ldots, m_d)$. So in $j_1$ we find

$$(\ell + 1)A_3 - (2 + m_2 + 2m_1(1 - m_2) + 2x_1 + 2m_d(1 - m_1) - r_1(A_3 - 1))$$
$$= (\ell_1 + 1)A_3 - (2 + m_2^* + 2m_1^*(1 - m_2^*) + 2x_1^* + 2m_d^*(1 - m_1^*) - r_1^*(A_3 - 1))$$
$$\implies ((\ell - \ell_1) + (r_1 - r_1^*))\frac{A_3}{2} + (x_1^* - x_1) = \frac{r_1 - r_1^*}{2}.$$

Since $((\ell - \ell_1) + (r_1 - r_1^*))\frac{A_3}{2} + (x_1^* - x_1) \in \mathbb{Z}$ and $|r_1 - r_1^*| \leq 1$ by construction, we must have $r_1 = r_1^*$. This leaves us with
$$\ell - \ell_1 = \frac{2(x_1 - x_1^*)}{A_3}.$$

Following from $\ell - \ell_1 \in \mathbb{Z}$ and $2|x_1 - x_1^*| \leq A_3 - 4$ by construction, we see it must be the case $x_1 = x_1^*$ and so $\ell = \ell_1$.

To conclude our treatment of the major cases and subcases in which we did not reach a contradiction, we will make use of the following assumptions we were able to deduce in said cases: $(m_1, \ldots, m_d) = (m_1^*, \ldots, m_d^*)$, $x_1 = x_1^*$, $r_1 = r_1^*$, and $\ell = \ell_1$.

Applying our results to $j_3$ and $j_k$ for $4 \leq k \leq d$ when $d \geq 4$, the equalities in $j_3$ and $j_k$
$$s_1 + m_1 + r_1 = s_1^* + m_1^* + r_1^*$$
$$s_{k-2} + m_{k-1} + r_1 = s_{k-2}^* + m_{k-1}^* + r_1^*$$

imply $s_1 = s_1^*$ and for all $4 \leq k \leq d$ when $d \geq 4$, $s_{k-2} = s_{k-2}^*$. Since $s_1 = s_1^*$, our assumptions in $j_2$ yield
$$2 + \gamma + (1 + s_1 + m_d - r_1) + 4x_2 = 2 + \gamma_1 + (1 + s_1^* + m_d^* - r_1^*) + 4x_2^*$$

$$\implies x_2 - x_2^* = \frac{\gamma_1 - \gamma}{4}.$$

Given $x_2 - x_2^* \in \mathbb{Z}$ and $|\gamma_1 - \gamma| \leq 1$ by construction, we must have $\gamma = \gamma_1$ and so $x_2 = x_2^*$. Thus, $C_1 = C_{\ell,\gamma,t,p_1,s_1,\ldots,p_{d-2},s_{d-2}} = C_{\ell_1,\gamma_1,t_1,p_1^*,s_1^*,\ldots,p_{d-2}^*,s_{d-2}^*} = C_2$.

**Set Comparison Case 3:** We will consider when edges are shared from moves with the same orientation via set 3.

<u>Major Case 1</u>: Focusing on moves along $j_1$ and $j_2$, $m_1 = m_d$ and $m_1^* = m_d^*$ during these moves. Applying this to $j_1$, we have

$$(\ell+1)A_3 + (-2 + m_1 + 2m_d(1 - m_1) - (A_3 - 2)(1 - m_d))$$
$$= (\ell_1 + 1)A_2 + (-2 + m_1^* + 2m_d^*(1 - m_1^*) - (A_3 - 2)(1 - m_d^*))$$
$$\implies (\ell - \ell_1) + (m_d - m_d^*) = \frac{m_1^* - m_1}{A_3}.$$

Following from $(\ell - \ell_1) + (m_d - m_d^*) \in \mathbb{Z}$ and $|m_1^* - m_1| \leq 1$ with $A_3 \geq 4$ by construction, it must be the case $m_1 = m_1^*$. During $j_1$ moves, $m_1 = \cdots = m_d$ and $m_1^* = \cdots = m_d^*$. During $j_2$ moves, $m_1 = m_3 = \cdots m_d$, $m_2 = 1 - m_1$, $m_1^* = m_3^* = \cdots = m_d^*$, and $m_2^* = 1 - m_1^*$. Hence, in both cases our last result implies $(m_1, \ldots, m_d) = (m_1^*, \ldots, m_d^*)$ and so $\ell = \ell_1$ by the equality above.

<u>Major Case 2</u>: Focusing on moves along $j_k$ for $3 \leq k \leq d$, we see that during all of these moves $(m_1, \ldots, m_d) = (m_1^*, \ldots, m_d^*)$ and so in $j_1$ we have

$$(\ell+1)A_3 + (-2 + m_1 + 2m_d(1 - m_1) - (A_3 - 2)(1 - m_d))$$
$$= (\ell_1 + 1)A_2 + (-2 + m_1^* + 2m_d^*(1 - m_1^*) - (A_3 - 2)(1 - m_d^*))$$
$$\implies \ell = \ell_1.$$

To conclude our treatment of these major cases, we will work with the following assumptions that we deduced in all major cases: $(m_1, \ldots, m_d) = (m_1^*, \ldots, m_d^*)$ and $\ell = \ell_1$.

Inspecting $j_k$ for $3 \leq k \leq d$, we observe

$$p_{k-2} + A_{k+1}s_{k-2} + m_k + 2x_k = p_{k-2}^* + A_{k+1}s_{k-2}^* + m_k^* + 2x_k^*$$
$$\implies s_{k-2} = s_{k-2}^*.$$

Hence, for all $3 \leq k \leq d$, $s_{k-2} = s_{k-2}^*$. In particular, $s_1 = s_1^*$, implying in $j_2$

$$m_2 + \gamma + 4x_2 + 2(1 - m_2)m_1 + (s_1 + 2(1 - m_1)m_d) = m_2^* + \gamma_1 + 4x_2^* + 2(1 - m_2^*)m_1^* + (s_1^* + 2(1 - m_1^*)m_d^*)$$
$$\implies x_2 - x_2^* = \frac{\gamma_1 - \gamma}{4}.$$

Since $x_2 - x_2^* \in \mathbb{Z}$ and $|\gamma_1 - \gamma| \leq 1$ by construction, we must have $\gamma = \gamma_1$ and so $x_2 = x_2^*$. Consequently, $C_1 = C_{\ell,\gamma,t,p_1,s_1,\ldots,p_{d-2},s_{d-2}} = C_{\ell_1,\gamma_1,t_1,p_1^*,s_1^*,\ldots,p_{d-2}^*,s_{d-2}^*} = C_2$.

**Set Comparison Case 4:** We will consider when edges are shared from moves with the same orientation via set 4.

<u>Major Case 1</u>: Focusing on moves along $j_1$, we observe that during these moves $m_1 = m_3 = \cdots = m_d$, $m_2 = 1 = m_2^*$, and $m_1^* = m_3^* = \cdots = m_d^*$. Applying these observations to $j_1$, we get

$$(\ell+1)A_3 - (m_1 + m_1(1-m_2) + 2(1-m_1)m_d) = (\ell_1+1)A_3 - (m_1^* + m_1^*(1-m_2^*) + 2(1-m_1^*)m_d^*)$$

$$\implies \ell - \ell_1 = \frac{m_1 - m_1^*}{A_3}.$$

Given $\ell - \ell_1 \in \mathbb{Z}$ and $|m_1 - m_1^*| \leq 1$ by construction, it must be the case $m_1 = m_1^*$, meaning $(m_1, \ldots, m_d) = (m_1^*, \ldots, m_d^*)$, and so $\ell = \ell_1$.

<u>Major Case 2</u>: Focusing on moves along $j_k$ for $2 \leq k \leq d$, note that during all such moves $(m_1^*, \ldots, m_d^*) = (m_1, \ldots, m_d)$. Applying this to $j_1$, we have

$$(\ell+1)A_3 - (m_1 + m_1(1-m_2) + 2(1-m_1)m_d) = (\ell_1+1)A_3 - (m_1^* + m_1^*(1-m_2^*) + 2(1-m_1^*)m_d^*)$$

$$\implies \ell = \ell_1.$$

To conclude our treatment of all the major cases, we will work with the following assumptions we deduced in all major cases: $(m_1, \ldots, m_d) = (m_1^*, \ldots, m_d^*)$ and $\ell = \ell_1$.

From our assumptions, the equalities $j_3$ and $j_k$ for $4 \leq k \leq d$

$$s_1 + m_1 m_3 = s_1^* + m_1^* m_3^*$$

$$s_{k-2} + m_{k-1} m_k = s_{k-2}^* + m_{k-1}^* m_k^*$$

imply $s_1 = s_1^*$ and $s_{k-2} = s_{k-2}^*$ for all $4 \leq k \leq d$. Since we have $s_1 = s_1^*$, in $j_2$ we obtain

$$2 + m_1 + (s_1 + (1-m_1)m_d + m_d) + \gamma + m_2(1-m_1) + 4x_2$$
$$= 2 + m_1^* + (s_1^* + (1-m_1^*)m_d^* + m_d^*) + \gamma_1 + m_2^*(1-m_1^*) + 4x_2^*$$

$$\implies x_2 - x_2^* = \frac{\gamma_1 - \gamma}{4}.$$

Noting that $x_2 - x_2^* \in \mathbb{Z}$ and $|\gamma_1 - \gamma| \leq 1$ by construction, we see $\gamma = \gamma_1$ and so $x_2 = x_2^*$. Thus, $C_1 = C_{\ell,\gamma,t,p_1,s_1,\ldots,p_{d-2},s_{d-2}} = C_{\ell_1,\gamma_1,t_1,p_1^*,s_1^*,\ldots,p_{d-2}^*,s_{d-2}^*} = C_2$.

Hence, for $\alpha_1 < \alpha^* < \alpha_2$ when $d \geq 3$, edges are shared from moves with the same orientation if and only if $C_1 = C_2$.

**$\alpha^*$−Case 3:** Let $\alpha^* = \alpha_2$. Then, $A_1 = 1$ and $A_2 = 1$, meaning only sets 1 and 3 are active. Note that $t = 0 = t_1$, $A_3 \geq 4$ is even, and for all $3 \leq k \leq d$, $p_{k-2} = 0 = p_{k-2}^*$, $A_{k+1} = 1$ and $x_k = 0 = x_k^*$.

**Set Comparison Case 1:** We will consider when edges are shared from moves with the same orientation via set 1.

<u>Major Case 1</u>: Focusing on moves along $j_1$, we see that $m_1 = m_3 = \cdots = m_d$, $m_2 = 1 - m_1$, $m_1^* = m_3^* = \cdots = m_d^*$, and $m_2^* = 1 - m_1^*$. Consequently, the above in $j_1$ yields

$$\ell A_3 + (m_1 + 2x_1 + 2m_d(1 - m_1)) = \ell_1 A_3 + (m_1^* + 2x_1^* + 2m_d^*(1 - m_1^*))$$
$$\implies (\ell - \ell_1)\frac{A_3}{2} + (x_1 - x_1^*) = \frac{m_1^* - m_1}{2}.$$

Since $(\ell - \ell_1)\frac{A_3}{2} + (x_1 - x_1^*) \in \mathbb{Z}$ and $|m_1^* - m_1| \leq 1$ with $A_3 \geq 4$ by construction, $m_1 = m_1^*$, meaning $(m_1, \ldots, m_d) = (m_1^*, \ldots, m_d^*)$. So we are left with

$$\ell - \ell_1 = \frac{2(x_1 - x_1^*)}{A_3}.$$

Given $\ell - \ell_1 \in \mathbb{Z}$ and $2|x_1 - x_1^*| \leq A_3 - 4$ by construction, it follows that $x_1 = x_1^*$ and so $\ell = \ell_1$.

<u>Major Case 2:</u> Focusing on moves along $j_k$ for $2 \leq k \leq d$, observe that during these moves $(m_1^*, \ldots, m_d^*) = (m_1, \ldots, m_d)$ and so in $j_1$ we find

$$\ell A_3 + (m_1 + 2x_1 + 2m_d(1 - m_1)) = \ell_1 A_3 + (m_1^* + 2x_1^* + 2m_d^*(1 - m_1^*))$$
$$\implies \ell - \ell_1 = \frac{2(x_1 - x_1^*)}{A_3}.$$

Since $\ell - \ell_1 \in \mathbb{Z}$ and $2|x_1 - x_1^*| \leq A_3 - 4$ by construction, it follows that $x_1 = x_1^*$ and so $\ell = \ell_1$.

To conclude our treatment of all the major cases, we will work with the following assumptions we deduced in all major cases: $(m_1, \ldots, m_d) = (m_1^*, \ldots, m_d^*)$, $x_1 = x_1^*$, and $\ell = \ell_1$.

Now, our assumptions in $j_2$ imply

$$m_2 + \gamma + 2x_2 = m_2^* + \gamma_1 + 2x_2^*$$
$$\implies x_2 - x_2^* = \frac{\gamma_1 - \gamma}{2}.$$

Following from $x_2 - x_2^* \in \mathbb{Z}$ and $|\gamma_1 - \gamma| \leq 1$ by construction, it must be the case $\gamma = \gamma_1$ and so $x_2 = x_2^*$. Inspecting $j_k$ for $3 \leq k \leq d$, observe

$$p_{k-2} + A_{k+1}s_{k-2} + m_k + 2x_k = p_{k-2}^* + A_{k+1}s_{k-2}^* + m_k^* + 2x_k^*$$
$$\implies s_{k-2} = s_{k-2}^*.$$

Hence, for all $3 \leq k \leq d$, $s_{k-2} = s_{k-2}^*$ and so $C_1 = C_{\ell,\gamma,t,p_1,s_1,\ldots,p_{d-2},s_{d-2}} = C_{\ell_1,\gamma_1,t_1,p_1^*,s_1^*,\ldots,p_{d-2}^*,s_{d-2}^*} = C_2$.

**Set Comparison Case 2:** We will consider when edges are shared from moves with the same orientation via set 3.

<u>Major Case 1:</u> Focusing on moves along $j_1$ and $j_2$, we see that during these moves $m_1 = m_d$ and $m_1^* = m_d^*$, giving us in $j_1$

$$(\ell + 1)A_3 + (-2 + m_1 + 2m_d(1 - m_1)) = (\ell_1 + 1)A_3 + (-2 + m_1^* + 2m_d^*(1 - m_1^*))$$
$$\implies \ell - \ell_1 = \frac{m_1^* - m_1}{A_3}.$$

Since $\ell - \ell_1 \in \mathbb{Z}$ and $|m_1^* - m_1| \leq 1$ with $A_3 \geq 4$ by construction, we must have $m_1 = m_1^*$, meaning $(m_1, \ldots, m_d) = (m_1^*, \ldots, m_d^*)$, and so $\ell = \ell_1$.

**Major Case 2:** Focusing on moves along $j_k$ for $3 \leq k \leq d$, we see that during all such moves $(m_1, \ldots, m_d) = (m_1^*, \ldots, m_d^*)$, giving us in $j_1$

$$(\ell+1)A_3 + (-2 + m_1 + 2m_d(1-m_1)) = (\ell_1+1)A_3 + (-2 + m_1^* + 2m_d^*(1-m_1^*))$$
$$\implies \ell = \ell_1.$$

To conclude our treatment of all the major cases, we will work with the following assumptions we were able to deduce in all major cases: $(m_1, \ldots, m_d) = (m_1^*, \ldots, m_d^*)$ and $\ell = \ell_1$.

Our assumptions in $j_2$ yield

$$m_2 + \gamma + 2x_2 + 2(1-m_2)m_1 + 2(1-m_1)m_d = m_2^* + \gamma_1 + 2x_2^* + 2(1-m_2^*)m_1^* + 2(1-m_1^*)m_d^*$$
$$\implies x_2 - x_2^* = \frac{\gamma_1 - \gamma}{2}.$$

Since $x_2 - x_2^* \in \mathbb{Z}$ and $|\gamma_1 - \gamma| \leq 1$ by construction, we see $\gamma = \gamma_1$ and so $x_2 = x_2^*$. Inspecting $j_k$ for $3 \leq k \leq d$, observe

$$p_{k-2} + A_{k+1}s_{k-2} + m_k + 2x_k = p_{k-2}^* + A_{k+1}s_{k-2}^* + m_k^* + 2x_k^*$$
$$\implies s_{k-2} = s_{k-2}^*.$$

Thus, for all $3 \leq k \leq d$, $s_{k-2} = s_{k-2}^*$ and so $C_1 = C_{\ell,\gamma,t,p_1,s_1,\ldots,p_{d-2},s_{d-2}} = C_{\ell_1,\gamma_1,t_1,p_1^*,s_1^*,\ldots,p_{d-2}^*,s_{d-2}^*} = C_2$.

**$\alpha^*$−Case 4:** Let $\alpha_2 < \alpha^* \leq \alpha_d$. Then, sets 1 and 3 only are active as $A_1 = 1$ and $A_2 = 1$. Note that $t = 0 = t_1$, $A_3 \geq 4$ is even, and that when $\alpha_{z-1} < \alpha^* \leq \alpha_z$, for all $z+1 \leq k \leq d$ with $3 \leq z \leq d$ it is the case $p_{k-2} = 0 = p_{k-2}^*$, $A_{k+1} = 1$ and $x_k = 0 = x_k^*$.

**Set Comparison Case 1:** We will consider when edges are shared from moves with the same orientation via set 1.

Observe that the argument required here is precisely that made in set comparison case 1 of $\alpha^*$−case 3 modulo the justification for why $p_{k-2} = p_{k-2}^*$, $x_k = x_k^*$, and $s_{k-2} = s_{k-2}^*$ when $\alpha_{z-1} < \alpha^* \leq \alpha_z$ for all $3 \leq k \leq z$ with $3 \leq z \leq d$.

Note that we thus have at this point $\ell = \ell_1$, $(m_1, \ldots, m_d) = (m_1^*, \ldots, m_d^*)$, $x_2 = x_2^*$, $r_2 = r_2^*$, and $\gamma = \gamma_1$ during all moves. So let $\alpha_{z-1} < \alpha^* \leq \alpha_z$ with $3 \leq z \leq d$. Then, for all $z+1 \leq k \leq d$, $p_{k-2} = 0 = p_{k-2}^*$ and $x_k = 0 = x_k^*$, which by the argument made in set comparison case 1 of the previous $\alpha^*$−case implies $s_{k-2} = s_{k-2}^*$ for each $z+1 \leq k \leq d$.

Now, for all components $j_k$ with $3 \leq k \leq z$, we see

$$p_{k-2} + 2A_{k+1}s_{k-2} + m_k + 2x_k = p_{k-2}^* + 2A_{k+1}s_{k-2}^* + m_k^* + 2x_k^*$$
$$\implies A_{k+1}(s_{k-2} - s_{k-2}^*) + (x_k - x_k^*) = \frac{p_{k-2}^* - p_{k-2}}{2}.$$

Following from $A_{k+1}(s_{k-2} - s^*_{k-2}) + (x_k - x^*_k) \in \mathbb{Z}$ and $|p^*_{k-2} - p_{k-2}| \leq 1$ by construction, we have $p_{k-2} = p^*_{k-2}$ for all $3 \leq k \leq z$, leaving us with

$$s_{k-2} - s^*_{k-2} = \frac{x^*_k - x_k}{A_{k+1}}.$$

Given $s_{k-2} - s^*_{k-2} \in \mathbb{Z}$ and $|x^*_k - x_k| \leq A_{k+1} - 1$ with $A_{k+1} \geq 2$ by construction, it follows that $x_k = x^*_k$ and $s_{k-2} = s^*_{k-2}$ for all $3 \leq k \leq z$. Thus, for all $\alpha_2 < \alpha^* \leq \alpha_d$, we conclude $C_1 = C_{\ell,\gamma,t,p_1,s_1,\ldots,p_{d-2},s_{d-2}} = C_{\ell_1,\gamma_1,t_1,p^*_1,s^*_1,\ldots,p^*_{d-2},s^*_{d-2}} = C_2$.

**Set Comparison Case 2:** We will consider when edges are shared from moves with the same orientation via set 3.

<u>Major Case 1</u>: Focusing on moves along $j_1$, we see $m_1 = m_3 = \cdots = m_d$, $m_2 = 1 - m_1$, $m^*_1 = m^*_3 = \cdots = m^*_d$ and $m^*_2 = 1 - m^*_1$ during these moves, giving us in $j_1$

$$(\ell + 1)A_3 + (-2 + m_1 + 2m_d(1 - m_1)) = (\ell_1 + 1)A_3 + (-2 + m^*_1 + 2m^*_d(1 - m^*_1))$$
$$\implies \ell - \ell_1 = \frac{m^*_1 - m_1}{A_3}.$$

Following from $\ell - \ell_1 \in \mathbb{Z}$ and $|m^*_1 - m_1| \leq 1$ with $A_3 \geq 4$ by construction, we have $m_1 = m^*_1$, meaning $(m_1, \ldots, m_d) = (m^*_1, \ldots, m^*_d)$, and so $\ell = \ell_1$.

Applying what we have so far to $j_2$, we get

$$m_2 + \gamma + 2x_2 + 2(1 - m_2)m_1 + 2((1 - m_1)m_d - (1 - R_1)((1 - m_2)m_1 + m_d(1 - m_1))$$
$$= m^*_2 + \gamma_1 + 2x^*_2 + 2(1 - m^*_2)m^*_1 + 2((1 - m^*_1)m^*_d - (1 - R^*_1)((1 - m^*_2)m^*_1 + m^*_d(1 - m^*_1))$$
$$\implies (x_2 - x^*_2) + (R_1 - R^*_1)((1 - m_2)m_1 + m_d(1 - m_1)) = \frac{\gamma_1 - \gamma}{2}.$$

Since $(x_2 - x^*_2) + (R_1 - R^*_1)((1 - m_2)m_1 + m_d(1 - m_1)) \in \mathbb{Z}$ and $|\gamma_1 - \gamma| \leq 1$ by construction, it must be the case $\gamma = \gamma_1$. This leaves us with

$$x_2 - x^*_2 = (R^*_1 - R_1)((1 - m_2)m_1 + m_d(1 - m_1)).$$

We now case on $R^*_1 - R_1$:

<u>Case 1:</u> Let $R_1 = R^*_1$. Then, our equality above implies $x_2 = x^*_2$, meaning $r_2 = r^*_2$.

<u>Case 2:</u> Let $R_1 \neq R^*_1$. Then, without loss of generality, let $R_1 = 0$ and $R^*_1 = 1$. From our last equality in $j_2$ prior to these cases, we now see that $x_2 - x^*_2 = m_1$. From this, it follows that either $(r_2 - r^*_2) - m_1 = 0$ or $(r_2 - r^*_2) - m_1 = -2$. In particular, we see that $m_1 = |r_2 - r^*_2|$ in both cases, a result we will use without further mention.

From our assumption $R_1 = 0$, we know that there exists a smallest integer $3 \leq w \leq d$ such that $x_w \neq (1 - r_2)(A_{w+1} - 1)$. Note that this implies $\alpha_{w-1} < \alpha^* \leq \alpha_d$. Similarly, our assumption $R^*_1 = 1$ implies $x^*_z = (1 - r^*_2)(A_{z+1} - 1)$ for all $3 \leq z \leq d$.

To complete our treatment of case 2, we will require the following result we will prove using induction:

We will prove that for every $3 \leq w \leq d$ such that $w$ is the smallest integer with $x_w \neq (1-r_2)(A_{w+1}-1)$, it is the case that for all $3 \leq z \leq w-1$, $x_z = (1-r_2)(A_{z+1}-1)$. Note that for all $3 \leq k \leq d$, $(r_2^* - r_2)m_k = r_2^* - r_2$ by our previous deduction in $j_2$. We proceed by induction on $w$:

<u>Base Case 1:</u> Let $w = 3$. Then, the statement is vacuously true as there does not exist an integer $z$ such that $3 \leq z \leq 2$.

<u>Base Case 2:</u> Let $w = 4$. Then, $z = 3$ and so our assumptions for this subcase give us in $j_3$ that

$$p_1 + 2A_4 s_1 + (1 - 2r_2\Psi_0)m_3 + 2x_3 + 2\eta_{\alpha_3}(1 - R_1)((1 - r_2)X_0 - r_2\Psi_0)(1 - m_3)m_d$$
$$= p_1^* + 2A_4 s_1^* + (1 - 2r_2^*\Psi_0^*)m_3^* + 2x_3^* + 2\eta_{\alpha_3}(1 - R_1^*)((1 - r_2^*)X_0^* - r_2^*\Psi_0^*)(1 - m_3^*)m_d^*$$
$$\implies A_4(s_1 - s_1^*) + (r_2^* - r_2)m_3 + (x_3 - x_3^*) + (1 - 2r_2)(1 - m_3)m_d = \frac{p_1^* - p_1}{2}.$$

Since $A_4(s_1 - s_1^*) + (r_2^* - r_2)m_3 + (x_3 - x_3^*) + (1 - 2r_2)(1 - m_3)m_d \in \mathbb{Z}$ and $|p_1^* - p_1| \leq 1$ by construction, it must be the case $p_1 = p_1^*$. Recalling that $x_3^* = (1 - r_2^*)(A_4 - 1)$ since $R_1^* = 1$, we are left with

$$(s_1 - s_1^*) - (1 - r_2^*) = -\frac{x_3 + (1 - r_2)}{A_4}.$$

Following from $(s_1 - s_1^*) - (1 - r_2^*) \in \mathbb{Z}$ and $1 - r_2 \leq x_3 + (1 - r_2) \leq A_4 - r_2$, we see that we must have $x_3 = (1 - r_2)(A_4 - 1)$.

<u>Induction Step:</u> Let $3 \leq w \leq d$ and $d \geq 3$ be such that $3 \leq w + 1 \leq d$ so that the following makes sense to consider. Otherwise, we will have shown that $x_z = (1-r_2)(A_{z+1}-1)$ for all $3 \leq z \leq d$ when $w = d$, which is not possible as $R_1 = 0$ by assumption. Further assume that for all $3 \leq z \leq y - 1$ with $3 \leq y \leq w$, $x_z = (1 - r_2)(A_{z+1} - 1)$. We will show that our statement holds for $3 \leq w + 1 \leq d$ by showing that $x_w = (1 - r_2)(A_{w+1} - 1)$ since our induction hypothesis grants us the above for all $3 \leq z \leq w - 1$. Note that our assumptions on $x_z$ and $x_{z^*}^*$, with $x_{z^*}^* = (1 - r_2^*)(A_{z^*+1} - 1)$ for all $3 \leq z^* \leq d$ since $R_1^* = 1$, imply $(1-r_2)X_{w-3} = (1-r_2)$, $r_2\Psi_{w-3} = r_2$, $(1-r_2^*)X_{w-3}^* = (1-r_2^*)$, and $r_2^*\Psi_{w-3}^* = r_2^*$. The corresponding assumptions and results we have so far, along with the observations $\eta_{w+1} = 1$ since $w + 1 \leq d$ and $\eta_{\alpha_w} = \eta_{\alpha_{w-1}} = 1$ since $\alpha^* > \alpha_w > \alpha_{w-1}$, yield in $j_w$

$$p_{w-2} + 2A_{w+1}s_{w-2} + (1 - 2r_2\Psi_{w-3})m_w + 2x_w + 2\eta_{w+1}\eta_{\alpha_w}(1 - R_1)((1-r_2)X_{w-3} - r_2\Psi_{w-3})(1 - m_w)m_d$$
$$= p_{w-2}^* + 2A_{w+1}s_{w-2}^* + (1 - 2r_2^*\Psi_{w-3}^*)m_w^* + 2x_w^* + 2\eta_{w+1}\eta_{\alpha_w}(1 - R_1^*)((1-r_2^*)X_{w-3}^* - r_2^*\Psi_{w-3}^*)(1 - m_w^*)m_d^*$$
$$\implies A_{w+1}(s_{w-2} - s_{w-2}^*) + (r_2^* - r_2)m_w + (x_w - x_w^*) + (1 - 2r_2)(1 - m_w)m_d = \frac{p_{w-2}^* - p_{w-2}}{2}.$$

Since $A_{w+1}(s_{w-2} - s_{w-2}^*) + (r_2^* - r_2)m_w + (x_w - x_w^*) + (1 - 2r_2)(1 - m_w)m_d \in \mathbb{Z}$ and $|p_{w-2}^* - p_{w-2}| \leq 1$, it must be the case $p_{w-2} = p_{w-2}^*$. Hence, observing that $x_w^* = (1 - r_2^*)(A_{w+1} - 1)$ since $R_1^* = 1$, we are left with

$$(s_{w-2} - s_{w-2}^*) - (1 - r_2^*) = -\frac{x_w + (1 - r_2)}{A_{w+1}}.$$

Given $(s_{w-2} - s_{w-2}^*) - (1 - r_2^*) \in \mathbb{Z}$ and $1 - r_2 \leq x_w + (1 - r_2) \leq A_{w+1} - r_2$ by construction, it follows that $x_w = (1 - r_2)(A_{w+1} - 1)$.

By the Principle of Strong Mathematical Induction, it follows that for every $3 \leq w \leq d$, it is the case that $x_z = (1-r_2)(A_{z+1}-1)$ for every $3 \leq z \leq w-1$.

From this, we know that $x_z = (1-r_2)(A_{z+1}-1)$ for all $3 \leq z \leq w-1$ with $3 \leq w \leq d$. Given the definitions of $j_d$ and $j_k$ for $3 \leq k \leq d-1$ differ slightly, we case on whether $w = d$:

<u>Subcase 1:</u> Let $w = d$. Then, $x_z = (1-r_2)(A_{z+1}-1)$ for all $3 \leq z \leq w-1$ and $x_d \neq (1-r_2)(A_{d+1}-1)$. Hence, noting that our assumptions on $x_z$ and $x_z^*$ for $3 \leq z \leq w-1$ imply $r_2\Psi_{w-3} = r_2$ and $r_2^*\Psi_{w-3}^* = r_2^*$, in $j_w$ we have

$$p_{w-2} + 2A_{w+1}s_{w-2} + (1-2r_2\Psi_{w-3})m_w + 2x_w = p_{w-2}^* + 2A_{w+1}s_{w-2}^* + (1-2r_2^*\Psi_{w-3}^*)m_w^* + 2x_w^*$$

$$\implies (s_{w-2} - s_{w-2}^*) - (1-r_2^*) = -\frac{x_w + (1-r_2)}{A_{w+1}}.$$

Given $(s_{w-2} - s_{w-2}^*) - (1-r_2^*) \in \mathbb{Z}$ and $x_w \neq (1-r_1)(A_{w+1}-1)$, it follows that $-\frac{x_w+(1-r_2)}{A_{w+1}} \notin \mathbb{Z}$, giving us a contradiction as the above equality is not possible. Thus, since $3 \leq w \leq d$ and $w = d$, we see there does not exist a smallest integer $w$ such that $x_w \neq (1-r_2)(A_{w+1}-1)$.

<u>Subcase 2:</u> Let $3 \leq w < d$. Then, $x_z = (1-r_2)(A_{z+1}-1)$ for all $3 \leq z \leq w-1$ and $x_w \neq (1-r_2)(A_{w+1}-1)$. Hence, noting that our assumptions on $x_z$ and $x_z^*$ for $3 \leq z \leq w-1$ imply $r_2\Psi_{w-3} = r_2$, $(1-r_2)X_{w-3} = (1-r_2)$, $r_2^*\Psi_{w-3}^* = r_2^*$ and $(1-r_2^*)X_{w-3}^* = (1-r_2^*)$, in $j_w$ we get

$$p_{w-2} + 2A_{w+1}s_{w-2} + (1-2r_2\Psi_{w-3})m_w + 2x_w + 2\eta_{\alpha_w}(1-R_1)((1-r_2)X_{w-3} - r_2\Psi_{w-3})(1-m_w)m_d$$
$$= p_{w-2}^* + 2A_{w+1}s_{w-2}^* + (1-2r_2^*\Psi_{w-3}^*)m_w^* + 2x_w^* + 2\eta_{\alpha_w}(1-R_1)((1-r_2^*)X_{w-3}^* - r_2^*\Psi_{w-3}^*)(1-m_w^*)m_d^*$$

$$\implies (s_{w-2} - s_{w-2}^*)A_{w+1} + (r_2^* - r_2)m_w + (x_w - x_w^*) = \frac{p_{w-2}^* - p_{w-2}}{2}.$$

Given $(s_{w-2} - s_{w-2}^*)A_{w+1} + (r_2^* - r_2)m_w + (x_w - x_w^*) \in \mathbb{Z}$ and $|p_{w-2}^* - p_{w-2}| \leq 1$ by construction, it follows that $p_{w-2} = p_{w-2}^*$. This leaves us with

$$(s_{w-2} - s_{w-2}^*) - (1-r_2^*) = -\frac{x_w + (1-r_2)}{A_{w+1}}.$$

Since $(s_{w-2} - s_{w-2}^*) - (1-r_2^*) \in \mathbb{Z}$ and $x_w \neq (1-r_2)(A_{w+1}-1)$, we have $-\frac{x_w+(1-r_2)}{A_{w+1}} \notin \mathbb{Z}$, giving us a contradiction as the above equality is not possible. The above contradicts our assumption that $3 \leq w \leq d-1$ was the smallest integer such that $x_w \neq (1-r_2)(A_{w+1}-1)$.

Further, since $w$ was the arbitrary smallest and all components $j_k$ for $3 \leq k \leq d-1$ are defined as above, the above holds for all such $w$. So there is no such $3 \leq w \leq d-1$. Now, if $w = d$, then by subcase 1 we conclude $w \neq d$ as we get a contradiction. Having exhausted all possible $3 \leq w \leq d$, we conclude that there does not exist a smallest integer $3 \leq w \leq d$ such that $x_w \neq (1-r_2)(A_{w+1}-1)$. This means that the scenario we have considered never occurs.

<u>Major Case 2:</u> Focusing on moves along $j_2$, we see that which $j_2$ moves in $C_1$ have the same orientation as those in $C_2$ depends on $R_1$ and $R_1^*$. Hence, we proceed by casing on $R_1$ and $R_1^*$:

<u>Case 1</u>: Let $R_1 = 0$. Then, stair-casing/normal $j_2$ moves are taking place in $C_1$. We now further case on $R_1^*$:

<u>Subcase 1</u>: Let $R_1^* = 0$. Then, stair-casing/normal $j_2$ moves are taking place in $C_2$. Hence, $(m_1, \ldots, m_d) = (m_1^*, \ldots, m_d^*)$ and so in $j_1$ we get

$$(\ell + 1)A_3 + (-2 + m_1 + 2m_d(1 - m_1)) = (\ell_1 + 1)A_3 + (-2 + m_1^* + 2m_d^*(1 - m_1^*))$$
$$\implies \ell = \ell_1.$$

So $\ell = \ell_1$ and $(m_1, \ldots, m_d) = (m_1^*, \ldots, m_d^*)$. Applying our results to $j_2$, we obtain

$$m_2 + \gamma + 2x_2 + 2(1 - m_2)m_1 + 2((1 - m_1)m_d - (1 - R_1)((1 - m_2)m_1 + m_d(1 - m_1)))$$
$$= m_2^* + \gamma_1 + 2x_2^* + 2(1 - m_2^*)m_1^* + 2((1 - m_1^*)m_d^* - (1 - R_1^*)((1 - m_2^*)m_1^* + m_d^*(1 - m_1^*)))$$
$$\implies x_2 - x_2^* = \frac{\gamma_1 - \gamma}{2}.$$

Following from $x_2 - x_2^* \in \mathbb{Z}$ and $|\gamma_1 - \gamma| \leq 1$ by construction, it must be the case that $\gamma = \gamma_1$ and so $x_2 = x_2^*$, meaning $r_2 = r_2^*$.

<u>Subcase 2</u>: Let $R_1^* = 1$. Then, $C_2$ is being defined by column transitions along $j_2$, and so only $C_1$'s $j_2$ move $(m_1, \ldots, m_d) = (0, \ldots, 0)$ has the same orientation as the $j_2$ moves in $C_2$. Noting that during these moves $m_1^* = \cdots = m_d^*$, in $j_1$ we find

$$(\ell + 1)A_3 + (-2 + m_1 + 2m_d(1 - m_1)) = (\ell_1 + 1)A_3 + (-2 + m_1^* + 2m_d^*(1 - m_1^*))$$
$$\implies \ell - \ell_1 = \frac{m_1^*}{A_3}.$$

Since $\ell - \ell_1 \in \mathbb{Z}$ and $m_1^* \in \{0, 1\}$ with $A_3 \geq 4$ by construction, it follows that $m_1^* = 0$, meaning $(m_1, \ldots, m_d) = (m_1^*, \ldots, m_d^*)$, and $\ell = \ell_1$. Applying what we have so far to $j_2$, we see

$$m_2 + \gamma + 2x_2 + 2(1 - m_2)m_1 + 2((1 - m_1)m_d - (1 - R_1)((1 - m_2)m_1 + m_d(1 - m_1)))$$
$$= m_2^* + \gamma_1 + 2x_2^* + 2(1 - m_2^*)m_1^* + 2((1 - m_1^*)m_d^* - (1 - R_1^*)((1 - m_2^*)m_1^* + m_d^*(1 - m_1^*)))$$
$$\implies x_2 - x_2^* = \frac{\gamma_1 - \gamma}{2}.$$

Given $x_2 - x_2^* \in \mathbb{Z}$ and $|\gamma_1 - \gamma| \leq 1$ by construction, it must be the case $\gamma = \gamma_1$ and so $x_2 = x_2^*$, meaning $r_2 = r_2^*$. Now since $R_1 = 0$, we know there exists a smallest $3 \leq w \leq d$ such that $x_w \neq (1 - r_2)(A_{w+1} - 1)$.

To complete our treatment of this subcase, we will require the following result we will prove by way of induction:

We will prove that for every $3 \leq w \leq d$ such that $w$ is the smallest integer with $x_w \neq (1 - r_2)(A_{w+1} - 1)$, it is the case that for all $3 \leq z \leq w - 1$, $x_z = (1 - r_2)(A_{z+1} - 1)$. We proceed by induction on $w$:

<u>Base Case 1</u>: Let $w = 3$. Then, the statement is vacuously true as there does not exist an integer $z$ such that $3 \leq z \leq 2$.

<u>Base Case 2</u>: Let $w = 4$. Then, $z = 3$ and so our assumptions for this subcase give us in $j_3$ that

$$p_1 + 2A_4 s_1 + (1 - 2r_2\Psi_0)m_3 + 2x_3 + 2\eta_{\alpha_3}(1 - R_1)((1 - r_2)X_0 - r_2\Psi_0)(1 - m_3)m_d$$
$$= p_1^* + 2A_4 s_1^* + (1 - 2r_2^*\Psi_0^*)m_3^* + 2x_3^* + 2\eta_{\alpha_3}(1 - R_1^*)((1 - r_2^*)X_0^* - r_2^*\Psi_0^*)(1 - m_3^*)m_d^*$$
$$\implies A_4(s_1 - s_1^*) + (r_2^* - r_2)m_3 + (x_3 - x_3^*) + (1 - 2r_2)(1 - m_3)m_d = \frac{p_1^* - p_1}{2}.$$

Since $A_4(s_1 - s_1^*) + (r_2^* - r_2)m_3 + (x_3 - x_3^*) + (1 - 2r_2)(1 - m_3)m_d \in \mathbb{Z}$ and $|p_1^* - p_1| \leq 1$ by construction, it must be the case $p_1 = p_1^*$. Recalling that $x_3^* = (1 - r_2^*)(A_4 - 1)$ since $R_1^* = 1$, we are left with

$$(s_1 - s_1^*) - (1 - r_2^*) = -\frac{x_3 + (1 - r_2)}{A_4}.$$

Following from $(s_1 - s_1^*) - (1 - r_2^*) \in \mathbb{Z}$ and $1 - r_2 \leq x_3 + (1 - r_2) \leq A_4 - r_2$, we see that we must have $x_3 = (1 - r_2)(A_4 - 1)$.

<u>Induction Step</u>: Let $3 \leq w \leq d$ and $d \geq 3$ be such that $3 \leq w + 1 \leq d$ so that the following makes sense to consider. Otherwise, we will have shown that $x_z = (1 - r_2)(A_{z+1} - 1)$ for all $3 \leq z \leq d$ when $w = d$, which is not possible as $R_1 = 0$ by assumption. Further assume that for all $3 \leq z \leq y - 1$ with $3 \leq y \leq w$, $x_z = (1 - r_2)(A_{z+1} - 1)$. We will show that our statement holds for $3 \leq w + 1 \leq d$ by showing that $x_w = (1 - r_2)(A_{w+1} - 1)$ since our induction hypothesis grants us the above for all $3 \leq z \leq w - 1$. Note that $r_2 = r_2^*$ and that our assumptions on $x_z$ and $x_{z^*}^*$, with $x_{z^*}^* = (1 - r_2^*)(A_{z^*+1} - 1)$ for all $3 \leq z^* \leq d$ since $R_1^* = 1$, imply $(1 - r_2)X_{w-3} = (1 - r_2)$, $r_2\Psi_{w-3} = r_2$, $(1 - r_2^*)X_{w-3}^* = (1 - r_2^*)$, and $r_2^*\Psi_{w-3}^* = r_2^*$. The corresponding assumptions and results we have so far, along with the observations $\eta_{w+1} = 1$ since $w + 1 \leq d$ and $\eta_{\alpha_w} = \eta_{\alpha_{w-1}} = 1$ since $\alpha^* > \alpha_w > \alpha_{w-1}$, yield in $j_w$

$$p_{w-2} + 2A_{w+1}s_{w-2} + (1 - 2r_2\Psi_{w-3})m_w + 2x_w + 2\eta_{w+1}\eta_{\alpha_w}(1 - R_1)((1 - r_2)X_{w-3} - r_2\Psi_{w-3})(1 - m_w)m_d$$
$$= p_{w-2}^* + 2A_{w+1}s_{w-2}^* + (1 - 2r_2^*\Psi_{w-3}^*)m_w^* + 2x_w^* + 2\eta_{w+1}\eta_{\alpha_w}(1 - R_1^*)((1 - r_2^*)X_{w-3}^* - r_2^*\Psi_{w-3}^*)(1 - m_w^*)m_d^*$$
$$\implies A_{w+1}(s_{w-2} - s_{w-2}^*) + (r_2^* - r_2)m_w + (x_w - x_w^*) + (1 - 2r_2)(1 - m_w)m_d = \frac{p_{w-2}^* - p_{w-2}}{2}.$$

Since $A_{w+1}(s_{w-2} - s_{w-2}^*) + (r_2^* - r_2)m_w + (x_w - x_w^*) + (1 - 2r_2)(1 - m_w)m_d \in \mathbb{Z}$ and $|p_{w-2}^* - p_{w-2}| \leq 1$, it must be the case $p_{w-2} = p_{w-2}^*$. Hence, observing that $x_w^* = (1 - r_2^*)(A_{w+1} - 1)$ since $R_1^* = 1$, we are left with

$$(s_{w-2} - s_{w-2}^*) - (1 - r_2^*) = -\frac{x_w + (1 - r_2)}{A_{w+1}}.$$

Given $(s_{w-2} - s_{w-2}^*) - (1 - r_2^*) \in \mathbb{Z}$ and $1 - r_2 \leq x_w + (1 - r_2) \leq A_{w+1} - r_2$ by construction, we must have $x_w = (1 - r_2)(A_{w+1} - 1)$.

By the Principle of Strong Mathematical Induction, it follows that for every $3 \leq w \leq d$, it is the case that $x_z = (1 - r_2)(A_{z+1} - 1)$ for every $3 \leq z \leq w - 1$.

From this, we know that $x_z = (1 - r_2)(A_{z+1} - 1)$ for all $3 \leq z \leq w - 1$. Given the definitions of $j_d$ and $j_k$ for $3 \leq k \leq d - 1$ differ slightly, we case on whether $w = d$:

<u>Sub-subcase 1:</u> Let $w = d$. Then, $x_d \neq (1-r_2)(A_{d+1}-1)$ and so in $j_w$ we see

$$p_{w-2} + 2A_{w+1}s_{w-2} + (1-2r_2\Psi_{w-3})m_w + 2x_w = p^*_{w-2} + 2A_{w+1}s^*_{w-2} + (1-2r^*_2\Psi^*_{w-3})m^*_w + 2x^*_w$$

$$\implies A_{w+1}(s_{w-2} - s^*_{w-2}) + (x_w - x^*_w) = \frac{p^*_{w-2} - p_{w-2}}{2}.$$

Since $A_{w+1}(s_{w-2} - s^*_{w-2}) + (x_w - x^*_w) \in \mathbb{Z}$ and $|p^*_{w-2} - p_{w-2}| \leq 1$ by construction, it must be the case $p_{w-2} = p^*_{w-2}$. This leaves us with

$$(s_{w-2} - s^*_{w-2}) - (1 - r^*_2) = -\frac{x_w + (1 - r_2)}{A_{w+1}}.$$

Given $(s_{w-2} - s^*_{w-2}) - (1 - r^*_2) \in \mathbb{Z}$ and $x_w \neq (1-r_2)(A_{w+1}-1)$, it follows that $-\frac{x_w + (1-r_2)}{A_{w+1}} \notin \mathbb{Z}$, giving us a contradiction as the above equality is not possible. Thus, since $3 \leq w \leq d$ and $w = d$, we see there does not exist a smallest integer $w$ such that $x_w \neq (1-r_2)(A_{w+1}-1)$.

<u>Sub-subcase 2:</u> Let $3 \leq w < d$. Then, $x_z = (1-r_2)(A_{z+1}-1)$ for all $3 \leq z \leq w-1$ and $x_w \neq (1-r_2)(A_{w+1}-1)$. Hence, noting that $r_2 = r^*_2$ and our assumptions on $x_z$ and $x^*_z$ for $3 \leq z \leq w-1$ imply $r_2\Psi_{w-3} = r_2$, $(1-r_2)X_{w-3} = (1-r_2)$, $r^*_2\Psi^*_{w-3} = r^*_2$ and $(1-r^*_2)X^*_{w-3} = (1-r^*_2)$, in $j_w$ we get

$$p_{w-2} + 2A_{w+1}s_{w-2} + (1-2r_2\Psi_{w-3})m_w + 2x_w + 2\eta_{\alpha_w}(1-R_1)((1-r_2)X_{w-3} - r_2\Psi_{w-3})(1-m_w)m_d$$
$$= p^*_{w-2} + 2A_{w+1}s^*_{w-2} + (1-2r^*_2\Psi^*_{w-3})m^*_w + 2x^*_w + 2\eta_{\alpha_w}(1-R^*_1)((1-r^*_2)X^*_{w-3} - r^*_2\Psi^*_{w-3})(1-m^*_w)m^*_d$$

$$\implies (s_{w-2} - s^*_{w-2})A_{w+1} + (r^*_2 - r_2)m_w + (x_w - x^*_w) = \frac{p^*_{w-2} - p_{w-2}}{2}.$$

Given $(s_{w-2} - s^*_{w-2})A_{w+1} + (r^*_2 - r_2)m_w + (x_w - x^*_w) \in \mathbb{Z}$ and $|p^*_{w-2} - p_{w-2}| \leq 1$ by construction, it follows that $p_{w-2} = p^*_{w-2}$. This leaves us with

$$(s_{w-2} - s^*_{w-2}) - (1 - r^*_2) = -\frac{x_w + (1-r_2)}{A_{w+1}}.$$

Since $(s_{w-2} - s^*_{w-2}) - (1-r^*_2) \in \mathbb{Z}$ and $x_w \neq (1-r_2)(A_{w+1}-1)$, we have $-\frac{x_w+(1-r_2)}{A_{w+1}} \notin \mathbb{Z}$, giving us a contradiction as the above equality is not possible. The above contradicts our assumption that $3 \leq w \leq d-1$ was the smallest integer such that $x_w \neq (1-r_2)(A_{w+1}-1)$.

Further, since $w$ was the arbitrary smallest and all components $j_k$ for $3 \leq k \leq d-1$ are defined as above, the above holds for all such $w$. So there is no such $3 \leq w \leq d-1$. Now, if $w = d$, then by sub-subcase 1 we conclude $w \neq d$ as we get a contradiction. Having exhausted all possible $3 \leq w \leq d$, we conclude that there does not exist a smallest integer $3 \leq w \leq d$ such that $x_w \neq (1-r_2)(A_{w+1}-1)$. This means that the scenario we have considered never occurs.

<u>Case 2:</u> Let $R_1 = 1$. Noting that the subcase with $R_1 = 1$ and $R^*_1 = 0$ is symmetric to that of subcase 2 of case 1, we see that we get a contradiction there by the same argument. So let $R^*_1 = 1$. Then, during these $j_2$ moves $m_1 = \cdots = m_d$ and $m^*_1 = \cdots = m^*_d$, giving us in $j_1$

$$(\ell + 1)A_3 + (-2 + m_1 + 2m_d(1 - m_1)) = (\ell_1 + 1)A_3 + (-2 + m^*_1 + 2m^*_d(1 - m^*_1))$$

$$\implies \ell - \ell_1 = \frac{m_1^* - m_1}{A_3}.$$

Given $\ell - \ell_1 \in \mathbb{Z}$ and $|m_1^* - m_1| \leq 1$ with $A_3 \geq 4$ by construction, it follows that $m_1 = m_1^*$, meaning $(m_1, \ldots, m_d) = (m_1^*, \ldots, m_d^*)$, and so $\ell = \ell_1$. Now, in $j_2$ we have

$$m_2 + \gamma + 2x_2 + 2(1-m_2)m_1 + 2((1-m_1)m_d - (1-R_1)((1-m_2)m_1 + m_d(1-m_1)))$$
$$= m_2^* + \gamma_1 + 2x_2^* + 2(1-m_2^*)m_1^* + 2((1-m_1^*)m_d^* - (1-R_1^*)((1-m_2^*)m_1^* + m_d^*(1-m_1^*)))$$
$$\implies x_2 - x_2^* = \frac{\gamma_1 - \gamma}{2}.$$

Noting that $x_2 - x_2^* \in \mathbb{Z}$ and $|\gamma_1 - \gamma| \leq 1$ by construction, we get $\gamma = \gamma_1$ and so $x_2 = x_2^*$, meaning $r_2 = r_2^*$.

<u>Major Case 3:</u> Focusing on moves along $j_k$ for $3 \leq k \leq d$, we proceed by casing on the three general states defined by set 3:

<u>Case 1:</u> Let $R_1 = 0$ and $R_{k-1} = 1$. This means $C_1$ is stair-casing, so $x_z = (1-r_2)(A_{z+1} - 1)$ for all $3 \leq z \leq k-1$ and that there exists at least one integer $k \leq y \leq d$ such that $x_y \neq (1-r_2)(A_{y+1} - 1)$. This last deduction then implies that there exists a smallest integer $k \leq w \leq d$ such that $x_w \neq (1-r_2)(A_{w+1} - 1)$. Since this is to hold for all $r_2 \in \{0, 1\}$, it must be the case $A_{w+1} \geq 2$, which is if and only if $\alpha_{w-1} < \alpha^* \leq \alpha_d$. Hence, for this case, assume $\alpha_{w-1} < \alpha^* \leq \alpha_d$.

<u>Subcase 1:</u> Let $R_1^* = 0$ and $R_{k-1}^* = 1$. Then, $C_2$ is stair-casing, meaning $x_z^* = (1-r_2^*)(A_{z+1} - 1)$ for all $3 \leq z \leq k-1$ and there exists at least one integer $k \leq y \leq d$ such that $x_y^* \neq (1-r_2^*)(A_{y+1} - 1)$. Since we require $C_1$ and $C_2$ to have $j_k$ moves with the same orientation, we see that this is if and only if $r_2 = r_2^*$.

Consequently, during these $j_k$ moves for $3 \leq k \leq d$, we see $m_1 = \cdots = m_{k-1}$, $m_k = \cdots = m_d = 1 - m_1$, $m_1^* = \cdots = m_{k-1}^*$, and $m_k^* = \cdots = m_d^* = 1 - m_1^*$. Applying this to $j_1$, we get

$$(\ell+1)A_3 + (-2 + m_1 + 2m_d(1-m_1)) = (\ell_1+1)A_3 + (-2 + m_1^* + 2m_d^*(1-m_1^*))$$
$$\implies \ell - \ell_1 = \frac{m_d^* - m_d}{A_3}.$$

Since $\ell - \ell_1 \in \mathbb{Z}$ and $|m_d^* - m_d| \leq 1$ with $A_3 \geq 4$ by construction, it follows that $m_d = m_d^*$, meaning $(m_1, \ldots, m_d) = (m_1^*, \ldots, m_d^*)$, and so $\ell = \ell_1$. Our results and assumptions in $j_2$ imply

$$m_2 + \gamma + 2x_2 + 2(1-m_2)m_1 + 2((1-m_1)m_d - (1-R_1)((1-m_2)m_1 + m_d(1-m_1)))$$
$$= m_2^* + \gamma_1 + 2x_2^* + 2(1-m_2^*)m_1^* + 2((1-m_1^*)m_d^* - (1-R_1^*)((1-m_2^*)m_1^* + m_d^*(1-m_1^*)))$$
$$\implies x_2 - x_2^* = \frac{\gamma_1 - \gamma}{2}.$$

Following from $x_2 - x_2^* \in \mathbb{Z}$ and $|\gamma_1 - \gamma| \leq 1$ by construction, it must be the case $\gamma = \gamma_1$ and so $x_2 = x_2^*$. The last deduction implies $r_2 = r_2^*$, which is consistent with our initial assumption.

<u>Subcase 2:</u> Since $\alpha_{w-1} < \alpha^* \leq \alpha_d$, let $R_1^* = 0$ and $R_{k-1}^* = 0$ for all $4 \leq k \leq d$. We exclude $j_3$ since it does not have any normal moves given $R_2 = 1$ always. So $C_2$ is being defined by normal moves, and during all $j_k$ moves here for $4 \leq k \leq d$ when $d \geq 4$, $m_1 = \cdots = m_{k-1}$, $m_k = \cdots = m_d = 1 - m_1$, $m_1^* = \cdots = m_{k-1}^*$ and $m_k^* = \cdots = m_d^* = 1 - m_1^*$. Applying this to $j_1$, we have

$$(\ell+1)A_3 + (-2 + m_1 + 2m_d(1-m_1)) = (\ell_1+1)A_3 + (-2 + m_1^* + 2m_d^*(1-m_1^*))$$

$$\implies \ell - \ell_1 = \frac{m_d^* - m_d}{A_3}.$$

Given $\ell - \ell_1 \in \mathbb{Z}$ and $|m_d^* - m_d| \leq 1$ with $A_3 \geq 4$ by construction, it follows that $m_d = m_d^*$, meaning $(m_1, \ldots, m_d) = (m_1^*, \ldots, m_d^*)$, and so $\ell = \ell_1$. Our results and assumptions in $j_2$ imply

$$m_2 + \gamma + 2x_2 + 2(1-m_2)m_1 + 2((1-m_1)m_d - (1-R_1)((1-m_2)m_1 + m_d(1-m_1)))$$
$$= m_2^* + \gamma_1 + 2x_2^* + 2(1-m_2^*)m_1^* + 2((1-m_1^*)m_d^* - (1-R_1^*)((1-m_2^*)m_1^* + m_d^*(1-m_1^*)))$$

$$\implies x_2 - x_2^* = \frac{\gamma_1 - \gamma}{2}.$$

Since $x_2 - x_2^* \in \mathbb{Z}$ and $|\gamma_1 - \gamma| \leq 1$ by construction, it follows that $\gamma = \gamma_1$ and so $x_2 = x_2^*$, meaning $r_2 = r_2^*$.

To complete our treatment of this subcase, we will require the following result which we will prove by induction on $w^*$:

We will show that for every smallest $3 \leq w^* \leq k-1$ such that $x_{w^*}^* \neq (1-r_2^*)(A_{w^*+1} - 1)$ during a given $j_k$ move for $4 \leq k \leq d$, it is the case that for every $3 \leq z \leq w^* - 1$, $x_z^* = (1-r_2^*)(A_{z+1} - 1)$. Note that $k$ and $d$ are such that what follows makes sense to consider.

Base Case 1: Let $w^* = 3$. Then, the statement holds vacuously as there does not exist an integer $z$ such that $3 \leq z \leq 2$.

Base Case 2: Let $w^* = 4$. Then, $z = 3$ and so by our results and assumptions from this subcase along with the observation $\eta_{\alpha_3} = 1$ since $\alpha^* > \alpha_{w-1} > \alpha_{w^*-1}$, in $j_3$ we obtain

$$p_1 + 2A_4 s_1 + (1 - 2r_2 \Psi_0) m_3 + 2x_3 + 2\eta_{\alpha_3}(1-R_1)((1-r_2)X_0 - r_2\Psi_0)(1-m_3)m_d$$
$$= p_1^* + 2A_4 s_1^* + (1 - 2r_2^* \Psi_0^*) m_3^* + 2x_3^* + 2\eta_{\alpha_3}(1-R_1^*)((1-r_2^*)X_0^* - r_2^*\Psi_0^*)(1-m_3^*)m_d^*$$

$$\implies A_4(s_1 - s_1^*) + (x_3 - x_3^*) = \frac{p_1^* - p_1}{2}.$$

Since $A_4(s_1 - s_1^*) + (x_3 - x_3^*) \in \mathbb{Z}$ and $|p_1^* - p_1| \leq 1$ by construction, it must be the case that $p_1 = p_1^*$. This leaves us with

$$s_1 - s_1^* = \frac{x_3^* - x_3}{A_4}.$$

Given $s_1 - s_1^* \in \mathbb{Z}$ and $|x_3^* - x_3| \leq A_4 - 1$ by construction, we must have $x_3^* = x_3 = (1-r_2)(A_4 - 1) = (1-r_2^*)(A_4 - 1)$.

Induction Step: Assume that $3 \leq w^* \leq k-1$ with $4 \leq k \leq d$ and $d \geq 4$ are such that $3 \leq w^* + 1 \leq d$ so that the following makes sense to consider. Otherwise, we will have shown that $x_z^* = (1-r_2^*)(A_{z+1}-1)$ for all $3 \leq z \leq k-1$ when $w^* = k-1$, which is not possible as $R_{k-1}^* = 0$ by assumption. Further, assume that for all $3 \leq z \leq y^* - 1$ with $3 \leq y^* \leq w^*$, $x_z^* = (1 - r_2^*)(A_{z+1} - 1)$. We will show that our statement holds for $3 \leq w^* + 1 \leq k-1$ by showing that $x_{w^*}^* = (1-r_2^*)(A_{w^*+1} - 1)$ since our induction hypothesis already grants us the above for all $3 \leq z \leq w^* - 1$. Note that $r_2 = r_2^*$, and that our assumptions on $x_z$ and $x_z^*$ for $3 \leq z \leq w^* - 1$ imply $(1-r_2)X_{w^*-3} = (1-r_2)$, $r_2 \Psi_{w^*-3} = r_2$, $(1-r_2^*)X_{w^*-3}^* = (1-r_2^*)$ and $r_2^*\Psi_{w^*-3}^* = r_2^*$. Applying the pertinent assumptions and results we have so far along with the observations $\eta_{w^*+1} = 1$ since $w^* + 1 \leq k - 1 < d$ and $\eta_{\alpha_{w^*}} = 1 = \eta_{\alpha_{w^*-1}}$ since $\alpha^* > \alpha_{w-1} > \alpha_{w^*-1}$, in $j_{w^*}$ we have

$$p_{w^*-2}+2A_{w^*+1}s_{w^*-2}+(1-2r_2\Psi_{w^*-3})m_{w^*}+2x_{w^*}+2\eta_{w^*+1}\eta_{\alpha_{w^*}}(1-R_1)((1-r_2)X_{w^*-3}-r_2\Psi_{w^*-3})(1-m_{w^*})m_d$$

$$= p_{w^*-2}^*+2A_{w^*+1}s_{w^*-2}^*+(1-2r_2^*\Psi_{w^*-3}^*)m_{w^*}^*+2x_{w^*}^*+2\eta_{w^*+1}\eta_{\alpha_{w^*}}(1-R_1^*)((1-r_2^*)X_{w^*-3}^*-r_2^*\Psi_{w^*-3}^*)(1-m_{w^*}^*)m_d^*$$

$$\implies A_{w^*+1}(s_{w^*-2} - s_{w^*-2}^*) + (x_{w^*} - x_{w^*}^*) = \frac{p_{w^*-2}^* - p_{w^*-2}}{2}.$$

Since $A_{w^*+1}(s_{w^*-2} - s_{w^*-2}^*) + (x_{w^*} - x_{w^*}^*) \in \mathbb{Z}$ and $|p_{w^*-2}^* - p_{w^*-2}| \leq 1$, it must be the case $p_{w^*-2} = p_{w^*-2}^*$. Note that $x_{w^*} = (1-r_2)(A_{w^*+1} - 1)$ since $R_{k-1} = 1$. Hence, this leaves us with

$$(s_{w^*-2} - s_{w^*-2}^*) + (1 - r_2) = \frac{x_{w^*}^* + (1 - r_2^*)}{A_{w^*+1}}.$$

Given $(s_{w^*-2} - s_{w^*-2}^*) + (1 - r_2) \in \mathbb{Z}$ and $1 - r_2^* \leq x_{w^*}^* + (1 - r_2^*) \leq A_{w^*+1} - r_2^*$ by construction, it follows that $x_{w^*}^* = (1 - r_2^*)(A_{w^*+1} - 1)$.

By the Principle of Strong Mathematical Induction, it follows that for every $3 \leq w^* \leq k - 1$, it is the case that $x_z^* = (1 - r_2^*)(A_{z+1} - 1)$ for every $3 \leq z \leq w^* - 1$.

Now, we see that during all $j_k$ moves for $4 \leq k \leq d$, it is the case that the smallest $3 \leq w^* \leq k - 1$ such that $x_{w^*}^* \neq (1 - r_2^*)(A_{w^*+1} - 1)$ satisfies $x_z^* = (1 - r_2)(A_{z+1} - 1)$ for every $3 \leq z \leq w^* - 1$. By the same argument made in the induction step, in $j_{w^*}$ we get

$$p_{w^*-2}+2A_{w^*+1}s_{w^*-2}+(1-2r_2\Psi_{w^*-3})m_{w^*}+2x_{w^*}+2\eta_{w^*+1}\eta_{\alpha_{w^*}}(1-R_1)((1-r_2)X_{w^*-3}-r_2\Psi_{w^*-3})(1-m_{w^*})m_d$$

$$= p_{w^*-2}^*+2A_{w^*+1}s_{w^*-2}^*+(1-2r_2^*\Psi_{w^*-3}^*)m_{w^*}^*+2x_{w^*}^*+2\eta_{w^*+1}\eta_{\alpha_{w^*}}(1-R_1^*)((1-r_2^*)X_{w^*-3}^*-r_2^*\Psi_{w^*-3}^*)(1-m_{w^*}^*)m_d^*$$

$$\implies A_{w^*+1}(s_{w^*-2} - s_{w^*-2}^*) + (x_{w^*} - x_{w^*}^*) = \frac{p_{w^*-2}^* - p_{w^*-2}}{2}.$$

Since $A_{w^*+1}(s_{w^*-2} - s_{w^*-2}^*) + (x_{w^*} - x_{w^*}^*) \in \mathbb{Z}$ and $|p_{w^*-2}^* - p_{w^*-2}| \leq 1$, it must be the case $p_{w^*-2} = p_{w^*-2}^*$. Note that $x_{w^*} = (1-r_2)(A_{w^*+1} - 1)$ since $R_{k-1} = 1$. Hence, this leaves us with

$$(s_{w^*-2} - s_{w^*-2}^*) + (1 - r_2) = \frac{x_{w^*}^* + (1 - r_2^*)}{A_{w^*+1}}.$$

Given $(s_{w^*-2} - s_{w^*-2}^*) + (1 - r_2) \in \mathbb{Z}$ and $x_{w^*}^* + (1 - r_2^*) \neq (1 - r_2^*)A_{w^*+1}$ by construction and our assumption for this subcase, it follows that $\frac{x_{w^*}^* + (1-r_2^*)}{A_{w^*+1}} \notin \mathbb{Z}$. Consequently, we have a contradiction as the above equality is not possible. Further, since $3 \leq w^* \leq k - 1$ is the arbitrary smallest integer assumed to satisfy $x_{w^*}^* \neq (1 - r_2^*)(A_{w^*+1} - 1)$, it follows that this holds for all such $w^*$, meaning this scenario does not ever occur as there does not exist such a $w^*$ under these assumptions.

<u>Subcase 3:</u> Let $R_1^* = 1$ and note that here we treat $j_k$ moves for all $3 \leq k \leq d$. Then, $C_2$ is being defined by column transitions, and so $x_z^* = (1 - r_2^*)(A_{z+1} - 1)$ for all $3 \leq z \leq d$. In particular, there exists a smallest integer $k \leq w \leq d$ such that $x_w \neq (1 - r_2)(A_{w+1} - 1)$. Observing that $m_1 = \cdots = m_{k-1}$, $m_k = \cdots = m_d = 1 - m_1$, $m_1^* = \cdots = m_{k-1}^*$ and $m_k^* = \cdots = m_d^* = 1 - m_1^*$, in $j_1$ we find

$$(\ell+1)A_3 + (-2 + m_1 + 2m_d(1-m_1)) = (\ell_1+1)A_3 + (-2 + m_1^* + 2m_d^*(1-m_1^*))$$
$$\implies \ell - \ell_1 = \frac{m_d^* - m_d}{A_3}.$$

Noting that $\ell - \ell_1 \in \mathbb{Z}$ and $|m_d^* - m_d| \leq 1$ with $A_3 \geq 4$ by construction, it follows that $m_d = m_d^*$, meaning $(m_1, \ldots, m_d) = (m_1^*, \ldots, m_d^*)$, and so $\ell = \ell_1$. Our results and assumptions in $j_2$ imply

$$m_2 + \gamma + 2x_2 + 2(1-m_2)m_1 + 2((1-m_1)m_d - (1-R_1)((1-m_2)m_1 + m_d(1-m_1)))$$
$$= m_2^* + \gamma_1 + 2x_2^* + 2(1-m_2^*)m_1^* + 2((1-m_1^*)m_d^* - (1-R_1^*)((1-m_2^*)m_1^* + m_d^*(1-m_1^*)))$$
$$\implies (x_2 - x_2^*) - m_d = \frac{\gamma_1 - \gamma}{2}.$$

Following from $(x_2 - x_2^*) - m_d \in \mathbb{Z}$ and $|\gamma_1 - \gamma| \leq 1$ by construction, it must be the case $\gamma = \gamma_1$. Noting that this implies $x_2 - x_2^* = m_d$ with $m_d \in \{0,1\}$, we see

$$\max\{r_2, r_2^*\} - \min\{r_2, r_2^*\} = m_d$$
$$\implies \left(\frac{r_1 + r_2^*}{2} + \frac{|r_2 - r_2^*|}{2}\right) - \left(\frac{r_2 + r_2^*}{2} - \frac{|r_2 - r_2^*|}{2}\right) = m_d$$
$$\implies |r_2 - r_2^*| = m_d.$$

To conclude our treatment of this subcase, we will need the following result that is analogous to that required in the previous subcase.

We wish to show that for every smallest integer $k \leq w \leq d$ such that $x_w \neq (1-r_2)(A_{w+1}-1)$, it is the case that $x_z = (1-r_2)(A_{z+1}-1)$ for all $k \leq z \leq w-1$ as we already have $x_z = (1-r_2)(A_{z+1}-1)$ for all $3 \leq z \leq k-1$ by our case assumption $R_1 = 0$ and $R_{k-1} = 1$. We proceed by induction on $w$:

<u>Base Case 1:</u> Let $w = k$. Then, the statement is vacuously true as there does not exist an integer $z$ such that $k \leq z \leq k-1$.

<u>Base Case 2:</u> Let $w = k+1$. Then, $z = k$ and so noting that $k \leq z < w$ and $\alpha^* > \alpha_{w-1} > \alpha_{z-1}$, it follows that $\eta_{\alpha_{z-1}} = 1 = \eta_{\alpha_k}$ and $\eta_{k+1} = 1$. Since for every $3 \leq y \leq k-1$, we know $x_y = (1-r_2)(A_{y+1}-1)$ and $x_y^* = (1-r_2^*)(A_{y+1}-1)$, observe that $r_2\Psi_{k-3} = r_2$ and $r_2^*\Psi_{k-3}^* = r_2^*$. Hence, in $j_k$ we get

$$p_{k-2} + 2A_{k+1}s_{k-2} + (1 - 2r_2\Psi_{k-3})m_k + 2x_k + 2\eta_{k+1}\eta_{\alpha_k}(1-R_1)((1-r_2)X_{k-3} - r_2\Psi_{k-3})(1-m_k)m_d$$
$$= p_{k-2}^* + 2A_{k+1}s_{k-2}^* + (1 - 2r_2^*\Psi_{k-3}^*)m_k^* + 2x_k^* + 2\eta_{k+1}\eta_{\alpha_k}(1-R_1^*)((1-r_2^*)X_{k-3}^* - r_2^*\Psi_{k-3}^*)(1-m_k^*)m_d^*$$
$$\implies A_{k+1}(s_{k-2} - s_{k-2}^*) + (r_2^* - r_2)m_k + (x_k - x_k^*) + (1-2r_2)(1-m_k)m_d = \frac{p_{k-2}^* - p_{k-2}}{2}.$$

Given $A_{k+1}(s_{k-2} - s_{k-2}^*) + (r_2^* - r_2)m_k + (x_k - x_k^*) + (1-2r_2)(1-m_k)m_d \in \mathbb{Z}$ and $|p_{k-2}^* - p_{k-2}| \leq 1$, it must be the case $p_{k-2} = p_{k-2}^*$. Hence, observing that $x_k^* = (1-r_2^*)(A_{k+1}-1)$ since $R_1^* = 1$, we are left with

$$(s_{k-2} - s_{k-2}^*) - (1-r_2^*) = -\frac{x_k + (1-r_2)}{A_{k+1}}.$$

Noting that $(s_{k-2} - s^*_{k-2}) - (1-r^*_2) \in \mathbb{Z}$ and $1 - r_2 \leq x_k + (1 - r_2) \leq A_{k+1} - r_2$ by construction, it follows that $x_k = (1 - r_2)(A_{k+1} - 1)$.

Induction Step: Let $k \leq w \leq d$ with $3 \leq k \leq d$ and $d \geq 3$ be such that $3 \leq w + 1 \leq d$ so that the following makes sense to consider. Otherwise, we will have shown that $x_z = (1 - r_2)(A_{z+1} - 1)$ for all $k \leq z \leq d$ when $w = d$, which is not possible as $R_1 = 0$ and $R_{k-1} = 1$ by assumption. Further, assume that for every $k \leq y \leq w$ and $k \leq z \leq y - 1$, $x_z = (1 - r_2)(A_{z+1} - 1)$. We will show our statement holds for $k \leq w + 1 \leq d$ by showing $x_w = (1 - r_2)(A_{w+1} - 1)$ as our induction hypothesis already grants us $x_z = (1 - r_2)(A_{z+1} - 1)$ for all $k \leq z \leq w - 1$. Hence, observing that $x^*_{z^*} = (1 - r^*_2)(A_{z^*+1} - 1)$ for all $3 \leq z^* \leq d$, $r_2 \Psi_{w-3} = r_2$, $r^*_2 \Psi^*_{w-3} = r^*_2$ and $\eta_{\alpha_{w-1}} = 1 = \eta_{\alpha_w}$ since $\alpha^* > \alpha_w > \alpha_{w-1}$, in $j_w$ we have

$$p_{w-2} + 2A_{w+1}s_{w-2} + (1 - 2r_2\Psi_{w-3})m_w + 2x_w + 2\eta_{w+1}\eta_{\alpha_w}(1 - R_1)((1 - r_2)X_{w-3} - r_2\Psi_{w-3})(1 - m_w)m_d$$
$$= p^*_{w-2} + 2A_{w+1}s^*_{w-2} + (1 - 2r^*_2\Psi^*_{w-3})m^*_w + 2x^*_w + 2\eta_{w+1}\eta_{\alpha_w}(1 - R^*_1)((1 - r^*_2)X^*_{w-3} - r^*_2\Psi^*_{w-3})(1 - m^*_w)m^*_d$$
$$\implies A_{w+1}(s_{w-2} - s^*_{w-2}) + (r^*_2 - r_2)m_w + (x_w - x^*_w) + (1 - 2r_2)(1 - m_w)m_d = \frac{p^*_{w-2} - p_{w-2}}{2}.$$

Since $A_{w+1}(s_{w-2} - s^*_{w-2}) + (r^*_2 - r_2)m_w + (x_w - x^*_w) + (1 - 2r_2)(1 - m_w)m_d \in \mathbb{Z}$ and $|p^*_{w-2} - p_{w-2}| \leq 1$, it must be the case $p_{w-2} = p^*_{w-2}$. Hence, observing that $x^*_w = (1 - r^*_2)(A_{w+1} - 1)$ since $R^*_1 = 1$, we are left with

$$(s_{w-2} - s^*_{w-2}) - (1 - r^*_2) = -\frac{x_w + (1 - r_2)}{A_{w+1}}.$$

Following from $(s_{w-2} - s^*_{w-2}) - (1 - r^*_2) \in \mathbb{Z}$ and $1 - r_2 \leq x_w + (1 - r_2) \leq A_{w+1} - r_2$ by construction, we conclude that we must have $x_w = (1 - r_2)(A_{w+1} - 1)$.

Thus, by the Principle of Strong Mathematical Induction, we see that for every smallest $k \leq w \leq d$ such that $x_w \neq (1 - r_2)(A_{w+1} - 1)$, it is the case that for every $k \leq z \leq w - 1$ with $k \leq w \leq d$, $x_z = (1 - r_2)(A_{z+1} - 1)$.

Now, we have that during all $j_k$ moves for $3 \leq k \leq d$, it is the case that the smallest $k \leq w \leq d$ such that $x_w \neq (1 - r_2)(A_{w+1} - 1)$ satisfies $x_z = (1 - r_2)(A_{z+1} - 1)$ for all $k \leq z \leq w - 1$ by the induction argument above. Further, $x_z = (1 - r_2)(A_{z+1} - 1)$ for all $3 \leq z \leq k - 1$ by our assumption $R_{k-1} = 1$. Hence, the same argument made in the induction step applies with $\alpha^* > \alpha_{w-1}$, $x^*_{z^*} = (1-r^*_2)(A_{z^*+1}-1)$ for all $3 \leq z^* \leq d$ since $R^*_1 = 1$, and $m_k = \cdots = m_w = \cdots = m_d = |r_2 - r^*_2|$. Given the slight difference in the component definition between $j_y$ for $k \leq y \leq d-1$ and $j_d$, we case on whether $w = d$:

Sub-subcase 1: Let $w = d$. Then, in this case, we obtain

$$p_{w-2} + 2A_{w+1}s_{w-2} + (1 - 2r_2\Psi_{w-3})m_w + 2x_w = p^*_{w-2} + 2A_{w+1}s^*_{w-2} + (1 - 2r^*_2\Psi^*_{w-3})m^*_w + 2x^*_w$$
$$\implies A_{w+1}(s_{w-2} - s^*_{w-2}) + (r^*_2 - r_2)m_w + (x_w - x^*_w) = \frac{p^*_{w-2} - p_{w-2}}{2}.$$

Following from $A_{w+1}(s_{w-2} - s^*_{w-2}) + (r^*_2 - r_2)m_w + (x_w - x^*_w) \in \mathbb{Z}$ and $|p^*_{w-2} - p_{w-2}| \leq 1$ by construction, we must have $p_{w-2} = p^*_{w-2}$. This leaves us with

$$(s_{w-2} - s^*_{w-2}) - (1 - r^*_2) = -\frac{x_w + (1 - r_2)}{A_{w+1}}.$$

Since $(s_{w-2} - s^*_{w-2}) - (1 - r^*_2) \in \mathbb{Z}$ and $x_w \neq (1 - r_2)(A_{w+1} - 1)$ by construction and our case assumption, $-\frac{x_w + (1-r_2)}{A_{w+1}} \notin \mathbb{Z}$, giving us a contradiction as the above equality is not possible. It now follows that this scenario is impossible as the smallest $w$ is also the largest it could be, in this case, and so there does not exist such a $w$ by the contradiction above.

Sub-subcase 2: Let $k \leq w < d$. Then, $x_z = (1 - r_2)(A_{z+1} - 1)$ for all $3 \leq z \leq w - 1$ and $x_w \neq (1 - r_2)(A_{w+1} - 1)$, meaning $\alpha_{w-1} < \alpha^* \leq \alpha_d$ and $\eta_{\alpha_{w-1}} = 1$ by the same reasoning as before. Since $R^*_1 = 1$, we know $x^*_{z^*} = (1 - r^*_2)(A_{z^*+1} - 1)$ for all $3 \leq z^* \leq d$. Now, observe that our assumptions on $x_z$ and $x^*_z$ for $3 \leq z \leq w - 1$ imply $(1 - r_2)X_{w-3} = (1 - r_2)$, $r_2\Psi_{w-3} = r_2$, $(1 - r^*_2)X^*_{w-3} = (1 - r^*_2)$, and $r^*_2\Psi^*_{w-3} = r^*_2$. Recalling that $m_w = |r_2 - r^*_2|$, $R_1 = 0$ and $R^*_1 = 1$, all of the above in $j_w$ yield

$$p_{w-2} + 2A_{w+1}s_{w-2} + (1 - 2r_2\Psi_{w-3})m_w + 2x_w + 2\eta_{w+1}\eta_{\alpha_w}(1 - R_1)((1 - r_2)X_{w-3} - r_2\Psi_{w-3})(1 - m_w)m_d$$
$$= p^*_{w-2} + 2A_{w+1}s^*_{w-2} + (1 - 2r^*_2\Psi^*_{w-3})m^*_w + 2x^*_w + 2\eta_{w+1}\eta_{\alpha_w}(1 - R^*_1)((1 - r^*_2)X^*_{w-3} - r^*_2\Psi^*_{w-3})(1 - m^*_w)m^*_d$$

$$\implies A_{w+1}(s_{w-2} - s^*_{w-2}) + (r^*_2 - r_2) + (x_w - x^*_w) = \frac{p^*_{w-2} - p_{w-2}}{2}.$$

Since $A_{w+1}(s_{w-2} - s^*_{w-2}) + (r^*_2 - r_2) + (x_w - x^*_w) \in \mathbb{Z}$ and $|p^*_{w-2} - p_{w-2}| \leq 1$ by construction, we must have $p_{w-2} = p^*_{w-2}$. So we now get

$$(s_{w-2} - s^*_{w-2}) - (1 - r^*_2) = -\frac{x_w + (1 - r_2)}{A_{w+1}}.$$

Following from $(s_{w-2} - s^*_{w-2}) - (1 - r^*_2) \in \mathbb{Z}$ by construction and $x_w + (1 - r_2) \neq (1 - r_2)A_{w+1}$ by assumption, we have $-\frac{x_w + (1-r_2)}{A_{w+1}} \notin \mathbb{Z}$, giving us a contradiction as the above equality is not possible. The above contradicts our assumption that $k \leq w \leq d - 1$ was the smallest integer such that $x_w \neq (1 - r_2)(A_{w+1} - 1)$.

Further, since $w$ was the arbitrary smallest and all components $j_k$ for $3 \leq k \leq d - 1$ are defined as above, the above holds for all such $w$. So there is no such $k \leq w \leq d - 1$. Now, if $w = d$, then by sub-subcase 1 we conclude $w \neq d$ as we get a contradiction. Having exhausted all possible $k \leq w \leq d$, we conclude that there does not exist a smallest integer $k \leq w \leq d$ such that $x_w \neq (1 - r_2)(A_{w+1} - 1)$. This means that the scenario we have considered never occurs.

Case 2: Focusing on moves along $j_k$ for $4 \leq k \leq d$ when $d \geq 4$ since $j_3$ is not subject to normal moves, either $\alpha_2 < \alpha^* \leq \alpha_{k-1}$ or let $\alpha_{k-1} < \alpha^* \leq \alpha_d$ with $R_1 = 0$ and $R_{k-1} = 0$. Then, $C_1$ is being defined by normal moves and so, in the case $\alpha_{k-1} < \alpha^* \leq \alpha_d$, there exists a smallest integer $3 \leq w \leq k - 1$ such that $x_w \neq (1 - r_2)(A_{w+1} - 1)$.

Subcase 1: Assume that $\alpha_2 < \alpha^* \leq \alpha_{k-1}$ if this is already being assumed and otherwise assume $\alpha_{k-1} < \alpha^* \leq \alpha_d$, $R^*_1 = 0$ and $R^*_{k-1} = 0$. Now, observing that during all such moves $m_1 = \cdots = m_{k-1}$, $m_k = \cdots = m_d = 1 - m_1$ and $(m^*_1, \ldots, m^*_d) = (m_1, \ldots, m_d)$, in $j_1$ we have

$$(\ell + 1)A_3 + (-2 + m_1 + 2m_d(1 - m_1)) = (\ell_1 + 1)A_3 + (-2 + m^*_1 + 2m^*_d(1 - m^*_1))$$
$$\implies \ell = \ell_1.$$

Hence, we know $\ell = \ell_1$ and so our assumptions in $j_2$ yield

$$m_2 + \gamma + 2x_2 + 2(1-m_2)m_1 + 2((1-m_1)m_d - (1-R_1)((1-m_2)m_1 + m_d(1-m_1)))$$
$$= m_2^* + \gamma_1 + 2x_2^* + 2(1-m_2^*)m_1^* + 2((1-m_1^*)m_d^* - (1-R_1^*)((1-m_2^*)m_1^* + m_d^*(1-m_1^*)))$$
$$\implies (x_2 - x_2^*) + (R_1 - R_1^*)m_d = \frac{\gamma_1 - \gamma}{2}.$$

Since $(x_2 - x_2^*) + (R_1 - R_1^*)m_d \in \mathbb{Z}$ and $|\gamma_1 - \gamma| \leq 1$ by construction, it follows that $\gamma = \gamma_1$ and so $x_2 - x_2^* = (R_1^* - R_1)m_d$. If $R_1 = R_1^*$, then $x_2 = x_2^*$, meaning $r_2 = r_2^*$.

To treat the first of the following sub-subcases, without loss of generality let $R_1 = 0 \neq 1 = R_1^*$. Since $R_1^* = 1$, note that this implies $x_{z^*}^* = (1-r_2^*)(A_{z^*+1} - 1)$ for all $3 \leq z^* \leq d$. Then, $x_2 - x_2^* = m_d$, meaning $|r_2 - r_2^*| = m_d$. We now case on $\alpha^*$:

<u>Sub-subcase 1:</u> Let $\alpha_2 < \alpha^* \leq \alpha_{k-1}$. Then, for all $k \leq z \leq d$, $p_{z-2} = 0 = p_{z-2}^*$, $x_z = 0 = x_z^*$, $A_{z+1} = 1$, and consequently in $j_z$ we find

$$p_{z-2} + A_{z+1}s_{z-2} + m_z + 2x_z = p_{z-2}^* + A_{z+1}s_{z-2}^* + m_z^* + 2x_z^*$$
$$\implies s_{z-2} = s_{z-2}^*.$$

Thus, for all $k \leq z \leq d$, $s_{z-2} = s_{z-2}^*$. Now, for $3 \leq y \leq k-1$, let $\alpha_{y-1} < \alpha^* \leq \alpha_y$. Then, for all $y+1 \leq z \leq k-1$, the argument above implies $p_{z-2} = 0 = p_{z-2}^*$, $x_z = 0 = x_z^*$, $A_{z+1} = 1$, and $s_{z-2} = s_{z-2}^*$.

Since $R_1 = 0$ and we know $x_z = 0 = x_z^*$ for all $y+1 \leq z \leq d$ since $\alpha_{y-1} < \alpha^* \leq \alpha_y$, there exists a smallest integer $3 \leq w \leq y$ such that $x_w \neq (1-r_2)(A_{w+1} - 1)$. We will show via induction on $w$ that it must then be the case that for all $3 \leq z \leq w-1$, $x_z = (1-r_2)(A_{z+1} - 1)$:

<u>Base Case 1:</u> Let $w = 3$. Then, the statement holds vacuously as there does not exist an integer $z$ such that $3 \leq z \leq 2$.

<u>Base Case 2:</u> Let $w = 4$. Then, $z = 3$ and so our subcase results and assumptions in $j_3$ yield

$$p_1 + 2A_4 s_1 + (1 - 2r_2\Psi_0)m_3 + 2x_3 + 2\eta_{\alpha_3}(1-R_1)((1-r_2)X_0 - r_2\Psi_0)(1-m_3)m_d$$
$$= p_1^* + 2A_4 s_1^* + (1 - 2r_2^*\Psi_0^*)m_3^* + 2x_3^* + 2\eta_{\alpha_3}(1-R_1^*)((1-r_2^*)X_0^* - r_2^*\Psi_0^*)(1-m_3^*)m_d^*$$
$$\implies A_4(s_1 - s_1^*) + (r_2^* - r_2)m_3 + (x_3 - x_3^*) + (1 - 2r_2)(1 - m_3)m_d = \frac{p_1^* - p_1}{2}.$$

Since $A_4(s_1 - s_1^*) + (r_2^* - r_2)m_3 + (x_3 - x_3^*) + (1 - 2r_2)(1 - m_3)m_d \in \mathbb{Z}$ and $|p_1^* - p_1| \leq 1$ by construction, it must be that $p_1 = p_1^*$. This leaves us with

$$(s_1 - s_1^*) - (1 - r_2^*) = -\frac{x_3 + (1 - r_2)}{A_4}.$$

Given $(s_1 - s_1^*) - (1 - r_2^*) \in \mathbb{Z}$ and $1 - r_2 \leq x_3 + (1-r_2) \leq A_4 - r_2$ by construction, it follows that $x_3 = (1-r_2)(A_4 - 1)$.

<u>Induction Step:</u> Let $3 \leq w \leq y$ with $3 \leq y \leq k-1$, $4 \leq k \leq d$ and $d \geq 4$ be such that $3 \leq w+1 \leq y$ so that the following makes sense to consider. Otherwise, we will have shown that $x_z = (1-r_2)(A_{z+1}-1)$ for all $3 \leq z \leq y$ when $w = y$, which is not possible as $R_1 = 0$ and $x_z = (1-r_2)(A_{z+1} - 1) = 0$ since $A_{z+1} = 1$ for all $y+1 \leq z \leq d$ by assumption. Further, assume that for every $3 \leq y^* \leq w$ and $3 \leq z \leq y^*-1$, $x_z = (1-r_2)(A_{z+1}-1)$. We will show our statement holds for $3 \leq w+1 \leq y$ by showing $x_w = (1-r_2)(A_{w+1} - 1)$ as our induction hypothesis already grants us $x_z = (1-r_2)(A_{z+1} - 1)$ for all $3 \leq z \leq w-1$. Hence, observing that $x^*_{z^*} = (1-r^*_2)(A_{z^*+1} - 1)$ for all $3 \leq z^* \leq d$, $r_2\Psi_{w-3} = r_2$, $(1-r_2)X_{w-3} = (1-r_2)$, $r^*_2\Psi^*_{w-3} = r^*_2$, $(1-r^*_2)X^*_{w-3} = (1-r^*_2)$ and $\eta_{\alpha_{w-1}} = 1 = \eta_{\alpha_w}$ since $\alpha^* > \alpha_w > \alpha_{w-1}$, in $j_w$ we have

$$p_{w-2} + 2A_{w+1}s_{w-2} + (1-2r_2\Psi_{w-3})m_w + 2x_w + 2\eta_{w+1}\eta_{\alpha_w}(1-R_1)((1-r_2)X_{w-3} - r_2\Psi_{w-3})(1-m_w)m_d$$
$$= p^*_{w-2} + 2A_{w+1}s^*_{w-2} + (1-2r^*_2\Psi^*_{w-3})m^*_w + 2x^*_w + 2\eta_{w+1}\eta_{\alpha_w}(1-R^*_1)((1-r^*_2)X^*_{w-3} - r^*_2\Psi^*_{w-3})(1-m^*_w)m^*_d$$
$$\implies A_{w+1}(s_{w-2} - s^*_{w-2}) + (r^*_2 - r_2)m_w + (x_w - x^*_w) + (1-2r_2)(1-m_w)m_d = \frac{p^*_{w-2} - p_{w-2}}{2}.$$

Since $A_{w+1}(s_{w-2} - s^*_{w-2}) + (r^*_2 - r_2)m_w + (x_w - x^*_w) + (1-2r_2)(1-m_w)m_d \in \mathbb{Z}$ and $|p^*_{w-2} - p_{w-2}| \leq 1$, it must be the case $p_{w-2} = p^*_{w-2}$. Consequently, recalling that $x^*_w = (1-r^*_2)(A_{w+1} - 1)$ since $R^*_1 = 1$, we are left with

$$(s_{w-2} - s^*_{w-2}) - (1 - r^*_2) = -\frac{x_w + (1-r_2)}{A_{w+1}}.$$

Following from $(s_{w-2} - s^*_{w-2}) - (1 - r^*_2) \in \mathbb{Z}$ and $1 - r_2 \leq x_w + (1-r_2) \leq A_{w+1} - r_2$ by construction, we have that $x_w = (1-r_2)(A_{w+1} - 1)$.

Thus, by the Principle of Strong Mathematical Induction, we have that for every smallest $3 \leq w \leq y$ such that $x_w \neq (1-r_2)(A_{w+1}-1)$, it is the case that for every $3 \leq z \leq w-1$, $x_z = (1-r_2)(A_{z+1}-1)$.

Now, given $3 \leq w \leq y$ and $x_z = (1-r_2)(A_{z+1}-1)$ for all $3 \leq z \leq w-1$, we case on whether $w = y$ as the definitions of $j_y$ and $j_{y^*}$ for $3 \leq y^* < y$ differ slightly:

$w$−Case 1: Let $w = y$. Then, $\eta_{\alpha_w} = 0$ and $\eta_{\alpha_{w-1}} = 1$ since $\alpha_{y-1} < \alpha^* \leq \alpha_y$ and so in $j_w$ we find

$$p_{w-2} + 2A_{w+1}s_{w-2} + (1 - 2r_2\Psi_{w-3})m_w + 2x_w = p^*_{w-2} + 2A_{w+1}s^*_{w-2} + (1 - 2r^*_2\Psi^*_{w-3})m^*_w + 2x^*_w$$
$$\implies A_{w+1}(s_{w-2} - s^*_{w-2}) + (r^*_2 - r_2)m_w + (x_w - x^*_w) = \frac{p^*_{w-2} - p_{w-2}}{2}.$$

Since $A_{w+1}(s_{w-2} - s^*_{w-2}) + (r^*_2 - r_2)m_w + (x_w - x^*_w) \in \mathbb{Z}$ and $|p^*_{w-2} - p_{w-2}| \leq 1$, it must be the case $p_{w-2} = p^*_{w-2}$. Hence, observing that $x^*_w = (1-r^*_2)(A_{w+1} - 1)$ since $R^*_1 = 1$, we are left with

$$(s_{w-2} - s^*_{w-2}) - (1 - r^*_2) = -\frac{x_w + (1-r_2)}{A_{w+1}}.$$

Given $(s_{w-2} - s^*_{w-2}) - (1 - r^*_2) \in \mathbb{Z}$ and $x_w \neq (1-r_2)(A_{w+1}-1)$ by construction and our case assumption, it follows that $-\frac{x_w + (1-r_2)}{A_{w+1}} \notin \mathbb{Z}$, giving us a contradiction as the above equality is not possible. Since $3 \leq w \leq y$ and $w = y$, there does not exist $3 \leq w \leq y$ such that $x_w \neq (1-r_2)(A_{w+1} - 1)$. Hence, we conclude that the scenario we have considered does not occur.

$w$−Case 2: Let $3 \leq w < y$. Then, $\eta_{\alpha_w} = 1 = \eta_{\alpha_{w-1}}$ since $\alpha^* > \alpha_{y-1} \geq \alpha_w > \alpha_{w-1}$ and so by an analogous argument as that made in the induction step, in $j_w$ we have

$$p_{w-2}+2A_{w+1}s_{w-2}+(1-2r_2\Psi_{w-3})m_w+2x_w+2\eta_{w+1}\eta_{\alpha_w}(1-R_1)((1-r_2)X_{w-3}-r_2\Psi_{w-3})(1-m_w)m_d$$
$$= p^*_{w-2}+2A_{w+1}s^*_{w-2}+(1-2r^*_2\Psi^*_{w-3})m^*_w+2x^*_w+2\eta_{w+1}\eta_{\alpha_w}(1-R^*_1)((1-r^*_2)X^*_{w-3}-r^*_2\Psi^*_{w-3})(1-m^*_w)m^*_d$$
$$\implies (s_{w-2}-s^*_{w-2})A_{w+1}+(r^*_2-r_2)+(x_w-x^*_w) = \frac{p^*_{w-2}-p_{w-2}}{2}.$$

Since $(s_{w-2}-s^*_{w-2})A_{w+1}+(r^*_2-r_2)+(x_w-x^*_w) \in \mathbb{Z}$ and $|p^*_{w-2}-p_{w-2}| \leq 1$ by construction, we must have $p_{w-2} = p^*_{w-2}$. Noting that $\eta_{w+1} = 1$ since $w < y$, we get

$$(s_{w-2}-s^*_{w-2})-(1-r^*_2) = -\frac{x_w+(1-r_2)}{A_{w+1}}.$$

Following from $(s_{w-2}-s^*_{w-2})-(1-r^*_2) \in \mathbb{Z}$ by construction and $x_w+(1-r_2) \neq (1-r_2)A_{w+1}$ by assumption, we have $-\frac{x_w+(1-r_2)}{A_{w+1}} \notin \mathbb{Z}$, giving us a contradiction as the above equality is not possible. The above contradicts our assumption that $3 \leq w < y$ was the smallest integer such that $x_w \neq (1-r_2)(A_{w+1}-1)$.

Further, since $w$ was the arbitrary smallest and all components $j_z$ for $3 \leq z < y$ are defined as above, the above holds for all $3 \leq w < y$ and so there is no such $w$. Now, if $w = y$, then by sub-subcase 1 we conclude $w \neq y$ as we get a contradiction. Having exhausted all possible $3 \leq w \leq y$, we conclude that there does not exist a smallest integer $3 \leq w \leq y$ such that $x_w \neq (1-r_2)(A_{w+1}-1)$. This means that the scenario we have considered never occurs.

<u>Sub-subcase 2:</u> Let $\alpha_{k-1} < \alpha^* \leq \alpha_d$, $R^*_1 = 0$, and $R^*_{k-1} = 0$. Since $\alpha_{k-1} < \alpha^* \leq \alpha_d$, note that $R_1 = 0$ and $R_{k-1} = 0$. Since $R_1 = R^*_1$, our last equality in $j_2$ implies $x_2 = x^*_2$, meaning $r_2 = r^*_2$.

<u>Subcase 2:</u> Let $R^*_1 = 1$, meaning $C_2$ is being defined by column transitions and for all $3 \leq z \leq d$, $x^*_z = (1-r^*_2)(A_{z+1}-1)$. By how we have defined column transitions to be dependent on $r^*_2$ to dictate the orientation of these moves, we must have $\alpha_{k-1} < \alpha^* \leq \alpha_d$. Hence, let $R_1 = 0$ and $R_{k-1} = 0$. Then, there exists a smallest integer $3 \leq w \leq k-1$ such that $x_w \neq (1-r_2)(A_{w+1}-1)$. Noting that $m_1 = \cdots = m_{k-1}$, $m_k = \cdots = m_d = 1-m_1$, $m^*_1 = \cdots = m^*_{k-1}$ and $m^*_k = \cdots = m^*_d = 1-m^*_1 = (1-r^*_2)m_d + r^*_2(1-m_d)$, in $j_1$ we get

$$(\ell+1)A_3 + (-2+m_1+2m_d(1-m_1)) = (\ell_1+1)A_3 + (-2+m^*_1+2m^*_d(1-m^*_1))$$
$$\implies \ell - \ell_1 = \frac{m^*_d - m_d}{A_3}.$$

Following from $\ell - \ell_1 \in \mathbb{Z}$ and $|m^*_d - m_d| \leq 1$ with $A_3 \geq 4$ by construction, we see that $m_d = m^*_d$. This implies

$$m^*_d - m_d = 0$$
$$\implies (-1)^{m_d}r^*_2 = 0.$$

The above gives us $r^*_2 = 0$. So we get $(m^*_1, \ldots, m^*_d) = (m_1, \ldots, m_d)$. Applying this to $j_2$, we find

$$m_2 + \gamma + 2x_2 + 2m_1(1-m_2) + 2((1-m_1)m_d - (1-R_1)((1-m_2)m_1 + m_d(1-m_1)))$$
$$= m^*_2 + \gamma_1 + 2x^*_2 + 2m^*_1(1-m^*_2) + 2((1-m^*_1)m^*_d - (1-R^*_1)((1-m^*_2)m^*_1 + m^*_d(1-m^*_1)))$$
$$\implies (x_2 - x^*_2) - m_d = \frac{\gamma_1 - \gamma}{2}.$$

Given $(x_2-x_2^*)-m_d \in \mathbb{Z}$ and $|\gamma_1-\gamma| \leq 1$ by construction, it follows that $\gamma = \gamma_1$ and so $x_2-x_2^* = m_d$. The last deduction implies $m_d = r_2$.

We will now show that for every smallest $3 \leq w \leq k-1$, it is the case that for all $3 \leq z \leq w-1$, $x_z = (1-r_2)(A_{z+1}-1)$. Note that $\eta_{\alpha_z} = 1$ since $\alpha^* > \alpha_{k-1} \geq \alpha_w > \alpha_{w-1} \geq \alpha_z$. We proceed by induction on $w$:

<u>Base Case 1</u>: Let $w = 3$. Then, the statement holds true vacuously as there does not exist an integer $z$ such that $3 \leq z \leq 2$.

<u>Base Case 2</u>: Let $w = 4$. Then, $z = 3$ and so in $j_3$ we have

$$p_1 + 2A_4 s_1 + (1-2r_2\Psi_0)m_3 + 2x_3 + 2\eta_{\alpha_3}(1-R_1)((1-r_2)X_0 - r_2\Psi_0)(1-m_3)m_d$$
$$= p_1^* + 2A_4 s_1^* + (1-2r_2^*\Psi_0^*)m_3^* + 2x_3^* + 2\eta_{\alpha_3}(1-R_1^*)((1-r_2^*)X_0^* - r_2^*\Psi_0^*)(1-m_3^*)m_d^*$$
$$\implies A_4(s_1 - s_1^*) - r_2 m_3 + (x_3 - x_3^*) + (1-2r_2)(1-m_3)m_d = \frac{p_1^* - p_1}{2}.$$

Since $A_4(s_1 - s_1^*) - r_2 m_3 + (x_3 - x_3^*) + (1-2r_2)(1-m_3)m_d \in \mathbb{Z}$ and $|p_1^* - p_1| \leq 1$ by construction, it must be that $p_1 = p_1^*$. This leaves us with

$$(s_1 - s_1^*) - (1 - r_2^*) = -\frac{x_3 + (1-r_2)}{A_4}.$$

Given $(s_1 - s_1^*) - (1 - r_2^*) \in \mathbb{Z}$ and $1-r_2 \leq x_3 + (1-r_2) \leq A_4 - r_2$ by construction, it follows that $x_3 = (1-r_2)(A_4-1)$.

<u>Induction Step</u>: Let $3 \leq w \leq k-1$ with $4 \leq k \leq d$ and $d \geq 4$ be such that $3 \leq w+1 \leq k-1$ so that the following makes sense to consider. Otherwise, we will have shown that $x_z = (1-r_2)(A_{z+1}-1)$ for all $3 \leq z \leq k-1$ when $w = k-1$, which is not possible as $R_{k-1} = 0$ by assumption. Further, assume that for every $3 \leq y^* \leq w$ and $3 \leq z \leq y^* - 1$, $x_z = (1-r_2)(A_{z+1}-1)$. We will show our statement holds for $3 \leq w+1 \leq k-1$ by showing $x_w = (1-r_2)(A_{w+1}-1)$ as our induction hypothesis already grants us $x_z = (1-r_2)(A_{z+1}-1)$ for all $3 \leq z \leq w-1$. Hence, observing that $x_{z^*}^* = (1-r_2^*)(A_{z^*+1}-1)$ for all $3 \leq z^* \leq d$, $r_2\Psi_{w-3} = r_2$, $(1-r_2)X_{w-3} = (1-r_2)$, $r_2^* = 0$, $m_d = r_2$ and $\eta_{\alpha_{w-1}} = 1 = \eta_{\alpha_w}$ since $\alpha^* > \alpha_w > \alpha_{w-1}$, in $j_w$ we have

$$p_{w-2} + 2A_{w+1}s_{w-2} + (1-2r_2\Psi_{w-3})m_w + 2x_w + 2\eta_{w+1}\eta_{\alpha_w}(1-R_1)((1-r_2)X_{w-3} - r_2\Psi_{w-3})(1-m_w)m_d$$
$$= p_{w-2}^* + 2A_{w+1}s_{w-2}^* + (1-2r_2^*\Psi_{w-3}^*)m_w^* + 2x_w^* + 2\eta_{w+1}\eta_{\alpha_w}(1-R_1^*)((1-r_2^*)X_{w-3}^* - r_2^*\Psi_{w-3}^*)(1-m_w^*)m_d^*$$
$$\implies A_{w+1}(s_{w-2} - s_{w-2}^*) - r_2 m_w + (x_w - x_w^*) + (1-2r_2)(1-m_w)m_d = \frac{p_{w-2}^* - p_{w-2}}{2}.$$

Since $A_{w+1}(s_{w-2} - s_{w-2}^*) - r_2 m_w + (x_w - x_w^*) + (1-2r_2)(1-m_w)m_d \in \mathbb{Z}$ and $|p_{w-2}^* - p_{w-2}| \leq 1$, it must be the case $p_{w-2} = p_{w-2}^*$. Hence, observing that $x_w^* = (1-r_2^*)(A_{w+1}-1)$ since $R_1^* = 1$, we are left with

$$(s_{w-2} - s_{w-2}^*) - (1 - r_2^*) = -\frac{x_w + (1-r_2)}{A_{w+1}}.$$

Following from $(s_{w-2} - s_{w-2}^*) - (1 - r_2^*) \in \mathbb{Z}$ and $1-r_2 \leq x_w + (1-r_2) \leq A_{w+1} - r_2$ by construction, we conclude that $x_w = (1-r_2)(A_{w+1}-1)$.

Thus, by the Principle of Strong Mathematical Induction, we have that for every smallest $3 \leq w \leq k-1$ such that $x_w \neq (1-r_2)(A_{w+1}-1)$, it is the case that for every $3 \leq z \leq w-1$, $x_z = (1-r_2)(A_{z+1}-1)$.

Now, during every $j_k$ move for $4 \leq k \leq d$, we see that for the smallest integer $3 \leq w \leq k-1$ such that $x_w \neq (1-r_2)(A_{w+1}-1)$, our induction argument above implies that for all $3 \leq z \leq w-1$, $x_z = (1-r_2)(A_{z+1}-1)$. By an analogous argument to that made in the induction step with $\eta_{\alpha_w} = 1 = \eta_{\alpha_{w-1}}$ since $\alpha^* > \alpha_{k-1} \geq \alpha_w > \alpha_{w-1}$, in $j_w$ we see

$$p_{w-2}+2A_{w+1}s_{w-2}+(1-2r_2\Psi_{w-3})m_w+2x_w+2\eta_{w+1}\eta_{\alpha_w}(1-R_1)((1-r_2)X_{w-3}-r_2\Psi_{w-3})(1-m_w)m_d$$
$$= p^*_{w-2}+2A_{w+1}s^*_{w-2}+(1-2r^*_2\Psi^*_{w-3})m^*_w+2x^*_w+2\eta_{w+1}\eta_{\alpha_w}(1-R^*_1)((1-r^*_2)X^*_{w-3}-r^*_2\Psi^*_{w-3})(1-m^*_w)m^*_d$$
$$\implies A_{w+1}(s_{w-2} - s^*_{w-2}) - r_2 m_w + (x_w - x^*_w) + (1-2r_2)(1-m_w)m_d = \frac{p^*_{w-2} - p_{w-2}}{2}.$$

Given $A_{w+1}(s_{w-2} - s^*_{w-2}) - r_2 m_w + (x_w - x^*_w) + (1-2r_2)(1-m_w)m_d \in \mathbb{Z}$ and $|p^*_{w-2} - p_{w-2}| \leq 1$, it must be the case $p_{w-2} = p^*_{w-2}$. Hence, observing that $x^*_w = (1-r^*_2)(A_{w+1}-1)$ since $R^*_1 = 1$, we are left with

$$(s_{w-2} - s^*_{w-2}) - (1 - r^*_2) = -\frac{x_w + (1-r_2)}{A_{w+1}}.$$

Since $(s_{w-2} - s^*_{w-2}) - (1-r^*_2) \in \mathbb{Z}$ and $x_w \neq (1-r_2)(A_{w+1}-1)$ by construction and our case assumption, it follows that $-\frac{x_w+(1-r_2)}{A_{w+1}} \notin \mathbb{Z}$, giving us a contradiction as the above equality is not possible. Further, since $3 \leq w \leq k-1$ is the arbitrary smallest integer assumed to satisfy $x_w \neq (1-r_2)(A_{w+1}-1)$, we have that this holds for all such $w$. Thus, this scenario never occurs as such a $w$ does not exist under these assumptions.

<u>Case 3:</u> Let $R_1 = 1$. Then, $C_1$ is being defined by column transitions along $j_k$ for $3 \leq k \leq d$, and so we must have $\alpha_{k-1} < \alpha^* \leq \alpha_d$. Further, we know that $x_z = (1-r_2)(A_{z+1}-1)$ for all $3 \leq z \leq d$. In this case, we only need to consider when $C_2$ is being defined by column transitions as the other subcases follow by symmetric arguments as before.

Hence, let $R^*_1 = 1$ so that $x^*_z = (1-r^*_2)(A_{z+1}-1)$ for all $3 \leq z \leq d$, and observe that during such moves $m_1 = \cdots = m_{k-1}$, $m_k = \cdots = m_d = 1 - m_1$ and $(m^*_1, \ldots, m^*_d) = ((1-|r^*_2 - r_2|)m_1 + |r^*_2 - r_2|(1-m_1), \ldots, (1-|r^*_2 - r_2|)m_d + |r^*_2 - r_2|(1-m_d))$. Applying this to $j_1$, we get

$$(\ell + 1)A_3 + (-2 + m_1 + 2m_d(1-m_1)) = (\ell_1 + 1)A_3 + (-2 + m^*_1 + 2m^*_d(1-m^*_1))$$
$$\implies \ell - \ell_1 = \frac{m^*_d - m_d}{A_3}.$$

Since $\ell - \ell_1 \in \mathbb{Z}$ and $|m^*_d - m_d| \leq 1$ with $A_3 \geq 4$ by construction, it follows that $m_d = m^*_d$, meaning $\ell = \ell_1$ and $(m_1, \ldots, m_d) = (m^*_1, \ldots, m^*_d)$. Note that this implies $r_2 = r^*_2$. In $j_2$, our last deduction yields

$$m_2 + \gamma + 2x_2 + 2m_1(1-m_2) + 2((1-m_1)m_d - (1-R_1)((1-m_2)m_1 + m_d(1-m_1)))$$
$$= m^*_2 + \gamma_1 + 2x^*_2 + 2m^*_1(1-m^*_2) + 2((1-m^*_1)m^*_d - (1-R^*_1)((1-m^*_2)m^*_1 + m^*_d(1-m^*_1)))$$
$$\implies x_2 - x^*_2 = \frac{\gamma_1 - \gamma}{2}.$$

Given $x_2 - x_2^* \in \mathbb{Z}$ and $|\gamma_1 - \gamma| \leq 1$ by construction, it must be the case $\gamma = \gamma_1$ and so $x_2 = x_2^*$, meaning $r_2 = r_2^*$ as we previously deduced.

To conclude our treatment of all cases where we did not arrive at a contradiction in this set comparison case, we would like to bring to the attention of the reader that in said cases we ended with the following results: $\ell = \ell_1$, $(m_1, \ldots, m_d) = (m_1^*, \ldots, m_d^*)$, $\gamma = \gamma_1$, $x_2 = x_2^*$ and $r_2 = r_2^*$. Further, note that during these cases, it was also the case that $R_1 = R_1^*$ and $R_{k^*-1} = R_{k^*-1}^*$ for $3 \leq k^* \leq d$ and $\alpha_2 < \alpha^* \leq \alpha_d$ as applicable for a given $j_k$ move for $1 \leq k \leq d$.

Since we showed that for all $\alpha_2 < \alpha^* \leq \alpha_d$ having equalities with $R_1 \neq R_1^*$ (and $R_{k-1} \neq R_{k-1}^*$ when applicable) would lead to contradictions, we may safely assume the following as necessary: $R_1 = R_1^*$, $R_{k^*-1} = R_{k^*-1}^*$ for $3 \leq k^* \leq d$, $\ell = \ell_1$, $(m_1, \ldots, m_d) = (m_1^*, \ldots, m_d^*)$, $\gamma = \gamma_1$, $x_2 = x_2^*$ and $r_2 = r_2^*$.

These assumptions make it so that we are considering the same type of move in $C_1$ and $C_2$ as applicable based on $k$ and $\alpha^*$. Since our assumptions are such that the moves along $j_k$ for $1 \leq k \leq d$ are indistinguishable in our equalities, we will case on $\alpha_{y-1} < \alpha^* \leq \alpha_y$ for $3 \leq y \leq d$ and consider any $j_k$ moves as applicable.

Letting $\alpha_{y-1} < \alpha^* \leq \alpha_y$ for $3 \leq y \leq d$, we case on whether $3 \leq z \leq y-1$ or $y \leq z \leq d$:

$z-$Case 1: Let $3 \leq z \leq y-1$. Then, for all $j_z$, our assumptions imply that $\eta_{\alpha_z} = 1 = \eta_{\alpha_{z-1}}$ since $\alpha^* > \alpha_{y-1} \geq \alpha_z > \alpha_{z-1}$ and $\eta_{z+1} = 1$ since $z \leq y-1 < d$.

We now set out to show that for all $3 \leq y \leq d$ and $3 \leq z \leq y-1$, our assumptions give us $x_z = x_z^*$, $p_{z-2} = p_{z-2}^*$, and $s_{z-2} = s_{z-2}^*$. We proceed by induction on $y$:

Base Case 1: Let $y = 3$. Then, the statement is vacuously true as there does not exist an integer $z$ such that $3 \leq z \leq 2$.

Base Case 2: Let $y = 4$. Then, $z = 3$, meaning $\eta_{\alpha_3} = 1 = \eta_{\alpha_2}$, and so our assumptions in $j_3$ yield

$$p_1 + 2A_4 s_1 + (1 - 2r_2 \Psi_0) m_3 + 2x_3 + 2\eta_{\alpha_3}(1 - R_1)((1 - r_2)X_0 - r_2 \Psi_0)(1 - m_3) m_d$$
$$= p_1^* + 2A_4 s_1^* + (1 - 2r_2^* \Psi_0^*) m_3^* + 2x_3^* + 2\eta_{\alpha_3}(1 - R_1^*)((1 - r_2^*)X_0^* - r_2^* \Psi_0^*)(1 - m_3^*) m_d^*$$
$$\implies A_4(s_1 - s_1^*) + (x_3 - x_3^*) = \frac{p_1^* - p_1}{2}.$$

Since $A_4(s_1 - s_1^*) + (x_3 - x_3^*) \in \mathbb{Z}$ and $|p_1^* - p_1| \leq 1$ by construction, it must be that $p_1 = p_1^*$. This leaves us with

$$s_1 - s_1^* = \frac{x_3^* - x_3}{A_4}.$$

Given $s_1 - s_1^* \in \mathbb{Z}$ and $|x_3^* - x_3| \leq A_4 - 1$ by construction, it follows that $x_3 = x_3^*$ and so $s_1 = s_1^*$.

Induction Step: Let $3 \leq y \leq d$ and $d \geq 3$ be such that $3 \leq y+1 \leq d$ so that the following makes sense to consider. Further, assume that for every $3 \leq y^* \leq y$ and $3 \leq z \leq y^*-1$, $p_{z-2} = p_{z-2}^*$, $x_z = x_z^*$, and $s_{z-2} = s_{z-2}^*$. We will show our statement holds for $3 \leq y+1 \leq d$ by showing $p_{y-2} = p_{y-2}^*$, $x_y = x_y^*$ and $s_{y-2} = s_{y-2}^*$ as our induction hypothesis already grants us $p_{z-2} = p_{z-2}^*$, $x_z = x_z^*$ and $s_{z-2} = s_{z-2}^*$ for all $3 \leq z \leq y-1$. Hence, observing that $r_2 \Psi_{y-3} = r_2^* \Psi_{y-3}^*$, $(1 - r_2)X_{y-3} = (1 - r_2^*)X_{y-3}^*$ and $\eta_{\alpha_{y-1}} = 1 = \eta_{\alpha_y}$ since $\alpha^* > \alpha_y > \alpha_{y-1}$, in $j_y$ we have

$$p_{y-2} + 2A_{y+1}s_{y-2} + (1 - 2r_2\Psi_{y-3})m_y + 2x_y + 2\eta_{y+1}\eta_{\alpha_y}(1 - R_1)((1 - r_2)X_{y-3} - r_2\Psi_{y-3})(1 - m_y)m_d$$

$$= p^*_{y-2} + 2A_{y+1}s^*_{y-2} + (1 - 2r^*_2\Psi^*_{y-3})m^*_y + 2x^*_y + 2\eta_{y+1}\eta_{\alpha_y}(1 - R^*_1)((1 - r^*_2)X^*_{y-3} - r^*_2\Psi^*_{y-3})(1 - m^*_y)m^*_d$$

$$\implies A_{y+1}(s_{y-2} - s^*_{y-2}) + (x_y - x^*_y) = \frac{p^*_{y-2} - p_{y-2}}{2}.$$

Since $A_{y+1}(s_{y-2} - s^*_{y-2}) + (x_y - x^*_y) \in \mathbb{Z}$ and $|p^*_{y-2} - p_{y-2}| \leq 1$, it must be the case $p_{y-2} = p^*_{y-2}$. Hence, we are left with

$$s_{y-2} - s^*_{y-2} = \frac{x^*_y - x_y}{A_{y+1}}.$$

Following from $s_{y-2} - s^*_{y-2} \in \mathbb{Z}$ and $|x^*_y - x_y| \leq A_{y+1} - 1$ by construction, we conclude that $x_y = x^*_y$ and so $s_{y-2} = s^*_{y-2}$.

Thus, by the Principle of Strong Mathematical Induction, we have that for every $3 \leq y \leq d$ and $3 \leq z \leq y - 1$, it is the case that $p_{z-2} = p^*_{z-2}$, $x_z = x^*_z$ and $s_{z-2} = s^*_{z-2}$ when $\alpha_{y-1} < \alpha^* \leq \alpha_y$.

Now, focusing on all applicable moves along $j_k$ for $3 \leq k \leq d$ and $\alpha_{y-1} < \alpha^* \leq \alpha_y$, our induction argument above gives us that for all $3 \leq y \leq d$ and $3 \leq z \leq y - 1$, $p_{z-2} = p^*_{z-2}$, $x_z = x^*_z$ and $s_{z-2} = s^*_{z-2}$.

$\underline{z\text{--Case 2}}$: Let $y \leq z \leq d$. Recalling from the previous $z$--case that $x_{z^*} = x^*_{z^*}$ for all $3 \leq z^* \leq y - 1$, we know in particular $\Psi_{y-3} = \Psi^*_{y-3}$. Since $\alpha_{y-1} < \alpha^* \leq \alpha_y$, it follows that $p_z = 0 = p^*_z$ and $x_z = 0 = x^*_z$ for all $y + 1 \leq z \leq d$. Observe that our observations above, along with $r_2 = r^*_2$, $\eta_{\alpha_{y-1}} = 1$ and $\eta_{\alpha_y} = 0$, in $j_y$ yield

$$p_{y-2} + 2A_{y+1}s_{y-2} + (1 - 2r_2\Psi_{y-3})m_y + 2x_y = p^*_{y-2} + 2A_{y+1}s^*_{y-2} + (1 - 2r^*_2\Psi^*_{y-3})m^*_y + 2x^*_y$$

$$\implies A_{y+1}(s_{y-2} - s^*_{y-2}) + (x_y - x^*_y) = \frac{p^*_{y-2} - p_{y-2}}{2}.$$

Since $A_{y+1}(s_{y-2} - s^*_{y-2}) + (x_y - x^*_y) \in \mathbb{Z}$ and $|p^*_{y-2} - p_{y-2}| \leq 1$, it must be the case $p_{y-2} = p^*_{y-2}$. Hence, we are left with

$$s_{y-2} - s^*_{y-2} = \frac{x^*_y - x_y}{A_{y+1}}.$$

Following from $s_{y-2} - s^*_{y-2} \in \mathbb{Z}$ and $|x^*_y - x_y| \leq A_{y+1} - 1$ by construction, we conclude that $x_y = x^*_y$ and so $s_{y-2} = s^*_{y-2}$. At this point, we have $p_{z-2} = p^*_{z-2}$, $x_z = x^*_z$ and $s_{z-2} = s^*_{z-2}$ for all $3 \leq z \leq y$.

Note that $p_{z-2} = 0 = p^*_{z-2}$, $A_{z+1} = 1$ and $x_z = 0 = x^*_z$ for all $y + 1 \leq z \leq d$ since $\alpha_{z-1} \geq \alpha_y \geq \alpha^*$. Hence, to conclude our treatment of this set comparison case, we must show $s_{z-2} = s^*_{z-2}$ for all $y + 1 \leq z \leq d$. Observing that $\eta_{\alpha_z} = 0 = \eta_{\alpha_{z-1}}$ since $\alpha_z > \alpha_{z-1} \geq \alpha_y > \alpha^*$, for $y + 1 \leq z \leq d$ in $j_z$ we obtain

$$p_{z-2} + A_{z+1}s_{z-2} + m_z + 2x_z = p^*_{z-2} + A_{z+1}s^*_{z-2} + m^*_z + 2x^*_z$$

$$\implies s_{z-2} = s^*_{z-2}.$$

Consequently, for all $y+1 \leq z \leq d$, $s_{z-2} = s^*_{z-2}$.

With this, for $d \geq 3$ and $\alpha_2 < \alpha^* \leq \alpha_d$, we have shown that $C_1$ and $C_2$ share an edge if and only if $(m_1, \ldots, m_d) = (m_1^*, \ldots, m_d^*)$ from the same start vertex, which is if and only if $C_1 = C_{\ell,\gamma,t,p_1,s_1,\ldots,p_{d-2},s_{d-2}} = C_{\ell_1,\gamma_1,t_1,p_1^*,s_1^*,\ldots,p_{d-2}^*,s_{d-2}^*} = C_2$.

Thus, for all $\alpha_1 \leq \alpha^* \leq \alpha_d$ with $d \geq 2$, $C_1$ and $C_2$ share an edge if and only if $C_1 = C_2$, proving Proposition 17.

∎

Combining the proofs of Propositions 15, 16 and 17, we have our proof of Theorem 14.

## 7.2 Proof of all Cycles being of the Same Length:

Here, we will show that all cycles defined by the edge set for the Lock-and-Key Decomposition are of the same length. We will assume the presence of the parameters $\ell, \gamma, t, p_1, s_1, \ldots, p_{d-2}, s_{d-2}$ as applicable based on the dimension $d$ of the torus of the form

$$C_{z_1} \square C_{z_d} \square C_{z_2} \square \cdots \square C_{z_{d-1}}.$$

Then, it suffices to show, for each $\alpha^*$–case of a given dimension case, that the cycles must all be of the same length given the number of edges defining the torus, the number of possible cycles, and our previous result that all cycles are edge-disjoint.

This will lead to the following result:

**Theorem 18** *Every cycle defined by the General Lock-and-Key Decomposition's set is of length $\frac{dz_d}{2} \prod_{k=0}^{d-2} A_{k+3}$, without distinguishing between partitioned and non-partitioned edges along each dimension.*

**Proof:** Let $d \in \mathbb{Z}^{\geq 2}$ and $\alpha_1 \leq \alpha^* \leq \alpha_d$. Note that the number of edges defining the torus of the form presented above is

$$d \prod_{k=1}^{d} z_k.$$

Given our last result established that the cycles defined by the Lock-and-Key decomposition are edge-disjoint, all we must show is that they are of the same length given the number of cycles defined and the number of edges in the torus to assign to each one.

**Dimension Case 1:** Let $d = 2$. Then, only the parameters $\ell$ and $\gamma$ are present. Note that $\eta_3 = 0$. We now further case on $\alpha^*$:

**$\alpha^*$−Case 1:** Let $\alpha^* = \alpha_1$. Then, $A_1 = 0$, $A_2 = 0$, and $A_3 = 2$. So we see that there are $\frac{z_1}{A_3} = \frac{z_1}{2}$ choices for $\ell$ and 2 for $\gamma$, meaning that all cycles must be of length

$$\frac{d}{\frac{z_1}{2} \cdot 2} \prod_{k=1}^{d} z_k = 2z_2.$$

**$\alpha^*$−Case 2:** Let $\alpha_1 < \alpha^* < \alpha_2$. Then, $A_1 = 1$, $A_2 = 0$, and $A_3 = 2(b_1+1)$. From this, we get that there are $\frac{z_1}{A_3} = \frac{z_1}{2(b_1+1)}$ choices for $\ell$ and 2 for $\gamma$, giving us that all cycles must be of length

$$\frac{d}{\frac{z_1}{2(b_1+1)} \cdot 2} \prod_{k=1}^{d} z_k = 2z_2(b_1+1).$$

**$\alpha^*$−Case 3:** Let $\alpha^* = \alpha_2$. Then, $A_1 = 1$, $A_2 = 1$, and $A_3 = z_1$. Noting that there are $\frac{z_1}{A_3} = 1$ choices for $\ell$ and 2 for $\gamma$, we have that all cycles are of length

$$\frac{d}{1 \cdot 2} \prod_{k=1}^{d} z_k = 2\left(\frac{z_1 z_2}{2}\right).$$

**Dimension Case 2:** Let $d \geq 3$. Then, the parameters $\ell, \gamma, t, p_1, s_1, \ldots, p_{d-2}$, and $s_{d-2}$ are all present as applicable based on $d$. We case on $\alpha^*$:

**$\alpha^*$−Case 1:** Let $\alpha^* = \alpha_1$. Then, $A_1 = 0$, $A_2 = 0$, $A_3 = 2$, $\eta_3 = 1$, and $\eta_{\alpha_k} = 0$ for all $2 \leq k \leq d-1$. Consequently, we see that there are $\frac{z_1}{A_3} = \frac{z_1}{2}$ choices for $\ell$, 2 for $\gamma$, 2 for $t$, 1 for each $p_k$ for all $1 \leq k \leq d-2$, $\frac{z_2}{2A_4} = \frac{z_2}{2}$ for $s_1$ and $\frac{z_{j+1}}{A_{j+3}} = z_{j+1}$ choices for each $s_j$ for all $2 \leq j \leq d-2$, we find that all cycles are of length

$$\frac{d \prod_{k=1}^{d} z_k}{\frac{z_1 z_2}{2^2} \cdot 2^2 \prod_{k=3}^{d-1} z_k} = dz_d.$$

**$\alpha^*$−Case 2:** Let $\alpha_1 < \alpha^* < \alpha_2$. Then, $A_1 = 1$, $A_2 = 0$, $A_3 = 2(b_1+1)$, $\eta_3 = 1$, and $\eta_{\alpha_k} = 0$ for all $2 \leq k \leq d-1$. Hence, it follows that there are $\frac{z_1}{A_3} = \frac{z_1}{2(b_1+1)}$ choices for $\ell$, 2 for $\gamma$, 1 for $t$, 1 for each $p_k$ for all $2 \leq k \leq d-2$ and $\frac{z_{j+1}}{A_{j+3}} = z_{j+1}$ for each $s_j$ for all $1 \leq j \leq d-2$, and so all cycles are of length

$$\frac{d \prod_{k=1}^{d} z_k}{\frac{z_1 z_2}{2(b_1+1)} \cdot 2 \prod_{k=3}^{d-1} z_k} = dz_d(b_1+1).$$

**$\alpha^*$−Case 3:** Let $\alpha^* = \alpha_2$. Then, $A_1 = 1$, $A_2 = 1$, $A_3 = z_1$, $\eta_3 = 1$, and $\eta_{\alpha_k} = 0$ for all $2 \leq k \leq d-1$. Now, we observe that there is $\frac{z_1}{A_3} = 1$ choice for $\ell$, 2 for $\gamma$, 1 for $t$, 1 for each $p_k$ for all $2 \leq k \leq d-2$, and $\frac{z_{j+1}}{A_{j+3}} = z_{j+1}$ for each $s_j$ for all $1 \leq j \leq d-2$. These observations tell us that all cycles are of length

$$\frac{d \prod_{k=1}^{d} z_k}{2 \prod_{k=2}^{d-1} z_k} = \frac{dz_1 z_d}{2}.$$

**α\*−Case 4:** Let $\alpha_{k-1} < \alpha^* \leq \alpha_k$ for $3 \leq k \leq d$. Then, $A_1 = 1 = A_2$, $A_3 = z_1$, $\eta_{\alpha_2} = \cdots = \eta_{\alpha_{k-1}} = 1$, and $\eta_{\alpha_{d-1}} = \cdots = \eta_{\alpha_k} = 0$. So we have that there is $\frac{z_1}{A_3} = 1$ choice for $\ell$, 2 for $\gamma$, 1 for $t$, 2 for each $p_j$ for all $1 \leq j \leq k-2$, 1 for each $p_j$ for all $k-1 \leq j \leq d-2$, $\frac{z_{j+1}}{2A_{j+3}} = 1$ for each $s_j$ for all $1 \leq j \leq k-3$, $\frac{z_{k-1}}{2A_{k+1}} = \frac{z_{k-1}}{2(b_{k-1}+1)}$ for $s_{k-2}$, and $\frac{z_{j+1}}{A_{j+3}} = z_{j+1}$ for each $s_j$ for all $k-1 \leq j \leq d-2$. The above tells us that all cycles must be of length

$$\frac{d \prod_{j=1}^{d} z_j}{2 \cdot 2^{k-2} \cdot \frac{z_{k-1}}{2(b_{k-1}+1)} \prod_{j=k}^{d-1} z_j} = \frac{dz_d(b_{k-1}+1)}{2^{k-2}} \prod_{j=1}^{k-2} z_j.$$

From the above and the definition of $A_{k+3}$ for $0 \leq k \leq d-2$, we conclude that all cycles defined by the Lock-and-Key decomposition's edge set definition are of length

$$\frac{dz_d}{2} \prod_{k=0}^{d-2} A_{k+3}.$$

Towards showing that the cycles in the decomposition of the underlying torus correspond to cycles decomposing the subdivided torus and are all of the same length with lengths of the desired form, we must show that the same number of edges are used along every dimension in defining a given cycle. This is an important consideration as otherwise having cycles of the same length decomposing an underlying torus would not translate to cycles decomposing the corresponding subdivided torus using cycles with all of the same length.

To see that our cycles indeed have the desired properties, we can observe that in all of our major sets the fundamental components of the cycles they define, illustrated in the discussion of the major sets of the Lock-and-Key decomposition, are defined by the same number of edges along every dimension. In particular, we mean that every edge belongs to a path of length $d$ with exactly one edge coming from each dimension and each path is a member of an edge-disjoint decomposition of the given cycle into paths as can be seen by parsing a given cycle by every $d$ edges. From this, we get that the resulting cycle has the same number of edges along every dimension, revealing another way in which the Lock-and-Key decomposition is symmetric.

Hence, dividing all $\alpha^*$-case lengths by the corresponding dimension $d$, we have that every cycle defined by the Lock-and-Key decomposition has

$$\frac{1}{d} \cdot \frac{dz_d}{2} \prod_{k=0}^{d-2} A_{k+3} = \frac{z_d}{2} \prod_{k=0}^{d-2} A_{k+3}$$

edges along each dimension.

Now, let $n \in \mathbb{Z}^+$ and $a = 2^{i_1} + \cdots + 2^{i_d}$ with $i_1 > \cdots > i_d = 0$. All that remains is to show that in all cases, the cycles decomposing the torus

$$C_{z_1} \square C_{z_d} \square C_{z_2} \square \cdots \square C_{z_{d-1}} = C_{2^{n2^{i_1+1}-i_1}} \square C_{2^{2n}} \square C_{2^{n2^{i_2+1}-i_2}} \square \cdots \square C_{2^{n2^{i_{d-1}+1}-i_{d-1}}}$$

are of length $a \cdot 2^\alpha$. Note that in each case we will be making use of Proposition 1 with each edge coming from $C_{z_j}$ counting as $2^{i_j}$ edges for all $1 \leq j \leq d$ to distinguish between partitioned and non-partitioned edges. Lastly, assume that $\alpha$ and $\alpha_k$ for $1 \leq k \leq d$ are as in the remark made following the Lock-and-Key Decomposition's edge set definition.

**Dimension Case 1:** Let $d = 2$. Then, $a = 2^{i_1} + 2^{i_2}$ and so we proceed by casing on $\alpha$:

**$\alpha$−Case 1:** Let $\alpha = \alpha_1 = 2n$. Recalling that all cycles are of length $2z_2$ in this case, we get that all cycles are of length

$$\sum_{j=1}^{2} z_2 2^{i_j} = a \cdot 2^{2n} = a \cdot 2^{\alpha}$$

when distinguishing between partitioned and non-partitioned edges.

**$\alpha$−Case 2:** Let $\alpha_1 < \alpha < \alpha_2$. Noting that all cycles are of length $2z_2(b_1 + 1)$ with $b_1 = 2^{\alpha-\alpha_1} - 1$, we find that all cycles are of length

$$\sum_{j=1}^{2} z_2(b_1 + 1) 2^{i_j} = a \cdot 2^{2n} \cdot 2^{\alpha-2n} = a \cdot 2^{\alpha}$$

when distinguishing between partitioned and non-partitioned edges.

**$\alpha$−Case 3:** Let $\alpha = \alpha_2 = 2an - (i_1 + 1)$. Observing that all cycles are of length $2(z_1 z_2/2)$ in this case, we deduce that all cycles are of length

$$\sum_{j=1}^{2} z_1 z_2 2^{i_j - 1} = \frac{a}{2} \cdot 2^{n 2^{i_1+1} - i_1} \cdot 2^{2n} = a \cdot 2^{2an-(i_1+1)} = a \cdot 2^{\alpha}$$

when distinguishing between partitioned and non-partitioned edges.

**Dimension Case 2:** Let $d \geq 3$. Then, $a = 2^{i_1} + \cdots + 2^{i_d}$. Below we case on $\alpha$ once again:

**$\alpha$−Case 1:** Let $\alpha = \alpha_1 = 2n$. Invoking the result from before that all cycles are of length $dz_d$ in this case, we have that all cycles are of length

$$\sum_{j=1}^{d} z_d 2^{i_1} = a \cdot 2^{2n} = a \cdot 2^{\alpha}$$

when distinguishing between partitioned and non-partitioned edges.

**$\alpha$−Case 2:** Let $\alpha_1 < \alpha < \alpha_2$. Given that all cycles are of length $dz_d(b_1 + 1)$ and $b_1 = 2^{\alpha-\alpha_1} - 1$, we have that all cycles are of length

$$\sum_{j=1}^{d} z_d(b_1 + 1) 2^{i_j} = a \cdot 2^{2n} \cdot 2^{\alpha-2n} = a \cdot 2^{\alpha}$$

when distinguishing between partitioned and non-partitioned edges.

**$\alpha$−Case 3:** Let $\alpha = \alpha_2 = 2(2^{i_1} + 1)n - (i_1 + 1)$. Following from our result that all cycles are of length $dz_1 z_d/2$, we conclude that all cycles are of length

$$\sum_{j=1}^{d} z_1 z_d 2^{i_j - 1} = \frac{a}{2} \cdot 2^{n 2^{i_1+1} - i_1} \cdot 2^{2n} = a \cdot 2^{2(2^{i_1}+1)n - (i_1+1)} = a \cdot 2^{\alpha}$$

when distinguishing between partitioned and non-partitioned edges.

**$\alpha$−Case 4:** Let $\alpha_{k-1} < \alpha \leq \alpha_k$ for all $3 \leq k \leq d$. Given all cycles are of length

$$\frac{d z_d (b_{k-1}+1)}{2^{k-2}} \prod_{j=1}^{k-2} z_j$$

and $b_{k-1} = 2^{\alpha - \alpha_{k-1}} - 1$ with $\alpha_{k-1} = \sum_{j=1}^{k-2} \left[n 2^{i_j+1} - i_j\right] + 2n - (k-2)$, we conclude that all cycles are of length

$$\sum_{j=1}^{d} \frac{z_d(b_{k-1}+1) 2^{i_j}}{2^{k-2}} \prod_{t=1}^{k-2} z_t = a \cdot 2^{2n} \cdot 2^{\alpha - \alpha_{k-1} - (k-2)} \cdot 2^{\sum_{t=1}^{k-2}\left[n 2^{i_t+1} - i_t\right]} = a \cdot 2^{\alpha}$$

when distinguishing between partitioned and non-partitioned edges.

Thus, we have that $Q_{2an}$ can be decomposed into cycles of length $a \cdot 2^{\alpha}$ for all

$$2n \leq \alpha \leq 2an - (d-1) - \sum_{k=1}^{d} i_k.$$

∎

**Open Problem:** The current question of great interest to us moving forward is hence the following:

Under what set of necessary-and-sufficient conditions can $Q_{2an}$ be decomposed into cycles of length $a \cdot 2^{\alpha}$ for

$$2an - (d-1) - \sum_{k=1}^{d} i_k < \alpha \leq 2an - \lceil \log_2(a) \rceil$$

when $d \geq 3$, where $n \geq 1$ and $a \geq 3$ is an odd integer with at least three powers of two in its binary representation. The author believes the explicit cycle decompositions obtained here using the longest cycles defined by the Lock-and-Key decomposition reveal a derivative decomposition method by which $Q_{2an}$ as an $a$-fold $Q_{2n}$ Cartesian product can be decomposed into cycles of the remaining lengths. This would be leveraging the fact that there are the same number of edges along each dimension of the torus since the Cartesian product is being taken against cycles all of the same length. Consequently, in contrast to the anchored product setting, there is more freedom in the kinds of symmetries one can introduce into the cycle definition as all edges are weighted equally, allowing for further continuation of the cycles.

## Declaration of Competing Interest

The author declares that they have no known competing financial interests or personal relationships that could have appeared to influence the work reported in this paper.

## Acknowledgements

The author would like to begin by thanking Carnegie Mellon University professor Dr. David Offner for his continual support throughout this research project and the constructive feedback in the presentation of this paper. This research project began in the summer of 2022 under the supervision of Dr. Offner as part of Carnegie Mellon University's Summer Undergraduate Research Fellowhsip (SURF) program, supported by the generosity of the Porges Family Fund, and was subsequently completed independently with the advising of Dr. Offner. In closing, the author wishes to extend their gratitude to Dr. Aris Winger for his mentorship in the process to starting this research project and meeting Dr. Offner.

## Data Availability

No data was used for the research described in the article.